\pdfoutput=1
\documentclass[a4paper,11pt]{amsart}

\usepackage{amsmath,amssymb,amsthm,mathtools}
\usepackage{enumitem}
\usepackage{etoolbox}
\usepackage[hidelinks]{hyperref}
\hypersetup{pdftitle={Finite-order independence: stability, lifting, and strong convergence},
  pdfauthor={Ethan Davis}}

\makeatletter
\patchcmd{\@setabstracta}{\box\abstractbox}{\unvbox\abstractbox}{}%
  {\errmessage{Unable to enable page breaks in the AMS abstract}}
\patchcmd{\@maketitle}{\@mkboth{\@nx\shortauthors}{\@nx\shorttitle}}%
  {\@mkboth{\MakeUppercase{\shortauthors}}{\shorttitle}}{}{}
\makeatother

\AtBeginDocument{%
  \setlength{\abovedisplayskip}{9pt plus 3pt minus 6pt}%
  \setlength{\belowdisplayskip}{9pt plus 3pt minus 6pt}%
  \setlength{\belowdisplayshortskip}{5.5pt plus 3pt minus 3pt}%
}

\numberwithin{equation}{section}
\newcommand{\eqtag}{\refstepcounter{equation}\tag{\theequation}}

\newtheorem{theorem}{Theorem}[section]
\newtheorem{lemma}[theorem]{Lemma}
\newtheorem{proposition}[theorem]{Proposition}
\newtheorem{corollary}[theorem]{Corollary}
\newtheorem{fact}[theorem]{Fact}
\newtheorem*{problemthirty}{Problem~30}
\newtheorem*{problemthirtyone}{Problem~31}
\newtheorem{introthm}{Theorem}

\newtheorem{introcor}[introthm]{Corollary}

\newcommand{\N}{\mathbb{N}}
\newcommand{\R}{\mathbb{R}}
\newcommand{\C}{\mathbb{C}}
\newcommand{\norm}[1]{\left\lVert #1\right\rVert}
\newcommand{\abs}[1]{\left\lvert #1\right\rvert}
\newcommand{\Ree}{\operatorname{Re}}
\newcommand{\Imm}{\operatorname{Im}}
\newcommand{\cU}{\mathcal U}
\newcommand{\ultraprod}{\prod_{\cU}}

\title[Finite-order independence]{Finite-order independence: stability, lifting, and strong convergence}
\author{Ethan Davis}
\address{Department of Mathematics, California Institute of Technology}
\email{erdavis@caltech.edu}
\thanks{The author was partially supported by NSF grant DMS-2104330.}
\date{}

\begin{document}
\raggedbottom
\begin{abstract}
We study finite-order independence under perturbation. Following Popa, we say that two unital von Neumann subalgebras $A,B$ of a tracial von Neumann algebra $(M,\tau)$ are $n$-independent if $\tau(a_1b_1\cdots a_rb_r)=0$ for $1\le r\le n$ and centered $a_i\in A$, $b_i\in B$. Thus, $1$-independence is orthogonality, and $n$-independence for every $n$ amounts to the usual notion of freeness. Our first perturbation theorem is for finite-dimensional $A,B$ in a $\mathrm{II}_1$ factor, with fixed isomorphism classes and a positive lower bound on the traces of minimal projections. If the moments $\tau(a_1b_1\cdots a_rb_r)$ are sufficiently small, uniformly over centered contractions, then there is a unitary $u$ near $1$ such that $uAu^*$ and $B$ are $n$-independent. Its operator-norm distance from $1$ is at most a constant times a power of the largest absolute value of these moments. The smallness threshold and the constant are uniform over the ambient factors, and the power depends only on $n$. Consequently, in tracial ultraproducts of $\mathrm{II}_1$ factors over countably cofinal ultrafilters, separable hyperfinite $n$-independent unital von Neumann subalgebras are contained in ultraproducts of finite-dimensional unital $^*$-subalgebras that are $n$-independent at every coordinate. The order-one case answers, for $\mathrm{II}_1$ factors, a lifting question posed by Houdayer and Ioana and listed by Kunnawalkam Elayavalli: unitaries generating orthogonal von Neumann subalgebras in an ultrapower have unitary representatives whose generated algebras are orthogonal at every coordinate. Moreover, for each fixed $n$, freely independent unitaries in an ultrapower of any tracial von Neumann algebra have unitary representatives generating $n$-independent algebras at every coordinate, a result obtained independently by Patchell and Shiner.

Using the lifting of orthogonality, we show that there is a separable non-Gamma $\mathrm{II}_1$ factor $M$ such that, for every countably cofinal ultrafilter $\mathcal U$ and all unitaries $u_1,u_2\in M^{\mathcal U}$ generating orthogonal von Neumann subalgebras, there are Haar unitaries $v_1,v_2\in M^{\mathcal U}$ with $[u_1,v_1]=[u_2,v_2]=[v_1,v_2]=0$, extending results of Chifan, Ioana, and Kunnawalkam Elayavalli as well as Houdayer and Ioana.

We also study order-one independence for pairs of finite partitions of $1$, allowing zero projections and fixing the two partition sizes. In any fixed tracial von Neumann algebra that is diffuse, of type~I, or the direct sum of a diffuse tracial von Neumann algebra and one scalar atom, a nearly independent pair is $L^2$-close to an independent pair whose projection traces may differ from the original ones. This fails for a suitable direct sum of a $\mathrm{II}_1$ factor and finitely many scalar atoms, giving a negative answer to the order-one lifting question for ultrapowers of this algebra.

Our second perturbation theorem is for matrix models. Suppose that finitely many tuples of matrices, each of fixed length, converge jointly and strongly to freely independent limits, and that, for each tuple and each $d$, the $^*$-polynomials of degree at most $d$ eventually span a space of the same dimension as in the limiting tracial GNS representation. We show that conjugation by unitaries tending to $1$ in operator norm then makes the joint moments agree, through degrees tending to infinity, with those of the free product of the marginal distributions at each dimension, preserving those distributions and the strong limit. We also prove a version with amalgamation over a fixed finite-dimensional unital $C^*$-algebra. For each $0<\gamma<1$, finitely many mutually independent normalized GUE matrices, grouped into $k\ge2$ families, can be conjugated by unitaries within operator-norm distance $N^{-\gamma}$ of $1$ so that, almost surely for all sufficiently large $N$, the joint moments of the families through degree $\lfloor c\log N\rfloor$, for some $c>0$, are those of freely independent families, each with the same distribution as before, with respect to the normalized trace. Logarithmic growth of the degree is optimal.

By the same perturbation method, for each $s\ge2$, we construct irreducible representations of $\mathbb F_s$ in $U(N)$ for every sufficiently large $N$ that send each generator to a conjugate of the cyclic shift and converge strongly to the left regular representation, as in the work of Collins and Male, but with the additional property that every nontrivial word of length at most a constant times $\log N$ has trace zero. Again, logarithmic growth is optimal. For independent Haar unitary matrices, we prove a quantitative Gaussian approximation for the traces of words whose number and lengths grow with the dimension, and we identify the limit laws of the largest of these traces.

For any finite real symmetric matrix $Q=(q_{ab})$ with $|q_{ab}|<1$, we construct independent centered self-adjoint Gaussian matrices that converge almost surely and strongly to the mixed $q$-Gaussian family with parameter matrix $Q$ along a suitable sequence of dimensions. Hence every such family generates an MF $C^*$-algebra. In particular, finite $q$-Gaussian families generate MF $C^*$-algebras for every $-1<q<1$, with $\operatorname{Ext}$ not a group when there are at least two generators, extending the range $|q|<\sqrt2-1$ in which Caspers and Chen constructed strongly convergent matrix models and obtained these conclusions.
\end{abstract}
\maketitle
\tableofcontents

\section{Introduction}

We ask when two unital von Neumann subalgebras whose mixed moments
approximate the free-product values can be perturbed within the ambient
algebra to make the agreement exact through a prescribed order.
For finite-dimensional subalgebras of a $\mathrm{II}_1$ factor, of fixed isomorphism classes and with the traces of their minimal projections bounded below, we show that a small unitary conjugation of one of them suffices when the mixed moments through that same order are sufficiently close to these values, with bounds that are uniform over all $\mathrm{II}_1$ factors (Theorem~\ref{introthm:stability}).

For tracial von Neumann algebras $(M_k,\tau_k)$, a mixed moment in their
tracial ultraproduct can vanish even though its coordinate values are
nonzero at every index. This raises a lifting question for freely independent
von Neumann subalgebras of the ultraproduct: can one choose
representatives satisfying the same mixed-moment identities at every
coordinate? We study finite-order versions of this question by correcting
small errors in the identities with small perturbations inside the
coordinate algebras. In tracial ultraproducts of $\mathrm{II}_1$ factors over countably cofinal ultrafilters, we lift the identities through any prescribed order for separable hyperfinite subalgebras (Theorem~\ref{introthm:lifting}). At order one, we thus answer, for $\mathrm{II}_1$ factors, a question of Houdayer and Ioana~\cite{houdayer-ioana} (Corollary~\ref{introcor:problem-30}), and we extend their theorem on commuting Haar unitaries in ultrapowers of a non-Gamma factor from pairs of unitaries generating $2$-independent von Neumann algebras to all pairs generating orthogonal ones (Theorem~\ref{introthm:non-gamma}). For ultrapowers of a suitable direct sum of a $\mathrm{II}_1$ factor and finitely many scalar atoms, the same question has a negative answer (Theorem~\ref{introthm:tracial}). For each fixed $n$, freely independent unitaries in an ultrapower of any tracial von Neumann algebra have unitary representatives generating $n$-independent von Neumann algebras at every coordinate (Corollary~\ref{cor:free-unitary-lifting}); this answers one reading of a question of Kunnawalkam Elayavalli~\cite[Problem~31]{kunnawalkam-elayavalli-open-problems} and was obtained independently by Patchell and Shiner~\cite[Theorem~C]{patchell-shiner}.

A similar question arises for asymptotically free matrix families~\cite{voiculescu-limit-laws}. Their
mixed moments approach the free-product values, but the identities need
not hold exactly at any finite dimension. We seek small perturbations of the matrices
that make these identities hold exactly through degrees that grow with the
dimension. The perturbations should preserve the
marginal distributions and tend to zero in operator norm, so that the
convergence of polynomial operator norms is retained. We obtain such perturbations for strongly convergent tuples of matrices under a condition on the dimensions of their polynomial spaces (Theorem~\ref{introthm:strong}). For independent GUE matrices, and for representations of free groups obtained from independent Haar unitary matrices, the perturbed matrices satisfy the identities through degrees of order $\log N$ in dimension $N$, and no faster growth is possible (Theorems~\ref{introthm:gue} and~\ref{introthm:clocks}).

By a separate construction, for every real symmetric matrix $Q$ with entries of absolute value less than $1$, we obtain independent Gaussian matrices that converge strongly, almost surely, to the mixed $q$-Gaussian family with parameter matrix $Q$. Hence the $C^*$-algebra generated by a finite $q$-Gaussian family is MF for every $-1<q<1$ (Theorem~\ref{introthm:mixed-q}); Caspers and Chen~\cite{caspers-chen} proved this for $|q|<\sqrt2-1$.

Let $(M,\tau)$ be a tracial von Neumann algebra, with faithful normal
tracial state $\tau$. For $n\ge1$ and unital von Neumann subalgebras
$A,B\subseteq M$, we say, following Popa~\cite{popa-kadison-singer,popa-independence}, that $A$ and $B$ are
$n$-independent if
\[
\tau(a_1b_1\cdots a_rb_r)=0\qquad(1\le r\le n)
\]
for centered $a_i\in A$ and $b_i\in B$. At order one this is
orthogonality in the sense of Popa~\cite[Lemma~2.1 and Definition~2.2]{popa-orthogonal-pairs}:
$\tau(ab)=\tau(a)\tau(b)$ for $a\in A$ and $b\in B$. In this case, we write $A\perp B$.
Requiring all orders gives free independence. Note that if $A$ and $B$ are
$n$-independent, then the trace of every alternating word with at most
$2n+1$
blocks agrees with its trace in the free product of $A$ and $B$ (Lemma~\ref{prop:truncated-free-product}). We
quantify the failure of $n$-independence by the defect
\[
\varepsilon_n(A,B):=
\sup\left\{
\abs{\tau(a_1b_1\cdots a_rb_r)}:
\begin{array}{l}
1\le r\le n,\quad a_i\in A,\ b_i\in B,\\
\tau(a_i)=\tau(b_i)=0,\\
\norm{a_i}_\infty,\norm{b_i}_\infty\le1\quad(1\le i\le r)
\end{array}
\right\}.
\]

\subsection*{The perturbation theorem and lifting}
In a $\mathrm{II}_1$ factor, our first perturbation theorem makes the
order-$n$ identities exact while preserving both marginal tracial states. 

\begin{introthm}[Theorem~\ref{thm:fixed-type-factor-no-loss}]
\label{introthm:stability}
Fix
finite-dimensional unital $C^*$-algebras $P,Q$, an order $n\ge1$, and
$\delta>0$. There are constants $\varepsilon_0>0$ and $C<\infty$, depending only on $n,P,Q,\delta$, with the following property. Let $A,B$ be unital $^*$-subalgebras of a
$\mathrm{II}_1$ factor $M$, $^*$-isomorphic to $P$ and $Q$, respectively, such that every minimal
projection of either algebra has trace at least $\delta$. If
$\varepsilon_n(A,B)\le\varepsilon_0$, then there is a unitary $u\in M$ such
that
$uAu^*$ and $B$ are $n$-independent and
\[
\norm{u-1}_\infty\le C\varepsilon_n(A,B)^{1/(2n^2+2n+1)}.
\]
\end{introthm}

Here only one subalgebra moves, and only the centered
alternating moments of order at most $n$ need be small. In an arbitrary tracial von Neumann algebra, we obtain such a unitary when $\varepsilon_{2n}(A,B)$ is small (Theorem~\ref{thm:fd-noncommutative-exactify}), and we lift separable abelian $2n$-independent subalgebras of a tracial ultraproduct over a countably cofinal ultrafilter to $n$-independent coordinate algebras (Theorem~\ref{thm:local-implies-hi}\textup{(i)}). Houdayer and Ioana~\cite[Theorem~5.1]{houdayer-ioana} proved this lifting at order one, and Patchell and Shiner~\cite[Theorem~C]{patchell-shiner} obtained it at every order independently of us. In $\mathrm{II}_1$ factors, we lift without decreasing the order.

Since $\varepsilon_0$ and $C$ do not depend on the factor, we can correct representatives in tracial
ultraproducts over countably cofinal ultrafilters $\mathcal U$, i.e., those for which some
sequence of sets in $\mathcal U$ has empty intersection. In particular, note that every
nonprincipal ultrafilter on $\N$ is countably cofinal.

\begin{introthm}[Theorem~\ref{thm:hyperfinite-lifting-cofinal}]
\label{introthm:lifting}
Let $n\ge1$, let $\mathcal U$ be a countably cofinal ultrafilter on a set $I$, and let $M_k$, $k\in I$, be $\mathrm{II}_1$ factors. If $A,B\subseteq\prod_{\mathcal U}M_k$ are separable
hyperfinite $n$-independent unital von Neumann subalgebras, then there are
finite-dimensional unital $^*$-subalgebras $C_k,D_k\subseteq M_k$,
$n$-independent for every $k\in I$, such that
\[
A\subseteq\prod_{\mathcal U}C_k,\qquad
B\subseteq\prod_{\mathcal U}D_k.
\]
\end{introthm}

The factors need not
be separable, and we do not assume that $W^*(A,B)$ is hyperfinite. If $A$ and $B$ are freely independent, then $C_k$ and $D_k$ can be chosen $r(k)$-independent with $r(k)\to\infty$ along $\mathcal U$ (Theorem~\ref{thm:hyperfinite-free-lifts-cofinal}). We obtain
the lifts by correcting finite-dimensional approximations to $A$ and $B$
and then making a diagonal choice along the ultrafilter. For countably complete ultrafilters, the conclusion holds for arbitrary tracial von Neumann algebras $M_k$, with hyperfinite $C_k,D_k$ (Theorem~\ref{thm:hyperfinite-lifting-complete}). Thus, for $\mathrm{II}_1$ factors, separable hyperfinite $n$-independent subalgebras lift over every ultrafilter.

\subsection*{Orthogonality and a question of Houdayer and Ioana}
Houdayer and
Ioana~\cite[Theorem~5.1 and the following paragraph]{houdayer-ioana}
obtained orthogonal coordinate von Neumann subalgebras from separable abelian
$2$-independent von Neumann subalgebras of a tracial ultraproduct and asked whether $1$-independence
suffices. Kunnawalkam Elayavalli's problem
list~\cite{kunnawalkam-elayavalli-open-problems} concerns ultraproducts of
$\mathrm{II}_1$ factors, and we formulate its lifting questions in this
setting. The list asks whether lifts of unitaries can preserve
orthogonality of their generated algebras at every coordinate.

\begin{problemthirty}[{Houdayer and Ioana~\cite[Theorem~5.1 and the following paragraph]{houdayer-ioana}; Kunnawalkam Elayavalli~\cite[Problem~30]{kunnawalkam-elayavalli-open-problems}, in our terminology}]
Let $N$ be a $\mathrm{II}_1$ factor and let $\mathcal U$ be a
nonprincipal ultrafilter on $\N$. Given unitaries $u,v\in N^{\mathcal U}$
such that $W^*(u)\perp W^*(v)$, can one choose unitary representatives
$u=(u_j)_{\mathcal U}$ and $v=(v_j)_{\mathcal U}$ satisfying
$W^*(u_j)\perp W^*(v_j)$ for every $j$?
\end{problemthirty}

This is the question of Houdayer and Ioana in the case of an ultrapower of a
$\mathrm{II}_1$ factor. Chifan, Ioana, and Kunnawalkam
Elayavalli~\cite[Lemma~3.1]{chifan-ioana-kunnawalkam-elayavalli} had already proved the case $u^2=v^3=1$.

\begin{introcor}[Corollary~\ref{cor:fixed-factor-problem-30}]
\label{introcor:problem-30}
Problem~30 has an affirmative answer.
\end{introcor}

Indeed, $W^*(u)$ and $W^*(v)$ are separable, abelian, and orthogonal. By Theorem~\ref{introthm:lifting} at order one, or by its abelian case, Theorem~\ref{thm:factor-no-loss-hierarchy}, we have $W^*(u)\subseteq\prod_{\mathcal U}C_j$ and $W^*(v)\subseteq\prod_{\mathcal U}D_j$ for finite-dimensional unital $^*$-subalgebras $C_j\perp D_j$ of $N$, and we take unitary representatives of $u$ and $v$ in these algebras. More generally, we thus obtain the separable hyperfinite case of the
lifting of orthogonality, which Patchell and Shiner~\cite{patchell-shiner} noted as open in general.

\subsection*{Non-Gamma factors}
Chifan, Ioana, and Kunnawalkam Elayavalli~\cite{chifan-ioana-kunnawalkam-elayavalli} and Houdayer and Ioana~\cite{houdayer-ioana} distinguished non-Gamma $\mathrm{II}_1$ factors up to elementary equivalence by finding, for suitable pairs of unitaries $u_1,u_2$ in the ultrapowers of a non-Gamma factor, commuting Haar unitaries $v_1,v_2$ with $[u_1,v_1]=[u_2,v_2]=0$. Using the order-one lifting, we find such $v_1,v_2$ for every pair generating orthogonal von Neumann algebras.

\begin{introthm}[Theorem~\ref{thm:houdayer-ioana-orthogonal-g}]
\label{introthm:non-gamma}
There exists a separable non-Gamma $\mathrm{II}_1$ factor $M$ with
the following property. If\/ $\mathcal U$ is a countably cofinal ultrafilter on a set $I$ and unitaries
$u_1,u_2\in M^{\mathcal U}$ satisfy $W^*(u_1)\perp W^*(u_2)$, then there are Haar unitaries
$v_1,v_2\in M^{\mathcal U}$ such that
\[
[u_1,v_1]=[u_2,v_2]=[v_1,v_2]=0.
\]
\end{introthm}

Chifan, Ioana, and Kunnawalkam Elayavalli~\cite[Theorem~B]{chifan-ioana-kunnawalkam-elayavalli} proved this
under the additional relations $u_1^2=u_2^3=1$, and Houdayer and Ioana~\cite[Theorem~G]{houdayer-ioana} for
$2$-independent $W^*(u_1)$ and $W^*(u_2)$. We use their construction unchanged, in the form of Houdayer and Ioana~\cite[Definition~5.7]{houdayer-ioana}: by~\cite[Corollary~4.3]{chifan-ioana-kunnawalkam-elayavalli}, each of its steps requires only orthogonality of the pair. The stronger assumptions were used only to lift $u_1,u_2$ to unitaries generating orthogonal von Neumann algebras in the factors of the construction, and by Theorem~\ref{thm:factor-no-loss-hierarchy} at order one, orthogonality alone suffices for this lifting.

\subsection*{The method}
To prove Theorem~\ref{introthm:stability}, we express finite-order independence of finite-dimensional unital $^*$-subalgebras as a finite system of real moment equations. When the first
subalgebra varies by unitary conjugation, the gradients of these equations
at a free pair have a positive-definite Gram matrix. The Gram entries
involve mixed moments through order $2n$. Hence, if $\varepsilon_{2n}(A,B)$
is sufficiently small, then, by comparison
with the free product, we obtain a
unitary $u$ such that $uAu^*$ and $B$ are $n$-independent and
$\norm{u-1}_\infty\le C\varepsilon_{2n}(A,B)$. This holds in an arbitrary tracial von
Neumann algebra, for copies of $P,Q$ as in Theorem~\ref{introthm:stability}, with the same lower bound $\delta$ on
the traces of minimal projections
(Theorem~\ref{thm:fd-noncommutative-exactify}).

To avoid this loss of order, one might try linearizing at a nearby
$n$-independent pair, but the moment derivative there need not
be surjective. For example, if $A$ and $B$ commute and are $1$-independent, then every gradient at
order one is a multiple of a commutator $[a,b]$ with $a\in A$ and $b\in B$, and hence vanishes.
However, in a $\mathrm{II}_1$ factor, an arbitrarily small preliminary conjugation can preserve
$n$-independence while making the derivative surjective. To construct
such conjugations, we work in a free-product extension containing a
unitary free from $W^*(A,B)$. By Popa's asymptotic freeness theorem~\cite[Theorem~2.1]{popa-free-independent-1995},
in the formulation of~\cite[Corollary~0.2]{popa-independence},
an ultrapower of the factor contains a unitary with the same tracial law as this one,
again free from $W^*(A,B)$. Correcting finitely many moments, we then obtain
the required identities inside the given factor
(Theorem~\ref{thm:finite-order-free-positioning}). For nearly $n$-independent pairs, we carry out this construction
quantitatively, bounding both the remaining moment errors and
the norm of a right inverse of the moment derivative. Finally, a small conjugation
removes these errors.

\subsection*{Arbitrary tracial coordinate algebras and a problem of Kunnawalkam Elayavalli}
For separable abelian unital von Neumann subalgebras
$A,B\subseteq\prod_{\mathcal U}M_k$ and a countably cofinal ultrafilter
$\mathcal U$, the argument without the preliminary conjugation works for arbitrary tracial von Neumann
algebras $M_k$, but it needs $2n$-independence of $A$ and $B$ to produce
$n$-independent coordinate algebras
(Theorem~\ref{thm:local-implies-hi}). We show that, in the question of Houdayer and Ioana, $1$-independence does not suffice for arbitrary tracial coordinate algebras, even for $\C^4$ with traces varying with the coordinate
(Proposition~\ref{prop:four-point-orthogonality-obstruction}), while for
$\mathrm{II}_1$-factor coordinates it does, by
Theorems~\ref{thm:factor-no-loss-hierarchy} and~\ref{thm:hyperfinite-lifting-complete}.

If the two separable abelian subalgebras are free, then
by a diagonal choice we obtain one pair of lifts whose independence orders tend to
infinity along a countably cofinal ultrafilter, even for arbitrary tracial
von Neumann coordinate algebras (Theorem~\ref{thm:free-coordinate-lifts}).

In Kunnawalkam Elayavalli's problem
list, $k$-independence means that every centered alternating
word of length $k$ has trace zero, where the length counts the letters, i.e., the blocks. To distinguish it from the
$n$-independence of this paper, which concerns all words $a_1b_1\cdots a_rb_r$ with $r\le n$, we call it exact
block-length $k$-independence. Note that orthogonality is exact block-length $2$-independence. The list also contains the following question.

\begin{problemthirtyone}[{Kunnawalkam Elayavalli~\cite[Problem~31]{kunnawalkam-elayavalli-open-problems}, in our terminology}]
Let $N$ be a $\mathrm{II}_1$ factor and let $\mathcal U$ be a
nonprincipal ultrafilter on $\N$. Do freely independent unitaries
$u,v\in N^{\mathcal U}$ have unitary representatives $u=(u_j)_{\mathcal U}$
and $v=(v_j)_{\mathcal U}$ for which $W^*(u_j)$ and $W^*(v_j)$ are exact
block-length $k$-independent for every $j$, for all $k$?
\end{problemthirtyone}

The quantifier over $k$ can be read in two ways. Requiring one pair of representatives to work for every $k$ amounts to
asking for free coordinate representatives. Under the other reading, which we call the fixed-length
interpretation, the representatives may depend on $k$. We show that the answer is then affirmative
for every tracial von Neumann algebra $N$: for each fixed $n$, the representatives can be chosen so that $W^*(u_j)$ and $W^*(v_j)$ are $n$-independent for every $j$, and hence exact block-length $k$-independent for every $k\le2n+1$
(Corollary~\ref{cor:free-unitary-lifting}).

However, stronger requirements can fail. In an ultrapower of the hyperfinite
$\mathrm{II}_1$ factor, a pair of free unitaries need not admit freely independent
coordinate representatives. This follows from Popa's asymptotic freeness
theorem~\cite[Theorem~2.1]{popa-free-independent-1995} and the nonamenability
of suitable free products. Boulanger, Curda, Harvey, Li, and Pi~\cite[Section~3.2]{boulanger-curda-harvey-li-pi}
gave this argument, for self-adjoint
elements, in proving that pairs of freely independent elements do not form a definable set.
We also show that one cannot in general require exact block-length $j$-independence
at coordinate $j$
(Proposition~\ref{prop:no-prescribed-diagonal}). Nor can convergence of the
independence orders along the ultrafilter always be replaced by ordinary
convergence: for some nonprincipal ultrafilter on $\N$, a pair of free
Haar unitaries in an ultrapower of the hyperfinite factor has no pair of unitary
representatives for which, for every fixed $n$, the two generated von
Neumann algebras are $n$-independent at all sufficiently large indices
(Proposition~\ref{prop:no-ordinary-eventual-lifts}).

Patchell and Shiner independently obtained the $2n$-to-$n$ lifting theorem
and the fixed-length solution of Problem~31~\cite[Theorem~C]{patchell-shiner}.
In the countably cofinal case, their theorem gives the same
$2n$-to-$n$ lifting conclusion as
Theorem~\ref{thm:local-implies-hi}\textup{(i)}.
For a nonprincipal ultrafilter on $\N$, their proof of
\cite[Theorem~A]{patchell-shiner} also constructs one pair of
unitary representatives of a given pair of free Haar unitaries in an ultrapower of a
${\rm II}_1$ factor, with independence orders tending to infinity
along the ultrafilter (by Proposition~\ref{prop:no-ordinary-eventual-lifts}, this cannot in general be
strengthened to ordinary convergence).
They also showed that, given finite sets $X,Y$ of centered elements in a
$\mathrm{II}_1$ factor and a prescribed finite order, there is a unitary
$u$ for which $uXu^*$ and $Y$ are independent to that
order, with alternating words of odd length included~\cite[Definition~2.1 and Theorem~B]{patchell-shiner}; the order-one case is due to Houdayer and Ioana~\cite[Theorem~H]{houdayer-ioana}. Independently of Patchell and Shiner, and by a slightly different method, we obtain such a unitary in Theorem~\ref{thm:finite-order-free-positioning}.

At order one, Houdayer and Ioana~\cite[Lemmas~5.2--5.4]{houdayer-ioana}
corrected the moments $\tau(x_iy_j)$ by small unitary conjugations,
using commutator gradients and inverting the nearly diagonal Gram
matrix of these gradients. Patchell and Shiner's paper and ours extend this method to
higher mixed moments, using cyclic gradients (expressed through
``generalized commutators'' in their paper) and Gram-matrix
inversion~\cite[Theorem~2.5]{patchell-shiner}. Their proof updates these commutators
at each step, while ours uses a contraction argument with the inverse of
the initial Gram matrix held fixed.

To obtain $n$-independence from a small defect of order $n$ in factors,
we also use the preliminary conjugation above
(Propositions~\ref{prop:quantitative-phase-marker} and~\ref{prop:dense-regular-exact-locus}).
This step goes beyond the results that the two papers share: in
$\mathrm{II}_1$ factors, it removes the loss of order and so gives Theorems~\ref{introthm:stability}
and~\ref{introthm:lifting}, and with them Corollary~\ref{introcor:problem-30}.

\subsection*{Partitions and the tracial version of the orthogonality question}
For partitions of $1$ in $\mathrm{II}_1$ factors, the lower bound on the traces of minimal projections
in Theorem~\ref{introthm:stability} can be removed if we allow both partitions to move. At order one, we
show that small $L^2$ perturbations of a nearly independent pair make it independent
while preserving the trace of every projection (Theorem~\ref{thm:finite-partition-stability}).
Chifan, Ioana, and Kunnawalkam
Elayavalli~\cite[Lemmas~3.1--3.2]{chifan-ioana-kunnawalkam-elayavalli}
proved the case of a projection and three mutually orthogonal projections
summing to $1$ in a $\mathrm{II}_1$ factor, without requiring all the
projection traces to remain fixed. At higher orders, we can also remove the lower bound if both
partitions and their traces are allowed to change slightly. For fixed order and
partition sizes, the perturbation bounds are independent of the ambient
$\mathrm{II}_1$ factor and of the traces of the individual projections
(Corollary~\ref{cor:factor-full-local-no-loss}). In both results, zero
projections are allowed, before and after
the correction.

Allowing $N$ in Problem~30 to be any tracial von Neumann algebra gives
what we call the tracial version of Problem~30. Note that the example on $\C^4$ of Proposition~\ref{prop:four-point-orthogonality-obstruction} does not answer it, since there the trace varies with the coordinate. We show that, for fixed $N$, this version has an
affirmative answer precisely when the equations
\[
\tau(p_iq_j)=\tau(p_i)\tau(q_j)
\]
are stable for partitions $\mathbf p=(p_i)$ and $\mathbf q=(q_j)$
of $1$ in $N$: for fixed partition sizes and any $\eta>0$,
sufficiently small errors in these equations can be corrected by
perturbing both partitions within $N$ with total $L^2$ movement less
than $\eta$. Here zero projections are allowed in both the original and the
corrected partitions, and the individual
projection traces need not be preserved
(Proposition~\ref{prop:fixed-ambient-stability-equivalence}).

\begin{introthm}[Corollary~\ref{cor:diffuse-scalar-problem-30}, Theorems~\ref{thm:finite-type-I-problem-30} and~\ref{thm:fixed-ambient-problem-30-counterexample}]
\label{introthm:tracial}
If a tracial von Neumann algebra $N$ is diffuse, of type~I, or the direct sum of a diffuse tracial von Neumann algebra and one scalar atom, then the tracial version of Problem~30 has an affirmative answer for $N$. For a suitable direct sum $N$ of a $\mathrm{II}_1$ factor and finitely many scalar atoms, it has a negative answer: for every nonprincipal ultrafilter $\mathcal U$ on $\N$, some unitaries $u,v\in N^{\mathcal U}$ with finite spectrum and $W^*(u)\perp W^*(v)$ have no unitary representatives $u=(u_j)_{\mathcal U}$ and $v=(v_j)_{\mathcal U}$ such that $W^*(u_j)\perp W^*(v_j)$ for every $j$.
\end{introthm}

When $N$ is not a factor, independence on each central summand need not
give independence in $N$, since an average of products of traces need
not be the product of the averages. Conversely, a pair can be independent in $N$ even though its restrictions
to the summands are not independent for their normalized traces.
We therefore study prescribed overlap matrices $(\tau(p_iq_j))_{i,j}$,
including targets other than the product matrices
$(\tau(p_i)\tau(q_j))_{i,j}$. On the $\mathrm{II}_1$-factor summands, we use finite-order independence
inside corners to split projections and to estimate the Gram matrices of
selected overlap gradients $i[p_i,q_j]$, which arise from unitary
conjugation of the first partition. We sum the Gram matrices, weighted by the central
traces, and apply the inverse estimate of
Lemma~\ref{lem:finite-coordinate-overlap-inverse} to correct the overlaps. Combining
this with separate arguments for the atomic part and the part with diffuse center, we
prove stability of the independence equations under $L^2$
perturbations in every fixed diffuse
tracial von Neumann algebra and in its direct sum with one scalar atom of
fixed weight $0<\lambda<1$
(Corollary~\ref{cor:diffuse-scalar-problem-30}).

By a separate compactness argument, we obtain stability
in every fixed type I tracial von Neumann algebra
(Theorem~\ref{thm:finite-type-I-problem-30}). When one partition is
binary, stability holds in every fixed tracial von Neumann algebra
(Corollary~\ref{cor:fixed-ambient-binary-problem-30}). All these stability
results allow both partitions and their traces to change.

However, there is no such stability theorem for all fixed tracial von Neumann
algebras. We construct a counterexample with finitely many scalar atoms
and one $\mathrm{II}_1$ factor, where even allowing both partitions and all their traces
to change does not give nearby independent partitions
(Theorem~\ref{thm:fixed-ambient-problem-30-counterexample}). The
construction starts from an overlap matrix $(\tau(p_iq_j))_{i,j}$
that is unstable in a fixed factor: there are pairs of partitions whose
overlap matrices converge to it, while their $L^2$ distance from every
pair with exactly this overlap matrix stays bounded away from zero. We
obtain this matrix from the nonlifting
theorem of Ioana and Spaas~\cite[Theorem~E]{ioana-spaas}. The scalar atoms force any nearby
independent pair to realize that matrix on the factor summand. Thus
stability of independence in factors need not survive adjoining finitely many scalar
summands, and the tracial version of Problem~30 has a negative answer.
More generally, if an overlap matrix is unstable in a fixed diffuse
tracial von Neumann algebra $D$, adjoining finitely many scalar atoms with suitable weights to
$D$ gives a fixed algebra for which the tracial version of Problem~30 fails
(Corollary~\ref{cor:universal-table-independence-equivalence}).

\subsection*{Matrix models}
For subalgebras of fixed finite-dimensional isomorphism classes with
traces of minimal projections bounded away from zero, we can also obtain $n$-independence from small defects of order $n$
in all sufficiently large matrix algebras. The proof passes
through a tracial ultraproduct of matrix algebras whose dimensions tend
to infinity
(Corollary~\ref{cor:large-matrix-fixed-type-no-loss}).

In matrix models, however, even a fixed number of matrices can generate algebras
of increasing dimension as the matrix size grows, so we seek exact
free-product identities for mixed moments only up to a polynomial degree,
which may grow with the dimension.
As in the work of Male~\cite{male-norm-polynomials} and of Collins and Male~\cite{collins-male}, strong
convergence of matrices means convergence of the normalized traces and of the operator norms of all
noncommutative $^*$-polynomials in the matrices. For families of independent normalized
matrices from the Gaussian unitary ensemble (GUE), Haagerup and
Thorbj{\o}rnsen~\cite[Theorem~A]{haagerup-thorbjornsen} proved the
almost-sure convergence of the norms.

Our second perturbation theorem is for strongly convergent tuples of matrices. The polynomial space of a tuple
through degree $d$ is the complex linear span of $1$ and all products
of at most $d$ factors chosen from the matrices of the tuple and their adjoints.

\begin{introthm}[Theorems~\ref{thm:gns-rank-strong-exactification} and~\ref{thm:multi-gns-rank-strong-exactification}]
\label{introthm:strong}
Let finitely many tuples of matrices in $M_{N_j}(\C)$, each of fixed length, converge jointly and strongly to free
limiting families in a tracial $C^*$-probability space.
Assume that, at every fixed degree, the
polynomial space of each tuple eventually has the same
dimension as the corresponding polynomial space in the limiting tracial
GNS representation. Then, for all sufficiently large $j$, there are unitaries in $U(N_j)$, one for each tuple, tending to $1$ in
operator norm, and integers $d_j\to\infty$ such that the conjugated tuples have the same joint strong limit and their joint moments of degree at most $d_j$ are the free-product moments determined by their marginal states at each matrix dimension. One of these unitaries can be taken to be $1$.
\end{introthm}

At each fixed degree, the dimension hypothesis ensures that polynomials
whose values form a basis in the limiting GNS space also have values
spanning the polynomial space in all sufficiently large matrices.
Note that strong convergence alone does not exclude additional
polynomial directions that disappear in the limit. The dimension hypothesis cannot be dropped: for some strongly convergent pairs of tuples with free limits and a faithful limiting trace, no unitaries tending to $1$ in operator norm make the joint moments agree with these free-product moments even through degree two
(Proposition~\ref{prop:strong-model-rank-obstructions}). When the limiting
trace is not faithful, polynomials of zero $L^2$-norm are identified with
zero in computing the limiting dimension. The correction still preserves
the supplied limiting $C^*$-norm, even when that norm is not determined by
the trace.

For the corresponding result with amalgamation over a fixed finite-dimensional unital
$C^*$-algebra $D$, we use all $^*$-polynomials in each tuple with coefficients in
$D$, including coefficients between generators. Suppose that the tuples are uniformly bounded in operator norm and,
together with $D$, converge in moments to families free over $D$, and that the
limiting trace is faithful on $D$. Suppose also that, at each fixed degree, these
polynomial spaces eventually have the corresponding limiting tracial
GNS dimensions. We conjugate the tuples by unitaries in $D'$ tending to $1$
in operator norm. These conjugations fix $D$ pointwise, preserve each marginal
$D$-valued law, and give exact amalgamated free-product identities for the
marginal laws at each dimension, through degrees tending to infinity. If the original
tuples converge strongly together with $D$, the corrected tuples have the same strong limit
(Theorem~\ref{thm:relative-gns-rank-exactification}).

\subsection*{GUE families}
Independent normalized GUE matrices, each forming one tuple, almost surely satisfy the hypotheses of Theorem~\ref{introthm:strong}: they converge strongly to free semicircular variables, and each of them has simple spectrum, which gives the dimension hypothesis (Corollary~\ref{cor:gue-growing-exactification}). For such matrices, we can ask more precisely how many identities can be
imposed for a given perturbation size.

\begin{introthm}[Theorem~\ref{thm:gue-relative-exactification} and Proposition~\ref{prop:relative-degree-ceiling}]
\label{introthm:gue}
For two independent normalized GUE matrices $G_N,H_N$ and every $0<\gamma<1$, there are $c>0$ and unitaries $u_N$, Borel functions of $(G_N,H_N)$, with $\norm{u_N-1}_\infty\le N^{-\gamma}$ such that, almost surely for all sufficiently large $N$, the joint moments of $u_NG_Nu_N^*$ and $H_N$ through degree $\lfloor c\log N\rfloor$ are those of freely independent variables with the same distributions as $G_N$ and $H_N$ with respect to $\operatorname{tr}_N$. On the other hand, almost surely for all sufficiently large $N$, no unitary $u\in U(N)$ gives these identities for $uG_Nu^*$ and $H_N$ through any degree larger than $4\log N/\log2+1$.
\end{introthm}

More generally, fix the numbers of GUE families and of generators in each, and a finite-dimensional unital $C^*$-algebra $D$. Let $N$ range over an unbounded set of dimensions, and fix unital $^*$-embeddings $D\to M_N(\C)$, with the normalized traces of minimal
projections in $D$ uniformly bounded away from zero. Assume that all generators are
mutually independent at each dimension. We show that, for every $0<\gamma<1$,
there are conjugating unitaries in $D'$, with one family fixed, whose
operator-norm distance from $1$ is always at most $N^{-\gamma}$.
Outside events of summable probability, they give exact amalgamated
identities through degree $\lfloor c\log N\rfloor$ for some $c>0$
(Theorem~\ref{thm:gue-relative-exactification}). If
the normalized traces restricted to $D$ converge, the correction preserves
the almost-sure strong limit of the tuple together with $D$.

When there are at least two nonempty families, logarithmic growth of the degree is optimal.
By counting, in a corner of $M_N(\C)$, the alternating words that are orthonormal in the free product, we show
that, almost surely and for all sufficiently large $N$, conjugation by
unitaries in $D'$ cannot give these identities through a degree larger
than $4\log N/\log s+1$, where $s\ge2$ is the number of families
(Proposition~\ref{prop:relative-degree-ceiling}).

The degree and the perturbation size cannot be chosen independently. Every
prescribed sublogarithmic degree can be attained with a deterministic
distance bound $N^{-1+o(1)}$ and summable failure probabilities
(Proposition~\ref{prop:gue-degree-motion-profile}). Conversely, if one of the
GUE families has at least two generators, or if $D\ne\C$, a deterministic bound of
$N^{-1+o(1)}$ can succeed with probability tending to one only at degrees
$o(\log N)$ (Theorem~\ref{thm:relative-degree-motion-lower-bound}).

For finitely many covariance equations, we determine the
limiting law of $N$ times the least distance of a correcting unitary from
$1$. Let $P_1,\ldots,P_k$ be real polynomials, linearly independent modulo
constants, let $G_N,H_N$ be independent GUE matrices, and let $S,T$ be
free standard semicircular variables. Then the minimum of
$N\norm{u-1}_\infty$ over unitaries making every covariance between
$uP_j(G_N)u^*$ and $H_N$ vanish, with $H_N$ fixed, converges in distribution to
$p(Z)$ (Corollary~\ref{cor:optcost-polynomial-cost}). Here $Z$ is a centered
Gaussian vector whose covariance is the Gram matrix of the centered
$P_j(S)$, and $p$ is the norm on $\R^k$ whose dual norm is
\[
z\longmapsto\tau\abs{i\Bigl[\sum_{j=1}^kz_jP_j(S),T\Bigr]}.
\]
Moreover, the minimum of $\norm{u-1}_\infty$ over unitaries correcting
two such equations simultaneously can exceed the larger of the two
minima obtained by correcting each equation separately
(Proposition~\ref{prop:optcost-genuine-joint-cost}).

\subsection*{Clock representations and traces of words}
By the same perturbation method, for each fixed $s\ge2$, we construct irreducible unitary representations of $\mathbb F_s$ in $U(N)$ that converge strongly to the left regular representation, as in the work of Collins and Male~\cite{collins-male}, but which moreover have trace zero on every nontrivial word of length at most a constant times $\log N$.

\begin{introthm}[Theorem~\ref{thm:irreducible-clock-strong-models} and Proposition~\ref{prop:exact-character-radius-ceiling}]
\label{introthm:clocks}
Fix $s\ge2$, and write $|g|$ for the reduced word length of $g\in\mathbb F_s$ in the standard generators. For every sufficiently large integer $N$, there is an irreducible unitary representation
$\pi_N:\mathbb F_s\to U(N)$ that sends each standard generator to a unitary conjugate of the cyclic shift and satisfies
\[
\operatorname{tr}_N(\pi_N(g))=0
\qquad
\left(1\le |g|\le
\left\lfloor\frac{\log N}{16\log(2s)}\right\rfloor\right).
\]
These representations converge strongly
to the left regular representation as $N\to\infty$. On the other hand, if a representation $\pi:\mathbb F_s\to U(N)$ satisfies $\operatorname{tr}_N(\pi(g))=0$ whenever $1\le|g|\le R$, then $R\le4\log N/\log(2s-1)+O_s(1)$.
\end{introthm}

As with GUE families, logarithmic growth is thus optimal.

To construct $\pi_N$, we start with
independent Haar unitary matrices and replace their eigenvalues, in
circular order, by the $N$th roots of unity, each occurring once. The
resulting unitaries, which we call clocks, are conjugate to the clock matrix
$\operatorname{diag}(1,e^{2\pi i/N},\ldots,e^{2\pi i(N-1)/N})$,
and hence to the cyclic shift. Collins and Male used this rounding in their proof of strong
convergence and showed that its operator-norm error tends to zero almost
surely~\cite[Lemma~3.3 and Section~3.4]{collins-male}. A unitary
$U\in U(N)$ satisfies $\operatorname{tr}_N(U^m)=0$ for $1\le m<N$ precisely when it is a
scalar multiple of a clock, so clocks agree with the regular character on the powers of a
generator through the largest possible degree.
We then make the mixed word traces vanish by small unitary conjugations,
and the conjugated unitaries are still clocks. Since both perturbations tend to zero in operator
norm, the corrected unitaries keep
the strong limit of Collins and Male~\cite[Theorem~1.4]{collins-male}. They still satisfy
Hastings's asymptotically optimal estimate for the norm of the uniform average of
conjugations by the generators and their inverses, as an operator on the traceless
matrices with the $L^2$-norm~\cite[p.~1]{hastings-quantum-expanders}; see also Hastings and Harrow~\cite[Theorem~6]{hastings-harrow}. Bordenave and Collins~\cite[Theorem~2]{bordenave-collins-compact-groups} also obtained the limit of these
norms, by a different method.

For the free product of two nontrivial finite groups, small conjugations likewise give unitary representations that converge strongly to the left regular representation and whose normalized traces vanish on all nontrivial elements of syllable length at most $2r_j+1$, where $r_j\to\infty$ (Corollary~\ref{cor:exact-character-strong-models}). Using finite quotients and the theorem of Bordenave and Collins~\cite[Theorem~9.1]{bordenave-collins-norm-polynomials} on the strong convergence of random permutation matrices tensored with deterministic unitaries, we also obtain finite-image representations of free groups, in every sufficiently large dimension $N$, that converge strongly and whose normalized traces vanish on all nontrivial words of length at most a constant times $\log N$ (Proposition~\ref{prop:finite-image-log-character-models}); these are reducible.

For the clock construction, we require moment estimates for traces of words in
independent Haar unitaries, with explicit dependence on word length.
From these estimates, we also obtain a quantitative Gaussian approximation
for word traces of the original independent Haar matrices
(Theorem~\ref{thm:haar-word-gaussian}). For fixed words, Gaussian limits were obtained by
Diaconis and Shahshahani~\cite{diaconis-shahshahani} (powers of one matrix),
R\u{a}dulescu~\cite{radulescu-haar-traces} (words in two matrices), and Mingo, \'Sniady, and
Speicher~\cite{mingo-sniady-speicher} (any finite number of matrices). For the traces of the first $m$ powers of a single Haar matrix, Johansson and Lambert~\cite{johansson-lambert} proved Gaussian approximation in total variation with $m$ growing with the dimension.

Fix $s\ge2$ and write
$\mathbf U_N$ for the original Haar tuple. Among the nontrivial cyclically reduced words of length at most $L_N$,
choose one representative from each class under
cyclic rotation and inversion. Let $k_N$ be the number of representatives
and $r(w)$ the largest integer $r\ge1$ such that $w=v^r$ for some
$v\in\mathbb F_s$. We show that if $L_N\to\infty$ and $L_N=o(N^{2/7})$, then the values
\[
\frac{|\operatorname{Tr}(w(\mathbf U_N))|^2}{r(w)}-\log k_N,
\]
computed with the unnormalized trace, form a point process converging to a Poisson
point process with intensity $e^{-x}\,dx$. Their maximum has limiting
distribution function $x\mapsto\exp(-e^{-x})$
(Corollary~\ref{cor:haar-word-full-ball-extremes}). In this regime,
with probability tending to one, every class attaining the largest
absolute value of the unnormalized trace is represented by a generator power
(Theorem~\ref{thm:haar-word-raw-extremes}). Since $L_N<N$ for large $N$, rounding each
Haar unitary to a clock makes the traces
of these generator powers zero.
With probability tending to one, the original Haar tuple can also be
perturbed by $o(1)$ in operator norm to clocks with trace zero at every
nontrivial group word in the range of Theorem~\ref{introthm:clocks}
(Proposition~\ref{prop:quantitative-haar-clock-exactification}).
Thus these fluctuation laws need not be preserved by perturbations
tending to zero in operator norm.

\subsection*{Mixed \texorpdfstring{$q$}{q}-Gaussian families}
The preceding strong models have free limits, possibly with
amalgamation. Finally, we construct strong Gaussian matrix models
for the $q$-Gaussian variables introduced by Bo\.{z}ejko and
Speicher~\cite{bozejko-speicher-brownian}, whose commutation relations were first considered by Frisch and
Bourret~\cite{frisch-bourret}, and for their mixed generalization
introduced by Speicher~\cite{speicher-generalized-statistics}.
In the mixed version, a matrix $Q=(q_{ab})$ supplies a deformation
parameter for each pair of variables, in place of the single parameter
$q$. These mixed $q$-Gaussian variables form a $q$-Gaussian family when
$q_{ab}=q$ for all $a,b$, and a free semicircular family when $Q=0$.
In general, the limiting variables need not be free.

\begin{introthm}[Theorem~\ref{thm:mixed-q-strong-gaussian-models}, Corollaries~\ref{cor:mixed-q-spectra-mf} and~\ref{cor:scalar-q-mf-ext}]
\label{introthm:mixed-q}
Let $Q=(q_{ab})_{a,b=1}^r$ be a real symmetric matrix with
$\max_{a,b}|q_{ab}|<1$. There are integers $D_m\uparrow\infty$ and, for each $m$, independent centered Gaussian self-adjoint matrices $X_1^{(m)},\ldots,X_r^{(m)}\in M_{D_m}(\C)$, all defined on one probability space, that converge strongly, almost surely,
to the mixed $q$-Gaussian family with parameter matrix $Q$. Hence the $C^*$-algebra
generated by the mixed
$q$-Gaussian family is
MF. In particular, the $C^*$-algebra generated by a finite $q$-Gaussian family is MF for every $-1<q<1$, and its
Brown--Douglas--Fillmore extension semigroup is not a group when the family
has at least two generators.
\end{introthm}

The matrices are mutually independent, but the entries of each matrix may be correlated. Here MF (matricial field) is in the sense of Blackadar and Kirchberg~\cite{blackadar-kirchberg}, and the extension semigroup is that of Brown, Douglas, and Fillmore~\cite{brown-douglas-fillmore-1973,brown-douglas-fillmore}. Haagerup and
Thorbj{\o}rnsen~\cite[Corollary~1]{haagerup-thorbjornsen} showed that the extension semigroup of the reduced $C^*$-algebra of a free group
of rank at least two is not a group.
For the assertion on the extension semigroup, we use
the uniqueness of the trace on these $q$-Gaussian $C^*$-algebras, proved
by Amrutam, Jekel, and Wasilewski~\cite{amrutam-jekel-wasilewski}. For $|q|<\sqrt2-1$, Caspers and
Chen~\cite{caspers-chen} constructed strongly convergent matrix models for $q$-Gaussian
families and obtained these conclusions. We reduce to families with zero diagonal parameters, for which
the Fock representation of the universal $C^*$-algebra of the commutation relations is faithful by a theorem of Kuzmin and
Pochekai~\cite[Theorem~6.2]{kuzmin-pochekai}, and apply the graph-product argument of Caspers and Chen to those families;
this reduction removes the restriction $|q|<\sqrt2-1$. Our convergence is for matrix coefficients of each fixed size; for each fixed degree, the results of Caspers and Chen also allow matrix coefficients whose dimension exceeds that of the random matrices. The nonquasidiagonality argument
is due to Caspers and Chen.

\subsection*{Organization}
The paper has three parts: stability and lifting (Part~I,
Sections~\ref{sec:local-correction}--\ref{sec:unstable-tables}), exact moment identities for matrix
models (Part~II, Sections~\ref{sec:scalar-models}--\ref{sec:gue-bounds}), and free groups and
$q$-Gaussian families (Part~III,
Sections~\ref{sec:clock-representations}--\ref{sec:mixed-q-strong-models}). The perturbation method of
Part~I is used again for matrices in Part~II and for free groups in Part~III.
For the fixed-length solution of Problem~31, it suffices to read
Section~\ref{sec:local-correction},
Subsection~\ref{subsec:lifting-finite-order}, and Corollary~\ref{cor:free-unitary-lifting}.
Readers interested mainly in matrix models may go directly from
Section~\ref{sec:higher-order-no-loss} to Section~\ref{sec:scalar-models}; the intervening sections
are not prerequisites. The construction of mixed $q$-Gaussian models in
Section~\ref{sec:mixed-q-strong-models} is independent of the moment
correction and uses only the GUE moment bound
\eqref{eq:relgauss-genus-majorant} from the earlier matrix sections.
When the parameters vanish between groups of variables, we apply
Theorem~\ref{thm:multi-gns-rank-strong-exactification} to these models in Corollary~\ref{cor:mixed-q-exact-free-blocks}.
The proofs in Section~\ref{subsec:haar-word-fluctuations} use the earlier sections only through
Lemmas~\ref{lem:group-word-gradient-gram} and~\ref{lem:haar-growing-word-moments} and their proofs.
Subsections~\ref{subsec:further-overlap-special-cases} and~\ref{subsec:optimal-correction-costs} are not used
in later sections.
In Appendix~\ref{app:second-order-geometry}, we give a second proof of
Corollary~\ref{thm:faithful-relative-partition-stability} that avoids Popa's asymptotic freeness theorem, yielding
another proof of Theorem~\ref{thm:finite-partition-stability} by the induction in its proof.

\subsection*{Conventions}
For von Neumann algebras, separability means separability of the predual.
All inclusions of von Neumann algebras are unital unless stated otherwise.
A \emph{partition of $1$} is a finite indexed family of mutually
orthogonal nonzero projections whose sum is $1$, unless zero
projections are explicitly allowed. More generally, a partition of a
projection $e$ is such a family with sum $e$.
A von Neumann algebra is \emph{diffuse} if it has no nonzero minimal projections.
We use $C^*(S)$ and, in an ambient von Neumann algebra, $W^*(S)$ for the
$C^*$-algebra and von Neumann algebra generated by $S$ together with the
ambient identity. In a corner $pMp$, that identity is $p$.
On $M_N(\C)$, $\operatorname{Tr}$ is the unnormalized trace and
$\operatorname{tr}_N:=N^{-1}\operatorname{Tr}$. We use $U(N)$ for the
unitary group of $M_N(\C)$, and $U(A)$ for the unitary group of a unital $C^*$-algebra $A$.
A normalized GUE matrix $G=(G_{ij})_{i,j=1}^N$ is a centered Gaussian self-adjoint matrix with
$\mathbb E G_{ij}G_{kl}=N^{-1}\delta_{il}\delta_{jk}$.
For a vector $y\in\R^d$, $|y|$ denotes its Euclidean norm, unless another norm is specified. We write $Z(N)$ for the center of a von Neumann algebra $N$. In $i[x,y]$, $e^{ih}$, and the factors $\frac i2$ and $\frac1{2i}$, the letter $i$ is the imaginary unit, also in formulas where $i$ is an index.

\subsection*{AI collaboration}
This work was carried out in collaboration with artificial intelligence agents,
including various models in OpenAI's Codex and GPT~Pro, and
Anthropic's Claude Code.

\subsection*{Acknowledgments}
This investigation began in January~2024, when the author posed a question
similar to the fixed-length interpretation of Problem~31 during Adrian
Ioana's talk at the UCLA Functional Analysis Seminar. Problem~31 was posed
independently by Kunnawalkam
Elayavalli~\cite[Problem~31]{kunnawalkam-elayavalli-open-problems}.
In September~2026, while this manuscript was in
preparation, the author became aware of an independent solution of the
fixed-length interpretation of Problem~31 by Greg Patchell and Austin
Shiner~\cite[Theorem~C]{patchell-shiner}.

The author thanks Sorin Popa for comments on a draft of this paper, Nicholas Boschert and Patrick Hiatt for early discussions
of these questions, and Patchell and Shiner for their cordial correspondence.

\part{Stability and lifting}
We first correct small defects of finite-order independence by small conjugations: in
Section~\ref{sec:local-correction} by a local inverse argument, and in Section~\ref{sec:higher-order-no-loss}
without loss of order in $\mathrm{II}_1$ factors. We then lift independence from tracial ultraproducts to
their coordinates (Section~\ref{sec:coordinate-lifting}) and ask whether nearly independent pairs of finite
partitions in a fixed tracial von Neumann algebra are close to independent pairs: they are in every diffuse
algebra and after adjoining one scalar atom, but not in a suitable direct sum of a $\mathrm{II}_1$ factor and
finitely many scalar atoms (Sections~\ref{sec:projection-tables}--\ref{sec:unstable-tables}).

\section{Comparison and local correction}
\label{sec:local-correction}

For finite-dimensional unital $^*$-subalgebras of a tracial von Neumann
algebra, we express finite-order independence as a finite system of real mixed-moment equations. We compute the
gradients of these equations in the free product and show that
their Gram matrix is positive definite. We then use these gradients to
find a small conjugation that makes the moments vanish exactly. To
control the Gram matrix, we need estimates for longer words than those
appearing in the equations. This is why, in Theorem~\ref{thm:fd-noncommutative-exactify}, we assume
a sufficiently small defect of order $2n$ to obtain independence of order $n$. In $\mathrm{II}_1$
factors, a preliminary
conjugation allows us to keep the same order, as we show in the next
section. We will also use these estimates for matrix models.

\subsection{Definitions and conventions}

Throughout, $I$ is a set and $\cU$ is an arbitrary ultrafilter on $I$. We call $\cU$
\emph{countably cofinal} (equivalently, countably incomplete) if there is a sequence $(X_m)_{m\ge1}$ in $\cU$ such that
$\bigcap_{m\ge1}X_m=\varnothing$. We call $\cU$ \emph{countably complete} if
$\bigcap_{m\ge1}X_m\in\cU$ for every sequence $(X_m)_{m\ge1}$ in $\cU$. Every ultrafilter satisfies exactly one of these two conditions.
On a countable index set, every countably complete ultrafilter is
principal.
If $(M_k,\tau_k)_{k\in I}$ are tracial von Neumann algebras, we write
$\ultraprod M_k:=\prod_{\cU}(M_k,\tau_k)$ for their tracial ultraproduct.
Two uniformly bounded families indexed by $I$ represent the same element when their difference tends to zero along $\cU$ in the $L^2$-norms induced by the traces.

For self-adjoint elements $g_1,\ldots,g_d$ in a tracial von Neumann algebra $(M,\tau)$,
their Gram matrix is the matrix of $L^2$ inner products
$G=(\tau(g_i g_j))_{i,j=1}^d$.

If $(M,\tau)$ is a tracial von Neumann algebra and $A\subseteq M$ is a unital von Neumann subalgebra, we write
$A\ominus \C 1:=\{a\in A:\tau(a)=0\}$.
For subsets $X,Y\subseteq M\ominus\C1$ and an integer $m\ge1$, we say that $X$ and $Y$ are
\emph{$m$-independent} if
\[
\tau(x_1y_1\cdots x_r y_r)=0
\]
for every $1\le r\le m$, $x_1,\dots,x_r\in X$, and $y_1,\dots,y_r\in Y$.
For unital von Neumann subalgebras $A,B\subseteq M$, we say, following Popa~\cite{popa-kadison-singer,popa-independence}, that $A$ and $B$ are $m$-independent when
$A\ominus\C1$ and $B\ominus\C1$ are $m$-independent. By multilinearity, replacing
$X,Y$ by their complex linear spans leaves the definition unchanged.
For $^*$-closed subspaces $E,F\subseteq M\ominus\C1$, however,
$m$-independence of the subspaces does not imply
$m$-independence of $W^*(E)$ and $W^*(F)$.
For unital von Neumann subalgebras $A,B\subseteq M$, define the defect by
\[
\varepsilon_m(A,B)
:=\max_{1\le r\le m}\ \sup
\left\{
\abs{\tau(a_1b_1\cdots a_rb_r)}:
\begin{array}{l}
 a_i\in A\ominus\C1,\ b_i\in B\ominus\C1,\\
 \norm{a_i}_\infty\le1,\ \norm{b_i}_\infty\le1
\end{array}
\right\}.
\]
Thus $A$ and $B$ are $m$-independent exactly when $\varepsilon_m(A,B)=0$.

\medskip
\noindent\emph{Finite-order convention.}
We count pairs of letters: $x_1y_1\cdots x_ry_r$ has pair length $r$
and block length $2r$. This is the convention of Houdayer and Ioana~\cite[Section~1, Terminology]{houdayer-ioana}.
Popa's convention for centered subsets~\cite[Section~3]{popa-kadison-singer}
also allows one or both outer letters $x_1,y_r$ of $x_1y_1\cdots x_ry_r$ to be
replaced by $1$, provided the resulting word is nonempty. For general
subsets, the two conventions differ: in $M_3(\C)$ with normalized trace,
the centered self-adjoint contractions
\[
x=\begin{pmatrix}-1&0&0\\0&0&0\\0&0&1\end{pmatrix},
\qquad
y=\frac12\begin{pmatrix}1&0&1\\0&-2&0\\1&0&1\end{pmatrix}
\]
satisfy $\tau(xy)=\tau(xyxy)=0$ but $\tau(xyx)=1/3$. However, for the
centered parts of unital von Neumann subalgebras, the two conventions agree
at each
order: order-$m$ independence in our sense implies that every centered
alternating word of block length at most $2m+1$ has trace zero
(Lemma~\ref{prop:truncated-free-product}).
We also use \emph{exact block-length $\ell$-independence}:
all centered alternating words having exactly $\ell$ blocks have trace zero, with no condition imposed at smaller lengths.
This is our reading of the convention in Problems~30--31 of~\cite{kunnawalkam-elayavalli-open-problems}.

A von Neumann algebra with separable predual is \emph{hyperfinite} if it has an increasing sequence of
finite-dimensional unital $^*$-subalgebras with weakly dense union. Every such algebra is
countably generated. Abelian von Neumann algebras with separable predual are hyperfinite.

For a tracial von Neumann algebra $(N,\tau)$ and integers $a,b\ge1$, let
$\mathcal P_{a,b}(N)$ denote the pairs of partitions of $1$ in $N$,
written $\mathbf p=(p_1,\ldots,p_a)$ and $\mathbf q=(q_1,\ldots,q_b)$,
with zero projections allowed in both partitions.
We retain the indices when
comparing two pairs of partitions. Write
\[
\alpha_i:=\tau(p_i),\qquad \beta_j:=\tau(q_j),\qquad
\Delta_{ij}(\mathbf p,\mathbf q):=\tau(p_iq_j)-\alpha_i\beta_j.
\]
We call $(\tau(p_iq_j))_{i,j}$ the \emph{overlap matrix} of the pair, and $\alpha=(\alpha_i)_i$ and
$\beta=(\beta_j)_j$ its \emph{trace vectors}. The trace vectors are the row and column sums of the overlap matrix.
Write $\Delta=(\Delta_{ij})$ for the covariance matrix, let $Z_{a,b}(N)$ be the set of pairs in $\mathcal P_{a,b}(N)$
for which $\Delta=0$, and let
\[
d_2((\mathbf p,\mathbf q),(\mathbf p',\mathbf q'))
:=\sum_i\norm{p_i-p_i'}_2+\sum_j\norm{q_j-q_j'}_2.
\]
The covariance matrix $\Delta$ has zero row and column sums:
\begin{equation}
\sum_j\Delta_{ij}=0,\qquad \sum_i\Delta_{ij}=0.
\label{eq:covariance-zero-row-column-sums}
\end{equation}
Thus its entries with $i<a$ and $j<b$ determine all the others.
We say that pairs in $Z_{a,b}(N)$ are \emph{independent}. Equivalently, the
unital $^*$-algebras they generate are $1$-independent. We call unitaries
$u,v$ independent when $W^*(u)$ and $W^*(v)$ are $1$-independent,
i.e., orthogonal. As in the Introduction, we write $A\perp B$ when $A$ and $B$ are orthogonal.

For partitions
\[
\mathbf p=(p_1,\ldots,p_a),
\qquad
\mathbf q=(q_1,\ldots,q_b)
\]
of $1$ in a tracial von Neumann algebra, allowing zero projections, set
\[
P_i:=p_i-\tau(p_i),
\qquad
Q_j:=q_j-\tau(q_j)
\]
and define their defect
\[
\operatorname{def}_n(\mathbf p,\mathbf q)
:=
\max_{\substack{1\le r\le n\\
i_1,\ldots,i_r\\j_1,\ldots,j_r}}
\abs{\tau(P_{i_1}Q_{j_1}\cdots P_{i_r}Q_{j_r})}.
\]
If
\[
A_{\mathbf p}:=\operatorname{span}\{p_1,\ldots,p_a\},
\qquad
B_{\mathbf q}:=\operatorname{span}\{q_1,\ldots,q_b\},
\]
then
\begin{equation}
\operatorname{def}_n(\mathbf p,\mathbf q)
\le \varepsilon_n(A_{\mathbf p},B_{\mathbf q})
\le (ab)^n\operatorname{def}_n(\mathbf p,\mathbf q).
\label{eq:labelled-defect-comparison}
\end{equation}
Indeed, the left inequality follows by taking the centered contractions $P_i,Q_j$.  Conversely, a centered contraction
$x\in A_{\mathbf p}$ has the form $x=\sum_i c_iP_i$, where $\abs{c_i}\le1$ whenever $p_i\ne0$. Terms with $p_i=0$ contribute nothing.
Expanding a centered contraction in the second algebra in the same way, we obtain at most $(ab)^r$ terms in a mixed moment of pair length
$r$. Thus, for fixed $a,b,n$, $\operatorname{def}_n(\mathbf p,\mathbf q)$ tends to zero if and only if
$\varepsilon_n(A_{\mathbf p},B_{\mathbf q})$ does.

\subsection{Free-product comparison}

We first compare the mixed moments of two subalgebras with their values in the free product,
with an error controlled by the defect.
An alternating word $W=c_1\cdots c_\ell$ is a sequence of blocks $c_t$
taken alternately from two subalgebras. A block that is a scalar still
counts as a block of the subalgebra it was taken from.
\begin{lemma}
\label{prop:truncated-free-product}
\label{lem:quantitative-truncated-comparison}
Let $(M,\tau)$ be a tracial von Neumann algebra, let $A_0,B_0\subseteq M$ be unital von Neumann
subalgebras, and let $m\ge1$. Set $\varepsilon=\varepsilon_m(A_0,B_0)$ and
let $(N,\varphi)=(A_0,\tau|_{A_0})*(B_0,\tau|_{B_0})$ be their free product von Neumann algebra. Define
\[
K_0=K_1=0,\qquad K_\ell=3^\ell K_{\ell-1}+2^{\ell+1}\quad(\ell\ge2).
\]
For every alternating word $W=c_1\cdots c_\ell$ in $A_0$ and $B_0$ with
$\ell\le2m+1$,
\[
|\tau(W)-\varphi(W)|
\le K_\ell\varepsilon\prod_{t=1}^\ell\norm{c_t}_\infty.
\]
For $m$-independent algebras, the two traces thus agree through block length $2m+1$.
If the blocks are contractions and at most $m$ come from each algebra,
the bound is $K_{2m}\varepsilon$.
\end{lemma}

\begin{proof}
We induct on $\ell$. The traces agree at lengths zero and one. Now we expand
$c_t=\tau(c_t)+c_t^\circ$, where
$\norm{c_t^\circ}_\infty\le2\norm{c_t}_\infty$. In each term containing a scalar
block, we remove the scalar blocks and multiply consecutive blocks from the
same algebra, preserving their order. The resulting alternating word is
shorter. By the induction hypothesis and
$\sum_{S\subsetneq\{1,\ldots,\ell\}}2^{|S|}=3^\ell-2^\ell$, the sum of these
trace differences is at most
\[
(3^\ell-2^\ell)K_{\ell-1}\varepsilon
\prod_t\norm{c_t}_\infty.
\]
For the fully centered term at even length, cyclicity and
$\ell/2\le m$ give the bound
$2^\ell\varepsilon\prod_t\norm{c_t}_\infty$.
At odd length $\ell\ge3$, we cyclically multiply the endpoint
blocks, which belong to the same algebra, and split the product into its scalar and centered parts. The
centered term has pair length $(\ell-1)/2\le m$, and centering increases
its norm bound by at most a factor of two. The scalar term has length $\ell-2$. Together they contribute at most
\[
2^\ell(2+K_{\ell-2})\varepsilon
\prod_t\norm{c_t}_\infty.
\]
The nondecreasing sequence $K_\ell$ satisfies
\[
(3^\ell-2^\ell)K_{\ell-1}
+2^\ell(2+K_{\ell-1})=K_\ell,
\]
so the induction holds for both even and odd lengths. The last assertion follows from
$\ell\le2m$.
\end{proof}

\begin{corollary}
\label{cor:filtered-truncated-comparison}
Fix integers $D\ge0$ and $m\ge1$.
Let $A_0,B_0$ be unital von Neumann subalgebras of a tracial von Neumann algebra $(M,\tau)$.
Let $\mathcal V_A(d)\subseteq A_0$ and $\mathcal V_B(d)\subseteq B_0$ be
the filtrations by total $^*$-degree of all noncommutative $^*$-polynomials in two fixed finite generating tuples, with
$\mathcal V_A(0)=\mathcal V_B(0)=\C1$ and
$\mathcal V_A(i)\mathcal V_A(j)\subseteq\mathcal V_A(i+j)$ and the analogous
property for the second filtration. If
$\mathcal V_A(D)\ominus\C1$ and $\mathcal V_B(D)\ominus\C1$ are
$m$-independent, then every joint noncommutative $^*$-polynomial of total $^*$-degree at most
$\min\{D,2m+1\}$ has the same trace as in the free product with the
tracial states induced on the two algebras.
\end{corollary}

\begin{proof}
We apply the preceding proof to each monomial. Note that removing scalar blocks,
multiplying blocks from the same algebra, and centering do not increase
total degree. Thus every centered block used in the induction belongs
to $\mathcal V_A(D)\ominus\C1$ or
$\mathcal V_B(D)\ominus\C1$. By hypothesis, the fully centered even alternating
moments of length at most $2m$ vanish, and the odd moments reduce to these and shorter
words exactly as before.
\end{proof}

For the local correction in Subsection~\ref{sec:pair-cyclic-coordinates}, we
use the following bases of the centered parts.

\begin{proposition}
\label{prop:fd-selfadjoint-basis}
Let $(M,\tau)$ be a tracial von Neumann algebra and let $A\subseteq M$ be a finite-dimensional unital von Neumann subalgebra $^*$-isomorphic to
\[
P=\bigoplus_{r=1}^{R}M_{d_r}(\C).
\]
Assume that every minimal projection of $A$ has trace at least $\delta>0$, and set $D:=\max_r d_r$. Then
$(A\ominus\C1)_{\rm sa}$ has a real $L^2$-orthogonal basis $\mathcal X$ such that every $x\in\mathcal X$ is a centered
self-adjoint contraction,
\[
\operatorname{span}_{\C}\mathcal X=A\ominus\C1,
\]
and
\[
\norm{x}_2^2\ge
\kappa(P,\delta):=\min\left\{1,2\delta^3,\frac{\delta}{2},\frac{2\delta}{D^2}\right\}.
\]
If $P=\C$, the basis is empty.
\end{proposition}

\begin{proof}
Choose matrix units $e_{ij}^{(r)}$ of $A$, let $z_r:=\sum_ie_{ii}^{(r)}$ be the corresponding
minimal central projections, and set
$\omega_r=\tau(z_r)$, $\mu_r=\omega_r/d_r\ge\delta$, and
$t_r=\sum_{s>r}\omega_s$. Take the central vectors
\[
c_r=t_rz_r-\omega_r\sum_{s>r}z_s\quad(1\le r<R),
\]
and, in the $r$th block, the vectors
\[
a_\ell^{(r)}=\frac{d_r-\ell}{d_r}e_{\ell\ell}^{(r)}
 -\frac1{d_r}\sum_{j>\ell}e_{jj}^{(r)}\quad(1\le\ell<d_r),
\qquad
x_{ij}^{(r)}=\frac{e_{ij}^{(r)}+e_{ji}^{(r)}}2,
\quad
y_{ij}^{(r)}=\frac{e_{ij}^{(r)}-e_{ji}^{(r)}}{2i}\quad(i<j).
\]
By the matrix-unit relations, these are mutually $L^2$-orthogonal
centered self-adjoint contractions, with
\[
\begin{aligned}
\norm{c_r}_2^2&=\omega_rt_r(\omega_r+t_r)\ge2\delta^3,\\
\norm{a_\ell^{(r)}}_2^2
 &=\mu_r\frac{(d_r-\ell)(d_r-\ell+1)}{d_r^2}
 \ge\frac{2\delta}{D^2},\\
\norm{x_{ij}^{(r)}}_2^2=\norm{y_{ij}^{(r)}}_2^2
 &=\frac{\mu_r}{2}\ge\frac\delta2.
\end{aligned}
\]
There are $(R-1)+\sum_r(d_r^2-1)=\dim_\C A-1$ such vectors, so they
form the required real basis and span $A\ominus\C1$ over $\C$.
\end{proof}

\subsection{Local correction}
\label{sec:pair-cyclic-coordinates}

To prove Theorem~\ref{thm:fd-noncommutative-exactify}, we compute the Gram matrix of the cyclic
gradients in the free product (Proposition~\ref{prop:noncommutative-free-gram}), linearize the moment equations under
conjugation of the first algebra (Proposition~\ref{prop:pair-linearization}), and then solve them by a contraction argument
(Lemma~\ref{lem:finite-dimensional-correction}). We begin by fixing the notation in the abelian case.

Fix $n\ge1$, a tracial von Neumann algebra $(M,\tau)$, and finite-dimensional abelian unital $^*$-subalgebras
$A_0=\bigoplus_{r=1}^a \C p_r$ and $B_0=\bigoplus_{s=1}^b \C q_s$ of $M$.

For these abelian algebras, the basis of Proposition~\ref{prop:fd-selfadjoint-basis}
takes the following form, also used by Houdayer and Ioana~\cite[proof of Theorem~5.1]{houdayer-ioana}:
\begin{align*}
x_r&:=\tau(p_{r+1}+\cdots+p_a)\,p_r-\tau(p_r)\,(p_{r+1}+\cdots+p_a),
\qquad 1\le r\le a-1,\\
y_s&:=\tau(q_{s+1}+\cdots+q_b)\,q_s-\tau(q_s)\,(q_{s+1}+\cdots+q_b),
\qquad 1\le s\le b-1.
\end{align*}
By that proposition, $X:=\{x_1,\dots,x_{a-1}\}\subseteq A_0\ominus
\C 1$ and $Y:=\{y_1,\dots,y_{b-1}\}\subseteq B_0\ominus \C 1$ are
$L^2$-orthogonal bases of self-adjoint contractions. They may be empty when $a=1$ or $b=1$, in which case the corresponding centered part is zero.

Note that the constructions below use only that $X$ and $Y$ are disjoint finite sets of self-adjoint letters,
so we can later apply them to other alphabets.
For $k\ge 1$, consider the formal alternating words in the alphabet $X\sqcup Y$:
\begin{align*}
\mathcal V_k^{xy}
&:=\{x_{i_1}y_{j_1}\cdots x_{i_k}y_{j_k}: 1\le i_r\le a-1,\ 1\le j_r\le b-1\},\\
\mathcal V_k^{yx}
&:=\{y_{j_1}x_{i_1}\cdots y_{j_k}x_{i_k}: 1\le i_r\le a-1,\ 1\le j_r\le b-1\}.
\end{align*}
We evaluate these words in $M$ when taking traces or $L^2$-inner products.
For a word $z_1z_2\cdots z_m$ of length $m\ge 1$, let
$\rho(z_1z_2\cdots z_m):=z_2\cdots z_m z_1$. Let $\sim$ be the equivalence relation on
$\mathcal V_k^{xy}\sqcup \mathcal V_k^{yx}$ generated by $\rho$ and adjoints.
Every $\sim$-class has a representative beginning with an $X$-letter,
since a word beginning with a $Y$-letter acquires this form after one rotation.
For each $k\le n$, choose a set of representatives $\mathcal S_k\subseteq \mathcal V_k^{xy}$, one from each $\sim$-class, and let
$\mathcal S_{\le n}:=\bigcup_{k=1}^n \mathcal S_k$. Split
$\mathcal S_{\le n}=\mathcal S_{\le n}^{\mathrm{sym}}\sqcup \mathcal S_{\le n}^{\mathrm{ns}}$, where $\mathcal S_{\le n}^{\mathrm{sym}}$ consists of the representatives $w$ for which $w^*$
belongs to the $\rho$-orbit of $w$. We call these representatives and their rotation
orbits symmetric, and the others nonsymmetric.

We record three facts about rotations, which are used in the proof of
Proposition~\ref{prop:noncommutative-free-gram}. Take a formal
word $w\in\mathcal V_k^{xy}$, of length $2k$. Its even rotations begin
in $X$, and its odd rotations begin in $Y$. Hence, if $\rho^j(w)=\rho^{j'}(w)$ as formal words, then
$j-j'$ is even, and the two rotations carry the same sign in
$\sum_{m=0}^{2k-1}(-1)^m\rho^m(w)$.
Moreover, the rotation orbit of a representative in $\mathcal S_{\le
n}^{\mathrm{ns}}$ is disjoint from that of its adjoint. Finally, distinct
representatives in $\mathcal S_{\le n}$ have disjoint orbits under rotation
and adjoint.

\begin{proposition}
\label{prop:noncommutative-free-gram}
Let $(A,\tau_A)$ and $(B,\tau_B)$ be tracial von Neumann algebras, and let
\[
\mathcal X=\{\xi_1,\ldots,\xi_p\}\subseteq A\ominus\C1,
\qquad
\mathcal Y=\{\eta_1,\ldots,\eta_q\}\subseteq B\ominus\C1
\]
be finite $L^2$-orthogonal families of nonzero self-adjoint contractions. Regard the two families as disjoint alphabets.
Fix $n\ge1$ and $\kappa\in(0,1]$, and assume
\[
\norm{z}_2^2\ge\kappa
\qquad(z\in\mathcal X\cup\mathcal Y).
\]
In the free product von Neumann algebra
\[
(N,\varphi):=(A,\tau_A)*(B,\tau_B),
\]
use the preceding choice of representatives for alternating formal words
through pair length $n$, now with alphabets $\mathcal X,\mathcal Y$.
As before, distinguish representatives whose adjoints are cyclic rotations
of themselves from those for which this fails. For a representative $w$ of pair length $k$, let
\[
S(w):=\sum_{m=0}^{2k-1}(-1)^m\rho^m(w),
\qquad
\Xi(w):=-\frac{i}{2}S(w).
\]
Set $g_w:=\Xi(w)$ for a symmetric orbit, and set
\[
g_w^{\Re}:=\frac{\Xi(w)+\Xi(w)^*}{2},
\qquad
g_w^{\Im}:=\frac{\Xi(w)-\Xi(w)^*}{2i}
\]
for a nonsymmetric orbit. These cyclic gradients form an $L^2(N,\varphi)$-orthogonal family of nonzero
centered self-adjoint elements. Its real Gram matrix satisfies
\[
G^{\mathrm{fr}}\ge\lambda_0 I,
\qquad
\lambda_0:=\frac18\kappa^{2n}.
\]
If one alphabet is empty, the family is empty and the assertion is vacuous.
\end{proposition}

\begin{proof}
Assume both alphabets are nonempty. The free-product Hilbert
space identifies a centered alternating word $z_1\cdots z_m$ with
$z_1\otimes\cdots\otimes z_m$. Distinct formal words are therefore
orthogonal, and
\[
\norm{z_1\cdots z_m}_2^2=\prod_{j=1}^m\norm{z_j}_2^2.
\]
For a representative $w$ of pair length $k\le n$, equal rotations in
$S(w)$ have the same sign. After they are collected, the coefficients are
nonzero integers, so
\[
\norm{S(w)}_2^2\ge\kappa^{2k}\ge\kappa^{2n}.
\]
The rotations and adjoints of distinct representatives span
orthogonal subspaces.

If $w^*=\rho^t(w)$, then $t$ is odd, so $S(w)^*=-S(w)$ and
$g_w=-iS(w)/2$ is self-adjoint with squared norm
$\norm{S(w)}_2^2/4$. Otherwise $S(w)$ and $S(w)^*$ are orthogonal
with equal norms. Hence $g_w^{\Re}$ and $g_w^{\Im}$ are orthogonal
and self-adjoint, each with squared norm $\norm{S(w)}_2^2/8$.
By traciality and the alternating signs, $\varphi(S(w))=0$.
Thus all the gradients are centered and their Gram matrix is at least
$(\kappa^{2n}/8)I$, as desired.
\end{proof}

The signed cyclic sum $S(w)$ records the derivative of a word moment under
conjugation of the first alphabet. In
Proposition~\ref{prop:pair-linearization} below, we identify that derivative
as $h\mapsto-2\varphi(h\Xi(w))$, and the real and
imaginary parts of the moment give the real moment equations. The choice of orbit
representatives removes duplicate equations, and a symmetric orbit
needs only its real equation.

Fix $n,p,q\ge1$ and a set $K$. For each $k\in K$, let $(M_k,\tau_k)$ be a tracial von Neumann algebra and let
\[
X_k=(x_{1,k},\ldots,x_{p,k}),
\qquad
Y_k=(y_{1,k},\ldots,y_{q,k})
\]
be lists of centered self-adjoint contractions. Regard their indices as two
disjoint alphabets, and choose the same formal representatives
$\mathcal S_{\le n}=\mathcal S_{\le n}^{\mathrm{sym}}\sqcup\mathcal S_{\le n}^{\mathrm{ns}}$ at every $k$, as above.
We write $\Xi_k(w)$ for the element $\Xi(w)$ of Proposition~\ref{prop:noncommutative-free-gram}, evaluated at $X_k,Y_k$ in $M_k$. Evaluate the corresponding real cyclic gradients in $M_k$ and order them as
\[
\Gamma_k=(\gamma_{1,k},\ldots,\gamma_{d,k}).
\]
For a formal word $w$, let $w_k(h)$ be obtained by conjugating every $X_k$-letter by $e^{ih}$ and leaving every $Y_k$-letter
fixed. Let $\Phi_k:(M_k)_{\rm sa}\to\R^d$ be the map whose coordinates at $h$ are
\begin{itemize}
\item $\Ree\tau_k(w_k(h))$ for each $w\in \mathcal S_{\le n}^{\mathrm{sym}}$, and
\item $\Ree\tau_k(w_k(h))$ and $\Imm\tau_k(w_k(h))$ for each $w\in \mathcal S_{\le n}^{\mathrm{ns}}$.
\end{itemize}
We call $\Phi_k$ the \emph{constraint map}. It has
$d=\abs{\mathcal S_{\le n}^{\mathrm{sym}}}+2\abs{\mathcal S_{\le n}^{\mathrm{ns}}}$ coordinates. On $\R^d$, we write $\norm{\cdot}$ for the Euclidean norm, and also for the operator norm of matrices acting on $\R^d$.

If $w$ is symmetric, then $w^*=\rho^s(w)$ for some odd $s$, and this relation persists after perturbation. Hence
\[
\tau_k(w_k(h))
=\tau_k(\rho^s(w_k(h)))
=\tau_k(w_k(h)^*)
=\overline{\tau_k(w_k(h))}
\]
for every $h=h^*$, so one real constraint suffices. As in the proof of
Proposition~\ref{prop:noncommutative-free-gram}, reindexing the signed cyclic sum gives
$\Xi_k(w)^*=\Xi_k(w)$ at $h=0$.

If $h=h^*$ and $\Phi_k(h)=0$, every alternating word through
pair length $n$ in the perturbed lists has trace zero. Indeed,
traciality and adjoints reduce its trace to that of a chosen
representative or its complex conjugate, and a symmetric
representative has real trace. For the complex spans
$E_k=\operatorname{span}_{\C}X_k$ and
$F_k=\operatorname{span}_{\C}Y_k$, multilinearity therefore gives
\begin{equation}
\tau_k(e^{ih}e_1e^{-ih}f_1\cdots e^{ih}e_re^{-ih}f_r)=0
\quad(1\le r\le n,\ e_i\in E_k,\ f_i\in F_k).
\label{eq:cyclic-zero-independence}
\end{equation}

\begin{proposition}\label{prop:pair-linearization}
For every $k$, the map $\Phi_k:(M_k)_{\rm sa}\to\R^d$ is real Fr\'echet differentiable at $0$ with respect to the $L^2$-norm, and
for $h\in (M_k)_{\rm sa}$,
\[
D\Phi_k(0)[h]
=\bigl(-2\tau_k(h\gamma_{1,k}),\dots,-2\tau_k(h\gamma_{d,k})\bigr).
\]
Moreover, with $R_k(h):=\Phi_k(h)-\Phi_k(0)-D\Phi_k(0)[h]$, one may take the constant
$C_{\mathrm{nl}}=2n^2\sqrt d$, independently of $k$, in
\[
\norm{R_k(h)-R_k(h')}
\le C_{\mathrm{nl}}(\norm{h}_2+\norm{h'}_2)\norm{h-h'}_2
\]
for all $h,h'\in (M_k)_{\rm sa}$. Consequently, $\norm{R_k(h)}\le C_{\mathrm{nl}}\norm{h}_2^2$.
\end{proposition}

\begin{proof}
For a representative $w$ of pair length $\ell\le n$, let $w_k(t;h)$
be the word obtained by conjugating its $X_k$-letters by $e^{ith}$,
and set
\[
\Xi_k^t(w;h)=-\frac i2\sum_{m=0}^{2\ell-1}(-1)^m\rho^m(w_k(t;h)).
\]
Let $R_{w,k}(h):=\tau_k(w_k(h))-\tau_k(w_k(0))+2\tau_k(h\Xi_k(w))$ be the remainder of the trace of $w$. By differentiation and traciality, we have
\[
\frac d{dt}\tau_k(w_k(t;h))=-2\tau_k(h\Xi_k^t(w;h)),
\qquad
R_{w,k}(h)=-2\int_0^1\tau_k\bigl(h(\Xi_k^t(w;h)-\Xi_k(w))\bigr)\,dt.
\]
Integrating the derivative of $s\mapsto e^{i(t-s)h}e^{ish'}$ on $[0,t]$
and applying the triangle inequality gives
$\norm{e^{ith}xe^{-ith}-e^{ith'}xe^{-ith'}}_2
\le2t\norm{x}_\infty\norm{h-h'}_2$.
By telescoping over the $\ell$ conjugated letters in each rotation, we have
\[
\norm{\Xi_k^t(w;h)-\Xi_k^t(w;h')}_2
 \le2\ell^2t\norm{h-h'}_2.
\]
Substituting this bound and its special case $h'=0$ into the integral and applying
Cauchy--Schwarz, we obtain
\[
|R_{w,k}(h)-R_{w,k}(h')|
\le2\ell^2(\norm h_2+\norm{h'}_2)\norm{h-h'}_2.
\]
Taking the prescribed real and imaginary coordinates gives the stated
derivative and the factor $\sqrt d$ in the remainder bound.
Setting $h'=0$, we obtain $L^2$-Fr\'echet differentiability. Note that all bounds are global
in $h,h'$, since the exponentials are unitary.
\end{proof}

\begin{lemma}
\label{lem:finite-dimensional-correction}
Let $E,F$ be finite-dimensional real normed spaces, with $|\cdot|$
denoting the norm on $F$. Let $T:E\to F$ be linear and let
$\Phi(h)=v+Th+R(h)$ be an $F$-valued map, with $R(0)=0$.
Let $J:F\to E$ be linear with $TJ=-I_F$ and
\[
\norm{Jy}\le A|y|,\qquad \norm{Jy}_*\le B|y|,
\]
where $A>0$, $B,C\ge0$, and $\norm{\cdot}_*$ is a second norm on $E$. Set $\eta=|v|$.
Suppose $\Phi$ is defined on the closed ball of radius $2A\eta$ and, there,
\[
|R(h)-R(k)|\le C(\norm h+\norm k)\norm{h-k},
\qquad 8A^2C\eta\le1.
\]
Then $\Phi(h)=0$ has a solution satisfying
\[
\norm h\le2A\eta,\qquad \norm h_*\le\frac32 B\eta.
\]
It is the limit of the iterates $h_0=0$, $h_{j+1}=J(v+R(h_j))$.
\end{lemma}

\begin{proof}
If $\eta=0$, take $h=0$. Otherwise set $r=2A\eta$. By the remainder bound with $k=0$, we have $|R(h)|\le C\norm
h^2$, and since $8A^2C\eta\le1$, we get $Cr^2\le\eta/2$. Hence
$\mathcal T(h)=J(v+R(h))$ maps the closed $r$-ball into the closed
ball of radius $3r/4$. On the closed $r$-ball, it has Lipschitz
constant at most $2ACr=4A^2C\eta\le1/2$.
By the contraction mapping theorem, the indicated iterates converge to a fixed point
$h=J(v+R(h))$. Applying $T$ yields $Th=-v-R(h)$, so $\Phi(h)=0$, and the fixed-point identity also gives
$\norm h_*\le B(\eta+Cr^2)\le3B\eta/2$.
\end{proof}

In matrix applications, suppose the data vary in a fixed real matrix space, $J$ and $v$
depend Borel measurably on those data, and $R$ depends jointly Borel measurably on the data and on $h$.
Then, on any Borel set where the hypotheses hold, the iterates and their limit are Borel measurable,
and extending the limit by zero outside that set gives a Borel choice of the self-adjoint correction.

\begin{theorem}
\label{thm:fd-noncommutative-exactify}
Fix $n\ge1$, finite-dimensional unital $C^*$-algebras
\[
P=\bigoplus_{r=1}^{R}M_{d_r}(\C),
\qquad
Q=\bigoplus_{s=1}^{S}M_{e_s}(\C),
\]
and $\delta>0$. There are constants
\[
\varepsilon_0=\varepsilon_0(n,P,Q,\delta)>0,
\qquad
C=C(n,P,Q,\delta)<\infty
\]
with the following property. Let $(M,\tau)$ be a tracial von Neumann algebra and let $A,B\subseteq M$ be finite-dimensional unital $^*$-subalgebras $^*$-isomorphic to
$P$ and $Q$, respectively. Assume that every minimal projection of $A$ and of $B$
has trace at least $\delta$. If
\[
\varepsilon:=\varepsilon_{2n}(A,B)\le\varepsilon_0,
\]
then there is $h=h^*\in M$ such that, for $u=e^{ih}$, the algebras $uAu^*$ and $B$ are $n$-independent and
\[
\norm{h}_2\le C\varepsilon,
\qquad
\norm{h}_\infty\le C\varepsilon,
\qquad
\norm{u-1}_\infty\le C\varepsilon.
\]

If neither $P$ nor $Q$ is scalar, set
\[
D_0:=\max\{d_1,\ldots,d_R,e_1,\ldots,e_S\},
\qquad
\kappa:=\min\left\{1,2\delta^3,\frac{\delta}{2},\frac{2\delta}{D_0^2}\right\}.
\]
Let $d=d(n,P,Q)$ be the number of coordinates of the constraint map for alphabets with
$\dim_{\C}P-1$ and $\dim_{\C}Q-1$ letters, let $(K_\ell)$ be the constants from
Lemma~\ref{lem:quantitative-truncated-comparison}, and set
\[
\lambda_0:=\frac18\kappa^{2n},
\qquad
K_G:=dn^2K_{4n},
\qquad
C_{\mathrm{nl}}:=2n^2\sqrt d,
\qquad
L_0:=(2\lambda_0)^{-1/2}.
\]
One may take
\[
\varepsilon_0
:=\min\left\{
\frac{\lambda_0}{2K_G},
\frac{\lambda_0}{4C_{\mathrm{nl}}\sqrt d}
\right\},
\qquad
C:=\max\left\{
2L_0\sqrt d,
\frac{3nd}{2\lambda_0}
\right\}.
\]
If $P=\C$ or $Q=\C$, one may take $h=0$, $C=0$, and arbitrary positive $\varepsilon_0$.
\end{theorem}

\begin{proof}
We omit the coordinate index $k$ from the notation above. If either algebra is scalar, take $h=0$. Otherwise choose the bases of
Proposition~\ref{prop:fd-selfadjoint-basis}, with squared $L^2$-norms at least
$\kappa$, and form the constraint map $\Phi$, the gradients $\Gamma=(\gamma_1,\ldots,\gamma_d)$,
and their Gram matrix $G$. By Proposition~\ref{prop:noncommutative-free-gram}, the Gram
matrix $G^{\mathrm{fr}}$ of the corresponding gradients in the
free product satisfies $G^{\mathrm{fr}}\ge\lambda_0 I$.
Note that the sum of the absolute values of the coefficients in each
gradient is at most $n$, and products of two gradient words have at most
$4n$ alternating contraction blocks and at most $2n$
blocks from either algebra. Hence, by Lemma~\ref{lem:quantitative-truncated-comparison},
we have
\[
|G_{ij}-G^{\mathrm{fr}}_{ij}|\le n^2K_{4n}\varepsilon,
\qquad
\norm{G-G^{\mathrm{fr}}}\le K_G\varepsilon\le\lambda_0/2.
\]
Thus $G\ge(\lambda_0/2)I$. Set $E=\operatorname{span}_{\R}\Gamma$ and define $L:\R^d\to E$ by
\[
L(y)=\frac12\sum_{\ell=1}^d(G^{-1}y)_\ell\gamma_\ell.
\]
The Gram identity and $\norm{\gamma_\ell}_\infty\le n$ imply
\[
D\Phi(0)L=-I,\qquad LD\Phi(0)=-I_E,\qquad
\norm{L(y)}_2\le L_0\norm y,\qquad
\norm{L(y)}_\infty\le\frac{n\sqrt d}{\lambda_0}\norm y.
\]
We apply Lemma~\ref{lem:finite-dimensional-correction} with $J=L$,
$A=L_0$, $B=n\sqrt d/\lambda_0$, $C=C_{\mathrm{nl}}$, and
$\eta=\norm{\Phi(0)}\le\sqrt d\,\varepsilon$.
By Proposition~\ref{prop:pair-linearization}, the remainder bound holds,
and the second restriction on $\varepsilon_0$ gives
$8L_0^2C_{\mathrm{nl}}\eta\le1$. The resulting zero satisfies
\[
\norm h_2\le2L_0\sqrt d\,\varepsilon,
\qquad
\norm h_\infty\le\frac{3nd}{2\lambda_0}\varepsilon.
\]
Since the complex spans of the bases are $A\ominus\C1$ and $B\ominus\C1$, we
obtain $n$-independence from \eqref{eq:cyclic-zero-independence}.
Finally, $\norm{e^{ih}-1}_\infty\le\norm h_\infty$, as desired.
\end{proof}

The $2n$-to-$n$ correction above may also be derived, with possibly different
constants, from the local correction theorem independently obtained by
Patchell and Shiner~\cite[Theorem~2.5]{patchell-shiner}.

\begin{corollary}
\label{thm:quantitative-local-exactify}
Fix $n,a,b\ge1$ and $\delta>0$. There are $\varepsilon_0=\varepsilon_0(n,a,b,\delta)>0$ and
$C=C(n,a,b,\delta)<\infty$ with the following property. Let $(M,\tau)$ be a tracial von Neumann algebra, and let
\[
A=\bigoplus_{r=1}^a\C p_r,
\qquad
B=\bigoplus_{s=1}^b\C q_s
\]
be finite-dimensional abelian unital von Neumann subalgebras of $M$ whose displayed minimal projections have trace at least $\delta$.
If $\varepsilon:=\varepsilon_{2n}(A,B)\le\varepsilon_0$, then there is $h=h^*\in M$ such that, for $w=e^{ih}$,
the algebras $wAw^*$ and $B$ are $n$-independent and
\[
\norm{h}_2,\ \norm{h}_\infty,\ \norm{w-1}_\infty\le C\varepsilon.
\]
For $a,b\ge2$, one may take the constants of Theorem~\ref{thm:fd-noncommutative-exactify} with
$P=\C^a$, $Q=\C^b$, $D_0=1$, $\kappa=2\delta^3$, and $d=d(n,a,b)$.
If $a=1$ or $b=1$, one may take arbitrary $\varepsilon_0>0$ and $C=0$.
\end{corollary}

\begin{proof}
We apply Theorem~\ref{thm:fd-noncommutative-exactify}. If $a,b\ge2$, then the trace bound implies $\delta\le1/2$, hence
\[
\min\{1,2\delta^3,\delta/2,2\delta\}=2\delta^3.
\]
If $a=1$ or $b=1$, take $h=0$.
\end{proof}

In the proof of Theorem~\ref{thm:pair-local-exactify}, we use the following observation. Fix $m\ge1$ and finite-dimensional unital $C^*$-algebras $P,Q$.
Consider unital $^*$-embeddings $\pi_t:P\to M_t$ and
$\sigma_t:Q\to M_t$ into tracial von Neumann algebras.
Convergence here means either ordinary sequential convergence
or convergence along a fixed ultrafilter.
Suppose that the induced states on $P$ and $Q$ converge in norm to $\alpha,\beta$
and that every fixed mixed moment through pair length $m$ converges
to the corresponding moment of an $m$-independent pair with these
marginal traces. Then
$\varepsilon_m(\pi_t(P),\sigma_t(Q))\to0$.
Indeed, after pulling back to $P,Q$ and centering with respect to the
induced states, the relevant multilinear moment forms converge
pointwise to zero. By finite dimensionality, this convergence is uniform on
products of the unit balls, at each of the finitely many
pair lengths.

\begin{theorem}\label{thm:pair-local-exactify}
Let $n\ge1$, let $M=\ultraprod M_k$ be a tracial ultraproduct of von Neumann algebras, and let
\[
A_0=\bigoplus_{r=1}^a\C p_r,
\qquad
B_0=\bigoplus_{s=1}^b\C q_s
\]
be $2n$-independent unital $^*$-subalgebras of $M$. For each $k$, let $A_k,B_k\subseteq M_k$ be generated by partitions
$(p_{r,k})_{r=1}^a$ and $(q_{s,k})_{s=1}^b$ of $1$, respectively, allowing zero projections, where
\[
p_r=(p_{r,k})_{\cU},
\qquad
q_s=(q_{s,k})_{\cU}.
\]
There are $S\in\cU$ and self-adjoint elements $h_k\in M_k$, $k\in S$, such that,
after setting $h_k:=0$ off $S$ and $w_k:=e^{ih_k}$ for every $k$, we have
\[
\norm{h_k}_2\longrightarrow_{\cU}0,
\qquad
\norm{h_k}_\infty\longrightarrow_{\cU}0,
\qquad
\norm{w_k-1}_\infty\longrightarrow_{\cU}0,
\]
and $w_kA_k w_k^*$ and $B_k$ are $n$-independent for every $k\in S$.
\end{theorem}

\begin{proof}
If $a=1$ or $b=1$, take $h_k=0$. Otherwise, on a
$\cU$-large set, the lifted minimal projections have traces at least
half the smallest limiting trace. By the observation before the theorem, with $m=2n$,
we have $\varepsilon_{2n}(A_k,B_k)\to_{\cU}0$.
Hence we may apply Corollary~\ref{thm:quantitative-local-exactify} on some
$S\in\cU$ and obtain $n$-independence with
\[
\|h_k\|_2,\ \|h_k\|_\infty,\ \|e^{ih_k}-1\|_\infty
\le C\varepsilon_{2n}(A_k,B_k).
\]
Set $h_k=0$ off $S$.
\end{proof}

\section{Preserving the order of independence in factors}
\label{sec:higher-order-no-loss}

In Theorem~\ref{thm:fd-noncommutative-exactify}, we assumed a small defect
of order $2n$ to obtain independence of order $n$, because the entries of
the Gram matrix at the given pair are mixed moments of order up to $2n$.
Linearizing at a nearby $n$-independent pair does not resolve this, since
the gradients there can be linearly dependent. For example, if $A$ and $B$
commute and are $1$-independent, every gradient at order one is a
multiple of a commutator $[x,y]$ and vanishes. In a $\mathrm{II}_1$
factor, we first make a small conjugation after which the Gram matrix is
bounded below by a power of the size of that conjugation. This allows us
to assume only a small defect of order $n$ in Theorem~\ref{thm:fixed-type-factor-no-loss} below.
We also show that finite-order independence can be reached by conjugation in every ${\rm II}_1$ factor
(Theorem~\ref{thm:finite-order-free-positioning}). For partitions at order one, we remove in
Subsection~\ref{sec:order-one-no-loss} the lower bound on the projection traces assumed in
Theorem~\ref{thm:fixed-type-factor-no-loss}, by allowing both partitions to move.

\begin{theorem}[Power-law stability]
\label{thm:fixed-type-factor-no-loss}
Fix $n\ge1$, finite-dimensional unital $C^*$-algebras $P,Q$, and $\delta>0$.  There are constants
\[
\varepsilon_0=\varepsilon_0(n,P,Q,\delta)>0,
\qquad
C=C(n,P,Q,\delta)<\infty
\]
that work in every ${\rm II}_1$ factor $(M,\tau)$. Fix unital $^*$-subalgebras
$A,B\subseteq M$ and $^*$-isomorphisms $P\cong A$ and $Q\cong B$.
Suppose that every minimal projection of $A$ and $B$ has trace at
least $\delta$ and that
\[
\varepsilon:=\varepsilon_n(A,B)\le\varepsilon_0.
\]
Then there is $u\in U(M)$ such that
\[
\norm{u-1}_\infty\le C\varepsilon^{1/(2n^2+2n+1)}
\]
and $uAu^*,B$ are $n$-independent.
If either $P$ or $Q$ is scalar, one may take $u=1$.
\end{theorem}

Working in a free product, we conjugate the first algebra by a unitary
with three spectral values and leave the second fixed. If the original
pair is $n$-independent, the conjugated pair remains $n$-independent and
its moment derivative is surjective (Lemma~\ref{lem:free-marker-regularity}).
If the defect of order $n$ is only small, the centered alternating moments defining that defect remain small after conjugation,
and the Gram matrix of the moment gradients is bounded below by a power of the size of the conjugation
(Proposition~\ref{prop:quantitative-phase-marker}).
The calculation uses
higher moments involving the auxiliary unitary, while the only smallness assumption
on the original pair is on its defect of order $n$.
Using Theorem~\ref{thm:finite-order-free-positioning} below, we realize the
moments used in this calculation inside the given factor. By
Lemma~\ref{lem:finite-dimensional-correction}, we then make the defining moments of $n$-independence vanish. Taking the size of the first
conjugation to be a power of the defect gives the exponent in the theorem.
We give the proof after Proposition~\ref{prop:dense-regular-exact-locus}.

To obtain finite-order independence by conjugation, we first perturb representatives of free families in an ultraproduct.

\begin{proposition}
\label{prop:finite-list-local-exactify}
Let $\omega$ be an ultrafilter on a set $J$, let $(N_k,\tau_k)_{k\in J}$ be tracial von Neumann algebras, and let
$N:=\prod_{\omega}(N_k,\tau_k)$, with ultraproduct trace $\tau=(\tau_k)_\omega$. Fix $p,q,n\ge1$.
For $1\le i\le p$ and $1\le j\le q$, let
\[
a_{i,k},b_{j,k}\in (N_k)_{\rm sa}
\]
be contractions satisfying $\tau_k(a_{i,k})=\tau_k(b_{j,k})=0$, and write
$a_i=(a_{i,k})_\omega$ and $b_j=(b_{j,k})_\omega$. Assume that
$a_1,\dots,a_p$ are nonzero and pairwise $L^2$-orthogonal, that $b_1,\dots,b_q$ have the same properties, and that the two
families lie in freely independent unital von Neumann subalgebras of $N$.

Then there are $S\in\omega$ and $h_k=h_k^*\in N_k$, $k\in S$, such that, after setting $h_k=0$ off $S$,
\[
\norm{h_k}_2\longrightarrow_\omega0,
\qquad
\norm{h_k}_\infty\longrightarrow_\omega0,
\qquad
\norm{e^{ih_k}-1}_\infty\longrightarrow_\omega0,
\]
and
\[
\tau_k\bigl(
e^{ih_k}a_{i_1,k}e^{-ih_k}b_{j_1,k}\cdots
e^{ih_k}a_{i_r,k}e^{-ih_k}b_{j_r,k}
\bigr)=0
\]
for every $k\in S$, $1\le r\le n$, and all choices of indices. Equivalently, the complex spans of the two perturbed lists are
$n$-independent for every $k\in S$.
\end{proposition}

\begin{proof}
We form the real constraint maps $\Phi_k:(N_k)_{\rm sa}\to\R^d$ and the Gram matrices
$G_k:=(\tau_k(\gamma_{s,k}\gamma_{t,k}))_{s,t=1}^d$ of their real cyclic gradients $\gamma_{1,k},\ldots,\gamma_{d,k}$.
Set $\kappa=\min\{\|a_i\|_2^2,\|b_j\|_2^2\}>0$. By freeness in $N$ and
Proposition~\ref{prop:noncommutative-free-gram}, we have
\[
G_k\longrightarrow_\omega G\ge\tfrac18\kappa^{2n}I,\qquad
\eta_k:=\|\Phi_k(0)\|\longrightarrow_\omega0.
\]
Choose $C_G<\infty$ and an $\omega$-large set on which
$\|G_k^{-1}\|\le C_G$. Then, using the map $L$ built from the Gram inverse
in the proof of
Theorem~\ref{thm:fd-noncommutative-exactify}, we apply
Lemma~\ref{lem:finite-dimensional-correction} on the real span of the
gradients. After shrinking this set, we obtain
\[
\Phi_k(h_k)=0,\qquad
\|h_k\|_2\le\sqrt{C_G}\eta_k,\qquad
\|h_k\|_\infty\le\tfrac34n\sqrt d\,C_G\eta_k.
\]
By \eqref{eq:cyclic-zero-independence}, the complex spans of the perturbed
lists are $n$-independent. Finally, we extend $h_k$ by zero outside this set
and use
$\|e^{ih_k}-1\|_\infty\le\|h_k\|_\infty$.
\end{proof}

\begin{theorem}[Finite-order independence by conjugation]
\label{thm:finite-order-free-positioning}
Let $(M,\tau)$ be a ${\rm II}_1$ factor and let $n\ge1$.
\begin{enumerate}[label=(\roman*)]
\item If $E,F\subseteq M\ominus\C1$ are finite-dimensional $^*$-closed complex linear subspaces, then there is a unitary
$u\in M$ such that $uEu^*$ and $F$ are $n$-independent.
\item Let $\mathbf x=(x_1,\dots,x_p)$ and $\mathbf y=(y_1,\dots,y_q)$ be finite tuples in $M$, set
$A:=W^*(\mathbf x)$ and $B:=W^*(\mathbf y)$, and let
\[
(P,\varphi):=(A,\tau|_A)*(B,\tau|_B)
\]
be their free product von Neumann algebra, with canonical embeddings $\iota_A:A\to P$ and $\iota_B:B\to P$. There is a unitary $u\in M$
such that
\[
\tau\bigl(R(u\mathbf xu^*,\mathbf y)\bigr)
=\varphi\bigl(R(\iota_A(\mathbf x),\iota_B(\mathbf y))\bigr)
\]
for every noncommutative $^*$-polynomial $R$ of total $^*$-degree at most $2n+1$.
\end{enumerate}
\end{theorem}

Patchell and Shiner independently obtained finite-order independence
by conjugation for finite centered sets~\cite[Definition~2.1 and Theorem~B]{patchell-shiner}, in Popa's convention for
subsets~\cite[Section~3]{popa-kadison-singer}.
One may also obtain part~\textup{(i)} from their theorem by
multilinearity, and part~\textup{(ii)} by applying it to finite bases
for the spaces of centered values of all $^*$-polynomials of degree at most
$2n+1$ in each tuple, then using
Corollary~\ref{cor:filtered-truncated-comparison}.

\begin{proof}
For \textup{(i)}, if $E=\{0\}$ or $F=\{0\}$, take $u=1$.
Otherwise choose $L^2$-orthogonal self-adjoint contractive bases $(e_i)$ of $E$
and $(f_j)$ of $F$.
Fix a nonprincipal ultrafilter $\omega$ on $\N$.
By Popa's asymptotic freeness theorem~\cite[Theorem~2.1]{popa-free-independent-1995},
in the formulation of~\cite[Corollary~0.2]{popa-independence},
applied in $M^\omega$ to $N_1=W^*(E)$, $N_2=W^*(F)$ over $B_1=B_2=\C1$,
there is a unitary $v\in M^\omega$ making $vN_1v^*$ free from
$N_2$. Both $N_1$ and $N_2$ have separable
predual, even when $M$ does not.
Choose unitary representatives $v_k\in M$ and apply
Proposition~\ref{prop:finite-list-local-exactify} to
$(v_ke_iv_k^*)_i,(f_j)_j$. For any $k$ in the $\omega$-large set provided by
that proposition, the unitary $u=e^{ih_k}v_k$ has the property required in
\textup{(i)}, by multilinearity.

For \textup{(ii)}, we apply \textup{(i)} to the centered images of all
marginal $^*$-polynomials of degree at most $d=2n+1$, which are
finite-dimensional $^*$-closed spaces. By
Corollary~\ref{cor:filtered-truncated-comparison}, applied with $D=d$ and
$m=n$ to the polynomial filtrations of the conjugated first tuple and of the
second tuple, we obtain the required identities.
Conjugation preserves the first marginal trace, so the free-product
trace of the conjugated pair is the one in the statement.
\end{proof}

For finite $X,Y\subseteq M\ominus\C1$, applying part~(i) to $\operatorname{span}_{\C}(X\cup X^*)$ and $\operatorname{span}_{\C}(Y\cup Y^*)$,
we obtain a unitary $u$ for which $uXu^*$ and $Y$ are $n$-independent. At $n=1$, this recovers a theorem of Houdayer and Ioana~\cite[Theorem~H]{houdayer-ioana},
and our proof follows the same route as theirs. Similarly, by applying part~(i) to the
centered parts of two finite-dimensional unital von Neumann subalgebras of a ${\rm
II}_1$ factor, we can put them in $n$-independent position.

\subsection{Regularizing conjugation}

For an $n$-independent pair of finite-dimensional unital $^*$-subalgebras,
we seek a small conjugation of the first that preserves $n$-independence
and makes the moment gradients linearly independent. For a pair with a small
defect, the same
conjugation keeps the defect small, and a quantitative lower bound for the Gram matrix of the gradients gives the power law in
Theorem~\ref{thm:fixed-type-factor-no-loss}.

Fix ordered matrix units in
$P=\bigoplus_{r=1}^R M_{d_r}(\C)$ and
$Q=\bigoplus_{s=1}^S M_{e_s}(\C)$, distinguishing repeated isomorphic
summands. For each tracial state $\alpha$ on $P$, let
$\mathcal X_P(\alpha)$ be the list given by the formulas in the proof of
Proposition~\ref{prop:fd-selfadjoint-basis}, with the central vectors
first and the other vectors in a fixed order. Only the central
coefficients depend on $\alpha$, so this list varies entrywise
continuously. When $\alpha$ is faithful, $\mathcal X_P(\alpha)$ is a real
$L^2(P,\alpha)$-orthogonal basis of $(P\ominus\C1)_{\rm sa}$ consisting of
centered self-adjoint contractions, and its complex span is $P\ominus\C1$.
If every minimal projection has $\alpha$-trace at least $\delta$,
each basis vector has squared $L^2$-norm at least $\kappa(P,\delta)$.
Define $\mathcal X_Q(\beta)$ analogously. For unital $^*$-embeddings
$\pi:P\to M$ and $\sigma:Q\to M$ into a tracial von Neumann algebra $(M,\tau)$, write
\[
X_\pi:=\pi\bigl(\mathcal X_P(\tau\circ\pi)\bigr),\qquad
Y_\sigma:=\sigma\bigl(\mathcal X_Q(\tau\circ\sigma)\bigr).
\]

We apply the construction of
Subsection~\ref{sec:pair-cyclic-coordinates} to these bases, regarded as disjoint alphabets. Let
$\Gamma=(\gamma_1,\ldots,\gamma_d)$ denote the family of real cyclic gradients through pair length $n$, formed from
$X_\pi$ and $Y_\sigma$ using the fixed matrix units. We call an $n$-independent pair $\pi(P),\sigma(Q)$
\emph{order-$n$ regular} when
\[
G:=\bigl(\tau(\gamma_s\gamma_t)\bigr)_{s,t=1}^d
\]
is positive definite.  Equivalently, the derivative at zero of the constraint map obtained by conjugating the first
algebra is onto $\R^d$.  If either $P$ or $Q$ is scalar, then $d=0$ and regularity is vacuous.

The following lemma is the analogue of
Proposition~\ref{prop:finite-list-local-exactify} for regular pairs:
regularity supplies the lower Gram bound, so the limiting pair need only be $n$-independent, not free.
We use it, together with Proposition~\ref{prop:dense-regular-exact-locus}, in the proof of
Corollary~\ref{cor:large-matrix-fixed-type-no-loss}.

\begin{lemma}
\label{lem:regular-limit-exactification}
Fix $n\ge1$ and $P,Q$ as above. Let $\omega$ be an ultrafilter on $J$, let $(M_k,\tau_k)_{k\in J}$ be tracial von Neumann algebras, and let
$\pi_k:P\to M_k$ and $\sigma_k:Q\to M_k$ be unital $^*$-embeddings.  Set
\[
\alpha_k:=\tau_k\circ\pi_k,
\qquad
\beta_k:=\tau_k\circ\sigma_k,
\qquad
\alpha:=\lim_\omega\alpha_k,
\qquad
\beta:=\lim_\omega\beta_k,
\]
and assume that each induced trace assigns every minimal projection a value at least $\delta>0$.  Let $\pi,\sigma$
be the induced embeddings into $\prod_\omega M_k$.  Suppose that $\pi(P),\sigma(Q)$ are $n$-independent and
order-$n$ regular relative to the bases obtained from the limiting traces
\[
\pi\bigl(\mathcal X_P(\alpha)\bigr),
\qquad
\sigma\bigl(\mathcal X_Q(\beta)\bigr).
\]
Then, on an $\omega$-large set, there are $h_k=h_k^*\in M_k$ such that
\[
\norm{h_k}_2\longrightarrow_\omega0,
\qquad
\norm{h_k}_\infty\longrightarrow_\omega0,
\]
and the unital $^*$-subalgebras
\[
e^{ih_k}\pi_k(P)e^{-ih_k},\qquad \sigma_k(Q)
\]
are $n$-independent.
\end{lemma}

\begin{proof}
If $P$ or $Q$ is scalar, take $h_k=0$. Otherwise the continuously
chosen bases $X_{\pi_k},Y_{\sigma_k}$ consist of centered contractions
and represent the limiting bases. Let $\Phi_k$ and $G_k$ be the real
constraint map and cyclic-gradient Gram matrix formed from $X_{\pi_k}$ and
$Y_{\sigma_k}$. By independence and regularity, we have
\[
\eta_k:=\|\Phi_k(0)\|\longrightarrow_\omega0,
\qquad G_k\longrightarrow_\omega G>0.
\]
Since $G_k^{-1}$ is uniformly bounded on an $\omega$-large set,
the correction argument in the proof of
Proposition~\ref{prop:finite-list-local-exactify} gives, after shrinking
that set, $\Phi_k(h_k)=0$ and $\|h_k\|_\infty\to_\omega0$.
The complex spans of these bases are $\pi_k(P)\ominus\C1$ and
$\sigma_k(Q)\ominus\C1$, so we obtain $n$-independence from
\eqref{eq:cyclic-zero-independence}, as desired.
\end{proof}

With the next lemma, we can compute the traces of words $c_0z_1c_1\cdots z_mc_m$ in a free product $C*P$,
with $c_j\in C$ and $z_i\in P$, from finitely many moments in $P$. We use it for conjugations by unitaries in $C$.

\begin{lemma}
\label{lem:finite-free-substitution}
Let $(C,\tau_C)$, $(P,\tau_P)$, and $(P',\tau_{P'})$ be tracial von Neumann algebras.  Suppose
\[
z_1,\ldots,z_m\in P,
\qquad
z'_1,\ldots,z'_m\in P'
\]
satisfy
\[
\tau_P(z_{i_1}\cdots z_{i_s})
=\tau_{P'}(z'_{i_1}\cdots z'_{i_s})
\]
for every increasing sequence $1\le i_1<\cdots<i_s\le m$, including the empty sequence.  Then, for every
$c_0,\ldots,c_m\in C$,
\[
\tau_{C*P}(c_0z_1c_1\cdots z_mc_m)
=\tau_{C*P'}(c_0z'_1c_1\cdots z'_mc_m),
\]
where both sides use the canonical free-product traces.

Moreover, for every $m\ge0$ there is a constant $L_m<\infty$, depending only on $m$, with the following property.  If all
$z_i,z_i',c_j$ are contractions and
\[
\abs{\tau_P(z_{i_1}\cdots z_{i_s})-\tau_{P'}(z'_{i_1}\cdots z'_{i_s})}
\le\eta
\]
for every increasing sequence, then
\[
\abs{\tau_{C*P}(c_0z_1c_1\cdots z_mc_m)
-\tau_{C*P'}(c_0z'_1c_1\cdots z'_mc_m)}
\le L_m\eta.
\]
\end{lemma}

\begin{proof}
We recursively center the maximal blocks in the two free-product words. Then
the fully centered alternating term has trace zero, and every other term has
fewer blocks after deleting scalar blocks and multiplying neighbors.
Every resulting $P$-monomial is of the form
$z_{i_1}\cdots z_{i_s}$ with $i_1<\cdots<i_s$, even when deleting
scalar blocks brings originally nonconsecutive blocks together.
Thus both traces are values of one universal polynomial
\[
F_m(\mathbf c,\mathbf z)=\sum_\nu a_\nu
       \mathbf c^{\alpha_\nu}\mathbf z^{\beta_\nu},
\]
where $\mathbf c$ consists of moments of ordered products of the
$c_j$, and $\mathbf z$ of the moments
$\tau_P(z_{i_1}\cdots z_{i_s})$ with $i_1<\cdots<i_s$. The trace in $C*P'$ is $F_m(\mathbf
c,\mathbf z')$ for the corresponding moments $\mathbf z'$ in $P'$. This proves the equality
assertion.

For contractions, all these variables lie in the unit ball of $\C$.
Telescoping each monomial therefore gives
\[
|F_m(\mathbf c,\mathbf z)-F_m(\mathbf c,\mathbf z')|
\le\left(\sum_\nu|a_\nu||\beta_\nu|_1\right)\eta.
\]
The sum in parentheses depends only on $m$ and can be computed from
the recursion, so we may take it as $L_m$.
\end{proof}

\begin{lemma}
\label{lem:free-marker-preservation}
Let $n\ge1$, and let $A,B$ be unital $n$-independent von Neumann subalgebras of a tracial von Neumann
algebra $(P,\tau)$. Given a tracial von Neumann algebra $(C,\tau_C)$, form
\[
(N,\varphi):=(C,\tau_C)*(P,\tau).
\]
For every unitary $u\in C$, the subalgebras $uAu^*$ and $B$ are $n$-independent in $N$.
\end{lemma}

\begin{proof}
Fix $1\le r\le n$ and centered letters $a_1,\ldots,a_r\in A$ and $b_1,\ldots,b_r\in B$.  Consider the ordered list
\[
z_1=a_1,z_2=b_1,\ldots,z_{2r-1}=a_r,z_{2r}=b_r.
\]
For any increasing subsequence, multiply consecutive selected letters from
the same algebra. This gives a word whose blocks alternate between $A$ and
$B$, with at most $2r\le2n$ blocks. By Lemma~\ref{prop:truncated-free-product},
which applies through $2n+1$ blocks, every such subsequence moment in $P$
equals the corresponding moment in
\[
M:=(A,\tau|_A)*(B,\tau|_B).
\]
Hence, by Lemma~\ref{lem:finite-free-substitution}, with the intervening copies of $u$ and $u^*$ as the $C$-letters, we have
\[
\varphi(ua_1u^*b_1\cdots ua_ru^*b_r)
=\tau_{C*M}(ua_1u^*b_1\cdots ua_ru^*b_r).
\]
In $C*M\cong C*A*B$, the algebra $B$ is free from $W^*(C,A)$, and $uAu^*\subseteq W^*(C,A)$, so the last trace is zero.
\end{proof}

Write the auxiliary unitary as $u=\tau_C(u)+v$, and suppose that
$v\ne0$. The vectors $v$ and $v^*$ mark the terms of greatest
reduced-word length in the moment gradients of $uAu^*$ and $B$. We prove that the gradients are linearly independent
by projecting a linear relation onto that length,
and then descending in length.

\begin{lemma}
\label{lem:free-marker-regularity}
Let $n\ge1$, let $(P_0,\tau)$ be a tracial von Neumann algebra, and let $A,B\subseteq P_0$ be
finite-dimensional unital von Neumann subalgebras which are $n$-independent. Let
\[
X\subseteq(A\ominus\C1)_{\rm sa},
\qquad
Y\subseteq(B\ominus\C1)_{\rm sa}
\]
be real $L^2$-orthogonal bases of nonzero self-adjoint contractions whose complex spans are $A\ominus\C1$ and
$B\ominus\C1$. Let $(C,\tau_C)$ be a tracial von Neumann algebra, and set
\[
(N,\varphi):=(C,\tau_C)*(P_0,\tau),
\qquad
u\in U(C),
\qquad
c:=\tau_C(u),
\qquad
v:=u-c.
\]
Assume that $u$ is not a scalar multiple of $1$, so that $v\ne0$. Then $uAu^*$ and $B$ are
$n$-independent, and the family of real cyclic gradients through pair length $n$, formed from $uXu^*$ and $Y$, has
a positive-definite Gram matrix.
\end{lemma}

\begin{proof}
The $n$-independence follows from Lemma~\ref{lem:free-marker-preservation}.
If $A=\C1$ or $B=\C1$, one of $X,Y$ is empty, so there are no mixed cyclic gradients. Otherwise, by order-one independence
and orthogonality of the two bases, $X\sqcup Y$ is linearly independent
in $L^2(P_0)\ominus\C1$.

Let $\Pi_\ell$ project onto reduced-word length $\ell$ in $L^2(C*P_0)$.
For $x\in X$, the expansion
\[
uxu^*=|c|^2x+cxv^*+\overline c\,vx+vxv^*
\]
has the unique length-three term $vxv^*$.
Write $\widehat w$ for the value of a formal alternating word $w$
after replacing each $X$-letter $x$ by $uxu^*$ and leaving the
$Y$-letters fixed. If $w$ has pair length $r$, its component of
reduced-word length $4r$ is
\[
\Pi_{4r}(\widehat w)=T(w),\qquad
T:\ x\mapsto v\otimes x\otimes v^*,\quad y\mapsto y,
\]
and every other term has length below $4r$.
Words beginning with an $X$-letter give tensors starting in $C$;
those beginning with a $Y$-letter give tensors starting in $P_0$.
Among tensors with the same starting factor, the $C$-slots have a fixed pattern
of $v$ and $v^*$. Hence, since $v\ne0$ and $X\sqcup Y$ is linearly
independent, the tensors of distinct formal rotations are linearly
independent.

By the coefficient calculation in the proof of
Proposition~\ref{prop:noncommutative-free-gram}, the formal
real gradients of each fixed pair length $r$ are linearly independent.
Since $T$ is injective, their components of reduced-word length $4r$ are still linearly
independent.
The lemma follows by projecting a relation onto length $4r$ for the largest
pair length present, and then descending in $r$.
\end{proof}

In the next proposition, and wherever we speak of the chosen bases of $A$ and $B$, we mean the bases $X_\pi$ and $Y_\sigma$ defined before Lemma~\ref{lem:regular-limit-exactification}, where $\pi:P\to A$ and $\sigma:Q\to B$ are the given $^*$-isomorphisms.

\begin{proposition}
\label{prop:quantitative-phase-marker}
Fix $n\ge1$, nonscalar finite-dimensional unital $C^*$-algebras $P,Q$, and $\delta>0$.  There are constants
\[
\varepsilon_*>0,
\qquad
0<C_0<\infty,
\qquad
c_G>0,
\]
depending only on $n,P,Q,\delta$, with the following property. Let $(P_0,\tau)$ be a tracial
von Neumann algebra, and let $A,B\subseteq P_0$ be unital $^*$-subalgebras with fixed
$^*$-isomorphisms $P\cong A$ and $Q\cong B$. Assume that every minimal projection of either
$A$ or $B$ has trace at least $\delta$, and set
\[
\varepsilon:=\varepsilon_n(A,B)\le\varepsilon_*.
\]
Let $X,Y$ be their chosen bases.  In
\[
(N,\varphi):=(C,\tau_C)*(P_0,\tau),
\qquad
C=\C^3,
\qquad
\tau_C(e_0)=\tau_C(e_+)=\tau_C(e_-)=\frac13,
\]
where $e_0,e_+,e_-$ are the minimal projections of $C$, set
\[
s_\theta:=e_0+e^{i\theta}e_++e^{-i\theta}e_-,
\qquad 0<\theta\le1.
\]
Let $\Phi_\theta$ and $G_\theta$ be the real constraint map and cyclic-gradient Gram matrix
through pair length $n$, formed from $s_\theta Xs_\theta^*$ and $Y$.  Then
\[
\norm{\Phi_\theta(0)}\le C_0\varepsilon,
\qquad
G_\theta\ge c_G\theta^{2n(n+1)}I.
\]
Moreover,
\[
\norm{s_\theta-1}_\infty=2\sin\frac\theta2.
\]
\end{proposition}

\begin{proof}
We begin with a preliminary bound for the Gram matrix of $X\sqcup Y$. Then we bound the constraint vector $\Phi_\theta(0)$, and
finally we bound $G_\theta$ from below by projecting onto the components of greatest reduced-word length, as in
the proof of Lemma~\ref{lem:free-marker-regularity}. Write $p:=\abs X$, $q:=\abs Y$, and
\[
\kappa:=\min\{\kappa(P,\delta),\kappa(Q,\delta)\}>0,
\qquad
\kappa_0:=\frac\kappa2.
\]
In the $L^2$ Gram matrix $G_Z$ of $Z:=X\sqcup Y$, the two diagonal blocks are diagonal matrices bounded below by $\kappa I$.
The off-diagonal entries have absolute value at most $\varepsilon$. Hence
\[
G_Z\ge\kappa_0I
\]
provided $\varepsilon\le\kappa/(2\sqrt{pq})$.  We include this restriction in $\varepsilon_*$.

For a constraint of pair length $r\le n$, we repeat the proof of
Lemma~\ref{lem:free-marker-preservation}, using
Lemmas~\ref{lem:quantitative-truncated-comparison} and
\ref{lem:finite-free-substitution}. Since every increasing subsequence of
its $2r$ letters has at most $r$ blocks from each algebra, the constraint is bounded by $L_{2r}K_{2n}\varepsilon$,
uniformly in $\theta$. Hence, if $d$ is
the number of real constraints, we may take
\[
C_0:=\sqrt d\,K_{2n}\max_{1\le r\le n}L_{2r}.
\]

Set $c_\theta:=\tau_C(s_\theta)=(1+2\cos\theta)/3$ and $v_\theta:=s_\theta-c_\theta$. Then
\[
\norm{v_\theta}_2^2=\norm{v_\theta^*}_2^2=1-c_\theta^2
=\frac49(1-\cos\theta)(2+\cos\theta)\ge\frac{\theta^2}{2}
\qquad(0<\theta\le1).
\]

Let $d_r$ be the number of real gradients of pair length $r$, and let
$T_r:\R^{d_r}\to L^2(N)_{\rm sa}$ send the standard basis to these gradients.
Let $\Pi_{4r}$ denote the projection onto reduced-word length $4r$ in $L^2(N)$.
Index the distinct formal tensors of length $4r$ by $\mathcal D_r$,
distinguishing their starting factors, and let $J_r^{\C}a$ be the
coefficient vector of $\Pi_{4r}T_ra$ in these tensors.
By the coefficient calculation in the proof of
Proposition~\ref{prop:noncommutative-free-gram}, the columns of $J_r^{\C}$
are orthogonal and have squared norm at least $1/8$. Indeed, equal rotations have the
same sign, and distinct representative
orbits have disjoint supports. The coefficients are rational
multiples of $1$ and $i$ depending only on $n,P,Q$.
For $x\in X$, the component of $s_\theta xs_\theta^*$ of reduced-word length
three is $v_\theta\otimes x\otimes v_\theta^*$, and since $x$ is self-adjoint,
taking adjoints leaves this component unchanged. Hence, in every rotation and
adjoint
of a word of pair length $r$, the component of length $4r$ has $v_\theta$
immediately before and $v_\theta^*$ immediately after each $X$-letter. For each
fixed starting factor, the tensor Gram matrix is therefore
$\norm{v_\theta}_2^{4r}$ times a principal submatrix of $G_Z^{\otimes2r}$.
Consequently,
\[
\norm{\Pi_{4r}T_ra}_2
\ge\frac{(\norm{v_\theta}_2^2\kappa_0)^r}{2\sqrt2}\norm a_2
\ge b_r\norm a_2,
\qquad
b_r:=\frac{(\kappa_0\theta^2/2)^r}{2\sqrt2}\le1.
\]

Now combine the maps by setting
\[
T(a_1,\ldots,a_n):=\sum_{r=1}^nT_ra_r.
\]
Every gradient has operator norm, hence $L^2$-norm, at
most $n$, so with $B_0:=n\sqrt d\ge1$ we have $\norm{T_r}\le B_0$.  Also
$\Pi_{4r}T_s=0$ for $s<r$.  For $a=(a_1,\ldots,a_n)$, projecting onto length
$4r$ gives
\[
b_r\norm{a_r}_2
\le\norm{Ta}_2+B_0\sum_{s>r}\norm{a_s}_2.
\]
Descending from $r=n$ yields
$\norm{a_r}_2\le(1+nB_0)^{n-r}\bigl(\prod_{s=r}^nb_s^{-1}\bigr)\norm{Ta}_2$
for every $r$, and therefore
\[
\norm a_2\le\sum_{r=1}^n\norm{a_r}_2
\le n(1+nB_0)^{n-1}(2\sqrt2)^n\Bigl(\frac{\kappa_0}{2}\Bigr)^{-n(n+1)/2}
\theta^{-n(n+1)}\norm{Ta}_2.
\]
Thus $\norm{Ta}_2\ge c_T\theta^{n(n+1)}\norm a_2$, where $c_T>0$ is the
reciprocal of the constant in the last display. Since $G_\theta=T^*T$, taking
$c_G=c_T^2$ proves the Gram estimate.  The formula for
$\norm{s_\theta-1}_\infty$ is immediate from its three spectral values.
\end{proof}

\emph{Transfer to the factor.}
We now realize the constraint vector and Gram matrix of Proposition~\ref{prop:quantitative-phase-marker}
inside a given factor.
Let $A,B$ be
unital $^*$-subalgebras of a ${\rm II}_1$ factor $(M,\tau)$,
with fixed $^*$-isomorphisms from the finite-dimensional $C^*$-algebras $P,Q$, and let
$X,Y$ be the bases chosen above. The pair need not be independent.
Choose mutually
orthogonal projections $e_0,e_+,e_-$ with
\[
e_0+e_++e_-=1,\qquad
\tau(e_0)=\tau(e_+)=\tau(e_-)=\tfrac13,
\]
and set $s_\theta=e_0+e^{i\theta}e_++e^{-i\theta}e_-$,
$0<\theta\le1$.
Applying Theorem~\ref{thm:finite-order-free-positioning}\textup{(ii)} at
order $4n$ to $(s_\theta)$ and $X\sqcup Y$, we obtain a unitary $z$. Then
the unitary $w=zs_\theta z^*$ satisfies
\[
\norm{w-1}_\infty=2\sin(\theta/2)\le\theta.
\]
Its joint $^*$-moments with $X\sqcup Y$ through total degree $8n+1$
agree with those in $W^*(s_\theta)*W^*(A\cup B)$, equipped with the
induced tracial states. Note that $A$ and $B$ have the same
joint law in this free product as in $M$. Since conjugation does not change the trace induced on $P$, the chosen
basis for the first algebra after conjugation is exactly $wXw^*$.
A conjugated word of pair length $r$ has degree
$4r$ in $w,w^*,X,Y$, and every gradient is a linear combination
of the cyclic rotations and adjoints of such a word. Thus constraints and gradients
have degree at most $4n$, and Gram products have degree at most $8n$.
The constraint vector and Gram matrix therefore agree with those in the
free product, i.e., with those of Proposition~\ref{prop:quantitative-phase-marker} for $P_0=W^*(A\cup B)$.

\begin{proposition}
\label{prop:dense-regular-exact-locus}
Fix $n\ge1$ and $P,Q$ as above. Let $(M,\tau)$ be a ${\rm II}_1$ factor, and let
\[
\pi:P\to M,
\qquad
\sigma:Q\to M
\]
be unital $^*$-embeddings whose images $A:=\pi(P)$ and $B:=\sigma(Q)$ are $n$-independent.  For every
$\zeta>0$, there is a unitary $w\in M$ such that
\[
\norm{w-1}_\infty<\zeta,
\]
the pair $wAw^*,B$ is $n$-independent, and it is order-$n$ regular relative to the chosen bases
\[
w\pi\bigl(\mathcal X_P(\tau\circ\pi)\bigr)w^*,
\qquad
\sigma\bigl(\mathcal X_Q(\tau\circ\sigma)\bigr).
\]
\end{proposition}

\begin{proof}
If $A=\C1$ or $B=\C1$, take $w=1$. Otherwise choose $\delta>0$ below the
traces of all minimal projections of $A$ and $B$, and $0<\theta\le1$ with $2\sin(\theta/2)<\zeta$.
By Proposition~\ref{prop:quantitative-phase-marker} at $\varepsilon=0$, the
constraint vector is zero and the Gram matrix is positive definite in the free product.
By the transfer to the factor described above, there is a unitary $w\in M$ with
$\|w-1\|_\infty=2\sin(\theta/2)<\zeta$ such that the constraint vector and Gram matrix formed from
$wXw^*$ and $Y$ are the same as in the free product, where $X,Y$ are the chosen bases of $A,B$.
Since $wXw^*$ and $Y$ are the bases in the statement, the pair $wAw^*,B$ is $n$-independent
and order-$n$ regular.
\end{proof}

\begin{proof}[Proof of Theorem~\ref{thm:fixed-type-factor-no-loss}]
If $\varepsilon=0$, take $u=1$. Otherwise let
$X,Y$ be the chosen bases for $A,B$ and set
\[
\nu:=2n(n+1),\qquad \alpha:=\frac1{\nu+1},\qquad
\theta:=\varepsilon^\alpha.
\]
With this choice, the bound $\theta$ for the preliminary conjugation
and the bound $\varepsilon/\theta^\nu$ for the final correction have
the same order. For $\varepsilon_0\le\min\{1,\varepsilon_*\}$,
we apply the transfer to the factor, described after
Proposition~\ref{prop:quantitative-phase-marker}, to obtain a unitary $w\in M$
with $\norm{w-1}_\infty\le\theta$ such that the constraint map $\Phi$ and
Gram matrix $G$ for $wXw^*,Y$ satisfy
\[
\eta:=\norm{\Phi(0)}\le C_0\varepsilon,
\qquad G\ge\lambda I,\qquad \lambda:=c_G\theta^\nu.
\]
Set $C_{\rm nl}:=2n^2\sqrt d$, where $d$ is the number of real constraints.
We repeat the Gram-inverse calculation in the proof of
Theorem~\ref{thm:fd-noncommutative-exactify}, using
Proposition~\ref{prop:pair-linearization} and
Lemma~\ref{lem:finite-dimensional-correction}, and obtain $h=h^*$ with $\Phi(h)=0$
and
\[
\norm h_\infty\le\frac{3n\sqrt d}{4\lambda}\eta
\quad\text{provided}\quad
\eta\le\frac{\lambda}{2C_{\rm nl}}.
\]
Since $\eta/\lambda\le(C_0/c_G)\varepsilon^\alpha$, it suffices to take
\[
\varepsilon_0\le
\min\left\{\varepsilon_*,1,
\left(\frac{c_G}{2C_{\rm nl}C_0}\right)^{\nu+1}\right\}.
\]
The complex spans of the bases $wXw^*$ and $Y$ are $wAw^*\ominus\C1$ and
$B\ominus\C1$, respectively. Thus \eqref{eq:cyclic-zero-independence}
gives $n$-independence of $uAu^*,B$, where $u=e^{ih}w$. Finally,
\[
\norm{u-1}_\infty\le\norm h_\infty+\theta
\le\left(1+\frac{3n\sqrt d\,C_0}{4c_G}\right)\varepsilon^\alpha.
\]
\end{proof}

The exponent $1/(2n^2+2n+1)$ comes from the Gram estimate and the descent
over reduced-word lengths. We do not know whether this exponent is optimal. Note that $\varepsilon_0$ and $C$ can be computed from the proofs,
but the positioning and correcting unitaries are chosen nonconstructively.

\begin{corollary}
\label{cor:fixed-type-kadison-kastler}
Fix $n\ge1$, finite-dimensional unital $C^*$-algebras $P,Q$, and $\delta>0$, and let $\varepsilon_0,C$ be supplied by
Theorem~\ref{thm:fixed-type-factor-no-loss}.  If $A,B$ satisfy the hypotheses of that theorem and
$\varepsilon:=\varepsilon_n(A,B)\le\varepsilon_0$, then there is a copy $A'\subseteq M$ of $P$
and a trace-preserving $^*$-isomorphism $A\to A'$ such that $A',B$ are $n$-independent and
\[
d_{\rm KK}(A,A')\le2C\varepsilon^{1/(2n^2+2n+1)}.
\]
Here $d_{\rm KK}$ denotes the Kadison--Kastler distance~\cite{kadison-kastler}, i.e., the Hausdorff
distance in operator norm between the closed unit balls.
\end{corollary}

\begin{proof}
For the unitary $u$ of Theorem~\ref{thm:fixed-type-factor-no-loss},
take $A'=uAu^*$ and use
$d_{\rm KK}(A,uAu^*)\le2\|u-1\|_\infty$.
\end{proof}

\begin{theorem}
\label{thm:factor-faithful-no-loss}
Fix $n,a,b\ge1$ and $\delta>0$.  There are constants
\[
\varepsilon_0=\varepsilon_0(n,a,b,\delta)>0,
\qquad
C=C(n,a,b,\delta)<\infty
\]
with the following property, uniformly over all ${\rm II}_1$ factors $(M,\tau)$.  Let
\[
\mathbf p=(p_1,\ldots,p_a),
\qquad
\mathbf q=(q_1,\ldots,q_b)
\]
be partitions of $1$ in $M$ satisfying
\[
\tau(p_i)\ge\delta,
\qquad
\tau(q_j)\ge\delta
\]
for all indices.  If
\[
\operatorname{def}_n(\mathbf p,\mathbf q)\le\varepsilon_0,
\]
then there is a unitary $u\in M$ such that
\[
\norm{u-1}_\infty
\le C\operatorname{def}_n(\mathbf p,\mathbf q)^{1/(2n^2+2n+1)}
\]
and the unital $^*$-algebras generated by $u\mathbf p u^*$ and $\mathbf q$ are $n$-independent.  In particular, only the first
partition moves and every projection retains its trace.  Moreover,
\[
\norm{up_iu^*-p_i}_\infty\le2\norm{u-1}_\infty
\qquad(1\le i\le a),
\]
and
\[
\sum_{i=1}^a\norm{up_iu^*-p_i}_2
\le2a\norm{u-1}_\infty.
\]
\end{theorem}

\begin{proof}
We apply Theorem~\ref{thm:fixed-type-factor-no-loss} with $P=\C^a$ and $Q=\C^b$.  By
\eqref{eq:labelled-defect-comparison}, it suffices to divide the threshold
in that theorem by $(ab)^n$. Its estimate becomes
\[
\norm{u-1}_\infty
\le C_{P,Q,\delta}(ab)^{n/(2n^2+2n+1)}
\operatorname{def}_n(\mathbf p,\mathbf q)^{1/(2n^2+2n+1)},
\]
and we absorb the constant factor into $C$. The two distance bounds follow from
\[
up_iu^*-p_i=(u-1)p_iu^*+p_i(u^*-1)
\]
and summation.
\end{proof}

\subsection{Order-one stability}
\label{sec:order-one-no-loss}

In Theorem~\ref{thm:factor-faithful-no-loss}, we corrected small errors in the moments by a
small unitary conjugation of the first partition, assuming a positive lower
bound on the projection traces. At order one, we can remove this
assumption by allowing both partitions to move in $L^2$, while preserving
every projection trace. The equations to be corrected are
$\tau(p_iq_j)=\tau(p_i)\tau(q_j)$.
In the proof, we merge projections whose traces are
small, apply induction, and then split the merged
projections to restore the original traces. For this splitting with the required trace pairings,
we use the next lemma and its corollary. A noncommutative Lyapunov result covering the lemma
is given in~\cite[Theorem~3.1]{akemann-weaver-automatic-convexity}. Noncommutative Lyapunov
theorems go back to Akemann and Anderson~\cite{akemann-anderson-lyapunov}.

\begin{lemma}
\label{lem:finite-functional-purification}
Let $D$ be a diffuse finite von Neumann algebra and let $\phi_1,\ldots,\phi_d$ be normal real-linear functionals on
$D_{\rm sa}$.  For every $0\le x\le1$, there is a projection $e\in D$ such that
\[
\phi_\ell(e)=\phi_\ell(x)
\qquad(1\le\ell\le d).
\]
\end{lemma}

\begin{proof}
By the Krein--Milman theorem, the nonempty weak-$*$ compact convex set
\[
K=\{0\le y\le1:\phi_\ell(y)=\phi_\ell(x)\text{ for all }\ell\}
\]
has an extreme point $y$. If $y$ is not a projection, then
$r=1_{[\epsilon,1-\epsilon]}(y)$ is nonzero for some $\epsilon>0$. By diffuseness,
$(rDr)_{\rm sa}$ is infinite-dimensional. We may therefore choose a nonzero
$h$ in this space and in the common kernel of the finitely many $\phi_\ell$, and scale it so that
$\norm h_\infty\le\epsilon$. Then $y\pm h\in K$, contradicting extremality.
Thus $y$ is the required projection.
\end{proof}

\begin{corollary}
\label{cor:proportional-projection-splitting}
Let $(D,\tau)$ be a diffuse tracial von Neumann algebra, let $r\in D$ be a projection, let $x_1,\ldots,x_d\in D$, and let
$s_1,\ldots,s_k\ge0$ satisfy $\sum_i s_i=1$.  There is a partition $(r_i)_{i=1}^k$ of $r$, with zero projections allowed, such that
\[
\tau(r_i)=s_i\tau(r),
\qquad
\tau(r_ix_\ell)=s_i\tau(rx_\ell)
\qquad(i,\ell).
\]
If $Z(D)$ is diffuse and $r$ is central, the $r_i$ may be chosen central.
\end{corollary}

\begin{proof}
The case $r=0$ is immediate.  We apply Lemma~\ref{lem:finite-functional-purification} in $rDr$ to $s_1r$, with the functionals given by
the trace and by the real and imaginary parts of the trace pairings with the $x_\ell$, and obtain $r_1$.  If $s_1<1$, we apply it in the remaining
corner to $s_2(r-r_1)/(1-s_1)$, and continue.  After each step, the functional values on the remaining projection equal their original values times the total remaining
weight. We choose the zero projection for a zero weight, and stop if no weight remains. For the central assertion, we perform
the same construction in the diffuse abelian von Neumann algebra $rZ(D)$.
\end{proof}

\begin{theorem}[Trace-preserving partition stability]
\label{thm:finite-partition-stability}
For all integers $a,b\ge1$ and every $\eta>0$, there is
\[
\varepsilon=\varepsilon(a,b,\eta)>0
\]
with the following property. Let $(M,\tau)$ be a ${\rm II}_1$ factor, and let
\[
1=\sum_{i=1}^a p_i,
\qquad
1=\sum_{j=1}^b q_j
\]
be partitions of $1$, possibly containing zero projections. If
\[
\max_{i,j}\abs{\Delta_{ij}(\mathbf p,\mathbf q)}<\varepsilon,
\]
then there are partitions
$\mathbf p'=(p'_1,\ldots,p'_a)$ and
$\mathbf q'=(q'_1,\ldots,q'_b)$ of $1$ in $M$, again possibly containing zero projections, such that
\[
\tau(p'_i)=\tau(p_i),
\qquad
\tau(q'_j)=\tau(q_j),
\qquad
\tau(p'_iq'_j)=\tau(p'_i)\tau(q'_j)
\]
for all indices, and
\[
\sum_i\norm{p_i-p'_i}_2+
\sum_j\norm{q_j-q'_j}_2<\eta.
\]
The constant $\varepsilon$ does not depend on $M$ or on the traces $\tau(p_i)$ and $\tau(q_j)$.
\end{theorem}

To prove this result, we use the following consequence of
Theorem~\ref{thm:factor-faithful-no-loss} for projection traces
bounded away from zero. In Appendix~\ref{app:second-order-geometry}, we give another proof of this consequence,
using second derivatives at points where the derivative of the moment map is not surjective.

\begin{corollary}
\label{thm:faithful-relative-partition-stability}
For all integers $a,b\ge2$ and all $\delta,\rho>0$, there is
\[
\varepsilon_{\rm rel}
=\varepsilon_{\rm rel}(a,b,\delta,\rho)>0
\]
with the following property. Let $(M,\tau)$ be a ${\rm II}_1$ factor and let
$\mathbf p,\mathbf q$ be partitions of $1$ into $a$ and $b$ projections, respectively. If
\[
\min_i\tau(p_i)\ge\delta,
\qquad
\min_j\tau(q_j)\ge\delta,
\qquad
\max_{i,j}\abs{\Delta_{ij}(\mathbf p,\mathbf q)}
<\varepsilon_{\rm rel},
\]
then there is $h=h^*\in M$ with $\norm{h}_\infty<\rho$ such that, for $w=e^{ih}$,
\[
\tau(wp_iw^*q_j)=\tau(p_i)\tau(q_j)
\qquad(1\le i\le a,\ 1\le j\le b).
\]
\end{corollary}

\begin{proof}
Here $\operatorname{def}_1(\mathbf p,\mathbf q)=\max_{i,j}|\Delta_{ij}|$.
Applying Theorem~\ref{thm:factor-faithful-no-loss} with $n=1$, we obtain the
required identities by conjugating only the first partition, preserving every
projection trace. By the uniform bound in that theorem, the conjugating unitary $w$
becomes arbitrarily close to $1$
as the defect tends
to zero. Choosing $\varepsilon_{\rm rel}$ sufficiently small and taking the
principal logarithm, we may therefore write $w=e^{ih}$ with
$\norm h_\infty<\rho$.
\end{proof}

\begin{proof}[Proof of Theorem~\ref{thm:finite-partition-stability}]
We induct on $a+b$. The case $a=1$ or $b=1$ is immediate. If the
conclusion fails for $(a,b)$, choose ${\rm II}_1$ factors
$(M_n,\tau_n)$ and pairs $(\mathbf p_n,\mathbf q_n)$ of partitions of $1$ with
\begin{equation}
\max_{i,j}|\Delta_{ij}(\mathbf p_n,\mathbf q_n)|\longrightarrow0,
\label{eq:boundary-vanishing-covariance}
\end{equation}
but at distance at least some $\eta>0$ from every independent pair of
partitions with the same traces at corresponding indices. Note that independent pairs with these prescribed traces exist in a
diffuse abelian von Neumann subalgebra. Passing to a subsequence, we may assume that the traces $\tau_n(p_{i,n})$ and $\tau_n(q_{j,n})$
converge for all $i,j$.
If all limiting traces are positive,
then by Corollary~\ref{thm:faithful-relative-partition-stability}, there is
a unitary $e^{ih_n}$ with $\norm{h_n}_\infty\to0$. Hence the total movement
is at most
$2a\norm{h_n}_\infty$, a contradiction.

Otherwise, after exchanging sides and renumbering, we may assume
$\alpha_{a,n}\to0$, where $\alpha_{i,n}:=\tau_n(p_{i,n})$. We replace the last two $p$-projections by their sum
$P_n=p_{a-1,n}+p_{a,n}$. After merging, the covariance defects still
tend to zero.
By induction and a diagonal choice, we obtain independent partitions
$(\widetilde p_{1,n},\ldots,\widetilde p_{a-2,n},\widetilde P_n)$ and
$(\widetilde q_{j,n})_j$ with the same traces as the merged partitions and total error $o(1)$.
Set $\beta_{j,n}=\tau_n(q_{j,n})$ and
$\gamma_n=\tau_n(\widetilde P_n)=\alpha_{a-1,n}+\alpha_{a,n}$.
If $\gamma_n>0$, we apply Corollary~\ref{cor:proportional-projection-splitting}
to $\widetilde P_n$, with weight $\alpha_{a,n}/\gamma_n$ and
the elements $\widetilde q_{j,n}$, to obtain $e_n\le\widetilde P_n$
satisfying
\[
\tau_n(e_n)=\alpha_{a,n},\qquad
\tau_n(e_n\widetilde q_{j,n})=\alpha_{a,n}\beta_{j,n}.
\]
If $\gamma_n=0$, take $e_n=0$ and the same identities hold.
We then replace $\widetilde P_n$ by the two projections $\widetilde P_n-e_n$ and $e_n$, and keep all other
corrected projections. This restores the trace at every original index and
preserves independence of the two partitions. Since
\[
\norm{e_n-p_{a,n}}_2^2
=2\alpha_{a,n}-2\tau_n(e_np_{a,n})\le2\alpha_{a,n},
\]
the total error is at most the merged error plus
$2\sqrt{2\alpha_{a,n}}$, hence tends to zero. This contradicts the
choice of the original partitions and completes the induction.
\end{proof}

\section{Coordinate lifting and its limits}
\label{sec:coordinate-lifting}

We apply the perturbation theorems to subalgebras of tracial
ultraproducts. For a countably cofinal ultrafilter, we lift increasing
finite-dimensional approximations, correct the coordinate algebras at
each fixed stage, and let the stage increase slowly enough that the
corrected representatives still define the original elements of the ultraproduct.
In Subsection~\ref{subsec:lifting-finite-order}, we show that for separable abelian subalgebras the order of independence is preserved over $\mathrm{II}_1$ factors.
For arbitrary tracial coordinate algebras, we assume twice the order
in the ultraproduct (Theorem~\ref{thm:local-implies-hi}).

In Subsection~\ref{subsec:lifting-hyperfinite},
we treat pairs of hyperfinite subalgebras, not assumed jointly hyperfinite, over factor coordinates. With a countably complete
ultrafilter, all finite-dimensional approximations can instead be
lifted at once (Subsection~\ref{subsec:lifting-countably-complete}).
In Subsection~\ref{subsec:lifting-free-inputs}, for
free abelian subalgebras, we construct lifts with independence orders tending to infinity along
the ultrafilter, and we show that those orders cannot in general be prescribed or
required to diverge in the ordinary sense. In the same subsection, we obtain free
coordinate lifts over factors when the two subalgebras jointly generate a
hyperfinite von Neumann algebra. Finally, in
Subsection~\ref{subsec:lifting-local-stability}, we show that lifting at a fixed
order $n$ is equivalent to local stability
of the $n$-independence equations for pairs of partitions of $1$ of any
fixed sizes, allowing zero projections. There we also answer Problem~30,
apply the order-one case to ultrapowers of a non-Gamma factor, and give
obstructions over nonfactor coordinates.

\subsection{Finite-order lifts}
\label{subsec:lifting-finite-order}

\begin{theorem}[Lifting without decreasing the order]
\label{thm:factor-no-loss-hierarchy}
Fix $n\ge1$.  Let $\cU$ be a countably cofinal ultrafilter on an arbitrary set $I$, let
$(M_k,\tau_k)_{k\in I}$ be arbitrary ${\rm II}_1$ factors, and set $M=\prod_{\cU}M_k$.  If
$A,B\subseteq M$ are separable abelian $n$-independent von Neumann subalgebras, then there are finite-dimensional abelian
unital $^*$-subalgebras $C_k,D_k\subseteq M_k$ such that
\[
A\subseteq\ultraprod C_k,
\qquad
B\subseteq\ultraprod D_k,
\]
and $C_k,D_k$ are $n$-independent for every $k\in I$.
\end{theorem}

We give the proof after some auxiliary results. To prove the theorem, we first lift the finite-dimensional approximations so that each coordinate partition refines the preceding one, and then correct each stage separately. The corrected partitions need not refine one another. Since the
corrections tend to zero in $L^2$, a diagonal choice still represents
every element of each fixed approximating algebra. At indices where
no correction is chosen, we take $C_k$ and $D_k$ to be scalar.

\begin{lemma}
\label{lem:nested-lifts-general}
Let $M=\prod_{\cU}(M_k,\tau_k)$ be a tracial ultraproduct of von Neumann algebras,
and let $p=(p_k)_{\cU}$, where every $p_k$ is a projection. Let $r\ge1$. If $e_1,\dots,e_r\le p$
are pairwise orthogonal projections with $\sum_{j=1}^r e_j=p$, then there are projections $e_{j,k}\le p_k$, possibly zero, such that, for
every $k$, $\sum_{j=1}^r e_{j,k}=p_k$, and $(e_{j,k})_{\cU}=e_j$ for $1\le j\le r$.

Moreover, let $A^{(1)}\subseteq A^{(2)}\subseteq \cdots \subseteq M$ be an increasing sequence of finite-dimensional abelian
von Neumann subalgebras. Then there exist finite-dimensional abelian von Neumann subalgebras $A_k^{(m)}\subseteq M_k$ such that
$A^{(m)}\subseteq \ultraprod A_k^{(m)}$ and $A_k^{(m)}\subseteq A_k^{(m+1)}$ for all $m$ and all $k\in I$. Write
$A^{(m)}=\bigoplus_{i=1}^{r_m}\C p_i^{(m)}$ with minimal projections $p_i^{(m)}$, and after reindexing at each stage write
$p_i^{(m)}=\sum_{j\in J_{m,i}}p_j^{(m+1)}$. The lifts may be chosen from coordinate partitions
$1=\sum_{i=1}^{r_m}p_{i,k}^{(m)}$, possibly containing zero projections, such that $(p_{i,k}^{(m)})_{\cU}=p_i^{(m)}$,
$A_k^{(m)}=\mathrm{span}\{p_{i,k}^{(m)}:p_{i,k}^{(m)}\ne 0\}$, and
$p_{i,k}^{(m)}=\sum_{j\in J_{m,i}}p_{j,k}^{(m+1)}$ for all $m,i$ and $k\in I$.
\end{lemma}

\begin{proof}
Choose projection representatives $f_{j,k}$ of $e_j$. Starting with
$r_{1,k}=p_k$, define for $j<r$
\[
e_{j,k}=1_{[1/2,1]}(r_{j,k}f_{j,k}r_{j,k}),\qquad
r_{j+1,k}=r_{j,k}-e_{j,k},
\]
and take $e_{r,k}=r_{r,k}$.
Inductively the compression represents $e_j$, and so does its spectral
cut. The resulting projections partition $p_k$ at every coordinate.

For the second assertion, we apply this construction to $A^{(1)}$ and then
inside each already lifted minimal projection at every refinement. This
preserves the prescribed indices and refinement identities, and taking spans
gives the required increasing coordinate algebras.
\end{proof}

\begin{fact}
\label{lem:cofinal-diagonal-choice}
Let $\cU$ be a countably cofinal ultrafilter on $I$, and let
$U_m\in\cU$ for $m\ge1$. There is a function
$\mu:I\to\N\cup\{0\}$ such that
\[
\{k:\mu(k)\ge m\}\in\cU\quad(m\ge1),
\qquad k\in U_{\mu(k)}\quad\text{whenever }\mu(k)\ge1.
\]
In particular, $\mu(k)\xrightarrow{k\to\cU}\infty$.
\end{fact}

Indeed, choose decreasing $S_m\in\cU$ with empty intersection and
$S_m\subseteq\bigcap_{j=1}^mU_j$, and take
$\mu(k)=\max(\{m\ge1:k\in S_m\}\cup\{0\})$. Houdayer and Ioana made this diagonal choice in Case~1 of the proof of~\cite[Theorem~5.1]{houdayer-ioana}, following the proof of~\cite[Lemma~2.2]{boutonnet-chifan-ioana}.

\begin{lemma}
\label{lem:cofinal-corrected-stage-assembly}
Let $\cU$ be countably cofinal, let $(M_k,\tau_k)_{k\in I}$ be tracial von Neumann algebras, set
$M:=\ultraprod M_k$, and let $A,B\subseteq M$ be unital von Neumann subalgebras. Let $(A^{(m)})_m$ and $(B^{(m)})_m$
be increasing finite-dimensional unital $^*$-subalgebras with weakly dense unions in $A,B$.
Suppose that, for every $k\in I$, $(A_k^{(m)})_m$ and $(B_k^{(m)})_m$ are increasing finite-dimensional unital von Neumann
subalgebras of $M_k$ such that $A^{(m)}\subseteq\ultraprod A_k^{(m)}$ and
$B^{(m)}\subseteq\ultraprod B_k^{(m)}$ for every $m$.

Let $U_m\in\cU$, let $\eta_m\ge0$ satisfy $\eta_m\to0$, and let $s_m\ge1$ be integers. Assume that there are unitaries
$v_k^{(m)}\in M_k$ such that, for every $k\in U_m$, one has $\norm{v_k^{(m)}-1}_2\le\eta_m$ and
$v_k^{(m)}A_k^{(m)}v_k^{(m)*},B_k^{(m)}$ are $s_m$-independent.

Then there are a function $\mu:I\to\N\cup\{0\}$ and finite-dimensional unital von Neumann subalgebras $C_k,D_k$ of $M_k$
such that $\mu(k)\to_{\cU}\infty$, $A\subseteq\ultraprod C_k$, and
$B\subseteq\ultraprod D_k$. Whenever $\mu(k)\ge1$, the pair $C_k,D_k$ is $s_{\mu(k)}$-independent. When $\mu(k)=0$,
both algebras are scalar.

In particular, if $s_m=s$ for every $m$, then $C_k,D_k$ are
$s$-independent for every $k$. If $s_m\to\infty$, define
$r(k):=s_{\mu(k)}$ when $\mu(k)\ge1$ and $r(k):=1$ when $\mu(k)=0$.
Then $r(k)\to_{\cU}\infty$ and $C_k,D_k$ are $r(k)$-independent for
every $k$. If all the algebras $A_k^{(m)},B_k^{(m)}$ are abelian, then
so are $C_k,D_k$.
\end{lemma}

\begin{proof}
Choose $\mu$ by Fact~\ref{lem:cofinal-diagonal-choice}. For
$\mu(k)\ge1$, set
\[
C_k=v_k^{(\mu(k))}A_k^{(\mu(k))}v_k^{(\mu(k))*},
\qquad D_k=B_k^{(\mu(k))},
\]
and take $C_k=D_k=\C1$ when $\mu(k)=0$. These algebras are finite
dimensional and have the required independence at every coordinate.

Now let $a=(a_k)_{\cU}\in A^{(m)}$ with bounded $a_k\in A_k^{(m)}$. We
conjugate $a_k$ by $v_k^{(\mu(k))}$ on $\{\mu\ge m\}$ and use zero
elsewhere. The resulting $c_k\in C_k$ satisfy
\[
\|c_k-a_k\|_2\le2\sup_j\|a_j\|_\infty\,\eta_{\mu(k)}
\longrightarrow_{\cU}0
\quad\text{on }\{\mu\ge m\}.
\]
Hence $A^{(m)}\subseteq\ultraprod C_k$. Without conjugation, the same argument gives
$B^{(m)}\subseteq\ultraprod D_k$. Since both coordinate ultraproducts are
weakly closed, they contain $A$ and $B$, respectively. Finally, the
assertions about the orders follow from $\mu\to_{\cU}\infty$.
\end{proof}

\begin{proof}[Proof of Theorem~\ref{thm:factor-no-loss-hierarchy}]
In $A$ and $B$, choose increasing sequences $(A^{(m)})_m$ and $(B^{(m)})_m$ of finite-dimensional abelian
unital $^*$-subalgebras with weakly dense unions. We lift their minimal
projections by Lemma~\ref{lem:nested-lifts-general} and let $A_k^{(m)}$ and
$B_k^{(m)}$ be the coordinate algebras spanned by the resulting partitions. At each fixed
stage $m$, the traces of the lifted minimal projections are bounded
below on a $\cU$-large set, and the defect $\operatorname{def}_n$ of the resulting pairs of partitions tends
to zero.
Hence, by Theorem~\ref{thm:factor-faithful-no-loss}, on a set $U_m\in\cU$
there is a unitary $v_k^{(m)}$ with
$\|v_k^{(m)}-1\|_\infty<2^{-m}$ such that
$v_k^{(m)}A_k^{(m)}v_k^{(m)*}$ and $B_k^{(m)}$ are $n$-independent.
Extending it by $1$ off $U_m$ and applying
Lemma~\ref{lem:cofinal-corrected-stage-assembly} with $\eta_m=2^{-m}$ and
$s_m=n$, we obtain the theorem.
\end{proof}

Unlike the proof of Theorem~\ref{thm:factor-no-loss-hierarchy}, the proofs of Proposition~\ref{prop:higher-order-local-global}
and Corollary~\ref{cor:fixed-ambient-binary-problem-30} require the following lemma, in which both partitions are corrected.

\begin{lemma}
\label{lem:cofinal-two-sided-labelled-assembly}
Fix $n\ge1$. Let $\cU$ be a countably cofinal ultrafilter on $I$, let $(M_k,\tau_k)_{k\in I}$ be tracial von Neumann algebras, set
$M=\prod_{\cU}(M_k,\tau_k)$, and let $A,B\subseteq M$ be separable abelian von Neumann subalgebras. Let
\[
A^{(m)}=\bigoplus_{i=1}^{a_m}\C p_i^{(m)},
\qquad
B^{(m)}=\bigoplus_{j=1}^{b_m}\C q_j^{(m)}
\]
be increasing finite-dimensional unital $^*$-subalgebras with weakly dense unions in $A,B$. Suppose their minimal projections are represented by
partitions $p_{i,k}^{(m)},q_{j,k}^{(m)}$, possibly containing zero projections, that preserve the indices and all refinement identities, as in Lemma~\ref{lem:nested-lifts-general}. For every $m$, assume
there are $U_m\in\cU$ and, for $k\in U_m$, partitions of $1$ with the same indices, again possibly containing zero projections,
\[
(\widetilde p_{i,k}^{(m)})_{i=1}^{a_m},
\qquad
(\widetilde q_{j,k}^{(m)})_{j=1}^{b_m}
\]
whose generated unital $^*$-algebras are $n$-independent and for which
\[
\eta_m:=
\sup_{k\in U_m}\left(
\sum_i\norm{\widetilde p_{i,k}^{(m)}-p_{i,k}^{(m)}}_2+
\sum_j\norm{\widetilde q_{j,k}^{(m)}-q_{j,k}^{(m)}}_2
\right)\longrightarrow0.
\]
Then there are finite-dimensional abelian unital $^*$-subalgebras $C_k,D_k\subseteq M_k$,
$n$-independent for every $k$, such that
\[
A\subseteq\ultraprod C_k,
\qquad
B\subseteq\ultraprod D_k.
\]
\end{lemma}

\begin{proof}
Choose $\mu$ by Fact~\ref{lem:cofinal-diagonal-choice}. When
$\mu(k)\ge1$, let $C_k,D_k$ be generated by the corrected partitions
at stage $\mu(k)$. When $\mu(k)=0$, take $C_k=D_k=\C1$.
Thus $C_k,D_k$ are $n$-independent at every index.

For a fixed $p_i^{(m)}$, let $J_{m,i}^r$ be the set of its refinement indices at
stage $r\ge m$. On $\{\mu\ge m\}$ the projection
$c_k=\sum_{\ell\in J_{m,i}^{\mu(k)}}\widetilde p_{\ell,k}^{(\mu(k))}$
belongs to $C_k$ and satisfies
\[
\|c_k-p_{i,k}^{(m)}\|_2
\le\sum_\ell\|\widetilde p_{\ell,k}^{(\mu(k))}
                 -p_{\ell,k}^{(\mu(k))}\|_2
\le\eta_{\mu(k)}\longrightarrow_{\cU}0.
\]
Set $c_k=0$ elsewhere. This sequence represents $p_i^{(m)}$, and the same
construction works for every $q_j^{(m)}$.
Since the coordinate ultraproducts are weakly closed, weak density gives
$A\subseteq\ultraprod C_k$ and $B\subseteq\ultraprod D_k$, as desired.
\end{proof}

For arbitrary tracial coordinate algebras, we can use Theorem~\ref{thm:pair-local-exactify}
in place of Theorem~\ref{thm:factor-faithful-no-loss} in the proof of
Theorem~\ref{thm:factor-no-loss-hierarchy} if we assume $2n$-independence
in the ultraproduct. This gives part~\textup{(i)} of the next theorem. For countably complete ultrafilters there is no loss of order (part~\textup{(ii)}; see also
Subsection~\ref{subsec:lifting-countably-complete}).

\begin{theorem}\label{thm:local-implies-hi}
Let $n\ge1$, let $I$ be a set, let $\cU$ be an ultrafilter on $I$, let
$(M_k,\tau_k)_{k\in I}$ be tracial von Neumann algebras, and let
$M:=\ultraprod M_k$, with trace $\tau=(\tau_k)_{\cU}$.
Let $A,B\subseteq M$ be separable abelian von Neumann subalgebras which are $2n$-independent.
\begin{enumerate}[label=\textup{(\roman*)}]
\item If $\cU$ is countably cofinal, there are finite-dimensional abelian von Neumann subalgebras
$C_k,D_k\subseteq M_k$ such that
\[
A\subseteq\ultraprod C_k,\qquad B\subseteq\ultraprod D_k,
\]
and $C_k,D_k$ are $n$-independent for every $k\in I$.
\item If $\cU$ is countably complete, there are abelian von Neumann subalgebras
$C_k,D_k\subseteq M_k$, not necessarily finite-dimensional, such that
\[
A\subseteq\ultraprod C_k,\qquad B\subseteq\ultraprod D_k,
\]
and $C_k,D_k$ are $2n$-independent for every $k\in I$.
\end{enumerate}
\end{theorem}

For $n=1$, part~\textup{(i)} is due to Houdayer and
Ioana~\cite[Theorem~5.1 and its proof]{houdayer-ioana}. Part~\textup{(i)}
was obtained independently by Patchell and
Shiner~\cite[Theorem~C]{patchell-shiner}.

\begin{proof}
In the countably cofinal case, we take finite-dimensional abelian stages in
$A$ and $B$ with weakly dense unions and lift them compatibly by
Lemma~\ref{lem:nested-lifts-general}. Then we
apply Theorem~\ref{thm:pair-local-exactify} at each stage $m$, choosing
a set $U_m\in\cU$ on which the correcting unitaries $v_k^{(m)}$ satisfy $\|v_k^{(m)}-1\|_2<2^{-m}$.
By Lemma~\ref{lem:cofinal-corrected-stage-assembly}, with $\eta_m=2^{-m}$
and $s_m=n$, we obtain \textup{(i)}.
For \textup{(ii)}, we apply
Theorem~\ref{thm:hyperfinite-lifting-complete}\textup{(i)} at order $2n$, requiring $C_k$
and $D_k$ to be abelian.
\end{proof}

\subsection{Hyperfinite von Neumann subalgebras}
\label{subsec:lifting-hyperfinite}

For hyperfinite von Neumann subalgebras of an ultraproduct of factors, we lift
the two algebras separately and then conjugate one of the
finite-dimensional images to obtain the required independence.
The matrix-unit lifts must respect the inclusions between successive
approximating algebras.

We need the following version, for varying coordinate factors, of the partial-isometry lemma stated by Hadwin and Li~\cite[Lemma~4.9]{hadwin-li}, who attribute it to Connes~\cite{connes-outer-conjugacy}.

\begin{lemma}
\label{lem:equal-trace-partial-isometry-representatives}
Let $\omega$ be an ultrafilter on $I$, let $(M_k,\tau_k)_{k\in I}$ be ${\rm II}_1$ factors, and let
$E_k,F_k\in M_k$ be projections satisfying
\[
\tau_k(E_k)=\tau_k(F_k)
\qquad(k\in I).
\]
Set $E=(E_k)_\omega$ and $F=(F_k)_\omega$.  If $V\in\prod_\omega M_k$ is a partial isometry satisfying
\[
V^*V=E,
\qquad
VV^*=F,
\]
then there are partial isometries $V_k\in M_k$ such that
\[
V_k^*V_k=E_k,
\qquad
V_kV_k^*=F_k
\qquad(k\in I),
\]
and $V=(V_k)_\omega$.
\end{lemma}

\begin{proof}
For a bounded representative $V=(X_k)_\omega$, set $A_k=F_kX_kE_k$ and
write $A_k=W_k|A_k|$, with right support $e_k\le E_k$ and left support
$f_k\le F_k$. Since $A_k$ represents $V$,
\[
\norm{A_k^*A_k-E_k}_2+\norm{A_kA_k^*-F_k}_2\longrightarrow_\omega0.
\]
Functional calculus gives $\norm{W_k-A_k}_2\to_\omega0$ and
$\tau_k(E_k-e_k)\to_\omega0$. Since the projections $E_k-e_k$ and $F_k-f_k$ have equal
trace, we may choose a partial isometry $U_k$ between them in the factor $M_k$.
Then $V_k=W_k+U_k$ has right support $E_k$ and left support
$F_k$ at every coordinate, and
$\norm{U_k}_2^2=\tau_k(E_k-e_k)\to_\omega0$. Hence $(V_k)_\omega=V$.
\end{proof}

The partial-isometry lemma stated in~\cite[Lemma~4.9]{hadwin-li}
concerns a constant separable factor. For prescribed right and left
supports, their traces must agree at every coordinate. However, this
hypothesis is absent from the printed statement. In Connes's original statement for
$\omega$-centralizing sequences, the supports are assumed equivalent at every
coordinate~\cite[Proposition~1.1.3(b)]{connes-outer-conjugacy}.

\emph{Compatible lifts of matrix units.}
The compatible extension
in~\cite[Lemma~4.10]{hadwin-li} uses this equal-trace case. To extend a
prescribed trace-preserving lift of the smaller algebra, we first lift the
required partition of one minimal diagonal projection in each matrix block
by Lemma~\ref{lem:nested-lifts-general}. Redistributing
subprojections from pieces with excess trace to pieces with deficient
trace makes all their traces exactly the prescribed ones. The
redistributed subprojections have traces tending to zero, so the pieces
still represent the same projections in the ultraproduct;
compare~\cite[proof of Theorem~4.4]{hadwin-li}. We then transport these pieces to
the other diagonal corners using the prescribed matrix units. This
gives a refinement compatible with the given lift. Finally, in each matrix
block of the larger algebra, we use
Lemma~\ref{lem:equal-trace-partial-isometry-representatives} to lift the
remaining partial isometries along a spanning tree connecting the
equal-trace diagonal projections, and define the remaining matrix units
algebraically. This justifies the compatible lifting argument
in~\cite[Lemma~4.10]{hadwin-li} for varying coordinate factors.

\begin{theorem}[Lifting hyperfinite $n$-independent subalgebras]
\label{thm:hyperfinite-lifting-cofinal}
Let $n\ge1$, let $\cU$ be a countably cofinal ultrafilter on $I$, let
$(M_k,\tau_k)_{k\in I}$ be ${\rm II}_1$ factors, and set $M:=\ultraprod M_k$.
Let $A,B\subseteq M$ be separable hyperfinite unital von Neumann subalgebras which are $n$-independent. Then there are
finite-dimensional unital von Neumann subalgebras $C_k,D_k\subseteq M_k$ such that
\[
A\subseteq\ultraprod C_k,
\qquad
B\subseteq\ultraprod D_k,
\]
and $C_k,D_k$ are $n$-independent for every $k\in I$.
\end{theorem}

\begin{proof}
Choose increasing sequences of finite-dimensional unital $^*$-subalgebras
$A^{(m)}\subseteq A$ and $B^{(m)}\subseteq B$ with weakly dense unions.
By Lemma~\ref{lem:equal-trace-partial-isometry-representatives}, the compatible extension above, and the
lifting argument in~\cite[Lemma~4.10 and Theorem~4.11]{hadwin-li}, there are trace-preserving
unital $^*$-homomorphisms $\pi_k:A\to M_k$, $\sigma_k:B\to M_k$ representing the
inclusions. Here we apply~\cite[Theorem~4.11]{hadwin-li} separately to the countably generated
hyperfinite von Neumann algebras $A$ and $B$, with factor coordinates and a
countably cofinal ultrafilter (called ``nontrivial'' there).
Set $A_k^{(m)}=\pi_k(A^{(m)})$ and $B_k^{(m)}=\sigma_k(B^{(m)})$.

For each fixed $m$, compactness of the unit balls of
$A^{(m)}\ominus\C1$ and $B^{(m)}\ominus\C1$ gives
$\varepsilon_n(A_k^{(m)},B_k^{(m)})\to_{\cU}0$.
The marginal traces are fixed and faithful, so by
Theorem~\ref{thm:fixed-type-factor-no-loss} there are correcting unitaries
at operator-norm distance $<2^{-m}$ from $1$ on a set $U_m\in\cU$.
Extending them by $1$ off $U_m$ and applying
Lemma~\ref{lem:cofinal-corrected-stage-assembly} with $\eta_m=2^{-m}$ and
$s_m=n$, we obtain the theorem.
\end{proof}

The factor assumption cannot be dropped
(Proposition~\ref{prop:four-point-orthogonality-obstruction}).

\subsection{Countably complete ultrafilters}
\label{subsec:lifting-countably-complete}

For a countably complete ultrafilter $\cU$ on $I$, no function
$I\to\N$ tends to infinity along $\cU$, so the diagonal choice of stages
is unavailable. It is also unnecessary, since countable completeness makes
all the scalar equalities below exact on one set in $\cU$, as in the
proof of~\cite[Lemma~2.3(2)]{boutonnet-chifan-ioana}. Houdayer and Ioana
applied this observation to abelian subalgebras at order one in Case~2 of
the proof of~\cite[Theorem~5.1]{houdayer-ioana}. We use it below for
hyperfinite subalgebras at every order of independence and for freeness.

\begin{fact}
\label{lem:countably-complete-exact-zero}
Let $\cU$ be countably complete on $I$, and let
$(z_k^{(j)})_{k\in I}$, $j\ge1$, be scalar families with
$z_k^{(j)}\xrightarrow{k\to\cU}0$ for every fixed $j$. Then
\[
S:=\bigcap_{j,q\ge1}\{k\in I:|z_k^{(j)}|<1/q\}\in\cU,
\qquad z_k^{(j)}=0\quad(k\in S,\ j\ge1).
\]
\end{fact}

We can therefore lift all finite-dimensional approximations at once
and take the von Neumann algebra they generate at each coordinate.
For freely independent von Neumann algebras, every mixed-moment identity is then preserved.

\begin{lemma}
\label{lem:complete-hyperfinite-lifts}
Let $\cU$ be countably complete, let $(M_k,\tau_k)_{k\in I}$ be tracial von Neumann algebras, and let
$A\subseteq\ultraprod M_k$ be a separable hyperfinite unital von Neumann subalgebra. If
\[
A^{(1)}\subseteq A^{(2)}\subseteq\cdots\subseteq A
\]
is an increasing sequence of finite-dimensional unital $^*$-subalgebras with weakly dense union, then there are $S\in\cU$ and,
for every $k\in S$, compatible trace-preserving unital $^*$-homomorphisms
\[
\pi_k^{(m)}:A^{(m)}\to M_k
\qquad(m\ge1)
\]
such that $a=(\pi_k^{(m)}(a))_{\cU}$ for every $m$ and every $a\in A^{(m)}$, after the coordinate families are defined
to be $0$ off $S$.
\end{lemma}

\begin{proof}
Choose complete systems of matrix units in every $A^{(m)}$, expressing each stage-$m$ matrix unit as a fixed
finite linear combination of stage-$(m+1)$ matrix units. Then choose bounded representatives of these countably many matrix units.
For these representatives, consider the $L^2$ errors in the matrix-unit identities
$(e_{ij}^{(r)})^*=e_{ji}^{(r)}$, $e_{ij}^{(r)}e_{kl}^{(s)}=\delta_{rs}\delta_{jk}e_{il}^{(r)}$, and
$\sum_{r,i}e_{ii}^{(r)}=1$ at each stage, and in the refinement identities. For each diagonal matrix
unit, also consider the absolute differences between its trace in $A$ and
the traces of its coordinate representatives. These form countably many
nonnegative scalar families, each tending to zero along $\cU$ because the
identities hold in the ultraproduct. By
Fact~\ref{lem:countably-complete-exact-zero}, there is one $S\in\cU$ on which all
these errors vanish.

For $k\in S$, the chosen representatives thus form compatible unital systems of matrix units with the original diagonal traces.
Sending each original matrix unit to its representative defines
compatible trace-preserving unital $^*$-homomorphisms $\pi_k^{(m)}$. Since each original diagonal projection has positive trace,
these homomorphisms are injective. By the choice of the matrix-unit representatives and linearity,
$\pi_k^{(m)}(a)$ represents $a$.
\end{proof}

\begin{theorem}
\label{thm:hyperfinite-lifting-complete}
\label{thm:hyperfinite-free-lifts-complete}
\label{prop:countably-complete-no-loss}
Let $\cU$ be countably complete, let $(M_k,\tau_k)_{k\in I}$ be arbitrary
tracial von Neumann algebras, and let $A,B\subseteq\ultraprod M_k$ be
separable hyperfinite unital von Neumann subalgebras.
\begin{enumerate}[label=\textup{(\roman*)}]
\item If $A,B$ are $n$-independent for some fixed $n\ge1$, there are
hyperfinite unital von Neumann subalgebras $C_k,D_k\subseteq M_k$, $n$-independent for every
$k\in I$, such that
\[
A\subseteq\ultraprod C_k,\qquad B\subseteq\ultraprod D_k.
\]
\item If $A,B$ are freely independent, the algebras $C_k,D_k$ can be
chosen freely independent for every $k\in I$.
\end{enumerate}
On either side, the coordinate algebras may be required to be abelian when
the corresponding original algebra is abelian, and finite dimensional when that algebra
is finite dimensional.
\end{theorem}

\begin{proof}
For each of $A$ and $B$, choose an increasing finite-dimensional chain with
weakly dense union, abelian when the algebra is abelian and constant when it
is finite dimensional.
By Lemma~\ref{lem:complete-hyperfinite-lifts}, there are compatible
trace-preserving coordinate $^*$-homomorphisms on one $S_0\in\cU$.
Choose centered self-adjoint bases at every stage and enumerate the
alternating moments of their coordinate images through pair length $n$ in
\textup{(i)}, or all
finite pair lengths in \textup{(ii)}.
By Fact~\ref{lem:countably-complete-exact-zero}, this countable list
vanishes on a single $S\in\cU$ with $S\subseteq S_0$.
Hence the coordinate stage algebras have the required independence.

For each $k\in S$, let $C_k,D_k$ be the von Neumann algebras
generated by the respective increasing unions, and take scalars off $S$.
The inclusions follow since the coordinate ultraproducts are weakly closed;
independence passes to the limit by bounded $L^2$ approximation from the
stages. Each of $C_k,D_k$ is hyperfinite because its finite-dimensional
stages form an increasing sequence. If the chain is abelian, so is the
generated algebra. Finally, if it is constant and finite dimensional, the
generated algebra is finite dimensional as well.
\end{proof}

The assertion of Theorem~\ref{thm:hyperfinite-lifting-complete} about finite-dimensional algebras has a trivial converse: finite-dimensional
coordinate algebras on one side force the corresponding original
algebra to be finite dimensional. Indeed, a tracial ultraproduct of
finite-dimensional $C^*$-algebras over a countably complete ultrafilter
is finite dimensional.

\subsection{Freely independent von Neumann algebras and growing orders}
\label{subsec:lifting-free-inputs}

If the two algebras are free, we may correct their lifts at every
finite order. For a countably cofinal ultrafilter, a diagonal argument
then gives one pair of lifts whose independence orders tend to infinity along the ultrafilter. In Corollary~\ref{cor:free-unitary-lifting}, we apply this, together with Theorem~\ref{thm:local-implies-hi}, to free unitaries and Problem~31.

\begin{theorem}\label{thm:free-coordinate-lifts}
Let $I$ be a set, let $\cU$ be an ultrafilter on $I$, let $(M_k,\tau_k)_{k\in I}$ be tracial von Neumann algebras, and let
$M:=\ultraprod M_k$, with trace $\tau=(\tau_k)_{\cU}$. Let $A,B\subseteq M$ be freely independent separable abelian von
Neumann subalgebras.
\begin{enumerate}[label=\textup{(\roman*)}]
\item If $\cU$ is countably cofinal, there are finite-dimensional abelian von Neumann subalgebras
$C_k,D_k\subseteq M_k$ and a function $r:I\to\N$ such that
\[
A\subseteq\ultraprod C_k,\qquad B\subseteq\ultraprod D_k,
\qquad r(k)\longrightarrow_{\cU}\infty,
\]
and $C_k,D_k$ are $r(k)$-independent for every $k\in I$.
\item If $\cU$ is countably complete, there are abelian von Neumann subalgebras
$C_k,D_k\subseteq M_k$, not necessarily finite-dimensional, such that
\[
A\subseteq\ultraprod C_k,\qquad B\subseteq\ultraprod D_k,
\]
and $C_k,D_k$ are freely independent for every $k\in I$.
\end{enumerate}
\end{theorem}

\begin{proof}
Suppose first that the ultrafilter is countably cofinal. We use the increasing finite-dimensional chains
and their compatible projection lifts from the proof of Theorem~\ref{thm:local-implies-hi}. Since the original pair is free, its
stage-$m$ algebras are $2m$-independent. Then we apply Theorem~\ref{thm:pair-local-exactify} at target order $m$ and restrict to a set in $\cU$ on which the correcting unitary is within $2^{-m}$
of $1$ in $L^2$. We extend the unitary by $1$ off that set and apply Lemma~\ref{lem:cofinal-corrected-stage-assembly} with
$\eta_m=2^{-m}$ and $s_m=m$. This proves \textup{(i)}.

For \textup{(ii)}, we apply
Theorem~\ref{thm:hyperfinite-lifting-complete}\textup{(ii)}, requiring $C_k$ and $D_k$
to be abelian.
\end{proof}

\begin{theorem}
\label{thm:hyperfinite-free-lifts-cofinal}
Let $\cU$ be countably cofinal, let $(M_k,\tau_k)_{k\in I}$ be ${\rm II}_1$ factors, and set $M:=\ultraprod M_k$.
Let $A,B\subseteq M$ be freely independent separable hyperfinite unital von Neumann subalgebras. Then there are
finite-dimensional unital von Neumann subalgebras $C_k,D_k\subseteq M_k$ and a function $r:I\to\N$ such that
\[
A\subseteq\ultraprod C_k,
\qquad
B\subseteq\ultraprod D_k,
\qquad
r(k)\longrightarrow_{\cU}\infty,
\]
and $C_k,D_k$ are $r(k)$-independent for every $k\in I$.
\end{theorem}

\begin{proof}
We repeat the proof of Theorem~\ref{thm:hyperfinite-lifting-cofinal},
requiring $m$-independence at stage $m$. By freeness, we have
$\varepsilon_m(A_k^{(m)},B_k^{(m)})\to_{\cU}0$.
Since the isomorphism classes and the positive traces of the minimal
projections are fixed at this stage, we obtain from Theorem~\ref{thm:fixed-type-factor-no-loss}
the required correction by a unitary at operator-norm distance
less than $2^{-m}$ from $1$ on a set in $\cU$. Finally, we extend the
unitary by $1$ elsewhere and apply
Lemma~\ref{lem:cofinal-corrected-stage-assembly} with $\eta_m=2^{-m}$
and $s_m=m$.
\end{proof}

Applying Theorems~\ref{thm:local-implies-hi} and~\ref{thm:free-coordinate-lifts} to the abelian von Neumann algebras generated by two unitaries gives the following.

\begin{corollary}
\label{cor:free-unitary-lifting}
Let $(N,\tau)$ be a tracial von Neumann algebra, let $\cU$ be a nonprincipal ultrafilter on $\N$, and let
$u,v\in N^\cU$ be freely independent unitaries.
\begin{enumerate}[label=\textup{(\roman*)}]
\item For every fixed $n\ge1$, there are unitary lifts $u=(u_j)_\cU$ and $v=(v_j)_\cU$ such that
$W^*(u_j)$ and $W^*(v_j)$ are $n$-independent for every $j\in\N$.
\item There are unitary lifts $u=(u_j)_\cU$, $v=(v_j)_\cU$ and positive integers $n(j)\to_\cU\infty$ such that
$W^*(u_j)$ and $W^*(v_j)$ are $n(j)$-independent for every $j\in\N$.
\end{enumerate}
\end{corollary}

\begin{proof}
The von Neumann algebras $A:=W^*(u)$ and $B:=W^*(v)$ are freely independent, separable and abelian,
and a nonprincipal ultrafilter on $\N$ is countably cofinal.
For \textup{(i)}, we apply Theorem~\ref{thm:local-implies-hi}\textup{(i)} to the constant family $M_j=N$, using freeness to
supply $2n$-independence.
For \textup{(ii)}, we apply Theorem~\ref{thm:free-coordinate-lifts}\textup{(i)} and set $n(j):=r(j)$.
In both cases, we choose bounded representatives in the respective abelian
coordinate von Neumann subalgebras and replace them by their polar parts, extended by
$1$ on their kernels. This gives unitary representatives $u=(u_j)_\cU$ and
$v=(v_j)_\cU$. The von Neumann subalgebras $W^*(u_j)$ and $W^*(v_j)$ inherit
the asserted independence.
\end{proof}

By part~\textup{(i)}, the fixed-length interpretation of Problem~31 has an affirmative answer.
By Proposition~\ref{prop:no-prescribed-diagonal}, one
cannot in general require exact block-length $j$-independence at coordinate
$j$. In particular, one cannot always choose freely independent
coordinate lifts, since these would satisfy every prescribed length. Boulanger, Curda, Harvey, Li, and Pi~\cite[Section~3.2]{boulanger-curda-harvey-li-pi} gave an earlier argument for self-adjoint elements, using Popa's asymptotic freeness theorem~\cite[Theorem~2.1]{popa-free-independent-1995} and the nonamenability of suitable free products.

Both parts may also be derived from the lifting theorem independently
obtained by Patchell and
Shiner~\cite[Theorem~C]{patchell-shiner}.
For part~\textup{(ii)}, one combines it with the diagonal selection
in their proof of~\cite[Theorem~A]{patchell-shiner},
which gives independence orders tending to infinity along~$\cU$.

For the next two propositions, on prescribed lengths and on eventual
independence at every fixed order, we require the following separation lemma.
It is a quantitative form of the fact that the hyperfinite ${\rm II}_1$
factor $R$ cannot contain two freely independent
copies of $\C^3$ whose minimal projections all have trace less than $1/2$.
More precisely, there is an $\varepsilon>0$ such that, for every $p$ and a suitable $N\ge p$, some $p$-independent
pair of unitaries whose cubes are $1$ stays at distance
at least $\varepsilon$ from every $N$-independent
pair of such unitaries.

\begin{lemma}
\label{lem:finite-order-separation}
Let $R$ be the hyperfinite ${\rm II}_1$ factor. Set
\[
\mathcal X:=\{(x,y)\in U(R)^2:x^3=y^3=1\},
\qquad
 d((x,y),(x',y')):=\norm{x-x'}_2+\norm{y-y'}_2,
\]
and let $\mathcal D_m\subseteq\mathcal X$ consist of the pairs for which
$W^*(x)$ and $W^*(y)$ are $m$-independent.
There is $\varepsilon>0$ such that, for every $p\ge1$, there are $N\ge p$
and $\xi=(x,y)\in\mathcal D_p$ with
$\operatorname{dist}(\xi,\mathcal D_N)\ge\varepsilon$ and with all three
spectral projections of each of $x,y$ having traces in $(1/4,1/2)$.
\end{lemma}

\begin{proof}
The space $\mathcal X$ is complete, and the sets $\mathcal D_m$ are closed
and decreasing. Indeed, since $x^3=y^3=1$, the centered parts of $W^*(x)$ and $W^*(y)$ are spanned by
$x-\tau(x),\ x^2-\tau(x^2)$ and $y-\tau(y),\ y^2-\tau(y^2)$, respectively. Hence a pair belongs to
$\mathcal D_m$ if and only if finitely many traces of products of these elements vanish.
Set $\zeta=e^{2\pi i/3}$.

Let $\mathcal K$ consist of the pairs whose three spectral projections on
each side all have trace $1/3$, and fix $0<\rho<1/12$. For
\[
e_r(x):=\frac13\sum_{s=0}^2\zeta^{-rs}x^s,\qquad0\le r\le2,
\]
we have $\norm{e_r(x)-e_r(x')}_2\le\norm{x-x'}_2$. Thus every pair in
$\mathcal L:=\{\xi\in\mathcal X:\operatorname{dist}(\xi,\mathcal K)\le\rho\}$
has all six spectral weights in $(1/4,1/2)$.
By Theorem~\ref{thm:finite-order-free-positioning}\textup{(i)}, applied to
the centered parts of two uniformly weighted copies of $\C^3$, we have
$\mathcal K\cap\mathcal D_m\ne\varnothing$ for every $m$. In contrast,
\begin{equation}
\mathcal L\cap\bigcap_{m\ge1}\mathcal D_m=\varnothing.
\label{eq:diagonal-empty-intersection}
\end{equation}
Indeed, by Lemma~\ref{prop:truncated-free-product}, such a pair would be free. Since all its spectral weights
$\alpha_r$ and $\beta_r$ are less than $1/2$, Dykema's theorem for tracial
free products~\cite[Theorem~2.3]{dykema-free-dimension}
gives an interpolated free group factor~\cite{dykema-interpolated,radulescu-interpolated}, with no atomic summands:
\[
(\C^3,\alpha)*(\C^3,\beta)\cong L(\mathbb F_s),
\qquad s=2-\sum_{r=1}^3\alpha_r^2-\sum_{r=1}^3\beta_r^2>1.
\]
This contradicts amenability of von Neumann subalgebras of $R$.

We claim that some $\varepsilon>0$ has the following property: for every
$p$ there are $N\ge p$ and $\xi\in\mathcal L\cap\mathcal D_p$ with
\begin{equation}
\operatorname{dist}(\xi,\mathcal D_N)\ge\varepsilon.
\label{eq:diagonal-uniform-separation}
\end{equation}
Otherwise, for every $\eta>0$ there is $p$ such that
\begin{equation}
\operatorname{dist}(\xi,\mathcal D_N)<\eta
\quad(\xi\in\mathcal L\cap\mathcal D_p,\ N\ge p).
\label{eq:diagonal-empty-interiors}
\end{equation}
Choose $\eta_q>0$ such that $\sum_q\eta_q<\rho$, and strictly increasing
$p_q$ for which \eqref{eq:diagonal-empty-interiors} holds with
$\eta=\eta_q$. Starting with
$\xi_1\in\mathcal K\cap\mathcal D_{p_1}$, choose recursively
\[
\xi_{q+1}\in\mathcal D_{p_{q+1}},\qquad
 d(\xi_q,\xi_{q+1})<\eta_q.
\]
The sum of the successive distances is less than $\rho$, so all
$\xi_q$ lie in $\mathcal L$.
They form a Cauchy sequence whose limit belongs to every closed
$\mathcal D_m$, contradicting \eqref{eq:diagonal-empty-intersection}.
\end{proof}

\begin{proposition}
\label{prop:no-prescribed-diagonal}
Let $R$ be the hyperfinite ${\rm II}_1$ factor, and let $\cU$ be a nonprincipal ultrafilter on $\N$ such that
\[
O:=\{j\in\N:j\text{ is odd}\}\in\cU.
\]
There are freely independent unitaries $u,v\in R^\cU$ with
$u^3=v^3=1$ such that, for every pair of unitary representatives
$u=(a_j)_\cU$ and $v=(b_j)_\cU$,
\[
\bigl\{j\in\N:W^*(a_j),W^*(b_j)
\text{ are exact block-length }j\text{-independent}\bigr\}\notin\cU.
\]
\end{proposition}

\begin{proof}
Let $\varepsilon,\mathcal X,d,\mathcal D_m$ be as in Lemma~\ref{lem:finite-order-separation}, and set $\zeta=e^{2\pi i/3}$.
Choose, recursively, increasing integers $p_t,N_t$ and pairs $\xi_t=(x_t,y_t)$ whose six spectral weights lie in $(1/4,1/2)$, such that
\begin{equation}
N_t\ge p_t,
\qquad
\xi_t\in\mathcal D_{p_t},
\qquad
\operatorname{dist}(\xi_t,\mathcal D_{N_t})\ge\varepsilon,
\qquad
p_{t+1}>N_t.                                                        \label{eq:diagonal-obstruction-schedule}
\end{equation}
Let $J_t:=2N_t+1$. On every odd coordinate $j\ge J_1$, define
\[
(z_j,w_j):=\xi_t
\quad\text{when }J_t\le j<J_{t+1},
\]
and set $z_j=w_j=1$ at all remaining coordinates. The unitaries
\[
u:=(z_j)_\cU,
\qquad
v:=(w_j)_\cU
\]
are freely independent in $R^\cU$. To see this, we center each coordinate polynomial using its coordinate trace. Since $O\in\cU$ and $p_t\to\infty$ along the odd coordinates, every centered alternating moment of a fixed even
length vanishes at all sufficiently large odd coordinates.
By Lemma~\ref{prop:truncated-free-product}, the same holds at odd lengths.

On odd coordinates, we use the following identity to deduce vanishing at all smaller lengths from vanishing at one length.
If $A\cong\C^3$ has a faithful trace and $A^\circ:=A\ominus\C1$, then
\begin{equation}
\operatorname{span}_{\C}(A^\circ A^\circ)=A.                       \label{eq:three-atom-centered-products}
\end{equation}
Indeed, if $e_1,e_2,e_3$ are the minimal projections and $\lambda_i:=\tau(e_i)>0$, the centered elements
$h_i:=e_1/\lambda_1-e_i/\lambda_i$ $(i=2,3)$ satisfy $h_2h_3=e_1/\lambda_1^2$ and
$h_i^2=e_1/\lambda_1^2+e_i/\lambda_i^2$. These products span all three minimal projections, proving \eqref{eq:three-atom-centered-products}.

Suppose that faithful copies $A,B\cong\C^3$ have vanishing centered
alternating moments at a fixed odd block length $K$. We first descend
through the odd lengths. If a centered alternating word $W$ of odd
length $\ell$ starts and ends in $A$, we use
\eqref{eq:three-atom-centered-products} in $B$ to write
$1=\sum_s b_sc_s$ with $b_s,c_s\in B^\circ$. Since
$\tau(W)=\sum_s\tau(c_sWb_s)$, vanishing at length $\ell+2$ implies
vanishing at length $\ell$. The same argument applies with $A$ and $B$
interchanged.

For a centered alternating word $W=a_1b_1\cdots a_rb_r$ of even length $2r$, write
$a_1=\sum_s y_sx_s$ with $x_s,y_s\in A^\circ$. Cyclicity gives
$\tau(W)=\sum_s\tau(x_sb_1a_2b_2\cdots a_rb_ry_s)$, so vanishing at
length $2r+1$ implies vanishing at length $2r$. Again, the same argument
applies with the algebras interchanged. It follows that exact block-length
$K$-independence implies vanishing at every block length at most $K$.
In particular,
\begin{equation}
\text{exact block-length }(2r+1)\text{-independence implies membership in }\mathcal D_r\quad(r\ge1). \label{eq:odd-block-independence-implication}
\end{equation}

Suppose, toward a contradiction, that $u=(a_j)_\cU$ and $v=(b_j)_\cU$ are unitary representatives and that
\[
S:=\{j:W^*(a_j),W^*(b_j)\text{ are exact block-length }j\text{-independent}\}\in\cU.
\]
Let $q:\mathbb T\to\{1,\zeta,\zeta^2\}$ choose a nearest cube root, and set $\widehat a_j:=q(a_j)$ and
$\widehat b_j:=q(b_j)$ by Borel functional calculus. Pointwise on $\mathbb T$, one has
$\abs{\lambda-q(\lambda)}\le\abs{\lambda^3-1}$.
Since $(a_j)_\cU=(z_j)_\cU$, $(b_j)_\cU=(w_j)_\cU$, and $z_j^3=w_j^3=1$, it follows that
\begin{equation}
d((\widehat a_j,\widehat b_j),(z_j,w_j))\longrightarrow_\cU0.       \label{eq:diagonal-representative-distance}
\end{equation}
Indeed, for the first unitary,
$\norm{\widehat a_j-z_j}_2\le\norm{a_j^3-1}_2+\norm{a_j-z_j}_2
\le4\norm{a_j-z_j}_2$; the second is identical.
The rounded unitaries generate unital $^*$-subalgebras of $W^*(a_j)$ and $W^*(b_j)$, so exact block-length $j$-independence persists for
$j\in S$. By \eqref{eq:diagonal-representative-distance} and the uniform lower weight bound for the pairs $\xi_t$, the spectral weights of the rounded unitaries are positive
on a $\cU$-large subset of $O\cap\{j\ge J_1\}$.

Intersect this set with $S$. If $j$ belongs to the block $J_t\le j<J_{t+1}$, then $j$ is odd and
$j\ge2N_t+1$. By \eqref{eq:odd-block-independence-implication}, we have
\[
(\widehat a_j,\widehat b_j)\in\mathcal D_{(j-1)/2}\subseteq\mathcal D_{N_t}.
\]
This contradicts \eqref{eq:diagonal-obstruction-schedule} and \eqref{eq:diagonal-representative-distance}. Hence $S\notin\cU$, as desired.
\end{proof}

The argument in Proposition~\ref{prop:no-prescribed-diagonal} depends on the ultrafilter containing the odd integers. Indeed, vanishing of the centered alternating moments at one even block length need not imply vanishing at smaller lengths. For example, for $m\ge2$, set
$\theta=\pi/(2m)$ and, in $M_2(\C)$, let
\[
s=\begin{pmatrix}1&0\\0&-1\end{pmatrix},
\qquad
t=\begin{pmatrix}\cos\theta&\sin\theta\\ \sin\theta&-\cos\theta\end{pmatrix}.
\]
The centered parts of $W^*(s)$ and $W^*(t)$ are spanned by $s$ and $t$, respectively, and
$\operatorname{tr}_2((st)^r)=\cos(r\theta)$. Hence the pair is exact block-length $2m$-independent but not exact block-length
$2$-independent.

We now show that the convergence of the independence orders along $\cU$ in
Theorem~\ref{thm:free-coordinate-lifts}\textup{(i)} cannot in general be
replaced by ordinary convergence to infinity.

\begin{proposition}
\label{prop:no-ordinary-eventual-lifts}
Let $R$ be the hyperfinite ${\rm II}_1$ factor. There are a nonprincipal
ultrafilter $\cU$ on $\N$ and freely independent Haar unitaries $u,v\in R^\cU$
such that, for every pair of unitary representatives
$u=(u_k)_\cU$, $v=(v_k)_\cU$, there is $n\ge1$ for which
\[
\{k\in\N:W^*(u_k),W^*(v_k)\text{ are not }n\text{-independent}\}
\]
is infinite.
\end{proposition}

\begin{proof}
We first construct freely independent unitaries $a,b\in R^{\cU}$ of
order three with the required property. We then realize them as
functions of freely independent Haar unitaries.
Let $\varepsilon,\mathcal X,d,\mathcal D_m$ be as in
Lemma~\ref{lem:finite-order-separation}. For each $p\ge1$, choose $N_p\ge p$
and
$\xi_p=(a_p,b_p)\in\mathcal D_p$, with all six spectral weights in
$(1/4,1/2)$, such that
\begin{equation}
\operatorname{dist}(\xi_p,\mathcal D_{N_p})\ge\varepsilon.
\label{eq:ordinary-eventual-separation}
\end{equation}
For nonprincipal ultrafilters $\mathcal V,\mathcal W$ on $\N$, define an ultrafilter $\cU$ on $\N^2$ by
\[
S\in\cU\quad\Longleftrightarrow\quad
\{p:\{k:(p,k)\in S\}\in\mathcal W\}\in\mathcal V,
\]
and identify $\N^2$ with $\N$ by a fixed bijection.
Let $a=(a_p)_{(p,k),\cU}$ and $b=(b_p)_{(p,k),\cU}$, using
representatives that are constant in $k$. These unitaries are free. Indeed,
for every fixed $n$, the pair $(a_p,b_p)$ is $n$-independent for all
sufficiently large $p$, and the centered first and second powers span
the centered parts of their generated unital $^*$-algebras.
Each unitary has three spectral projections of trace at least $1/4$.

Suppose that unitary representatives $c_{p,k},d_{p,k}$ of $a,b$ generate, for every fixed $n$,
$n$-independent von Neumann algebras at all but finitely many indices.
We round them to cube roots as in the proof of
Proposition~\ref{prop:no-prescribed-diagonal}. Then
\[
d((\widehat c_{p,k},\widehat d_{p,k}),\xi_p)
\le4\bigl(\|c_{p,k}-a_p\|_2+\|d_{p,k}-b_p\|_2\bigr)
\longrightarrow_\cU0.
\]
Rounding preserves independence. For each fixed $p$, the rounded
unitaries generate $N_p$-independent unital $^*$-algebras for all but finitely
many $k$. Hence, by \eqref{eq:ordinary-eventual-separation}, their distance from $\xi_p$ is at least
$\varepsilon$ on a $\mathcal W$-large subset of that row.
By the definition of $\cU$, the distance is at least $\varepsilon$
on a $\cU$-large set, contradicting its convergence to zero.

Set $M=R^\cU$, $A=W^*(a)$, and $B=W^*(b)$. Let $p_1,p_2,p_3$
be the spectral projections of $a$. In each $p_iMp_i$, choose a
separable diffuse abelian von Neumann algebra $E_i$ with unit $p_i$, and set
$E_0=\bigoplus_{i=1}^3E_i$. Then $E_0\supseteq A$.
Construct a separable diffuse abelian von Neumann algebra $F_0\supseteq B$
in the same way using the spectral projections of $b$.
By Popa's asymptotic freeness theorem~\cite[Theorem~2.1]{popa-free-independent-1995}, in the formulation of~\cite[Corollary~0.2 and proof of Corollary~4.4\textup{(2)}]{popa-independence},
applied with $P=A\vee B$ kept fixed, there is $w\in U(A'\cap
M)$ such that $E:=wE_0w^*$ and $P$ are free over $A$. Moreover,
$E\supseteq A$.
Amalgamated freeness gives the trace-preserving identifications
\[
W^*(E,P)\cong E*_A P\cong E*_A(A*B)\cong E*B,
\]
fixing $E$ and $B$. Thus $E$ and $B$ are freely independent.
Applying the same argument over $B$ to $F_0$ and $E\vee B$, keeping the latter fixed, we
obtain a separable diffuse abelian von Neumann algebra $F\supseteq B$
free from $E$.

Choose Haar unitaries $u$ and $v$ generating $E$ and $F$, respectively, so that
$a=f(u)$, $b=g(v)$, where $f,g$ take the three cube-root values on
consecutive arcs. The Haar measures of these arcs are the spectral
weights of $a$ and $b$, respectively.
For any unitary representatives $u_k,v_k$, their spectral laws
converge weakly along $\cU$ to Haar measure.
Since the arc boundaries are finite and Haar-null, continuous
approximation in $L^2$ gives
\[
(f(u_k))_\cU=a,\qquad (g(v_k))_\cU=b.
\]
For each $k$, the von Neumann algebras generated by $f(u_k),g(v_k)$ lie in those
of $u_k,v_k$, respectively. If the latter were eventually
$n$-independent for every fixed $n$, the former would be as well,
contradicting the choice of $a,b$.
\end{proof}

If the coordinate algebras are ${\rm II}_1$ factors, the ultrafilter is countably
cofinal, and $W^*(A,B)$ is separable and hyperfinite, then we can lift $W^*(A,B)$
by one trace-preserving
unital $^*$-homomorphism at each coordinate. Applying the same homomorphism to
both algebras preserves all their mixed moments, giving the following
sufficient condition for freely independent coordinate lifts.

\begin{theorem}
\label{thm:joint-hyperfinite-no-loss}
Let $\cU$ be countably cofinal, let $(M_k,\tau_k)_{k\in I}$ be
${\rm II}_1$ factors, and set $M:=\ultraprod M_k$.  Let
$A,B\subseteq M$ be freely independent separable abelian von Neumann
subalgebras.  If
\[
N:=W^*(A\cup B)
\]
is hyperfinite, then there are finite-dimensional abelian unital $^*$-subalgebras $C_k,D_k\subseteq M_k$ such that
\[
A\subseteq\ultraprod C_k,
\qquad
B\subseteq\ultraprod D_k,
\]
and $C_k,D_k$ are freely independent for every $k$.
\end{theorem}

\begin{proof}
Applying~\cite[Theorem~4.11]{hadwin-li} to the
countably generated hyperfinite von Neumann algebra $N$, with the compatible lifting
construction explained before
Theorem~\ref{thm:hyperfinite-lifting-cofinal},
we obtain trace-preserving unital $^*$-homomorphisms $\pi_k:N\to M_k$ representing its inclusion.
Choose increasing sequences of finite-dimensional abelian unital $^*$-subalgebras
$A^{(m)}\subseteq A$ and $B^{(m)}\subseteq B$ with weakly dense unions,
and choose $\mu\to_{\cU}\infty$ by
Fact~\ref{lem:cofinal-diagonal-choice} with $U_m=I$. For $\mu(k)\ge1$, set
\begin{equation}
 C_k=\pi_k(A^{(\mu(k))}),\qquad D_k=\pi_k(B^{(\mu(k))}),
 \label{eq:joint-hyperfinite-lifted-relations}
\end{equation}
and use scalars otherwise.
At each coordinate, the same trace-preserving homomorphism is applied
to both algebras, so their images are free. For fixed $a\in A^{(m)}$,
$\pi_k(a)\in C_k$ on $\{\mu\ge m\}$; the same holds for $B$.
Weak closure gives both ultraproduct inclusions.
\end{proof}

For freely independent abelian $A,B\ne\C1$, joint hyperfiniteness
forces $\dim A=\dim B=2$. Indeed, otherwise the two algebras contain
free copies of $\C^3$ and $\C^2$ with faithful traces, whose generated algebra
has a nonhyperfinite free-group-factor summand
by~\cite[Theorem~2.3]{dykema-free-dimension}. Thus the nontrivial free case
of Theorem~\ref{thm:joint-hyperfinite-no-loss} concerns two projections. If one algebra is scalar,
the other can be any separable abelian algebra.

If $A,B$ in Theorem~\ref{thm:joint-hyperfinite-no-loss} are only
$n$-independent, the same argument gives $n$-independent
finite-dimensional lifts.
In fact, by Theorem~\ref{thm:factor-no-loss-hierarchy}, this conclusion
holds without the joint hyperfiniteness assumption.

\subsection{Local stability and orthogonality}
\label{subsec:lifting-local-stability}

We now return to a fixed finite order and allow projection traces
to tend to zero. We show that lifting at order $n$ is equivalent to a local stability property of pairs of partitions, uniform in the projection traces (Proposition~\ref{prop:higher-order-local-global}). Combined with Theorem~\ref{thm:factor-no-loss-hierarchy}, this gives small $L^2$ perturbations of both partitions in $\mathrm{II}_1$ factors, with bounds independent of the projection traces (Corollary~\ref{cor:factor-full-local-no-loss}).

Let $\mathfrak C$ be a class of tracial von Neumann algebras.  We write
$\operatorname{LS}_n(\mathfrak C;a,b)$ for the following local statement: for every $\eta>0$ there is $\epsilon>0$,
uniform over $(N,\tau)\in\mathfrak C$ and over all trace vectors, such that, for every $(\mathbf p,\mathbf q)\in\mathcal P_{a,b}(N)$,
\[
\operatorname{def}_n(\mathbf p,\mathbf q)<\epsilon
\]
implies the existence of a pair $(\mathbf p',\mathbf q')\in\mathcal P_{a,b}(N)$ whose generated unital $^*$-algebras are
$n$-independent and for which
\begin{equation}
\sum_i\norm{p_i-p_i'}_2+\sum_j\norm{q_j-q_j'}_2<\eta.
\label{eq:higher-local-distance}
\end{equation}
The exact correction in this definition need not preserve the projection traces.

We write $\operatorname{Lift}_n(\mathfrak C)$ for the following property of lifting without decreasing the order.  Whenever $\cU$ is countably
cofinal, $M_k\in\mathfrak C$, and $A,B\subseteq\prod_{\cU}M_k$ are separable abelian von Neumann subalgebras which are $n$-independent,
there are finite-dimensional abelian unital $^*$-subalgebras $C_k,D_k\subseteq M_k$ such that
\[
A\subseteq\ultraprod C_k,
\qquad
B\subseteq\ultraprod D_k,
\]
and $C_k,D_k$ are $n$-independent for every $k$.

In this notation, Theorem~\ref{thm:factor-no-loss-hierarchy} states that
\[
\operatorname{Lift}_n(\{\text{${\rm II}_1$ factors}\})
\]
holds for every fixed finite $n$.

\begin{proposition}
\label{prop:higher-order-local-global}
Fix $n\ge1$, and let $\mathfrak C$ be any class of tracial von Neumann algebras.  Then
the following are equivalent:
\begin{enumerate}[label=\textup{(\roman*)},beginpenalty=10000]
\item $\operatorname{LS}_n(\mathfrak C;a,b)$ holds for every finite $a,b$;
\item $\operatorname{Lift}_n(\mathfrak C)$ holds.
\end{enumerate}
\end{proposition}

\begin{proof}
Assume \textup{(i)} and fix $\cU$, $(M_k)_{k\in I}$, and an $n$-independent
pair $A,B$ as in the definition of $\operatorname{Lift}_n(\mathfrak C)$.
Choose increasing chains of finite-dimensional abelian unital $^*$-subalgebras
$A^{(m)}\subseteq A$ and $B^{(m)}\subseteq B$ whose unions are weakly dense.  We use Lemma~\ref{lem:nested-lifts-general}
to lift the minimal projections, preserving their indices and refinement identities.
At each fixed stage, every moment of pair length at most $n$ in these projections converges to its value in the ultraproduct, so
\[
\operatorname{def}_n(\mathbf p_k^{(m)},\mathbf q_k^{(m)})
\longrightarrow_{\cU}0.
\]
We apply \textup{(i)} at stage $m$ with error $\eta_m\downarrow0$ on
a set $U_m\in\cU$.
By Lemma~\ref{lem:cofinal-two-sided-labelled-assembly}, at this fixed order $n$, there are
finite-dimensional abelian unital $^*$-subalgebras $C_k,D_k\subseteq M_k$, $n$-independent
for every $k$, whose respective ultraproducts
contain $A$ and $B$. This proves (ii).

Conversely, suppose (i) fails for some fixed partition sizes $a,b$. Choose $N_m\in\mathfrak C$ and pairs $(\mathbf p_m,\mathbf q_m)\in\mathcal P_{a,b}(N_m)$ whose defects tend to zero
but which remain at distance at least $\eta_0>0$ from every $n$-independent pair.  In any nonprincipal ultraproduct over $\N$,
the induced pair is $n$-independent.  We apply (ii), and then the first assertion of Lemma~\ref{lem:nested-lifts-general}
inside the coordinate algebras $C_m$ and $D_m$, to lift the two partitions in the ultraproduct exactly. The sum of the
$L^2$ distances between corresponding new and original projections tends to zero along
the ultrafilter, and the new projections generate $n$-independent unital $^*$-subalgebras of $C_m,D_m$, a contradiction.
\end{proof}

\begin{corollary}
\label{cor:factor-full-local-no-loss}
For every fixed $n,a,b\ge1$ and every $\eta>0$, there is
\[
\epsilon=\epsilon(n,a,b,\eta)>0,
\]
uniform over all ${\rm II}_1$ factors and all trace vectors, with the following property. If $M$ is such a factor and
$(\mathbf p,\mathbf q)\in\mathcal P_{a,b}(M)$ satisfies
\[
\operatorname{def}_n(\mathbf p,\mathbf q)<\epsilon,
\]
then there is a pair $(\mathbf p',\mathbf q')\in\mathcal P_{a,b}(M)$
whose generated unital $^*$-algebras are $n$-independent and which satisfies
\[
\sum_i\norm{p_i-p_i'}_2+
\sum_j\norm{q_j-q_j'}_2<\eta.
\]
The correction may change both partitions and their traces.
\end{corollary}

\begin{proof}
By Theorem~\ref{thm:factor-no-loss-hierarchy}, $\operatorname{Lift}_n(\{\text{${\rm II}_1$ factors}\})$ holds,
so the result follows from the converse direction of
Proposition~\ref{prop:higher-order-local-global}.
\end{proof}

Applying Theorem~\ref{thm:factor-no-loss-hierarchy} at order one, we obtain
an affirmative answer to Problem~30.

\begin{corollary}[Answer to Problem~30]
\label{cor:fixed-factor-problem-30}
Let $N$ be a ${\rm II}_1$ factor, let $\cU$ be a nonprincipal ultrafilter on $\N$, and let
$u,v\in N^\cU$ be unitaries with $W^*(u)\perp W^*(v)$. Then there are unitary representatives
\[
u=(u_j)_{\cU},
\qquad
v=(v_j)_{\cU}
\]
such that $W^*(u_j)$ and $W^*(v_j)$ are $1$-independent, equivalently orthogonal, for every $j\in\N$.
\end{corollary}

\begin{proof}
The von Neumann algebras $W^*(u),W^*(v)$ are $1$-independent, and a nonprincipal
ultrafilter on $\N$ is countably cofinal.
Hence, by Theorem~\ref{thm:factor-no-loss-hierarchy} at $n=1$, there are finite-dimensional
abelian unital $^*$-subalgebras $C_j,D_j\subseteq N$, orthogonal at every
coordinate, whose ultraproducts contain $u,v$.
Take bounded representatives in these algebras and replace them by
their polar parts, extended by $1$ on their kernels. The resulting
unitaries still represent $u,v$ and generate unital $^*$-subalgebras of
$C_j,D_j$, respectively.
\end{proof}

We combine Theorem~\ref{thm:factor-no-loss-hierarchy} with the construction of
Chifan, Ioana, and Kunnawalkam
Elayavalli~\cite[Definition~5.1 and Proposition~5.4]{chifan-ioana-kunnawalkam-elayavalli},
using Houdayer and Ioana's enumeration of pairs generating orthogonal
von Neumann algebras~\cite[Definition~5.7]{houdayer-ioana}.
Chifan, Ioana, and Kunnawalkam Elayavalli~\cite[Theorem~B]{chifan-ioana-kunnawalkam-elayavalli}
assumed $u_1^2=u_2^3=1$ in addition to orthogonality. Houdayer and
Ioana~\cite[Theorem~G]{houdayer-ioana} allowed arbitrary unitaries but assumed
$2$-independence. In the result below, we require only orthogonality.

\begin{theorem}
\label{thm:houdayer-ioana-orthogonal-g}
There exists a separable non-Gamma ${\rm II}_1$ factor $M$ with the following property.  If $\cU$ is a
countably cofinal ultrafilter on a set $I$ and unitaries $u_1,u_2\in M^\cU$ satisfy
\[
W^*(u_1)\perp W^*(u_2),
\]
then there are Haar unitaries $v_1,v_2\in M^\cU$ such that
\[
[u_1,v_1]=[u_2,v_2]=[v_1,v_2]=0.
\]
\end{theorem}

\begin{proof}
We start with a separable ${\rm II}_1$ factor $M_1$ with property~$(T)$
and use the construction of Houdayer and
Ioana~\cite[Definition~5.7, Corollary~5.6, and Proposition~5.9]{houdayer-ioana}.
It gives the separable non-Gamma factor $M=(\bigcup_mM_m)''$,
sequences $(u_1^{m,j},u_2^{m,j})_j$ that are $L^2$-dense among the pairs of unitaries in each $M_m$ that generate orthogonal von Neumann algebras,
and a bijection $\sigma:\N\to\N^2$ with $\sigma_1(\ell)\le\ell$.
At stage $\ell+1$, the two commuting Haar unitaries for the pair indexed by
$\sigma(\ell)$ come from~\cite[Corollary~4.3]{chifan-ioana-kunnawalkam-elayavalli}.

We follow the proof of~\cite[Proposition~5.3]{chifan-ioana-kunnawalkam-elayavalli}, as adapted by Houdayer and Ioana~\cite[Proposition~5.8]{houdayer-ioana}, with
Theorem~\ref{thm:factor-no-loss-hierarchy} at order one in place of~\cite[Theorem~5.1]{houdayer-ioana}.
For $Q_0=W^*(u_1,u_2)$, choose indices $n_i$ such that
$Q_0\subseteq\prod_\cU M_{n_i}$, using
\cite[Lemma~2.2]{boutonnet-chifan-ioana}.
By Theorem~\ref{thm:factor-no-loss-hierarchy} at order one, there are
orthogonal abelian von Neumann subalgebras $C_i,D_i\subseteq M_{n_i}$ with
$W^*(u_1)\subseteq\prod_\cU C_i$ and
$W^*(u_2)\subseteq\prod_\cU D_i$.
Take bounded representatives in these algebras and replace them by
their polar parts, extended by $1$ on their kernels. The resulting unitaries
$u_{1,i},u_{2,i}$ represent $u_1,u_2$ and generate orthogonal von Neumann
algebras.

Choose $f:I\to\N$ with $f(i)\to_\cU\infty$ and, by density, indices
$j_i$ with
\[
\|u_{1,i}-u_1^{n_i,j_i}\|_2+
\|u_{2,i}-u_2^{n_i,j_i}\|_2\le1/f(i).
\]
Set $\ell_i=\sigma^{-1}(n_i,j_i)$. Then $n_i=\sigma_1(\ell_i)\le\ell_i$.
The construction therefore gives Haar unitaries $v_{1,i},v_{2,i}\in M_{\ell_i+1}$
with
\[
[u_1^{n_i,j_i},v_{1,i}]=[u_2^{n_i,j_i},v_{2,i}]
=[v_{1,i},v_{2,i}]=0.
\]
Since $\|[u_{a,i},v_{a,i}]\|_2\le2/f(i)$,
the unitaries $v_a=(v_{a,i})_\cU$ are Haar and satisfy
the three required commutation relations.
\end{proof}

Note that orthogonality is strictly weaker than $2$-independence. Indeed, a trace-preserving copy of $M_2(\C)$ in a tracial von Neumann algebra contains
\[
Z=\begin{pmatrix}1&0\\0&-1\end{pmatrix},
\qquad
X=\begin{pmatrix}0&1\\1&0\end{pmatrix}.
\]
Then $W^*(Z)\perp W^*(X)$, but $\tau(ZXZX)=-1$, so the pair is not $2$-independent.

In the proof of Theorem~\ref{thm:houdayer-ioana-orthogonal-g}, we use that
each $M_{n_i}$ is a factor: over arbitrary tracial
coordinate von Neumann algebras, orthogonality cannot replace $2$-independence in
\cite[Theorem~5.1]{houdayer-ioana}
(Proposition~\ref{prop:four-point-orthogonality-obstruction}).

In the next two propositions, we give examples showing that independence of order one or two need
not lift to coordinate von Neumann algebras that are not factors. In both examples the underlying algebra stays fixed while the
tracial state varies with the coordinate. Thus the first example does not settle the tracial version of
Problem~30, where the algebra and its trace must both remain fixed. We return to fixed
tracial von Neumann algebras in
Sections~\ref{sec:fixed-ambient-beyond-factors}
and \ref{sec:unstable-tables}.

\begin{proposition}
\label{prop:four-point-orthogonality-obstruction}
For $k\in\N$, set $\theta_k=1/(16(k+1))$ and equip $M_k=\C^4$ with the faithful trace $\tau_k$ whose four atomic weights are
\[
\left(\frac14+\theta_k,\frac14-\theta_k,
\frac14-\theta_k,\frac14+\theta_k\right).
\]
Let $p_k=(1,1,0,0)$ and $q_k=(1,0,1,0)$. For every nonprincipal ultrafilter $\cU$ on
$\N$, the unital $^*$-algebras generated by $p=(p_k)_{\cU}$ and $q=(q_k)_{\cU}$ are $1$-independent in
$\prod_{\cU}(M_k,\tau_k)$. There are, however, no coordinate von Neumann subalgebras $C_k,D_k\subseteq M_k$
such that
\[
W^*(p)\subseteq\ultraprod C_k,
\qquad
W^*(q)\subseteq\ultraprod D_k,
\]
and $C_k,D_k$ are $1$-independent on a $\cU$-large set.
\end{proposition}

\begin{proof}
Here $\tau_k(p_k)=\tau_k(q_k)=1/2$ and
$\tau_k(p_kq_k)=1/4+\theta_k$, so $W^*(p)$ and $W^*(q)$
are orthogonal in the ultraproduct.
If such coordinate algebras $C_k,D_k$ existed, then by
Lemma~\ref{lem:nested-lifts-general} there would be projections $e_k\in
C_k$, $f_k\in D_k$ representing $p,q$.
Distinct projections in $\C^4$ have $L^2$ distance at least
$1/\sqrt8$, because every atomic weight is at least $1/8$.
Thus $e_k=p_k$ and $f_k=q_k$ on a $\cU$-large set,
where coordinate orthogonality would give
$0=\tau_k(p_kq_k)-\tau_k(p_k)\tau_k(q_k)=\theta_k>0$.
\end{proof}

\begin{proposition}
\label{prop:order-two-arbitrary-obstruction}
Let $\mathfrak T$ be the class of all tracial von Neumann algebras.  Then
\[
\operatorname{LS}_2(\mathfrak T;2,2)
\quad\text{and}\quad
\operatorname{Lift}_2(\mathfrak T)
\]
are both false.
\end{proposition}

\begin{proof}
Let $w_k=1/3-1/(24(k+1))$ and
$\tau_k=w_k\operatorname{tr}_2\oplus(1-w_k)\operatorname{tr}_2$
on $M_k=M_2(\C)\oplus M_2(\C)$. Take
\begin{equation}
p_k=(e_{11},e_{11}),\qquad
q_k=\left(e_{11},
\begin{pmatrix}
1/4&\sqrt3/4\\ \sqrt3/4&3/4
\end{pmatrix}\right).
\label{eq:order-two-obstruction-projections}
\end{equation}
Both have trace $1/2$. With $t_k=\tau_k(p_kq_k)$ and
$s_k=\tau_k(p_kq_kp_kq_k)$,
\begin{equation}
t_k=\frac12\left(w_k+\frac{1-w_k}{4}\right),\qquad
s_k=\frac12\left(w_k+\frac{1-w_k}{16}\right),
\label{eq:order-two-obstruction-moments}
\end{equation}
so
\begin{equation}
t_k\to\frac14,\qquad s_k-t_k+\frac1{16}\to0.
\label{eq:order-two-vanishing-defect}
\end{equation}
For projections $p,q$ of traces $\alpha,\beta$,
\begin{equation}
\begin{aligned}
&\tau((p-\alpha)(q-\beta)(p-\alpha)(q-\beta))\\
&=\tau(pqpq)+(4\alpha\beta-2\alpha-2\beta)\tau(pq)
+\alpha\beta^2+\alpha^2\beta-3\alpha^2\beta^2.
\end{aligned}
\label{eq:centered-binary-fourth-moment}
\end{equation}
When $W^*(p)$ and $W^*(q)$ are independent to order one, this
expression reduces to
$\tau(pqpq)-\alpha\beta(\alpha+\beta-\alpha\beta)$.
Thus the order-two defect of $(p_k,1-p_k)$ and $(q_k,1-q_k)$ tends
to zero.

Since both central weights exceed $1/4$, a projection whose rank differs
from one in at least one summand is at $L^2$ distance at least $1/\sqrt8$
from both $p_k$ and $q_k$.
For projections $e=(e_1,e_2)$ and $f=(f_1,f_2)$ of rank one in
each summand, set $\mu_r=\operatorname{Tr}(e_rf_r)\in[0,1]$.
Order-one and order-two independence require
\begin{equation}
w_k\mu_1+(1-w_k)\mu_2=\frac12,\qquad
w_k\mu_1^2+(1-w_k)\mu_2^2=\frac38.
\label{eq:order-two-exact-moment-constraints}
\end{equation}
By the first equation, $\mu_1=1/2+(1-w_k)d$ and $\mu_2=1/2-w_kd$
for some $d\in\R$. The interval constraints give $|d|\le1/(2(1-w_k))$, hence
\begin{equation}
w_k(1-w_k)d^2
\le\frac{w_k}{4(1-w_k)}<\frac18.
\label{eq:order-two-variance-obstruction}
\end{equation}
This contradicts the required variance $3/8-1/4=1/8$.
Every $2$-independent pair must therefore change the rank of at least one
projection in one summand. Its $L^2$ distance from the original pair is at
least $1/\sqrt8$, disproving $\operatorname{LS}_2(\mathfrak T;2,2)$.

For every nonprincipal $\cU$, the projections $p=(p_k)_\cU$ and
$q=(q_k)_\cU$ generate $2$-independent unital $^*$-algebras. Applying the converse
argument in Proposition~\ref{prop:higher-order-local-global} to the separation above, we
rule out $2$-independent coordinate von Neumann subalgebras representing them.
Thus $\operatorname{Lift}_2(\mathfrak T)$ fails.
\end{proof}

The two limiting projections in Proposition~\ref{prop:order-two-arbitrary-obstruction} generate a finite-dimensional, hence
hyperfinite, von Neumann algebra. However, the joint lifting argument in Theorem~\ref{thm:joint-hyperfinite-no-loss} does not apply, since the coordinate von Neumann algebras are not factors.

For the limiting projections, set $S=2p-1$ and $T=2q-1$. Then
\[
\tau((ST)^r)=\frac13+\frac23\cos\left(\frac{2\pi r}{3}\right),
\]
so the pair is not $3$-independent.

\section{Stability of prescribed overlap matrices}
\label{sec:projection-tables}

In Sections~\ref{sec:projection-tables}--\ref{sec:unstable-tables}, we study the
tracial version of Problem~30 for a fixed tracial von Neumann algebra $N$. By
Proposition~\ref{prop:fixed-ambient-stability-equivalence}, this is
equivalent to stability of the independence equations for pairs of
partitions of $1$ in $N$, allowing zero projections. A pair can be independent for the global trace although the overlap matrices of its restrictions to the central summands are not of product form. We therefore consider arbitrary prescribed values of
$\tau(p_iq_j)$.

We prove stability when one partition is binary. In
Subsection~\ref{subsec:further-overlap-special-cases}, we collect further special cases, including
overlap matrices supported on
graphs without four-cycles. In general, however, a feasible overlap matrix
need not be stable even in a fixed factor, as we
show in Section~\ref{sec:unstable-tables}. We also prove here Lemma~\ref{lem:finite-coordinate-overlap-inverse}, our main tool on the factor summands in Section~\ref{sec:fixed-ambient-beyond-factors}: when the
Gram matrix of a chosen family of commutators $i[p_i,q_j]$ is bounded
below, a small conjugation attains any prescribed nearby values of the corresponding
overlaps.

\subsection{Binary partitions and the overlap inverse lemma}

For a tracial von Neumann algebra $(D,\tau_D)$ and a matrix
$C\in M_{a\times b}(\R)$ with nonnegative entries summing to $1$, write
\[
\mathcal Y_D(C):=\{(\mathbf p,\mathbf q)\in\mathcal P_{a,b}(D):
\tau_D(p_iq_j)=C_{ij}\text{ for all }i,j\}.
\]
We say that the prescribed overlap matrix $C$ is stable in $D$ if every sequence $(\mathbf p^{(m)},\mathbf q^{(m)})$
in $\mathcal P_{a,b}(D)$ satisfying
\[
\max_{i,j}\abs{\tau_D(p_i^{(m)}q_j^{(m)})-C_{ij}}\longrightarrow0
\]
has distance from $\mathcal Y_D(C)$ tending to zero. If $D$ is
diffuse, every such $C$ is realized by splitting $1_D$ into pairwise
orthogonal projections of traces $C_{ij}$ and taking their row and
column sums.

For two projections in a tracial ultraproduct of ${\rm II}_1$ factors, the
unital von Neumann algebra they generate is
hyperfinite, and this allows trace-preserving lifts into the coordinate
factors
that retain any fixed finite collection of their joint moment relations.
For a projection $p$ and an arbitrary finite partition $(q_j)$, we use a
different argument. We first choose a positive contraction with the required
overlaps against the $q_j$, and then replace it by a projection with those
overlaps. The target traces and overlaps may vary in this
construction, which works in every diffuse tracial von Neumann algebra
but controls only the overlaps $\tau(pq_j)$, not higher joint moments.

In the next proposition, $\mathbf w$ denotes the tuple $(w_1,\ldots,w_s)$. We say that a von Neumann algebra of type~I has matrix fibers of size at most two if it is isomorphic to $A_1\oplus M_2(A_2)$ with $A_1,A_2$ abelian, possibly zero.

\begin{proposition}
\label{prop:two-projection-tracial-relation-stability}
Fix noncommutative $^*$-polynomials $w_1,\ldots,w_s$ in two self-adjoint variables and a continuous function
$\Phi:\C^s\to[0,\infty)$.  For projections $p,q$ in a tracial von Neumann algebra $(M,\tau)$,
define their defect
\[
D_{\Phi,\mathbf w}^{M}(p,q)
:=
\Phi\bigl(
\tau(w_1(p,q)),\ldots,\tau(w_s(p,q))
\bigr).
\]
For every $\eta>0$, there is
\[
\varepsilon=\varepsilon(\Phi,\mathbf w,\eta)>0
\]
such that, uniformly over all ${\rm II}_1$ factors $(M,\tau)$,
\[
D_{\Phi,\mathbf w}^{M}(p,q)<\varepsilon
\]
implies that there are projections $p',q'\in M$ satisfying
\[
D_{\Phi,\mathbf w}^{M}(p',q')=0,
\qquad
\norm{p-p'}_2+\norm{q-q'}_2<\eta.
\]
If $f_1,\ldots,f_k$ are continuous functions on $\C^s$, the choice
$\Phi=\max_\ell|f_\ell|$ gives uniform $L^2$-stability over
${\rm II}_1$ factors for the system
\[
f_\ell\bigl(\tau(w_1(p,q)),\ldots,\tau(w_s(p,q))\bigr)=0
\qquad(1\le\ell\le k).
\]
\end{proposition}

\begin{proof}
Suppose the assertion fails. Then we may choose pairs $(p_m,q_m)$ in ${\rm II}_1$ factors $(M_m,\tau_m)$ with defects
tending to zero but distance at least $\eta>0$ from the pairs satisfying the prescribed relation.
Take a tracial ultraproduct over a nonprincipal ultrafilter $\cU$ and
write $p,q$ for the projections they represent. By continuity, we have $D_{\Phi,\mathbf
w}^{N}(p,q)=0$ in the unital von Neumann algebra $N=W^*(p,q)$.

The von Neumann algebra $N$ is separable, and by Dixmier's description of two
projections~\cite{dixmier-position-relative} (see also~\cite{halmos-two-subspaces}), it is of type I with
matrix fibers of size at most two, hence hyperfinite.
Applying~\cite[Theorem~4.11]{hadwin-li}, as in
the proof of Theorem~\ref{thm:hyperfinite-lifting-cofinal}, we obtain trace-preserving
unital
$^*$-homomorphisms $\pi_m:N\to M_m$ representing the inclusion of $N$ in the ultraproduct.
Then $p'_m=\pi_m(p)$ and $q'_m=\pi_m(q)$ satisfy every prescribed moment relation
exactly at every coordinate and
\[
\norm{p_m-p'_m}_2+\norm{q_m-q'_m}_2\longrightarrow_{\cU}0,
\]
contradicting the separation.
\end{proof}

For a projection $p$ and a partition $\mathbf q=(q_j)_{j=0}^{b-1}$ of $1$ in a tracial von Neumann algebra $(M,\tau)$, allowing zero projections in $\mathbf q$, write
\[
\beta_j:=\tau(q_j),
\qquad
c_j:=\tau(pq_j).
\]
Then $(\boldsymbol\beta,\mathbf c)$ lies in the compact set
\[
\mathcal T_b
:=
\left\{(\boldsymbol\beta,\mathbf c)\in[0,1]^{2b}:
\sum_{j<b}\beta_j=1,\quad 0\le c_j\le\beta_j\text{ for every }j\right\}.
\]
The trace of $p$ is $\sum_jc_j$, and the complementary row of overlaps is $(\beta_j-c_j)_j$.
We measure the distance between $(p;\mathbf q)$ and a second such pair $(p';\mathbf q')$ by
\[
2\norm{p-p'}_2+\sum_{j<b}\norm{q_j-q'_j}_2,
\]
which is $d_2$ for the pairs of partitions $((p,1-p),\mathbf q)$ and $((p',1-p'),\mathbf q')$.

\begin{theorem}
\label{thm:diffuse-binary-table-stability}
Fix $b\ge1$.  For every $\eta>0$, there is $\varepsilon=\varepsilon(b,\eta)>0$ such that the following holds uniformly
over diffuse tracial von Neumann algebras $(D,\tau)$ and targets
$(\widehat{\boldsymbol\beta},\widehat{\mathbf c})\in\mathcal T_b$.
If $p$ is a projection in $D$, $\mathbf q$ is a partition of $1_D$ into $b$ projections, and
\[
\max_{j<b}
\left(
\abs{\tau(q_j)-\widehat\beta_j}
+\abs{\tau(pq_j)-\widehat c_j}
\right)<\varepsilon,
\]
then there are a projection $p'$ and a partition $\mathbf q'$ of $1_D$ such that
\[
\tau(q'_j)=\widehat\beta_j,
\qquad
\tau(p'q'_j)=\widehat c_j
\qquad(j<b),
\]
and
\[
2\norm{p-p'}_2+\sum_{j<b}\norm{q_j-q'_j}_2<\eta.
\]
Zero projections are allowed in both $\mathbf q$ and $\mathbf q'$.
\end{theorem}

\begin{proof}
We transfer subprojections from those $q_j$ with $\tau(q_j)>\widehat\beta_j$
to those with $\tau(q_j)<\widehat\beta_j$, obtaining a partition $\mathbf q'$ with
\[
\tau(q'_j)=\widehat\beta_j,\qquad
\norm{q_j-q'_j}_1=\norm{q_j-q'_j}_2^2
=\abs{\tau(q_j)-\widehat\beta_j}.
\]
Set $\beta_j=\widehat\beta_j$, $c_j=\tau(pq'_j)$, and
$d_j=\widehat c_j$. Then $\abs{c_j-d_j}<\varepsilon$.
Set $a_j=\sqrt{d_j/c_j}$ when $d_j<c_j$ and $a_j=1$ otherwise. Similarly,
set $b_j=\sqrt{(\beta_j-d_j)/(\beta_j-c_j)}$ when $d_j>c_j$
and $b_j=1$ otherwise. Define the positive contractions
\[
A=\sum_j a_jq'_j,\qquad B=\sum_j b_jq'_j,\qquad
x=ApA,\qquad y=1-B(1-x)B.
\]
Since $A,B$ commute with every $q'_j$,
\[
\tau(xq'_j)=\min\{c_j,d_j\},\qquad \tau(yq'_j)=d_j.
\]
For positive contractions $z,C$, H\"older's inequality gives
\[
\norm{z-CzC}_1
\le2\sqrt{\tau(z)\tau((1-C)^2z)}
\le2\sqrt{\tau((1-C)^2z)}.
\]
We apply this first to $(z,C)=(p,A)$ and then to $(1-x,B)$.
Writing $\delta_-:=\sum_j(c_j-d_j)^+$ and
$\delta_+:=\sum_j(d_j-c_j)^+$, we obtain
\[
\norm{y-p}_1\le2(\sqrt{\delta_-}+\sqrt{\delta_+})
\le2\sqrt{2b\varepsilon},
\]
since $(1-\sqrt{u/v})^2v\le v-u$ for $0\le u\le v$, $v>0$.
By Lemma~\ref{lem:finite-functional-purification}, applied to $y$ and the trace
pairings with $p,q'_0,\ldots,q'_{b-1}$, there is a projection $p'$ with the
required overlaps. Since $\sum_jq'_j=1$, these pairings
also give $\tau(p')=\tau(y)$. Hence
\[
\norm{p'-p}_2^2
=\tau((y-p)(1-2p))\le\norm{y-p}_1.
\]
Consequently, $2\norm{p-p'}_2+\sum_j\norm{q_j-q'_j}_2$ is at most
$2(2\sqrt{2b\varepsilon})^{1/2}+b\sqrt\varepsilon$.
Choosing $\varepsilon$ so that this is less than $\eta$ proves the
assertion uniformly in $D$ and $(\widehat{\boldsymbol\beta},\widehat{\mathbf c})$.
\end{proof}

In the proof of Theorem~\ref{thm:diffuse-binary-table-stability}, we
replaced a single positive contraction by a projection on which finitely many
normal functionals take the same values (Lemma~\ref{lem:finite-functional-purification}). However, this step
does not extend to partitions: replacing positive
contractions with sum $1$
by a partition of $1$ while preserving the values of finitely many normal functionals
can fail even in a diffuse factor.
To see this, let $S_1,S_2,S_3\in M_2(\C)$ be the rank-one projections onto
\[
(1,0),\qquad
(-1/2,\sqrt3/2),\qquad
(-1/2,-\sqrt3/2).
\]
The $S_i$ are linearly independent and sum to $\tfrac32$.  In $M_2(\C)\overline\otimes M$, with $M$ a ${\rm II}_1$ factor,
the positive contractions $\tfrac23S_i\otimes1$ sum to $1$.  We claim that there is no partition $(e_i)_{i=1}^3$ of $1$ with
\[
(\operatorname{id}\otimes\tau_M)(e_i)=\tfrac23S_i.
\]
Indeed, positivity and faithfulness of the expectation force $e_i\le S_i\otimes1$, so $e_i=S_i\otimes r_i$ for projections
$r_i\in M$.  Since $\sum_i e_i=1$ and the $S_i$ are linearly independent, we have $r_i=\tfrac23$, which is impossible.

We now ask when the overlaps of a projection $p$ with a finite partition $(q_j)$
can be corrected by conjugating the projection with a unitary close to
$1$ in operator norm, while leaving the partition unchanged. After omitting one redundant
overlap equation, the gradients of the remaining equations are the
commutators $i[p,q_j]$. When their Gram matrix is invertible, it gives
a right inverse for the derivative of the overlap map. We will also
use this construction on the factor summands of a general tracial von Neumann
algebra.

\begin{lemma}
\label{lem:finite-coordinate-overlap-inverse}
Let $(N,\tau)$ be a tracial von Neumann algebra, let $\mathbf p,\mathbf q$ be partitions of $1$
of sizes $a,b$, allowing zero projections, and let $z\in Z(N)$ be a central projection.  Fix a nonempty set $\mathcal J$ of index pairs, let
$d:=|\mathcal J|$, and set
\[
g_{ij}:=z\,i[p_i,q_j],
\qquad
G:=\bigl(\tau(g_{ij}g_{kl})\bigr)_{(i,j),(k,l)\in\mathcal J}.
\]
Suppose $G\ge\gamma I$ for some $\gamma>0$.  For real targets $v_{ij}$, define
\[
\delta:=\bigl(\tau(p_iq_j)-v_{ij}\bigr)_{(i,j)\in\mathcal J},
\qquad
K:=\frac{2\sqrt d}{\gamma},
\qquad
C:=4\sqrt d K^2.
\]
If $|\delta|\le(4C)^{-1}$, there is $h=h^*\in zN$ such that
\[
\tau(e^{ih}p_ie^{-ih}q_j)=v_{ij}\quad((i,j)\in\mathcal J),
\qquad
\norm{h}_\infty\le2K|\delta|,
\]
and
\[
\sum_i\norm{e^{ih}p_ie^{-ih}-p_i}_2\le4\sqrt a K|\delta|.
\]
\end{lemma}

\begin{proof}
For $h=h^*\in zN$, set
\[
\Phi(h)=\bigl(\tau(e^{ih}p_ie^{-ih}q_j)-v_{ij}\bigr)_{(i,j)\in\mathcal J},
\qquad
L(y)=\sum_{(i,j)\in\mathcal J}(G^{-1}y)_{ij}g_{ij}.
\]
Then $D\Phi(0)L=I$ and $\|L(y)\|_\infty\le K|y|$.
Differentiating the exponentials twice and using their unitarity gives
\[
\|D^2(e^{iH}p_ie^{-iH})[H_1,H_2]\|_\infty
\le4\|H_1\|_\infty\|H_2\|_\infty,
\]
so $\psi=\Phi\circ L=\delta+y+R(y)$ has $\|D^2\psi\|\le C$.
On the ball $|y|\le2|\delta|$, we have
$\abs{-\delta-R(y)}\le|\delta|+2C|\delta|^2\le2|\delta|$.
The map $y\mapsto-\delta-R(y)$ therefore maps this ball into itself
and has Lipschitz constant at most $2C|\delta|\le1/2$.
Taking its fixed point $y$, we obtain $h=L(y)$ with the asserted norm
bound. Finally,
\[
\sum_i\|e^{ih}p_ie^{-ih}-p_i\|_2
\le2\|h\|_\infty\sum_i\sqrt{\tau(p_i)}
\le4\sqrt a\,K|\delta|.
\]
\end{proof}

For a tracial von Neumann algebra $(M,\tau)$, let $\mathcal X(M)$ be the space of pairs consisting of a projection $p$ and a partition
$\mathbf q=(q_1,q_2,q_3)$ of $1$, allowing zero projections. We use the distance
\[
d_M((p;\mathbf q),(p';\mathbf q'))
:=2\norm{p-p'}_2+\sum_{j=1}^3\norm{q_j-q'_j}_2.
\]
The factor two counts both $p$ and $1-p$. Once the traces are fixed,
any two overlaps determine the third.

\begin{lemma}
\label{lem:binary-ternary-overlap-geometry}
Let $(M,\tau)$ be a tracial von Neumann algebra and $(p;\mathbf q)\in\mathcal X(M)$.  Set
\[
g_j:=i[p,q_j]\quad(j=2,3),
\qquad
w_{jk}:=\norm{q_jpq_k}_2^2\quad(j\ne k).
\]
The real Gram matrix of $g_2,g_3$ is
\[
\Gamma:=\bigl(\tau(g_jg_k)\bigr)_{j,k=2,3}
=2\begin{pmatrix}
w_{12}+w_{23}&-w_{23}\\
-w_{23}&w_{13}+w_{23}
\end{pmatrix},
\]
and
\[
\det\Gamma=4(w_{12}w_{13}+w_{12}w_{23}+w_{13}w_{23}).
\]
It is singular if and only if $[p,q_j]=0$ for at least one $j\in\{1,2,3\}$.  In that case,
the unital von Neumann subalgebra $N:=W^*(p,q_1,q_2,q_3)$ is a direct sum of von Neumann algebras generated by at most two projections. In particular, $N$ is a separable hyperfinite von Neumann algebra of type I
with matrix fibers of size at most two.
\end{lemma}

\begin{proof}
The subspaces $q_jMq_k$ are $L^2$-orthogonal and
$w_{jk}=w_{kj}\ge0$. Expanding $i[p,q_2]$ and $i[p,q_3]$ in these corners,
we obtain the stated Gram matrix and determinant. The determinant vanishes
exactly when at most one of $w_{12},w_{13},w_{23}$ is positive.
Equivalently, by faithfulness, there is an index $j$ with $q_jpq_k=0$ for both $k\ne j$, which means that $e=q_j$ commutes with $p$.

Such an $e$ is central in $N$. The corner $eN$ is generated by the one
projection $ep$, and $(1-e)N$ by the two projections $(1-e)p$ and any
$q_k$ with $k\ne j$, omitting zero corners. The remaining assertions
follow from Dixmier's description of two projections~\cite{dixmier-position-relative}, as in the proof of
Proposition~\ref{prop:two-projection-tracial-relation-stability}.
\end{proof}

In the setting of Lemma~\ref{lem:binary-ternary-overlap-geometry}, let
$r=1-q_1$. If $rpr=ar$ for some $a\in[0,1]$, then
\begin{equation}
\Gamma=2a(1-a)
\begin{pmatrix}
\tau(q_2)&0\\
0&\tau(q_3)
\end{pmatrix}.
\label{eq:scalar-compression-overlap-gram}
\end{equation}
In particular, $\Gamma$ is positive definite if $0<a<1$ and
$\tau(q_2),\tau(q_3)>0$.

If $\Gamma>0$, then by Lemma~\ref{lem:finite-coordinate-overlap-inverse} we
can prescribe sufficiently nearby values of
$\tau(pq_2),\tau(pq_3)$ by conjugating $p$ with a unitary close to $1$
in operator norm, while fixing $\mathbf q$ and all projection traces.
Here we apply the lemma with $N=M$, $z=1$, $\mathbf p=(p,1-p)$, and
$\mathcal J=\{(1,2),(1,3)\}$.

\begin{theorem}
\label{thm:binary-ternary-marginal-target-stability}
Let
\[
\Sigma_3:=\{\boldsymbol\beta\in[0,1]^3:\beta_1+\beta_2+\beta_3=1\},
\]
and fix a continuous map $T=(T_2,T_3):[0,1]\times\Sigma_3\to\R^2$. For a tracial von Neumann algebra $(M,\tau)$ and
$(p;\mathbf q)\in\mathcal X(M)$, define
\[
\delta_T(p;\mathbf q)
:=\max_{j=2,3}\abs{\tau(pq_j)-T_j(\tau(p),\tau(q_1),\tau(q_2),\tau(q_3))}.
\]
For every $\eta>0$, there is $\varepsilon=\varepsilon(T,\eta)>0$ such that, uniformly over all
${\rm II}_1$ factors $(M,\tau)$, $\delta_T(p;\mathbf q)<\varepsilon$ implies the existence of
$(p';\mathbf q')\in\mathcal X(M)$ satisfying
\[
\delta_T(p';\mathbf q')=0,
\qquad
d_M((p;\mathbf q),(p';\mathbf q'))<\eta.
\]
If the traces of $p$ and $\mathbf q$ are fixed in advance, we may also
require the corrected pair to have the same traces, with $\varepsilon$
depending also on these traces.
\end{theorem}

\begin{proof}
Suppose there is a counterexample sequence. After passing to a
subsequence, we may assume that its column traces and overlaps converge to some
$(\boldsymbol\beta,\mathbf c)\in\mathcal T_3$, where we index the coordinates of $\mathcal T_3$ by $1,2,3$. Set
$\alpha=\sum_jc_j$. Since $T$ is continuous, we have
$c_j=T_j(\alpha,\boldsymbol\beta)$ for $j=2,3$.
By Theorem~\ref{thm:diffuse-binary-table-stability}, the pairs can be
changed by a distance tending to zero so that their column traces
are exactly $\boldsymbol\beta$ and their overlaps are exactly
$\mathbf c$, uniformly over the varying factors. Their first
projection then has trace $\alpha$, so they satisfy the target
equations, a contradiction. When the original traces are prescribed,
they are also the limiting traces and hence are preserved by this
correction.
\end{proof}

Chifan, Ioana, and Kunnawalkam Elayavalli~\cite[Lemma~3.2]{chifan-ioana-kunnawalkam-elayavalli} proved the first assertion for $T_j(\alpha,\boldsymbol\beta)=\alpha\beta_j$. We sketch a second proof of the theorem, which explains the method of
Section~\ref{sec:fixed-ambient-beyond-factors}: where the commutator Gram
matrix is nondegenerate, a small conjugation corrects the overlaps, and
where it degenerates, the generated von Neumann algebra can be lifted. For
larger partitions, a degenerate Gram matrix no longer forces the generated
von Neumann algebra to be hyperfinite (see the
discussion before Lemma~\ref{lem:exact-markov-relabelling}). Take a
tracial ultraproduct, over a nonprincipal ultrafilter, of a counterexample sequence
$(p_m;\mathbf q_m)$, with limiting pair $(p;\mathbf q)$.
 By continuity of $T$, this pair satisfies the target
equations. Let $\Gamma$ be the Gram matrix of $i[p,q_2]$ and $i[p,q_3]$ there.

If $\Gamma>0$, let the targets at coordinate $m$ be $T_j$ evaluated at
the traces of $p_m$ and $\mathbf q_m$, for $j=2,3$. The overlap errors for these targets
tend to zero along the ultrafilter, and the coordinate Gram matrices tend to
$\Gamma$. By Lemma~\ref{lem:finite-coordinate-overlap-inverse}, we attain the
targets on an ultrafilter-large set by conjugating $p_m$ with unitaries
tending to $1$ in operator norm along the ultrafilter. These conjugations fix
$\mathbf q_m$ and all projection traces.

If $\Gamma$ is singular, then by
Lemma~\ref{lem:binary-ternary-overlap-geometry}, the unital von Neumann algebra generated
by the limiting projections is hyperfinite. Using the trace-preserving lifting argument in
the proof of Theorem~\ref{thm:hyperfinite-lifting-cofinal}, we then obtain coordinate representatives with the joint
moments of $(p;\mathbf q)$
at every index. Their $L^2$ distance from $(p_m;\mathbf q_m)$ tends
to zero along the ultrafilter, and if the traces of $p_m$ and $\mathbf q_m$ do not depend on $m$, the
representatives have the same traces. Both cases contradict the assumed
separation. In the singular case, the coordinate Gram matrices may
still be positive definite with least eigenvalues tending to zero.

\subsection{Further special cases}
\label{subsec:further-overlap-special-cases}

The results in this subsection are not used later. In a ${\rm II}_1$ factor, we perturb a projection,
keeping its overlaps with a partition, so that its commutators with all members of the partition but one have a
positive definite Gram matrix (Proposition~\ref{prop:binary-table-regularization}). In a diffuse algebra, we remove small
moments of a projection against finitely many centered self-adjoint elements by a small perturbation that
keeps its trace (Lemma~\ref{lem:conditioned-diffuse-balancing}). Overlap matrices whose support graph has
no four-cycle are stable in every diffuse algebra (Theorem~\ref{thm:four-cycle-free-table-stability}).

For a projection and a ternary partition,
we showed in Lemma~\ref{lem:binary-ternary-overlap-geometry} that a singular
commutator Gram matrix forces the generated unital von Neumann algebra
to be hyperfinite. For a partition of any finite size in a ${\rm II}_1$ factor, we show in the next proposition that the Gram matrix
can be made positive definite by perturbing the projection
while preserving all overlaps, provided that each
overlap lies strictly between zero and the corresponding partition trace.

\begin{proposition}
\label{prop:binary-table-regularization}
Let $(M,\tau)$ be a ${\rm II}_1$ factor, let $b\ge2$, and let $p$ be a projection and
$\mathbf q=(q_j)_{j=0}^{b-1}$ a partition of $1$ satisfying
\[
\beta_j:=\tau(q_j)>0,
\qquad
0<c_j:=\tau(pq_j)<\beta_j
\qquad(j<b).
\]
For every $\eta>0$, there is a projection $e\in M$ such that
\[
\norm{e-p}_2<\eta,
\qquad
\tau(eq_j)=c_j\quad(j<b),
\]
and
\[
\Gamma_e:=\bigl(\tau(i[e,q_j]\,i[e,q_k])\bigr)_{j,k=1}^{b-1}>0.
\]
More precisely, if $\beta_j\ge\beta_*>0$ and
$0<\delta\le c_j/\beta_j\le1-\delta$ for all $j$, then one may require
\[
\Gamma_e\ge\frac{s^2\delta^2\beta_*}{4(b-1)^2}I,
\qquad
s:=\min\{1/2,\eta^2\}.
\]
\end{proposition}

\begin{proof}
Let $\theta_j=c_j/\beta_j$ and $Y=\sum_j\theta_jq_j$.  Since $0<\theta_j<1$ for all $j$, there is
$\delta>0$ with $\delta\le Y\le1-\delta$; for the quantitative assertion, we take the $\delta$ given there.
For $0<s<1$, the positive contraction
\[
x_0:=(1-s)p+sY
\]
satisfies $\tau(x_0q_j)=c_j$ for every $j$ and $\tau(x_0)=\tau(p)$.
Moreover, $s\delta\le x_0\le1-s\delta$.

The Gram matrix will be positive definite if every off-diagonal block
$q_0eq_j$, $j>0$, is nonzero. We arrange this by adding off-diagonal terms to $x_0$
and preserving the trace pairings with these terms when we pass to a projection.
For each $j>0$, choose equal-trace subprojections of $q_0$ and $q_j$ with trace
\[
d_j:=\tfrac12\min(\beta_0,\beta_j),
\]
and a partial isometry $v_j\in q_0Mq_j$ between them.
Choose $t_j\in\R$ so that
\[
h_j:=e^{it_j}v_j+e^{-it_j}v_j^*
\qquad\text{satisfies}\qquad
\tau(ph_j)=0.
\]
The phase exists because $\tau(ph_j)=2\Ree(e^{it_j}\tau(pv_j))$.
Computing in the corners $q_jMq_k$ gives
\[
\norm{h_j}_\infty=1,\qquad
\tau(q_lh_j)=0,\qquad
\tau(h_jh_k)=0\ (j\ne k),\qquad
\tau(h_j^2)=2d_j.
\]
Note that the subprojections of $q_0$ need not be mutually orthogonal.
However, for $j\ne k$, any nonzero product contributing to $h_jh_k$ lies in
an
off-diagonal corner $q_jMq_k$ and has trace zero.

Set $\kappa=s\delta/[2(b-1)]$.  Since
$\norm{\kappa\sum_{j>0}h_j}_\infty\le s\delta/2$, the operator
\[
x:=x_0+\kappa\sum_{j>0}h_j
\]
is a positive contraction. We apply
Lemma~\ref{lem:finite-functional-purification} to $x$, preserving its trace
pairings with $p$, the $q_j$, and the $h_j$. Then the resulting projection $e$
has the prescribed overlaps.
Since $\sum_jq_j=1$, we also have $\tau(e)=\tau(p)$.
Using $\tau(ph_j)=0$, we obtain
\[
\begin{aligned}
\norm{e-p}_2^2
&=2\tau(p)-2\tau(xp)\\
&=2s\sum_j\beta_j\theta_j(1-\theta_j)
\le s/2.
\end{aligned}
\]
Thus $s=\min\{1/2,\eta^2\}$ gives the strict distance bound.  Moreover,
\[
\tau(eh_j)=\tau(xh_j)=2\kappa d_j.
\]
By Cauchy--Schwarz,
\[
2\kappa d_j
\le2\norm{q_0eq_j}_2\norm{v_j}_2
=2\norm{q_0eq_j}_2\sqrt{d_j},
\]
so $\norm{q_0eq_j}_2^2\ge\kappa^2d_j$.

For real coefficients $a_j$, set $a_0=0$. Expanding the commutators
in the mutually orthogonal corners $q_jMq_k$ gives
\[
\norm{\sum_{j>0}a_ji[e,q_j]}_2^2
=2\sum_{j<k}\norm{q_jeq_k}_2^2(a_j-a_k)^2
\ge2\kappa^2\min_{j>0}d_j\sum_{j>0}a_j^2.
\]
This proves positivity.  Since $d_j\ge\beta_*/2$, we also obtain the stated uniform lower bound.
\end{proof}

Note that the strict overlap inequalities are necessary. Indeed, if $c_j=0$ or
$c_j=\beta_j$, every projection $e$ with the prescribed overlaps
commutes with $q_j$, and $\Gamma_e$ is singular.

\begin{lemma}
\label{lem:conditioned-diffuse-balancing}
Let $(D,\tau)$ be a diffuse tracial von Neumann algebra.  Let $B_1,\ldots,B_d\in D$ be centered self-adjoint
elements with $\norm{B_j}_\infty\le L$ and
\[
G:=\bigl(\tau(B_jB_k)\bigr)_{j,k=1}^d\ge\gamma I,
\]
where $d\ge1$, $L>0$, and $\gamma>0$.  Let $p$ be a projection with
$\alpha:=\tau(p)\in[\delta,1-\delta]$, where $0<\delta\le1/2$, and set
\[
m:=\bigl(\tau(pB_j)\bigr)_{j=1}^d,
\qquad
K:=\frac{L\sqrt d}{\gamma}.
\]
If $\abs{m}\le\delta/(4K)$, there is a projection $e\in D$ such that
\[
\tau(e)=\alpha,\qquad
\tau(eB_j)=0\quad(1\le j\le d),\qquad
\norm{e-p}_2^2\le2K(2/\delta+1)\abs{m}.
\]
\end{lemma}

\begin{proof}
If $m=0$, take $e=p$.  Otherwise set
\[
w:=\sum_j(G^{-1}m)_jB_j,
\qquad
s:=\frac{2K\abs{m}}{\delta}\le\frac12.
\]
Then $\tau(w)=0$, $\tau(wB_j)=m_j$, and $\norm{w}_\infty\le K\abs{m}=s\delta/2$.
The operator
\[
x:=(1-s)p+s\alpha-(1-s)w
\]
is a positive contraction, since $(1-s)p+s\alpha$ lies between
$s\delta$ and $1-s\delta$, and $\|(1-s)w\|_\infty\le s\delta/2$.
Moreover, $\tau(x)=\alpha$ and $\tau(xB_j)=0$.
We apply Lemma~\ref{lem:finite-functional-purification} to $x$, preserving
its trace and its trace pairings with $B_1,\ldots,B_d,p$. Then the resulting
projection $e$ satisfies
\[
\begin{aligned}
\norm{e-p}_2^2
&=2(\alpha-\tau(xp))\\
&\le2\norm{x-p}_1\\
&\le2(s+K\abs{m})
=2K(2/\delta+1)\abs{m}.
\end{aligned}
\]
\end{proof}

Diffuseness cannot be omitted from
Lemma~\ref{lem:conditioned-diffuse-balancing}. In $M_3(\C)$ with
normalized trace, take $p=e_{11}$, the six
centered self-adjoint matrices
\[
e_{jk}+e_{kj},\qquad i(e_{jk}-e_{kj})
\qquad(1\le j<k\le3),
\]
and $B_\epsilon=\operatorname{diag}(\epsilon,1,-1-\epsilon)$, where $0<\epsilon<1/2$ tends to zero.
The first six are pairwise trace-orthogonal with squared $L^2$-norm $2/3$, and are orthogonal to $B_\epsilon$, while
\[
\tau(B_\epsilon^2)=\frac{\epsilon^2+1+(1+\epsilon)^2}{3}\ge\frac23.
\]
All seven matrices have operator norm at most $3/2$, and their real Gram matrix is at least $(2/3)I$.
The moments against $p$ are zero except for $\tau(pB_\epsilon)=\epsilon/3$.
A projection of trace $1/3$ with zero moment against each of the seven matrices would have all off-diagonal entries zero, hence would be one of the
$e_{jj}$.  However, none of the $e_{jj}$ has zero moment against $B_\epsilon$.

Suppose that any two distinct columns of the prescribed overlap matrix
have nonzero entries in at most one common row, i.e., its bipartite support graph has no four-cycle.
Then any partitions realizing the matrix commute, and we prove stability by lifting their common
refinement and adjusting the traces of its members. This condition is not necessary: when one partition
is binary, Theorem~\ref{thm:diffuse-binary-table-stability} gives stability for every support.

\begin{theorem}
\label{thm:four-cycle-free-table-stability}
Fix finite $a,b\ge1$ and a bipartite graph
$G\subseteq\{1,\ldots,a\}\times\{1,\ldots,b\}$ containing no four-cycle.  For every $\eta>0$, there is
$\varepsilon=\varepsilon(G,\eta)>0$ such that the following holds uniformly over all diffuse tracial von Neumann
algebras $(D,\tau_D)$ and all matrices $C\in M_{a\times b}(\R)$
with nonnegative entries summing to $1$ and $C_{ij}=0$ for $(i,j)\notin G$.
If $(\mathbf p,\mathbf q)\in\mathcal P_{a,b}(D)$ satisfies
\[
\max_{i,j}\abs{\tau_D(p_iq_j)-C_{ij}}<\varepsilon,
\]
then there is a commuting pair $(\mathbf p',\mathbf q')\in\mathcal P_{a,b}(D)$ with
\[
\tau_D(p'_iq'_j)=C_{ij}\qquad(i,j),
\qquad
d_2((\mathbf p,\mathbf q),(\mathbf p',\mathbf q'))<\eta.
\]

\end{theorem}

\begin{proof}
First suppose $p_iq_j=0$ off $G$. For distinct columns $j,k$, at most
one term in $\sum_iq_jp_iq_k=q_jq_k=0$ can be nonzero, since $G$ has no
four-cycle. That term also vanishes. Thus
$p_i=\sum_jq_jp_iq_j$, so all $p_i$ commute with all $q_j$.

Now suppose that uniformity fails, and take diffuse tracial von Neumann algebras $(D_m,\tau_m)$, targets $C^{(m)}$
vanishing off $G$, and pairs with overlap errors tending to zero but a fixed
positive distance from every commuting realization of $C^{(m)}$.
Passing to a subsequence, we may assume that $C^{(m)}\to C$. We then take an ultraproduct over
a nonprincipal ultrafilter $\cU$.
The limit has $p_iq_j=0$ off $G$, since
$\norm{p_iq_j}_2^2=\tau(p_iq_j)=0$ there. The products $r_{ij}=p_iq_j$
therefore form a common projection refinement, with zero projections retained. We lift $(r_{ij})$ to partitions $(r_{ij,m})$, also allowing zero projections,
by Lemma~\ref{lem:nested-lifts-general}. Their row and column sums
$\overline p_{i,m},\overline q_{j,m}$ approach the original pairs along $\cU$.

Write $c_{ij,m}=\tau_m(r_{ij,m})$ and
$s_m=\frac12\sum_{i,j}|c_{ij,m}-C_{ij}^{(m)}|\to_{\cU}0$.
By diffuseness, we can transfer subprojections of total trace $s_m$ from
those $r_{ij,m}$ with $c_{ij,m}>C_{ij}^{(m)}$ to those with
$c_{ij,m}<C_{ij}^{(m)}$. The resulting partition
$(r'_{ij,m})$ of $1$, allowing zero projections, has traces $C_{ij}^{(m)}$ and is unchanged outside a
projection $z_m$ of trace $s_m$. Its row and column sums
$p'_{i,m},q'_{j,m}$ commute, and
\[
d_2((\overline{\mathbf p}_m,\overline{\mathbf q}_m),
 (\mathbf p'_m,\mathbf q'_m))\le(a+b)\sqrt{s_m}\longrightarrow_{\cU}0.
\]
This contradicts the separation.
\end{proof}

\section{Independence in a fixed tracial von Neumann algebra}
\label{sec:fixed-ambient-beyond-factors}

We now keep both the von Neumann algebra and its trace fixed and ask whether
pairs of partitions of $1$ with small covariance defects can be perturbed
to independent pairs.

Central decomposition presents a difficulty even when all summands are
factors. The global overlaps are weighted sums of the overlaps on the
summands, but the product of the global traces of two projections need not be the weighted sum of the products of their traces on the summands. Independence can therefore result from
cancellation between summands. For example, let $N=M_1\oplus M_2$ be a
direct sum of ${\rm II}_1$ factors with trace weights $1/2$, and in each
summand take two projections of trace $1/2$, with overlap $1/2$ in $M_1$
and $0$ in $M_2$. The global binary partitions are independent, although
the covariances on the two summands are $1/4$ and $-1/4$. Correcting the two summands separately would move a pair that already
satisfies the global independence equations exactly.
We therefore prescribe the global overlaps instead. We split projections to control the
commutator Gram matrices, combine the resulting estimates with the
central trace weights, and then apply
Lemma~\ref{lem:finite-coordinate-overlap-inverse}.

For a tracial von Neumann algebra $(N,\tau)$, we say that \emph{partition stability} $\operatorname{PS}(N;a,b)$ holds if every sequence
$(\mathbf p^{(m)},\mathbf q^{(m)})$ in $\mathcal P_{a,b}(N)$ satisfying
\[
\max_{i,j}\abs{\Delta_{ij}(\mathbf p^{(m)},\mathbf q^{(m)})}\longrightarrow0
\]
also satisfies
\begin{equation}
d_2((\mathbf p^{(m)},\mathbf q^{(m)}),Z_{a,b}(N))\longrightarrow0.
\label{eq:fixed-ambient-sequential-stability}
\end{equation}

Equivalently, for fixed $N,a,b$ and every $\eta>0$, there is $\varepsilon=\varepsilon(N,a,b,\eta)>0$, which we call a \emph{modulus} of $\operatorname{PS}(N;a,b)$, such that
\begin{equation}
\max_{i,j}\abs{\Delta_{ij}(\mathbf p,\mathbf q)}<\varepsilon
\quad\Longrightarrow\quad
d_2((\mathbf p,\mathbf q),Z_{a,b}(N))<\eta.
\label{eq:fixed-ambient-stability-modulus}
\end{equation}

For a nonprincipal ultrafilter $\cU$ on $\N$, write $\operatorname{P30}(N,\cU)$ for the assertion that every
pair of unitaries $u,v\in N^\cU$ with $W^*(u)\perp W^*(v)$ has unitary representatives whose generated unital von Neumann algebras are orthogonal at
every coordinate. For a fixed tracial von Neumann algebra $N$, the tracial version of Problem~30, as in the Introduction, asks whether $\operatorname{P30}(N,\cU)$ holds for every nonprincipal ultrafilter $\cU$ on $\N$.

For a tracial von Neumann algebra $(D,\tau_D)$ and $0<\lambda<1$, write
$\C\oplus_\lambda D$ for $\C\oplus D$
with trace $(z,x)\mapsto\lambda z+(1-\lambda)\tau_D(x)$.

\begin{corollary}
\label{cor:diffuse-scalar-problem-30}
Let $D$ be any diffuse tracial von Neumann algebra.  Then $\operatorname{PS}(D;a,b)$ holds for every finite $a,b$,
as does $\operatorname{PS}(\C\oplus_\lambda D;a,b)$ for every fixed scalar weight $0<\lambda<1$.
For either choice of ambient von Neumann algebra $N$, every orthogonal pair
of separable abelian unital von Neumann subalgebras $A,B\subseteq N^\cU$, where $\cU$
is a nonprincipal ultrafilter on $\N$, admits finite-dimensional
abelian unital $^*$-subalgebras $C_k,D_k\subseteq N$ such that
$A\subseteq\prod_\cU C_k$, $B\subseteq\prod_\cU D_k$, and
$C_k\perp D_k$ for every $k$. In particular, every independent pair
of unitaries in $N^\cU$ has unitary representatives whose generated unital von Neumann algebras
are orthogonal at every coordinate.
\end{corollary}

We give the proof after
Theorem~\ref{thm:independent-atomic-central-assembly}.

\begin{proposition}
\label{prop:fixed-ambient-stability-equivalence}
For a fixed tracial von Neumann algebra $N$, the following are equivalent.
\begin{enumerate}[label=\textup{(\roman*)}]
\item $\operatorname{PS}(N;a,b)$ holds for every finite $a,b$.
\item For every nonprincipal ultrafilter $\cU$ on $\N$ and every orthogonal pair of separable abelian unital von Neumann
algebras $A,B\subseteq N^\cU$, there are finite-dimensional abelian unital $^*$-subalgebras $C_m,D_m\subseteq N$, orthogonal for every
$m$, such that
\[
A\subseteq\prod_\cU C_m,
\qquad
B\subseteq\prod_\cU D_m.
\]
\item $\operatorname{P30}(N,\cU)$ holds for every nonprincipal ultrafilter $\cU$ on $\N$.
\item $\operatorname{P30}(N,\cU_0)$ holds for one nonprincipal ultrafilter $\cU_0$ on $\N$.
\end{enumerate}
For fixed $N$, the tracial version of Problem~30 is therefore equivalent to stability of independence for pairs of partitions of $1$ in $N$, allowing zero projections.
\end{proposition}

\begin{proof}
The implication \textup{(i)}$\Rightarrow$\textup{(ii)} follows from
Proposition~\ref{prop:higher-order-local-global} with $n=1$ and
$\mathfrak C=\{N\}$, since the defect at $n=1$ is the covariance
defect defining $\operatorname{PS}(N;a,b)$. For \textup{(ii)}$\Rightarrow$\textup{(iii)},
we apply \textup{(ii)} to the unital von Neumann subalgebras $W^*(u),W^*(v)$ and round representatives of $u,v$ to
unitaries inside the coordinate abelian von Neumann algebras, as in
the proof of Corollary~\ref{cor:fixed-factor-problem-30}. The implication \textup{(iii)}$\Rightarrow$\textup{(iv)}
is immediate.

Assume \textup{(iv)} and suppose $\operatorname{PS}(N;a,b)$ fails. Then we
can choose pairs $(\mathbf p^{(m)},\mathbf q^{(m)})$ in $\mathcal P_{a,b}(N)$ and $\rho>0$ with
\begin{equation}
\max_{i,j}|\Delta_{ij}(\mathbf p^{(m)},\mathbf q^{(m)})|\longrightarrow0,
\qquad
d_2((\mathbf p^{(m)},\mathbf q^{(m)}),Z_{a,b}(N))\ge\rho.
\label{eq:fixed-ambient-separated-sequence}
\end{equation}
For distinct $\lambda_i\in\mathbb T$ and distinct $\mu_j\in\mathbb T$, set
$x_m=\sum_i\lambda_i p_i^{(m)}$ and $y_m=\sum_j\mu_j q_j^{(m)}$.
Their classes $x,y\in N^{\cU_0}$ are independent. By \textup{(iv)}, they
have unitary representatives $\widehat x_m,\widehat y_m$ whose generated
unital von Neumann algebras are orthogonal at every coordinate.

Applying Lemma~\ref{lem:nested-lifts-general} inside $\prod_{\cU_0}W^*(\widehat x_m)$ and
$\prod_{\cU_0}W^*(\widehat y_m)$, we obtain partitions in $W^*(\widehat x_m)$ and $W^*(\widehat y_m)$
representing the spectral projections of $x$ and $y$, respectively. These partitions are independent at every $m$ and approach
the original partitions in $L^2$ along $\cU_0$. This contradicts
\eqref{eq:fixed-ambient-separated-sequence} and proves \textup{(i)}.
\end{proof}

\subsection{Central decomposition}

If
$N$ is purely atomic, its unit ball is compact in $L^2$. If $Z(N)$ is
diffuse, we make the covariance vanish by replacing both partitions on a
small central projection
(Lemma~\ref{lem:atomless-center-central-replacement}).
Theorem~\ref{thm:atomless-center-problem-30} is a special case of Corollary~\ref{cor:diffuse-scalar-problem-30}. We prove it directly to obtain the explicit estimate \eqref{eq:central-replacement-distance}. We combine compactness on the atomic part with a coupling argument on the part with
diffuse center to prove stability in every type I tracial von Neumann algebra.

\begin{fact}
\label{prop:purely-atomic-diagonal-ultrapower}
Let $(N,\tau)$ be a purely atomic tracial von Neumann algebra. Its
operator-norm unit ball is compact in $L^2$,
and the diagonal map $N\to N^\cU$ is onto for every ultrafilter $\cU$.
In particular, $\operatorname{P30}(N,\cU)$ holds for every nonprincipal
ultrafilter $\cU$ on $\N$.
\end{fact}

Indeed, $N$ is a finite or countable product of matrix algebras with
trace weights $\lambda_r>0$ satisfying $\sum_r\lambda_r=1$.
On the unit ball, the blocks outside a finite set $F$ contribute at most
$4\sum_{r\notin F}\lambda_r$ to the squared $L^2$ distance.
Hence, by compactness in finitely many matrix blocks, the unit ball is
compact in $L^2$. Every bounded representative has an $L^2$ ultralimit in
$N$, so its ultrapower class is diagonal. Constant representatives of
independent unitaries give $\operatorname{P30}(N,\cU)$.

For a general tracial von Neumann algebra $(N,\tau)$, let
$T:N\to Z(N)$ be the canonical center-valued trace, normalized by
$T(z)=z$ for $z\in Z(N)$ and $\tau\circ T=\tau$.
When $Z(N)$ is diffuse, we can change a pair of partitions on a central
projection and leave them fixed on its complement. We choose the trace
of that projection so that the overlaps required on it are nonnegative.

\begin{lemma}
\label{lem:atomless-center-central-replacement}
Let $(N,\tau)$ be a tracial von Neumann algebra and assume that $Z(N)$ is diffuse. Let
$(\mathbf p,\mathbf q)\in\mathcal P_{a,b}(N)$ consist of nonzero projections, and write
\[
\begin{aligned}
\alpha_i&:=\tau(p_i)>0,
&
\beta_j&:=\tau(q_j)>0,\\
C_{ij}&:=\tau(p_iq_j),
&
\Delta_{ij}&:=C_{ij}-\alpha_i\beta_j.
\end{aligned}
\]
If every $\Delta_{ij}=0$, set $r=0$ and $t=0$. Otherwise set
\begin{equation}
r:=\max_{i,j}\frac{\Delta_{ij}^{+}}{\alpha_i\beta_j},
\qquad
t:=\frac{r}{1+r}.
\label{eq:central-replacement-weight}
\end{equation}
There are independent partitions $\mathbf p',\mathbf q'$ with the same scalar projection traces and
\begin{equation}
\sum_i\norm{p_i-p_i'}_2\le\sqrt{2at},
\qquad
\sum_j\norm{q_j-q_j'}_2\le\sqrt{2bt}.
\label{eq:central-replacement-distance}
\end{equation}
\end{lemma}

\begin{proof}
If $\Delta=0$, keep the original pair. Otherwise set
$c_{ij}=T(q_jp_iq_j)=T(p_iq_j)\ge0$. The vector measure
\begin{equation}
\nu(e)=(\tau(ec_{ij}))_{i,j}\qquad(e\text{ a central projection})
\label{eq:central-overlap-vector-measure}
\end{equation}
is atomless.
By the convexity of the range of an atomless vector measure~\cite[Theorem~I]{lyapunov-convexity}, there is a central projection $z$ with
\begin{equation}
\tau(zp_iq_j)=tC_{ij}\qquad(1\le i\le a,\ 1\le j\le b).
\label{eq:central-proportional-overlaps}
\end{equation}
Summing gives
\begin{equation}
\tau(z)=t,\qquad \tau(zp_i)=t\alpha_i,\qquad
\tau(zq_j)=t\beta_j.
\label{eq:central-proportional-margins}
\end{equation}
The residual matrix
\begin{equation}
H_{ij}=\alpha_i\beta_j-(1-t)C_{ij}
=t\alpha_i\beta_j-(1-t)\Delta_{ij}
\label{eq:central-residual-table}
\end{equation}
is nonnegative because $(1-t)\Delta_{ij}^{+}\le(1-t)r\alpha_i\beta_j
=t\alpha_i\beta_j$, and its row and column sums are
\begin{equation}
\sum_jH_{ij}=t\alpha_i,\qquad \sum_iH_{ij}=t\beta_j,
\qquad \sum_{i,j}H_{ij}=t.
\label{eq:central-residual-margins}
\end{equation}
Choose a central partition $z=\sum_{i,j}z_{ij}$, allowing zero projections, with
$\tau(z_{ij})=H_{ij}$ and set
\begin{equation}
p_i'=(1-z)p_i+\sum_jz_{ij},\qquad
q_j'=(1-z)q_j+\sum_iz_{ij}.
\label{eq:central-replacement-partitions}
\end{equation}
These are partitions with the original projection traces and
$\tau(p_i'q_j')=(1-t)C_{ij}+H_{ij}=\alpha_i\beta_j$.
Since $zp_i$ and $\sum_jz_{ij}$ both have trace $t\alpha_i$, we have
$\norm{p_i-p_i'}_2^2\le2t\alpha_i$, and similarly
$\norm{q_j-q_j'}_2^2\le2t\beta_j$.
By Cauchy--Schwarz, we obtain \eqref{eq:central-replacement-distance}.
\end{proof}

\begin{theorem}
\label{thm:atomless-center-problem-30}
Let $(N,\tau)$ be a tracial von Neumann algebra with diffuse center.
Then $\operatorname{PS}(N;a,b)$ holds for every finite $a,b$. Equivalently, if $\cU$ is a nonprincipal ultrafilter on $\N$ and
$A,B\subseteq N^\cU$ are orthogonal separable abelian unital von Neumann subalgebras, there
are finite-dimensional abelian unital $^*$-subalgebras $C_k,D_k\subseteq N$ with
$A\subseteq\prod_\cU C_k$, $B\subseteq\prod_\cU D_k$, and
$C_k\perp D_k$ for every $k$. In particular, $\operatorname{P30}(N,\cU)$ holds.
\end{theorem}

\begin{proof}
Suppose that $\operatorname{PS}(N;a,b)$ fails. Then we can choose a sequence
of partition pairs whose covariance defects tend to zero and whose
distances from $Z_{a,b}(N)$ are bounded below by a positive constant.
Passing to a subsequence, we may assume that the projection traces converge:
\begin{equation}
\tau(p_i^{(m)})\longrightarrow\alpha_i,\qquad
\tau(q_j^{(m)})\longrightarrow\beta_j.
\label{eq:atomless-center-limiting-margins}
\end{equation}
In each partition, we add all projections whose traces tend to zero to
one projection with positive limiting trace, and discard the emptied
indices. If we put zero projections back at these indices, the total
$L^2$ distance from the original partitions is at most
\[
2\sum_{\alpha_i=0}\sqrt{\tau(p_i^{(m)})}
+2\sum_{\beta_j=0}\sqrt{\tau(q_j^{(m)})}\longrightarrow0.
\]
The coarsened pair $(\widetilde{\mathbf p}^{(m)},\widetilde{\mathbf q}^{(m)})$
still has covariance defects tending to zero, and all its projection traces are eventually
at least some $\delta>0$. Hence we can apply
Lemma~\ref{lem:atomless-center-central-replacement}, with
\[
r_m\le\delta^{-2}\max_{i,j}
 |\Delta_{ij}(\widetilde{\mathbf p}^{(m)},\widetilde{\mathbf q}^{(m)})|
\longrightarrow0,
\qquad t_m=\frac{r_m}{1+r_m}\longrightarrow0.
\]
The correction to an independent pair has total $L^2$ distance at
most $\sqrt{2at_m}+\sqrt{2bt_m}$. Putting zero projections back at the discarded indices gives partitions of the original sizes that are independent and approach the original partitions in $L^2$, a contradiction. This proves
\eqref{eq:fixed-ambient-sequential-stability}. The remaining conclusions
follow from Proposition~\ref{prop:fixed-ambient-stability-equivalence}.
\end{proof}

Examples include $M_d(L^\infty(X))$ and $L^\infty(X)\,\overline\otimes\,R$ for an atomless probability space $X$, and more generally measurable fields of finite factors over such a space.

\begin{proposition}
\label{prop:atomless-center-global-overlap-convexity}
Let $(N,\tau)$ be a tracial von Neumann algebra with diffuse center, and let $T:N\to Z(N)$ be its center-valued trace. For
$\varepsilon\in\{0,1\}$, let
\[
(\mathbf p^{(\varepsilon)},\mathbf q^{(\varepsilon)})
\in\mathcal P_{a,b}(N)
\]
and assume that
\begin{equation}
T(p_i^{(0)})=T(p_i^{(1)})\quad(1\le i\le a),
\qquad
T(q_j^{(0)})=T(q_j^{(1)})\quad(1\le j\le b).
\label{eq:center-valued-marginal-agreement}
\end{equation}
Write
\[
C_{ij}^{(\varepsilon)}:=\tau(p_i^{(\varepsilon)}q_j^{(\varepsilon)}).
\]
For every $\theta\in[0,1]$, there is a pair
$(\mathbf r,\mathbf s)\in\mathcal P_{a,b}(N)$ with the same center-valued projection traces as the two given pairs and
\begin{equation}
\tau(r_is_j)
=\theta C_{ij}^{(0)}+(1-\theta)C_{ij}^{(1)}
\qquad(1\le i\le a,\ 1\le j\le b).
\label{eq:center-valued-overlap-convexity}
\end{equation}
\end{proposition}

\begin{proof}
Let $d_{ij}=T(p_i^{(0)}q_j^{(0)})-T(p_i^{(1)}q_j^{(1)})\in Z(N)_{\rm sa}$.
The vector measure
$\nu(e)=(\tau(e),(\tau(ed_{ij}))_{i,j})$ on the center is atomless.
By the convexity of the range of $\nu$~\cite[Theorem~I]{lyapunov-convexity}, there is a central projection $z$
with
\begin{equation}
\nu(z)=\theta\nu(1).
\label{eq:central-lyapunov-splitting}
\end{equation}
Set
\[
r_i=zp_i^{(0)}+(1-z)p_i^{(1)},\qquad
s_j=zq_j^{(0)}+(1-z)q_j^{(1)}.
\]
By bimodularity of $T$ and \eqref{eq:center-valued-marginal-agreement}, the partitions $\mathbf r$ and $\mathbf s$ have the required center-valued traces.
Finally,
$\tau(r_is_j)=C_{ij}^{(1)}+\tau(zd_{ij})$
is the required convex combination.
\end{proof}

To treat type I tracial von Neumann algebras, we regard partitions of $1$ as measurable fields of partitions of $1$ in matrix algebras, allowing zero projections in the fibers. With the following lemma, we can change the distribution of such a field while controlling its $L^2$ distance from the original one.

\begin{lemma}
\label{lem:atomless-compact-coupling-realization}
Let $(X,\mu)$ be a finite atomless measure space, let $(K,d)$ be a compact metric space, and let
$f:X\to K$ be measurable. Let $\nu$ be a finite Borel measure on $K$ with
$\nu(K)=\mu(X)$, and let $\pi$ be a finite Borel measure on $K\times K$
with first marginal $f_*\mu$ and second marginal $\nu$.
For every $\varepsilon>0$, there is a measurable
$g:X\to K$ such that
\[
g_*\mu=\nu
\]
and
\begin{equation}
\int_X d(f(x),g(x))^2\,d\mu(x)
\le
\int_{K\times K}d(s,t)^2\,d\pi(s,t)+\varepsilon.
\label{eq:atomless-coupling-realization-cost}
\end{equation}
\end{lemma}

\begin{proof}
If $\mu(X)=0$, take $g=f$. Otherwise, on an atomless measurable set $E$, we can realize every finite Borel measure
on $K$ of mass $\mu(E)$ as the law of a measurable map.

Choose a finite Borel partition $K=\bigsqcup_r A_r$ and points $a_r\in A_r$
so that
\[
|d(s,t)^2-d(a_r,t)^2|<\frac{\varepsilon}{2\mu(X)}
\quad(s\in A_r,\ t\in K).
\]
Set $E_r=f^{-1}(A_r)$ and $\nu_r(B)=\pi(A_r\times B)$.
The first marginal of $\pi$ gives $\nu_r(K)=\mu(E_r)$, so we can realize
$\nu_r$ by a map on $E_r$ and combine these maps into $g$. Then
$g_*\mu=\sum_r\nu_r=\nu$.
On each strip $A_r\times K$, the couplings $(f,g)_*\mu$ and $\pi$ have
the same second marginal $\nu_r$, hence the same integral of $d(a_r,t)^2$.
By the displayed bound, their integrals of $d(s,t)^2$ differ
by at most $\varepsilon$, proving
\eqref{eq:atomless-coupling-realization-cost}.
\end{proof}

\begin{theorem}
\label{thm:finite-type-I-problem-30}
Let $(N,\tau)$ be a type I tracial von Neumann algebra. Then
$\operatorname{PS}(N;a,b)$ holds for every finite $a,b$. Consequently, if $\cU$ is a nonprincipal ultrafilter on $\N$ and
$A,B\subseteq N^\cU$ are orthogonal separable abelian unital von Neumann subalgebras, there
are finite-dimensional abelian unital $^*$-subalgebras $C_k,D_k\subseteq N$ with
$A\subseteq\prod_\cU C_k$, $B\subseteq\prod_\cU D_k$, and
$C_k\perp D_k$ for every $k$. In particular, $\operatorname{P30}(N,\cU)$ holds.
\end{theorem}

\begin{proof}
Fix $a,b$ and write the homogeneous central decomposition as
\begin{equation}
z_dN\cong M_d(L^\infty(X_d,\mu_d)),\qquad
\tau|_{z_dN}=\int_{X_d}\operatorname{tr}_d(\,\cdot\,)\,d\mu_d,
\qquad \sum_d\mu_d(X_d)=1.
\label{eq:type-i-central-decomposition}
\end{equation}
Split $X_d$ into its atomless part $X_d^0$ and atoms $x_{d,r}$ of masses
$w_{d,r}$. Let $K_d:=\mathcal P_{a,b}(M_d(\C))$, a compact space. For $k=(\mathbf P,\mathbf Q)$
and $k'=(\mathbf P',\mathbf Q')$ in $K_d$, set
\begin{equation}
\rho_d(k,k')^2=\sum_i\norm{P_i-P_i'}_{2,\operatorname{tr}_d}^2
 +\sum_j\norm{Q_j-Q_j'}_{2,\operatorname{tr}_d}^2.
\label{eq:type-i-configuration-metric}
\end{equation}
Its squared diameter is at most $2(a+b)$, uniformly in $d$.

Consider a sequence of pairs $(\mathbf p^{(m)},\mathbf q^{(m)})$ with covariance
defects tending to zero, and suppose
that a subsequence stays a positive distance from the independent
pairs. On the summand $z_dN$, write the pair as a measurable field
$k_{d,m}:X_d\to K_d$, and set
$\nu_{d,m}=(k_{d,m})_*(\mu_d|_{X_d^0})$.
By compactness, we may pass to a diagonal subsequence on which
\begin{equation}
\nu_{d,m}\longrightarrow\nu_d\text{ weakly for every }d,
\qquad k_{d,m}(x_{d,r})\longrightarrow k_{d,r}\text{ for every }d,r.
\label{eq:type-i-subsequence-limits}
\end{equation}
The projection traces and overlaps $\operatorname{tr}_d(P_i)$,
$\operatorname{tr}_d(Q_j)$, $\operatorname{tr}_d(P_iQ_j)$ are bounded
continuous functions on $K_d$. Since the central trace weights are summable,
the global projection traces and overlaps also converge, say to
$\alpha_i,\beta_j$ and $C_{ij}$. The covariance defects tend to zero, so
\begin{equation}
C_{ij}=\alpha_i\beta_j.
\label{eq:type-i-limiting-product-table}
\end{equation}

For finite measures $\mu,\nu$ of equal mass on $K_d$, write
\[
W_2(\mu,\nu)^2
:=\inf_{\pi}\int_{K_d\times K_d}\rho_d(s,t)^2\,d\pi(s,t),
\]
where the infimum is over couplings with marginals $\mu$ and $\nu$.
Since $K_d$ is compact, weak convergence of measures of fixed mass gives
$W_2(\nu_{d,m},\nu_d)\to0$.
By Lemma~\ref{lem:atomless-compact-coupling-realization}, there are maps
\begin{equation}
k'_{d,m}:X_d^0\to K_d,\qquad
(k'_{d,m})_*(\mu_d|_{X_d^0})=\nu_d,
\label{eq:type-i-realized-limiting-law}
\end{equation}
with
\begin{equation}
\int_{X_d^0}\rho_d(k_{d,m},k'_{d,m})^2\,d\mu_d
\le W_2(\nu_{d,m},\nu_d)^2+2^{-m-d};
\label{eq:type-i-realization-cost}
\end{equation}
on a zero-mass part, we keep the original map. On every atom, we set
$k'_{d,m}(x_{d,r})=k_{d,r}$.
These measurable $K_d$-valued fields define pairs $(\mathbf p'^{(m)},\mathbf q'^{(m)})$ of partitions of $1$
in $N$ with exactly the limiting global traces and overlaps.
By \eqref{eq:type-i-limiting-product-table},
\begin{equation}
(\mathbf p'^{(m)},\mathbf q'^{(m)})\in Z_{a,b}(N)\quad\text{for every }m.
\label{eq:type-i-exact-configurations}
\end{equation}
Finally, let
\begin{equation}
S_m=\sum_d\int_{X_d^0}\rho_d(k_{d,m},k'_{d,m})^2\,d\mu_d
 +\sum_{d,r}w_{d,r}\rho_d(k_{d,m}(x_{d,r}),k_{d,r})^2.
\label{eq:type-i-aggregate-squared-distance}
\end{equation}
Every term tends to zero, and by the uniform diameter bound and the
summability of the masses, we may apply dominated convergence. Thus
$S_m\to0$, and
\[
d_2((\mathbf p^{(m)},\mathbf q^{(m)}),
 (\mathbf p'^{(m)},\mathbf q'^{(m)}))\le\sqrt{a+b}\,S_m^{1/2}\longrightarrow0.
\]
This contradiction proves $\operatorname{PS}(N;a,b)$. The other conclusions follow from Proposition~\ref{prop:fixed-ambient-stability-equivalence}.
\end{proof}

\subsection{Binary partitions and scalar extensions}

For pairs consisting of a projection and a finite partition, we
use compactness on the matrix-factor summands to reduce stability
in any fixed tracial von Neumann algebra to Theorem~\ref{thm:diffuse-binary-table-stability}
on the diffuse summand. We use the set $\mathcal T_b$, the indexing $\mathbf q=(q_j)_{j=0}^{b-1}$, and the distance between pairs $(p;\mathbf q)$ introduced before Theorem~\ref{thm:diffuse-binary-table-stability}. A system of relations among the traces and overlaps is stable in $N$ if pairs satisfying it up to errors tending to zero are at distance tending to zero from pairs satisfying it exactly.

\begin{proposition}
\label{prop:fixed-ambient-binary-table-stability}
Fix a tracial von Neumann algebra $(N,\tau)$, $b\ge1$, and
$(\boldsymbol\beta^*,\mathbf c^*)\in\mathcal T_b$.  If projections $p_m$ and partitions $\mathbf q_m$ of $1$ in $N$
satisfy
\[
\tau(q_{j,m})\longrightarrow\beta_j^*,
\qquad
\tau(p_mq_{j,m})\longrightarrow c_j^*
\qquad(j<b),
\]
then their distance from pairs with exactly these traces and overlaps tends to zero.
Zero projections are allowed in both the given and corrected partitions.
Consequently, every fixed finite system of continuous relations among the traces and overlaps is stable
in $N$.
\end{proposition}

\begin{proof}
Write $N=N_{\rm at}\oplus D$, where $N_{\rm at}$ is the sum of the
matrix-factor central summands and $D$ is the diffuse summand. By
Fact~\ref{prop:purely-atomic-diagonal-ultrapower}, the space of
pairs consisting of a projection and a partition of size $b$
in $N_{\rm at}$ is compact in $L^2$.
Suppose a sequence stays a positive distance from the pairs with the
prescribed traces and overlaps. Passing to a subsequence, we may assume that
its restrictions to $N_{\rm at}$ converge to $(p^{\rm at};\mathbf q^{\rm at})$. Replacing those restrictions by their limit changes the pairs by a distance tending to zero. If $D=0$, the limiting pair already has
the prescribed traces and overlaps, giving a contradiction.

Otherwise, set $t=\tau(1_D)>0$. Let
$\widehat\beta_j^{\rm at},\widehat c_j^{\rm at}$ be the column traces
and overlaps of $(p^{\rm at};\mathbf q^{\rm at})$, measured with the original trace $\tau$; set them equal
to zero when $N_{\rm at}=0$. With respect to the
normalized trace $\tau|_D/t$, the column traces and overlaps of the restrictions to $D$
converge to
\[
\beta_j^D=\frac{\beta_j^*-\widehat\beta_j^{\rm at}}t,\qquad
c_j^D=\frac{c_j^*-\widehat c_j^{\rm at}}t.
\]
Since $(\boldsymbol\beta^D,\mathbf c^D)\in\mathcal T_b$, by Theorem~\ref{thm:diffuse-binary-table-stability} we can change the restrictions to $D$ by a distance tending to zero so that they have exactly these traces and overlaps.
Together with $(p^{\rm at};\mathbf q^{\rm at})$, they have the prescribed traces and overlaps in $N$, contradicting
the assumed separation.

For the last assertion, compactness of $\mathcal T_b$ gives a convergent
subsequence of the trace and overlap vectors of any counterexample
sequence. By continuity, the limiting traces and overlaps satisfy the given relations, and by the first assertion there are nearby pairs with precisely those traces and overlaps, a contradiction.
\end{proof}

\begin{corollary}
\label{cor:fixed-ambient-binary-problem-30}
For every fixed tracial von Neumann algebra $N$ and every finite $b$,
\[
\operatorname{PS}(N;2,b)
\qquad\text{and}\qquad
\operatorname{PS}(N;b,2)
\]
hold. If $\cU$ is a nonprincipal ultrafilter on $\N$ and $u,v\in N^\cU$ are unitaries with $W^*(u)\perp W^*(v)$, one of them
having at most two spectral values, then they have unitary representatives whose generated unital von Neumann algebras are orthogonal at
every coordinate.
\end{corollary}

\begin{proof}
We apply Proposition~\ref{prop:fixed-ambient-binary-table-stability} to
\[
\Phi(\boldsymbol\beta,\mathbf c)
:=\max_{j<b}\abs{c_j-\Bigl(\sum_kc_k\Bigr)\beta_j}.
\]
Since $\tau(p)=\sum_jc_j$ and the second covariance row is the negative of the first, this is the binary covariance defect.  Swapping the
partitions gives $\operatorname{PS}(N;b,2)$.

For the unitary assertion, the conclusion is immediate if either $u$ or $v$ is scalar. After
exchanging $u,v$ if necessary, write
$u=\zeta_0p+\zeta_1(1-p)$. We use the unital von Neumann subalgebra $W^*(p)$ at every stage on the first side and increasing
finite-dimensional unital $^*$-subalgebras with weakly dense union in $W^*(v)$
on the second. By Lemma~\ref{lem:nested-lifts-general}, we can lift their
minimal projections while preserving every refinement.
At stage $r$, the partitions have sizes $(2,b_r)$ and their covariance
defect tends to zero along $\cU$. By the stability just proved, there are
independent partitions at total distance less than $2^{-r}$ on a $\cU$-large
set.
Applying Lemma~\ref{lem:cofinal-two-sided-labelled-assembly},
we then make a diagonal choice of stages tending to infinity along $\cU$. The resulting finite-dimensional abelian unital $^*$-subalgebras
are orthogonal at every coordinate, and
their respective ultraproducts contain $W^*(u)$ and $W^*(v)$.

At the selected stage, write $\widetilde p_m$ for the corrected
projection corresponding to $p$ at coordinate $m$, and take
$\zeta_0\widetilde p_m+\zeta_1(1-\widetilde p_m)$
as the representative of $u$ at that coordinate. On the second side, inside the
coordinate abelian von Neumann algebra, take the polar part of a bounded
representative of $v$ and extend it by $1$ on its kernel. This gives a unitary representative of $v$, as in the proof of
Corollary~\ref{cor:fixed-factor-problem-30}. At coordinates where no
stage is selected, take both unitaries to be scalar.
\end{proof}

For the von Neumann algebra $\C\oplus M$, where $M$ is a $\mathrm{II}_1$
factor and the scalar summand has a fixed trace weight, the global
covariance equations specify the overlaps needed in $M$. The estimates below
are uniform over the factors, and for partitions of sizes $(2,2)$ and $(2,3)$ we can also
preserve projection traces fixed in advance.

\begin{proposition}
\label{prop:scalar-factor-affine-reduction}
Fix $0<\lambda<1$ and a \({\rm II}_1\) factor $(M,\tau_M)$, and define the tracial von Neumann algebra $(N_{\lambda,M},\tau_N)$ by
\begin{equation}
N_{\lambda,M}:=\C\oplus M,
\qquad
\tau_N(s\oplus x):=\lambda s+(1-\lambda)\tau_M(x).
\label{eq:scalar-factor-ambient}
\end{equation}
Fix finite $a,b$ and indices $i_0\in\{1,\ldots,a\}$ and $j_0\in\{1,\ldots,b\}$. For
$(\mathbf p,\mathbf q)\in\mathcal P_{a,b}(M)$, set
\[
\alpha_i:=\tau_M(p_i),
\qquad
\beta_j:=\tau_M(q_j),
\qquad
c_{ij}:=\tau_M(p_iq_j)
\]
and
\begin{equation}
\begin{aligned}
F_{ij}^{i_0,j_0}(\mathbf p,\mathbf q)
&:=
\lambda\delta_{i,i_0}\delta_{j,j_0}+(1-\lambda)c_{ij}\\
&\quad-
\bigl(\lambda\delta_{i,i_0}+(1-\lambda)\alpha_i\bigr)
\bigl(\lambda\delta_{j,j_0}+(1-\lambda)\beta_j\bigr).
\end{aligned}
\label{eq:scalar-factor-affine-covariance}
\end{equation}
Then $\operatorname{PS}(N_{\lambda,M};a,b)$ is equivalent to the following assertion for every choice of $i_0,j_0$:
if a sequence $(\mathbf p^{(m)},\mathbf q^{(m)})$ in $\mathcal P_{a,b}(M)$ satisfies
\[
\max_{i,j}\abs{F_{ij}^{i_0,j_0}(\mathbf p^{(m)},\mathbf q^{(m)})}\longrightarrow0,
\]
then the pairs of partitions in $M$ approach the common zero set of the
functions in \eqref{eq:scalar-factor-affine-covariance}
with respect to $d_2$ in $M$.
\end{proposition}

\begin{proof}
For $(\mathbf p,\mathbf q)\in\mathcal P_{a,b}(M)$, set
\[
P_i:=\delta_{i,i_0}\oplus p_i,\qquad Q_j:=\delta_{j,j_0}\oplus q_j.
\]
Then $(\mathbf P,\mathbf Q)\in\mathcal P_{a,b}(N_{\lambda,M})$ has covariance matrix
\eqref{eq:scalar-factor-affine-covariance}, and every pair in $\mathcal P_{a,b}(N_{\lambda,M})$ has this form
for exactly one $(i_0,j_0)$. Two pairs with different $(i_0,j_0)$ are at distance at least $2\sqrt\lambda$
in $d_2$. For two pairs with the same $(i_0,j_0)$, i.e., the same scalar components, the distance $d_2$ in
$N_{\lambda,M}$ is $\sqrt{1-\lambda}$ times the distance $d_2$ of their components in $M$.

Assume $\operatorname{PS}(N_{\lambda,M};a,b)$, and let $(\mathbf p^{(m)},\mathbf q^{(m)})$ be a sequence in
$\mathcal P_{a,b}(M)$ with
\[
\max_{i,j}\abs{F_{ij}^{i_0,j_0}(\mathbf p^{(m)},\mathbf q^{(m)})}\longrightarrow0.
\]
The covariances of $(\mathbf P^{(m)},\mathbf Q^{(m)})$ tend to zero, so there are independent pairs in
$\mathcal P_{a,b}(N_{\lambda,M})$ whose distance from them tends to zero. For large $m$, these have the same
$(i_0,j_0)$, and their components in $M$ are common zeros of the functions in
\eqref{eq:scalar-factor-affine-covariance} whose distance from $(\mathbf p^{(m)},\mathbf q^{(m)})$ tends to zero.

Conversely, a sequence contradicting $\operatorname{PS}(N_{\lambda,M};a,b)$
has a subsequence with constant $(i_0,j_0)$. By the assertion for this $(i_0,j_0)$, its components in $M$
approach the common zero set of the functions in \eqref{eq:scalar-factor-affine-covariance}. Taking direct sums
with $\delta_{i,i_0}$ and $\delta_{j,j_0}$, we obtain independent pairs in $N_{\lambda,M}$
approaching the subsequence, a contradiction.
\end{proof}

Solving the exact equations \eqref{eq:scalar-factor-affine-covariance} for the factor overlap matrix gives
\begin{equation}
c_{ij}
=
\frac{
\bigl(\lambda\delta_{i,i_0}+(1-\lambda)\alpha_i\bigr)
\bigl(\lambda\delta_{j,j_0}+(1-\lambda)\beta_j\bigr)
-\lambda\delta_{i,i_0}\delta_{j,j_0}
}{1-\lambda}.
\label{eq:scalar-factor-required-table}
\end{equation}
These overlaps are generally different from $\alpha_i\beta_j$, so
Theorem~\ref{thm:finite-partition-stability} does not apply directly
to the covariance equations in
Proposition~\ref{prop:scalar-factor-affine-reduction}.
There is a further obstruction to using the regularizing conjugation
from Proposition~\ref{prop:dense-regular-exact-locus}
in a general tracial von Neumann algebra $N$: nonzero projections $e,f$ with
disjoint central supports satisfy $eNf=0$.

In Proposition~\ref{prop:two-projection-tracial-relation-stability},
we gave a direct argument for two binary partitions, including higher
joint moments, and in Theorem~\ref{thm:binary-ternary-marginal-target-stability},
we treated the binary--ternary case. For arbitrary finite partitions in
$N_{\lambda,M}$, we instead split projections and regroup
their pieces.

\begin{corollary}
\label{cor:binary-scalar-factor-affine-stability}
Fix $0<\lambda<1$ and set, for projections in a ${\rm II}_1$ factor,
\[
F_\lambda(p,q)
:=
\lambda+(1-\lambda)\tau(pq)
-\bigl(\lambda+(1-\lambda)\tau(p)\bigr)
 \bigl(\lambda+(1-\lambda)\tau(q)\bigr).
\]
For every $\eta>0$, there is $\varepsilon=\varepsilon(\lambda,\eta)>0$ such that, for every ${\rm II}_1$ factor
$(M,\tau)$ and projections $p,q\in M$,
\[
\abs{F_\lambda(p,q)}<\varepsilon
\]
implies the existence of projections $p',q'\in M$ satisfying
\[
F_\lambda(p',q')=0,
\qquad
\norm{p-p'}_2+\norm{q-q'}_2<\eta.
\]
Consequently,
\[
\operatorname{PS}(N_{\lambda,M};2,2)
\]
holds for every ${\rm II}_1$ factor $M$, with a modulus which, for fixed $\lambda$, is uniform in $M$.
\end{corollary}

\begin{proof}
The first assertion follows by applying
Proposition~\ref{prop:two-projection-tracial-relation-stability} to the
moment list $(X,Y,XY)$ and the continuous defect $|F_\lambda|$.
For the partition assertion, number each partition so that the
projection with scalar component $1$ comes first. The global covariance matrix has zero row and column sums and
is determined by its $(1,1)$ entry $F_\lambda(p,q)$. We correct $p,q$ by the
first assertion and keep the scalar components. The
distance $d_2$ is
\[
2\sqrt{1-\lambda}
\bigl(\norm{p-p'}_2+\norm{q-q'}_2\bigr).
\]
The conclusion then follows from
Proposition~\ref{prop:scalar-factor-affine-reduction}, with the same
uniformity over factors.
\end{proof}

\begin{corollary}
\label{cor:fixed-margin-binary-affine-stability}
Fix $0<\lambda<1$ and $\alpha,\beta\in[0,1]$, and set
\[
h_\lambda(\alpha,\beta)
:=(1-\lambda)\alpha\beta+\lambda(\alpha+\beta-1).
\]
Assume $h_\lambda(\alpha,\beta)\ge0$.  For every $\eta>0$, there is
\[
\varepsilon=\varepsilon(\lambda,\alpha,\beta,\eta)>0
\]
such that, uniformly over all ${\rm II}_1$ factors $(M,\tau)$, projections $p,q\in M$ satisfying
\[
\tau(p)=\alpha,
\qquad
\tau(q)=\beta,
\qquad
\abs{F_\lambda(p,q)}<\varepsilon
\]
admit projections $p',q'\in M$ satisfying
\[
\tau(p')=\alpha,
\qquad
\tau(q')=\beta,
\qquad
F_\lambda(p',q')=0,
\qquad
\norm{p-p'}_2+\norm{q-q'}_2<\eta.
\]
\end{corollary}

\begin{proof}
Let $h=h_\lambda(\alpha,\beta)$. The hypothesis and the definition
of $h$ give
\[
\max\{0,\alpha+\beta-1\}\le h\le\min\{\alpha,\beta\},
\]
so every ${\rm II}_1$ factor contains projections with traces $\alpha,\beta$ and overlap $h$. We apply
Proposition~\ref{prop:two-projection-tracial-relation-stability} to
$(X,Y,XY)$ with the continuous defect
\[
\Phi(a,b,c)=\max\{|a-\alpha|,|b-\beta|,|c-h|\}.
\]
For the given projections,
$|F_\lambda(p,q)|=(1-\lambda)|\tau(pq)-h|$.
If $\varepsilon_0$ is the proposition's threshold, the result follows with
$\varepsilon=(1-\lambda)\varepsilon_0$.
\end{proof}

In the next two corollaries, $\mathcal X(M)$ and the distance $d_M$ are those introduced before Lemma~\ref{lem:binary-ternary-overlap-geometry}, and $\Sigma_3$ is the simplex of Theorem~\ref{thm:binary-ternary-marginal-target-stability}.

\begin{corollary}
\label{cor:two-by-three-scalar-factor-affine-stability}
Fix $0<\lambda<1$, let $t=1-\lambda$, and, for a tracial von Neumann algebra $(M,\tau)$ and $(p;\mathbf q)\in\mathcal X(M)$, set
\[
A:=\lambda+t\tau(p),
\qquad
D_j:=\tau(pq_j)-A\tau(q_j)\quad(j=2,3).
\]
For every $\eta>0$, there is $\varepsilon=\varepsilon(\lambda,\eta)>0$, uniformly over all
${\rm II}_1$ factors, such that any $(p;\mathbf q)\in\mathcal X(M)$ satisfying $\max\{\abs{D_2},\abs{D_3}\}<\varepsilon$ lies at distance below $\eta$ in $d_M$ from some $(p';\mathbf q')\in\mathcal X(M)$ for which both $D_j$ vanish exactly.  Consequently,
\[
\operatorname{PS}(N_{\lambda,M};2,3)
\qquad\text{and}\qquad
\operatorname{PS}(N_{\lambda,M};3,2)
\]
hold for every ${\rm II}_1$ factor $M$, with a modulus uniform in $M$ for fixed $\lambda$.
\end{corollary}

\begin{proof}
We apply Theorem~\ref{thm:binary-ternary-marginal-target-stability} with
\[
T_j(\alpha,\boldsymbol\beta)=(\lambda+t\alpha)\beta_j
\qquad(j=2,3).
\]
For the second assertion, number each partition so that the projection with scalar component
$1$ comes first. Then \eqref{eq:scalar-factor-affine-covariance} gives
$F_{1j}^{1,1}=tD_j$ for $j=2,3$. Since all row and column sums are zero,
these entries determine all the others, so independence is equivalent to
$D_2=D_3=0$.
For pairs with the same scalar components, the distance $d_2$ in
$N_{\lambda,M}$ is $\sqrt t$ times $d_M$.
Hence, by Proposition~\ref{prop:scalar-factor-affine-reduction}, we obtain stability
for partitions of sizes $(2,3)$, with a modulus uniform over $M$ for each
fixed $\lambda$. Exchanging the partitions gives stability for sizes $(3,2)$.
\end{proof}

\begin{corollary}
\label{cor:fixed-margin-two-by-three-affine-stability}
Fix $0<\lambda<1$, $\alpha\in[0,1]$, and
$\boldsymbol\beta\in\Sigma_3$. For a ${\rm II}_1$ factor $(M,\tau)$ and $(p;\mathbf q)\in\mathcal X(M)$,
consider the partitions
\[
(1\oplus p,\ 0\oplus(1-p)),\qquad(1\oplus q_1,\ 0\oplus q_2,\ 0\oplus q_3)
\]
of $1$ in $N_{\lambda,M}$, allowing zero projections. There is such a $(p;\mathbf q)$ in a ${\rm II}_1$ factor with $\tau(p)=\alpha$,
$\tau(q_j)=\beta_j$, and these two partitions independent if and only if
\[
h_\lambda(\alpha,\beta_1)
:=(1-\lambda)\alpha\beta_1+\lambda(\alpha+\beta_1-1)\ge0.
\]
Assume this condition. For every $\eta>0$, there is
$\varepsilon=\varepsilon(\lambda,\alpha,\boldsymbol\beta,\eta)>0$ such that, for every ${\rm II}_1$ factor
$(M,\tau)$, the following holds. If $(p;\mathbf q)\in\mathcal X(M)$ satisfies $\tau(p)=\alpha$ and
$\tau(q_j)=\beta_j$, and every covariance of its two partitions has absolute value less than $\varepsilon$,
then there is $(p';\mathbf q')\in\mathcal X(M)$ with the same traces, with independent partitions, and with
$d_M((p;\mathbf q),(p';\mathbf q'))<\eta$.
\end{corollary}

\begin{proof}
Let $t=1-\lambda$, $A=\lambda+t\alpha$, and $c_j=\tau(pq_j)$.  Independence of the two partitions requires
\[
c_1=h_\lambda(\alpha,\beta_1),\qquad c_2=A\beta_2,\qquad c_3=A\beta_3.
\]
Necessity follows from $c_1\ge0$. Conversely, $\sum_jc_j=\alpha$, and all six entries
$c_j,\beta_j-c_j$ are nonnegative, since
\[
\beta_1-c_1=(1-\alpha)(\lambda+t\beta_1),
\qquad
\beta_j-c_j=(1-A)\beta_j\quad(j=2,3).
\]
Because the factor is diffuse, we can choose six projections, possibly zero,
forming a partition of $1$ with traces $c_j$ and $\beta_j-c_j$, and their
row and column sums realize the overlap matrix. Then by the trace-preserving clause of
Theorem~\ref{thm:binary-ternary-marginal-target-stability}, with the target
map from the proof of
Corollary~\ref{cor:two-by-three-scalar-factor-affine-stability}, we
obtain the required perturbation.
Indeed, $F_{1j}^{1,1}=tD_j$ for $j=2,3$, so a sufficiently small
covariance defect gives the target defect required by that theorem.
\end{proof}

The hyperfiniteness argument in
Lemma~\ref{lem:binary-ternary-overlap-geometry} does not extend directly
to larger partitions. For a projection $p$ and a partition $(q_j)$, consider
the graph with an edge between distinct $j$ and $k$ when $q_jpq_k\ne0$.
When this graph is disconnected, a connected component may still contain three or more vertices,
and the von Neumann algebra generated by $p$ and the
corresponding projections $q_j$ need not be hyperfinite.

We instead split each projection and redistribute its pieces according
to a stochastic matrix. Using Theorem~\ref{thm:finite-order-free-positioning} inside the corners, we choose the pieces so that we can compute both the new overlaps and the fourth moments that occur in the commutator Gram matrix.

\begin{lemma}
\label{lem:exact-markov-relabelling}
Let $(M,\tau)$ be a ${\rm II}_1$ factor, and let
$(p_i)_{0\le i<a}$ and $(q_j)_{0\le j<b}$ be partitions of $1$.
Zero projections are allowed throughout this lemma. Set
\[
\alpha_i:=\tau(p_i),
\qquad
C_{ij}:=\tau(p_iq_j).
\]
Let $S=(s_{ki})_{0\le k,i<a}$ have nonnegative entries with
$\sum_k s_{ki}=1$ for every $i$. There are partitions
$(f_{ki})_{0\le k<a}$ of the projections $p_i$ such that, on setting $p'_k:=\sum_i f_{ki}$,
\[
\tau(f_{ki}q_j)=s_{ki}C_{ij},
\qquad
\tau(p'_k)=\sum_i s_{ki}\alpha_i,
\qquad
\tau(p'_kq_j)=\sum_i s_{ki}C_{ij}.
\]
With the same choice of projections, the following identity holds for every real array $(x_{kj})$.
Define
\[
B_k:=\sum_jx_{kj}q_j,
\qquad
\overline B_i:=\sum_k s_{ki}B_k,
\]
\[
V_i(B):=\tau(p_iB^2)-\frac{\tau(p_iB)^2}{\alpha_i}\quad(\alpha_i>0),
\]
and set $V_i=0$ when $\alpha_i=0$.  Then
\[
\norm{\sum_k i[p'_k,B_k]}_2^2
=\norm{\sum_i i[p_i,\overline B_i]}_2^2
  +\sum_{i,k,l}s_{ki}s_{li}V_i(B_k-B_l).
\] Moreover,
\[
\sum_k\norm{p'_k-p_k}_2^2=2-2\sum_i s_{ii}\alpha_i.
\]
In particular, if $s_{ii}\ge1-s$ for all $i$, where $0\le s\le1$, then
\[
\sum_k\norm{p'_k-p_k}_2\le\sqrt{2as}.
\]
\end{lemma}

\begin{proof}
Set $f_{ki}=0$ when $\alpha_i=0$, and construct the other partitions
successively in increasing order of $i$. In the ${\rm II}_1$ factor
$M_i=p_iMp_i$, with unit $p_i$ and trace $\tau_i=\tau/\alpha_i$, choose a partition
$(e_{ki})_k$ with $\tau_i(e_{ki})=s_{ki}$. Let $\mathbf Y_i$ consist
of the real and imaginary parts of
\[
\begin{array}{ll}
p_iq_jp_i, & \text{all }j,\\
p_iq_jp_{i'}q_hp_i, & i'>i,\ \text{all }j,h,\\
p_iq_jf_{li'}q_hp_i, & i'<i,\ \text{all }l,j,h.
\end{array}
\]
Each real or imaginary part in $\mathbf Y_i$ counts as one entry,
regardless of its degree in the original projections.
We apply Theorem~\ref{thm:finite-order-free-positioning}\textup{(ii)} at
$n=2$ in $M_i$, and set $f_{ki}=u_ie_{ki}u_i^*$ for its corner
unitary $u_i$.
By the degree-two identities, we have $\tau(f_{ki}q_j)=s_{ki}C_{ij}$.
For $\alpha_i>0$, the free fourth-moment formula in the normalized
corner, multiplied by $\alpha_i$, gives
\[
\begin{aligned}
\tau(f_{ki}q_jf_{li}q_h)
={}&s_{ki}s_{li}\tau(p_iq_jp_iq_h)\\
&+(\delta_{kl}s_{ki}-s_{ki}s_{li})\frac{C_{ij}C_{ih}}{\alpha_i}.
\end{aligned}
\]
For $i<i'$ with $\alpha_i\alpha_{i'}>0$, we apply the imposed identities
first in $M_{i'}$ and then in $M_i$, using cyclicity, to obtain
\[
\begin{aligned}
\tau(f_{ki}q_jf_{li'}q_h)
&=s_{li'}\tau(f_{ki}q_jp_{i'}q_h)\\
&=s_{ki}s_{li'}\tau(p_iq_jp_{i'}q_h).
\end{aligned}
\]
The case $i>i'$ follows by cyclicity.
These identities have
degree at most four in the partition and the entries of $\mathbf Y_i$,
so the theorem's degree-five conclusion suffices.

For any partition $(e_r)$ of $1$ and self-adjoint elements $(D_r)$,
\[
\norm{\sum_r i[e_r,D_r]}_2^2
=2\sum_r\tau(e_rD_r^2)-2\Ree\sum_{r,t}\tau(e_rD_re_tD_t).
\]
We apply this to $(f_{ki},B_k)_{k,i}$ and $(p_i,\overline B_i)_i$.
By the preceding formulas, the terms from distinct original corners cancel. The difference of the two squared norms is
\[
2\sum_i\left(\sum_k s_{ki}V_i(B_k)-V_i(\overline B_i)\right)
=\sum_{i,k,l}s_{ki}s_{li}V_i(B_k-B_l).
\]
Here $V_i$ is $\alpha_i$ times the variance for the state
$B\mapsto\tau(p_iB)/\alpha_i$ when $\alpha_i>0$, and zero otherwise.
Since $V_i$ is a quadratic form and $\sum_ks_{ki}=1$, the last equality holds.

Since $f_{ki}\le p_i$, the family $(p'_k)_k$ is a partition of $1$, $p'_kp_k=f_{kk}$, and
\[
\sum_k\norm{p'_k-p_k}_2^2=2-2\sum_k s_{kk}\alpha_k.
\]
The distance bound then follows from Cauchy--Schwarz.
\end{proof}

\begin{corollary}
\label{cor:markov-absorption}
In the setting of Lemma~\ref{lem:exact-markov-relabelling}, let $0\le\theta\le1$.  There is a partition
$(p'_i)_i$ of $1$, allowing zero projections, such that
\[
\tau(p'_iq_j)=\theta C_{ij}\quad(i>0),
\qquad
\sum_i\norm{p'_i-p_i}_2\le\sqrt{2a(1-\theta)}.
\]
For the same partition and every real array $(x_{ij})_{i,j>0}$, set $B_0=0$ and
$B_i=\sum_{j>0}x_{ij}q_j$.  Then
\[
\norm{\sum_{i>0}i[p'_i,B_i]}_2^2
=\theta^2\norm{\sum_{i>0}i[p_i,B_i]}_2^2
 +2\theta(1-\theta)\sum_{i>0}V_i(B_i).
\]
For every $i$ with $\alpha_i>0$ and every $B=\sum_{j>0}x_jq_j$ with real $x_j$, one has the variance bound
\[
V_i(B)\ge\frac{C_{i0}}{\alpha_i}\sum_{j>0}C_{ij}x_j^2.
\]
\end{corollary}

\begin{proof}
We apply Lemma~\ref{lem:exact-markov-relabelling} with $s_{00}=1$, $s_{ii}=\theta$ and $s_{0i}=1-\theta$ for $i>0$, and all other entries zero.  Then
$\overline B_0=0$ and $\overline B_i=\theta B_i$ for $i>0$.  The two ordered pairs $(0,i)$ and $(i,0)$ in the variance formula give
$2\theta(1-\theta)V_i(B_i)$, proving the identity.
The distance bound follows from $s_{ii}\ge\theta$.
If $\alpha_i>0$, then, for real $x_j$ with $x_0=0$, positivity and $\sum_jC_{ij}=\alpha_i$ give
\[
V_i\left(\sum_jx_jq_j\right)
=\frac1{2\alpha_i}\sum_{j,h}C_{ij}C_{ih}(x_j-x_h)^2
\ge\frac{C_{i0}}{\alpha_i}\sum_{j>0}C_{ij}x_j^2.
\]
\end{proof}

\begin{corollary}
\label{cor:all-type-scalar-factor-affine-stability}
Fix $0<\lambda<1$ and finite $a,b\ge1$.  For every $\eta>0$, there is
$\varepsilon=\varepsilon(\lambda,a,b,\eta)>0$ such that, for every
${\rm II}_1$ factor $M$, every choice of indices $i_0,j_0$ for the scalar
summand, and every $(\mathbf p,\mathbf q)\in\mathcal P_{a,b}(M)$,
\[
\max_{i,j}|F_{ij}^{i_0,j_0}(\mathbf p,\mathbf q)|<\varepsilon
\]
implies the existence of $(\mathbf p',\mathbf q')\in\mathcal P_{a,b}(M)$ with
\[
F_{ij}^{i_0,j_0}(\mathbf p',\mathbf q')=0\quad\text{for all }i,j,
\qquad
\sum_i\norm{p_i-p'_i}_2+\sum_j\norm{q_j-q'_j}_2<\eta.
\]
Here $F_{ij}^{i_0,j_0}$ is the covariance defined in \eqref{eq:scalar-factor-affine-covariance}.  Consequently
$\operatorname{PS}(N_{\lambda,M};a,b)$ holds with a modulus uniform in $M$ and in all trace vectors.
\end{corollary}

\begin{proof}
Suppose the first assertion fails.  Choose a sequence of factors and partition pairs with defects tending to zero and
distance at least $\eta$ from the solutions of the covariance equations.  Passing to a subsequence, we may assume that $(i_0,j_0)$ is constant, and we number both partitions from $0$ so that $i_0=j_0=0$.  Set $t:=1-\lambda$ and, in each coordinate,
\[
\alpha_i:=\tau_M(p_i),\quad \beta_j:=\tau_M(q_j),\quad c_{ij}:=\tau_M(p_iq_j),
\qquad B_j:=\lambda\delta_{j0}+t\beta_j.
\]
Passing to a further subsequence, we may assume that the projection traces converge.
In each partition, the projection with index $0$ has global trace at
least $\lambda$. We add to it every projection whose global trace tends to
zero, and temporarily omit the emptied indices. This changes the factor
partitions by $o(1)$ in $L^2$. After correction, we put zero projections
back at the omitted indices. If only one projection remains
in either partition, the resulting global partitions are independent.
Otherwise, for all retained $i,j>0$, both $\alpha_i$ and $B_j$ have
positive limits, and
\[
F_{ij}=t(c_{ij}-\alpha_iB_j)\quad(i>0,\ \text{all }j).
\]
In particular,
\[
\lim\frac{c_{i0}}{\alpha_i}=\lim B_0\ge\lambda,
\qquad
\lim c_{ij}=\lim(\alpha_iB_j)>0\quad(i,j>0).
\]

We fix $0<\theta<1$ so close to $1$ that $\sqrt{2a(1-\theta)}<\eta/3$ and apply
Corollary~\ref{cor:markov-absorption} in each original factor. Write
$\widehat p_i$ for the projections in the resulting first partition.
The second partition is unchanged. For $i>0$,
\[
\widehat\alpha_i=\theta\alpha_i,
\qquad
\widehat c_{ij}=\theta c_{ij},
\qquad
\widehat F_{ij}=\theta F_{ij}.
\]  The commutator identity and the variance bound in that corollary give
\[
\norm{\sum_{i,j>0}x_{ij}i[\widehat p_i,q_j]}_2^2
\ge2\theta(1-\theta)\sum_{i,j>0}\frac{c_{i0}}{\alpha_i}c_{ij}x_{ij}^2.
\]
By the preceding positive limits, these Gram matrices have a common lower bound $\gamma>0$ at all sufficiently large
coordinates.  We apply Lemma~\ref{lem:finite-coordinate-overlap-inverse} in
each factor with $z=1$, $\mathcal J=\{(i,j):i,j>0\}$,
and targets $v_{ij}:=\widehat\alpha_iB_j$.  Its residual vector tends to zero, so conjugation by $e^{ih}$, with
$\norm{h}_\infty\to0$, makes these covariances zero while preserving all projection traces.
Since all row and column sums are zero, the entire covariance matrix then vanishes.  This last correction changes the
partitions by $o(1)$.

Altogether the pairs move by less than $\eta$ for large indices, a contradiction. The last assertion follows
from Proposition~\ref{prop:scalar-factor-affine-reduction} and the distance identity in its proof.
\end{proof}

\subsection{Combining central summands}

In Corollary~\ref{cor:all-type-scalar-factor-affine-stability},
we treated one scalar atom and one ${\rm II}_1$ factor. We extend this to other fixed tracial von Neumann algebras by combining corrections on their central summands. On every
diffuse summand, we split and regroup the projections with the same
stochastic matrices, so that cancellation between summands is retained. We sum
the factor Gram matrices with the central trace weights and apply
Lemma~\ref{lem:finite-coordinate-overlap-inverse} in the whole von Neumann
algebra. The resulting unitary is $1$ off the ${\rm II}_1$-factor summands,
and conjugating the first partition by it corrects the selected overlaps.
We remove any remaining error by replacing the partitions on a small central
projection in the part with diffuse center.

For a tracial von Neumann algebra $N$, let $z_{\rm at}$ be the
sum of the minimal central projections $z$ for which $zN$ is finite
dimensional. Then $N_{\rm at}:=z_{\rm at}N$ is a finite or countable
direct sum of matrix algebras, and $(1-z_{\rm at})N$ is diffuse. Write
$\lambda=\tau(z_{\rm at})$ and, when $\lambda>0$, use the normalized trace
$\tau_{\rm at}=\lambda^{-1}\tau$ on $N_{\rm at}$. We write $\Delta^{N}_{ij}$ for the quantities $\Delta_{ij}$ of Section~\ref{sec:local-correction} computed in $(N,\tau)$, and $\Delta^{N_{\rm at}}_{ij}$ for those computed in $(N_{\rm at},\tau_{\rm at})$.

\begin{theorem}
\label{thm:independent-atomic-central-assembly}
Fix a tracial von Neumann algebra $(N,\tau)$ and finite $a,b$.  Suppose that $(\mathbf p^{(m)},\mathbf q^{(m)})\in\mathcal P_{a,b}(N)$ satisfy
\[
\max_{i,j}\abs{\Delta_{ij}^{N}(\mathbf p^{(m)},\mathbf q^{(m)})}\longrightarrow0
\]
and, if $\lambda>0$, also
\[
\max_{i,j}\abs{\Delta_{ij}^{N_{\rm at}}
(z_{\rm at}\mathbf p^{(m)},z_{\rm at}\mathbf q^{(m)})}\longrightarrow0.
\]
Then there are pairs in $\mathcal P_{a,b}(N)$ whose total $L^2(N)$
distance from the given pairs tends to zero, which are independent in $N$,
and, when $\lambda>0$, whose restrictions to $N_{\rm at}$ are independent for $\tau_{\rm at}$.
\end{theorem}

\begin{proof}
Suppose the conclusion fails with separation $\rho>0$. When $\lambda>0$, by
Fact~\ref{prop:purely-atomic-diagonal-ultrapower}, $\mathcal P_{a,b}(N_{\rm at})$ is compact in $L^2$, so we may pass to
a subsequence on which the restrictions to $N_{\rm at}$ converge.
Their limit is independent. We replace the restrictions to $N_{\rm at}$ by this limit, changing the pairs by $o(1)$. Denote the projection trace vectors of the limit for $\tau_{\rm at}$ by $\mu,\nu$
and its weighted overlap matrix by
$C^{\rm at}=\lambda\mu\nu^{\mathsf T}$. This settles the case $\lambda=1$.

Next, we pass to a subsequence on which both global trace vectors
converge, and we merge projections whose global traces tend to
zero into projections with positive limiting trace, so that $\tau(p_i^{(m)})$ and $\tau(q_j^{(m)})$ have
positive limits for all remaining $i,j$. A projection whose trace tends to zero vanishes on $N_{\rm at}$, where
the pairs are now constant, so the merging does not change them there.
At the end, we put zero projections back at the discarded indices.
If either partition has only one nonzero member, it is already
independent of the other, so we may omit this case. When $\lambda>0$, we number the
rows and columns so that $\mu_0,\nu_0>0$. When $\lambda=0$, we set
$C^{\rm at}=0$, $\mu_i=\delta_{i0}$, and $\nu_j=\delta_{j0}$.

\emph{Splitting.} We decompose
\[
N=N_{\rm at}\oplus N_0\oplus\bigoplus_rM_r,\qquad
\lambda+w_0+\sum_rw_r=1,
\]
where $N_0$ is a von Neumann algebra with diffuse center and each $M_r$ is a ${\rm II}_1$
factor with trace weight $w_r>0$. Choose $s_p,s_q\in(0,1)$ with
$\sqrt{2as_p}+\sqrt{2bs_q}<\rho/3$, and set
\[
S=(1-s_p)I+s_p\mu\mathbf1^{\mathsf T},\qquad
T=(1-s_q)I+s_q\nu\mathbf1^{\mathsf T}.
\]
Here $I$ is the identity matrix and $\mathbf1$ is the column vector with all entries $1$, of size $a$ for $S$
and $b$ for $T$.
In each factor, we split each column projection according to $T$,
preserving its overlaps with every $p_i$ in the same proportions,
and regroup the pieces to form the new columns. We then apply
Lemma~\ref{lem:exact-markov-relabelling} to the row projections with $S$.
In $N_0$, we split each projection by multiplying it by the members of a
central partition of $1_{N_0}$. By Corollary~\ref{cor:proportional-projection-splitting},
we choose the central pieces so that the overlaps with the current
projections on the other side split in the required proportions.
The pair on $N_{\rm at}$ stays fixed. Let $C,A,B$ denote the global overlap
matrix and the row and column trace vectors before these splittings,
and set $F=C-AB^{\mathsf T}$. Since $S\mu=\mu$ and $T\nu=\nu$,
\[
C'=SCT^{\mathsf T},\qquad A'=SA,\quad B'=TB,\qquad
F'=SFT^{\mathsf T}\longrightarrow0.
\]
Using the diagonal entries of $S,T$ and summing with the trace
weights, we apply Cauchy--Schwarz to bound the total $L^2$ distance by
$\sqrt{2as_p}+\sqrt{2bs_q}$.
Denote the resulting pair by
$(\widehat{\mathbf p},\widehat{\mathbf q})$.

\emph{A Gram bound.}
In the factor $M_r$, let $\tau_r$ be the normalized trace and set
$\alpha_i=\tau_r(p_i)$ and
$\bar C_{ij}=\tau_r(p_i\widehat q_j)$. For a real array $(x_{kj})_{k,j>0}$, set
$B_k=\sum_{j>0}x_{kj}\widehat q_j$ and $B_0=0$.
Since $T_{0j}\ge s_q\nu_0$ for every $j$, by the variance bound in
Corollary~\ref{cor:markov-absorption} we have
\[
\bar C_{i0}=\sum_jT_{0j}\tau_r(p_iq_j)\ge s_q\nu_0\alpha_i,\qquad
V_i(B_k)\ge s_q\nu_0\sum_{j>0}\bar C_{ij}x_{kj}^2.
\]
In the variance identity of Lemma~\ref{lem:exact-markov-relabelling},
we keep the terms indexed by $(0,k)$ and $(k,0)$. Using
$s_{0i}\ge s_p\mu_0$, we obtain
\[
\left\|\sum_{k>0}i[\widehat p_k,B_k]\right\|_{2,r}^2
\ge2s_ps_q\mu_0\nu_0\sum_{k,j>0}
\tau_r(\widehat p_k\widehat q_j)x_{kj}^2.
\] Let $z_{\rm II}$ be the sum of
the central projections supporting the $M_r$, and set
$C^{\rm II}_{ij}=\sum_rw_r\tau_r(\widehat p_i\widehat q_j)$.
Multiplying the factor estimates by $w_r$ and summing gives
\begin{equation}
\left\|\sum_{k,j>0}x_{kj}z_{\rm II}i[\widehat p_k,\widehat q_j]\right\|_2^2
\ge2s_ps_q\mu_0\nu_0\sum_{k,j>0}C^{\rm II}_{kj}x_{kj}^2.
\label{eq:weighted-factor-gram-bound}
\end{equation}
Here both the norm and the overlaps $C^{\rm II}_{ij}$ are measured
with the global trace $\tau$.

\emph{The overlaps to correct.}
Passing to a subsequence, we may assume that the trace vectors and overlaps
on $N_0$ and on $z_{\rm II}N$ converge, omitting any zero summand.
When $w_0>0$, let $C_m^0$ be the overlap matrix on $N_0$ with its
normalized trace, and set
$\mathcal J_0=\{(i,j):i,j>0,\ \lim_m w_0C_{ij,m}^0>0\}$.
When $w_0=0$, set $\mathcal J_0=\varnothing$.
Let $\mathcal J_{\rm II}$ be its complement among the pairs with
$i,j>0$. We claim that $C^{\rm II}_{ij,m}$ has a positive limit for
every $(i,j)\in\mathcal J_{\rm II}$. Indeed, at such a pair, the weighted
overlap on $N_0$ tends to zero. If $C^{\rm at}_{ij}=0$, then, since $F'\to0$ and $A'_iB'_j$ has a positive limit, the global overlap
$C'_{ij}$ has a positive limit, and hence so does its part $C^{\rm II}_{ij,m}$ on the factor summands.

If $C^{\rm at}_{ij}>0$, then $\mu_i\nu_j>0$. Set $W=\sum_rw_r$.
When $W>0$, the two splittings give
\[
C^{\rm II}_{ij,m}\ge s_ps_q\mu_i\nu_jW>0.
\]
If $W=0$, then $w_0>0$, since the case $\lambda=1$ has already been
settled. The same splittings on $N_0$ give
$C^0_{ij,m}\ge s_ps_q\mu_i\nu_j>0$, which would put $(i,j)$ in
$\mathcal J_0$. Thus no pair belongs to $\mathcal J_{\rm II}$ when $W=0$, and the claim follows.

\emph{Conjugation.}
We apply Lemma~\ref{lem:finite-coordinate-overlap-inverse} in $N$, taking
$z=z_{\rm II}$, to the overlaps indexed by $\mathcal J_{\rm II}$ and
the global targets $A_i'B_j'$.
By \eqref{eq:weighted-factor-gram-bound} and the claim, the Gram matrix in that lemma is bounded below by a positive multiple of the identity for large $m$. Hence there is a single self-adjoint
$h_m\in z_{\rm II}N$ with $\|h_m\|_\infty\to0$ such that conjugating the
first partition by $e^{ih_m}$ corrects these overlaps. Take $h_m=0$ when
$\mathcal J_{\rm II}$ is empty.

\emph{Central replacement.}
After this conjugation, the remaining global covariance matrix $R_m$ tends to zero, has zero row and column sums, and satisfies $R_{ij,m}=0$ for all $(i,j)\in\mathcal J_{\rm II}$. If $R_m=0$, in particular if $w_0=0$, no further correction is needed.

Since $R_{i0,m}=-\sum_{j>0}R_{ij,m}$, a nonzero $R_{i0,m}$ with
$i>0$ requires $(i,j)\in\mathcal J_0$ for some $j>0$. Hence the
corresponding row projection in $N_0$ has positive limiting trace
with respect to the normalized trace on $N_0$. The same argument
applies to columns.
Every entry at which $R_m$ can be nonzero has positive limiting
overlap on $N_0$. Indeed, for entries in $\mathcal J_0$ this follows from
its definition, and for row zero and column zero we use
\[
C^0_{i0,m}\ge s_q\nu_0\alpha^0_{i,m},\qquad
C^0_{0j,m}\ge s_p\mu_0\beta^0_{j,m},\qquad
C^0_{00,m}\ge s_ps_q\mu_0\nu_0,
\]
where $\alpha^0_{i,m}$ and $\beta^0_{j,m}$ are the row and column
sums of $C_m^0$, respectively.

For large indices with $R_m\ne0$, set
\[
\kappa_m=\sqrt{\max_{i,j}|R_{ij,m}|},\qquad
H_m=C_m^0-\frac{R_m}{w_0\kappa_m}.
\]
Then $H_m$ is nonnegative and has the same row and column sums as $C_m^0$.
Choose a central projection in $N_0$ of normalized trace $\kappa_m$
that cuts every current overlap in the proportion $\kappa_m$.
Inside this corner, choose a central partition, allowing zero projections,
whose traces, measured in $N_0$, are $\kappa_mH_{ij,m}$, and replace the two partitions on this corner by its row and column sums. By Corollary~\ref{cor:proportional-projection-splitting}, both choices are possible.
The new overlap matrix on $N_0$ is
\[
(1-\kappa_m)C_m^0+\kappa_mH_m=C_m^0-R_m/w_0.
\]
This cancels the global covariance and preserves the trace vectors and the
pair on $N_{\rm at}$. The total $L^2$ distance is at most
$(a+b)\sqrt{w_0\kappa_m}=o(1)$. Together with the earlier estimates, this
contradicts the assumed separation $\rho$ once zero projections are put back at
the discarded indices.
\end{proof}

\begin{proof}[Proof of Corollary~\ref{cor:diffuse-scalar-problem-30}]
The purely atomic part is absent for $D$ and is scalar for
$\C\oplus_\lambda D$, so the additional independence hypothesis is
automatic.  The result then follows from
Theorem~\ref{thm:independent-atomic-central-assembly} and
Proposition~\ref{prop:fixed-ambient-stability-equivalence}.
\end{proof}

The construction in
Theorem~\ref{thm:independent-atomic-central-assembly} uses independence on the atomic
summand. Indeed, if
\[
C^{\rm at}=\lambda\mu\nu^{\mathsf T}+K,
\]
where $K\ne0$ has zero row and column sums, then the same two splittings on the diffuse summand give
\[
F'=SFT^{\mathsf T}+\bigl[1-(1-s_p)(1-s_q)\bigr]K,
\]
since $SKT^{\mathsf T}=(1-s_p)(1-s_q)K$. Thus $F'$ does not tend to zero when $F$ does.

If every prescribed overlap matrix were stable in every fixed
diffuse tracial von Neumann algebra, independence would be stable in every fixed tracial
von Neumann algebra: after passing to a subsequence, we could replace the
partitions on the atomic summand by their limits and correct the
remaining overlaps on the diffuse summand. For each fixed diffuse summand, we prove this implication and its converse
in Corollary~\ref{cor:universal-table-independence-equivalence}.

Cancellation can also occur between a scalar atom and a matrix factor. In the von Neumann algebra
$N=\C\oplus_\lambda M_2(\C)$, take the binary partitions generated by
$1\oplus p$ and $1\oplus q$, where $p,q\in M_2(\C)$ are rank-one projections satisfying
\[
\operatorname{tr}_2(pq)=\frac{1-\lambda}{4}.
\]
Both generating projections have trace $(1+\lambda)/2$, while
\[
\lambda+(1-\lambda)\operatorname{tr}_2(pq)
=\frac{(1+\lambda)^2}{4}.
\]
Thus the global binary partitions are independent, although the
factor covariance is $-\lambda/4$.
However, by Theorem~\ref{thm:finite-type-I-problem-30}, stability
still holds in this fixed type I von Neumann algebra.

If two pairs in $\mathcal P_{a,b}(N)$ agree off a central projection $z$, then corresponding projections differ by at most
$\sqrt{\tau(z)}$ in $L^2$, and their total distance is at most
$(a+b)\sqrt{\tau(z)}$. Thus changing the partitions on central
summands whose total trace tends to zero costs $o(1)$ in $L^2$.

\section{Instability and the tracial version of Problem~30}
\label{sec:unstable-tables}

We now construct a fixed tracial von Neumann algebra, the direct sum of
finitely many scalar atoms and one ${\rm II}_1$ factor, in which
independence of pairs of partitions is not stable. In the factor itself this independence is stable by
Corollary~\ref{cor:diffuse-scalar-problem-30},
but the factor admits a prescribed overlap matrix that is not stable.

We obtain this matrix from a homomorphism into an ultrapower that
cannot be lifted to coordinate homomorphisms. We then choose the
scalar weights so that any nearby independent pair in the direct sum
must realize the prescribed matrix on the factor. No such pair exists,
even if both partitions and all their traces are allowed to change.

We use the set $\mathcal Y_D(C)$ and the stability of a prescribed overlap matrix from Section~\ref{sec:projection-tables}, and $\operatorname{PS}(N;a,b)$ and $\operatorname{P30}(N,\cU)$ from Section~\ref{sec:fixed-ambient-beyond-factors}.

\begin{theorem}[A counterexample to the tracial version of Problem~30]
\label{thm:fixed-ambient-problem-30-counterexample}
There is a fixed tracial von Neumann algebra consisting of finitely many scalar atoms
and one ${\rm II}_1$ factor such that, for every nonprincipal ultrafilter
$\cU$ on $\N$, its ultrapower contains independent finite-spectrum unitaries
with no independent unitary representatives. More precisely, let
$M_1,L,n,T,\kappa$ be as in Proposition~\ref{prop:fixed-factor-table-obstruction}, and set
\[
t=\frac1{2L},\qquad E_{ij}=\frac1{n^2}-tT_{ij},
\qquad N=\C^{n^2}\oplus M_1.
\]
Equip $N$ with the faithful normal tracial state
\[
\tau_N((z_{ij}),x)=\sum_{i,j}E_{ij}z_{ij}+t\tau_{M_1}(x).
\]
There is a sequence of partition pairs in $\mathcal P_{n,n}(N)$ whose covariance defects
tend to zero and whose distances from $Z_{n,n}(N)$ are at least
\[
\rho=\min\left\{\frac1{\sqrt8\,L^2},\frac{\kappa}{\sqrt{2L}}\right\}>0
\]
at every coordinate.
\end{theorem}

We give the proof after
Corollary~\ref{cor:universal-table-independence-equivalence}. In the construction of
Proposition~\ref{prop:fixed-factor-table-obstruction}, zero entries of the prescribed
overlap matrix force zero products. By Lemma~\ref{lem:rectangular-partition-decoding},
we recover four partitions of one
projection, possibly with zero members, such that every member of the first two commutes with every
member of the last two. Coordinate representatives of the constructed pair that realize this overlap
matrix exactly would give commuting coordinate lifts of
those four partitions. After transport to a fixed corner, these lifts would contradict the choice of
the partitions obtained from the theorem of Ioana and Spaas.

\begin{lemma}
\label{lem:equal-trace-unitary-transport}
Let $(M,\tau)$ be a finite factor, and let $P,Q\in M$ be projections with $\tau(P)=\tau(Q)$.  There is a unitary
$w\in M$ such that
\[
wPw^*=Q,
\qquad
\norm{w-1}_2\le2\norm{P-Q}_2.
\]
\end{lemma}

\begin{proof}
Set $a=QP+(1-Q)(1-P)$. Then
\[
aP=Qa,\qquad a-1=(Q-P)(2P-1),\qquad
\norm{a-1}_2=\norm{P-Q}_2.
\]
The polar partial isometry of $a$ intertwines $P,Q$ on its supports.
Since $\tau(P)=\tau(Q)$, we can extend it across the missing supports, separately in the
$P,Q$ corners and in their complements, to a unitary $w$ with $wP=Qw$
and $w|a|=a$. Finally,
\[
\norm{w-a}_2^2=\norm{1-|a|}_2^2
\le\norm{1-a}_2^2
\]
by $\Ree\tau(a)\le\tau(|a|)$. The triangle inequality then gives
$\norm{w-1}_2\le2\norm{P-Q}_2$, as desired.
\end{proof}

\begin{lemma}
\label{lem:rectangular-partition-decoding}
Let $M$ be a finite von Neumann algebra, let $L\ge1$ be an integer, and let
\[
(e_{ac},f_{bd}),
\qquad
(g_{bc},h_{ad})
\]
be partitions of $1_M$, with all indices in $\{0,\ldots,L-1\}$.
Zero projections are allowed in both the given and the constructed partitions. Suppose
\[
\begin{array}{ll}
e_{ac}g_{bc'}=0 & (c\ne c'),\\
e_{ac}h_{a'd}=0 & (a\ne a'),\\
f_{bd}g_{b'c}=0 & (b\ne b'),\\
f_{bd}h_{ad'}=0 & (d\ne d').
\end{array}
\]
Set
\[
E:=\sum_{a,c}e_{ac},\qquad F:=1-E,
\qquad
G:=\sum_{b,c}g_{bc},\qquad H:=1-G,
\]
and let $J\subseteq M$ be the unital von Neumann subalgebra generated by the two partitions.  Then
\[
\Omega:=1-(E-G)^2\in Z(J).
\]
Set
\[
z:=1_{[1/4,3/4]}(\Omega),
\qquad
e_0:=Ez,\qquad f_0:=Fz,\qquad \Theta:=e_0Ge_0.
\]
If $z\ne0$, the element
\[
v:=\bigl[\Theta(e_0-\Theta)\bigr]^{-1/2}e_0Gf_0
\]
is a partial isometry with $vv^*=e_0$ and $v^*v=f_0$. The inverse
square root is taken in $e_0Me_0$ and has operator norm at most
$4/\sqrt3$.  For $z=0$, set $v=0$.  The four families
\[
\begin{aligned}
A_a^\sharp&:=\sum_c e_{ac}z,&
B_b^\sharp&:=v\left(\sum_d f_{bd}z\right)v^*,\\
C_c^\sharp&:=\sum_a e_{ac}z,&
D_d^\sharp&:=v\left(\sum_b f_{bd}z\right)v^*
\end{aligned}
\]
are partitions of $e_0$ satisfying
\[
[A_a^\sharp,C_c^\sharp]=[A_a^\sharp,D_d^\sharp]
=[B_b^\sharp,C_c^\sharp]=[B_b^\sharp,D_d^\sharp]=0
\]
for all indices $a,b,c,d$.
\end{lemma}

\begin{proof}
By direct computation, we have $0\le\Omega\le1$,
$[\Omega,E]=[\Omega,G]=0$, and
\[
E\Omega E=EGE,\quad F\Omega F=FHF,\quad
G\Omega G=GEG,\quad H\Omega H=HFH.
\]
The zero-product relations imply that $EGE$ commutes with each $\sum_a e_{ac}$ and
$EHE$ with each $\sum_c e_{ac}$. Since $EGE+EHE=E$, both operators
commute with both families of sums. For each $a,c$, the product of the
corresponding row and column sums is $e_{ac}$, so both operators commute
with $e_{ac}$. The same argument in the other three corners proves
$\Omega\in Z(J)$.

The case $z=0$ is immediate. Otherwise work in $zMz$, denote the
compressed operators by the same letters, and write $1$ for $z$.
For $K=EGF$ and $D_0=FGF$, since $\tfrac14\le\Omega\le\tfrac34$, we have
\[
KK^*=\Theta(E-\Theta)\ge\tfrac3{16}E,\qquad
K^*K=D_0(F-D_0)\ge\tfrac3{16}F.
\]
This proves the support identities for $v$ and the inverse-square-root bound.
Set
\[
G_b=\sum_c g_{bc},\quad G_c'=\sum_b g_{bc},\qquad
H_a=\sum_d h_{ad},\quad H_d'=\sum_a h_{ad},
\]
and define $u_G=\Omega^{-1/2}EG$ and
$u_H=(1-\Omega)^{-1/2}EH$, with inverses in $zMz$.
Their right supports are $G,H$, respectively, and both left
supports are $E$. The zero-product relations give
$EH_a=A_a^\sharp H$ and $EG_c'=C_c^\sharp G$, hence
\[
A_a^\sharp=u_HH_au_H^*,\qquad C_c^\sharp=u_GG_c'u_G^*.
\]
With $F_b=\sum_d f_{bd}$, we have $GF_b=G_bF$ and
\[
EGF_bGE=EG_bE(1-\Omega),\qquad
B_b^\sharp=vF_bv^*=\Omega^{-1}EG_bE=u_GG_bu_G^*.
\]
The same calculation with $H$, using $EHF=-EGF$, gives
$D_d^\sharp=u_HH_d'u_H^*$. Hence $B_b^\sharp$ commutes with $C_c^\sharp$, and
$A_a^\sharp$ with $D_d^\sharp$. The other two commutation relations
follow directly from the $e$ and $f$ partitions. All four families are partitions of $E=e_0$.
\end{proof}

\begin{proposition}
\label{prop:fixed-factor-table-obstruction}
There exist a ${\rm II}_1$ factor $M$, an integer $L\ge3$, and $\kappa>0$ with the following property.
Set $M_1=M_2(M)$, a ${\rm II}_1$ factor with normalized trace $\tau_{M_1}=\operatorname{tr}_2\otimes\tau_M$, and $n=2L^2$.
Index rows and columns by $I=\{0,1\}\times\mathbb Z_L\times\mathbb Z_L$, and define
\[
T_{(s,x,y),(r,u,v)}=
\begin{cases}
\dfrac1{4L^3},&s=r\text{ and }y=v+1,\\[1mm]
\dfrac1{4L^3},&s\ne r\text{ and }x=u,\\[1mm]
0,&\text{otherwise},
\end{cases}
\]
where the coordinates in $\mathbb Z_L$ are read modulo $L$.  Every row and column has sum $1/n$, and $T_{ii}=0$.
There are pairs $(\mathbf p^{(m)},\mathbf q^{(m)})\in\mathcal P_{n,n}(M_1)$ such that
\[
\max_{i,j}\abs{\tau_{M_1}(p_i^{(m)}q_j^{(m)})-T_{ij}}\longrightarrow0,
\qquad
d_2((\mathbf p^{(m)},\mathbf q^{(m)}),\mathcal Y_{M_1}(T))\ge\kappa
\]
for every $m$.  In particular, a fixed prescribed overlap matrix need not be stable even in a fixed diffuse factor.
\end{proposition}

\begin{proof}
Ioana and Spaas~\cite[Theorem~E and Section~5.2]{ioana-spaas} gave a ${\rm
II}_1$ factor $M$, a nonprincipal ultrafilter
$\omega$ on $\N$, and a trace-preserving unital $^*$-homomorphism of von Neumann algebras
$\Phi:L(\mathbb F_2\times\mathbb F_2)\to M^\omega$
whose restriction to $\mathbb F_2\times\mathbb F_2$ cannot be
represented by a sequence of homomorphisms $\mathbb F_2\times\mathbb F_2\to U(M)$,
even when no condition is imposed on the traces of their images.
Choose unitary representatives of the images $U_1,U_2,V_1,V_2$
of the four canonical group generators under $\Phi$. Consider the set of
unitary quadruples $(u_1,u_2,v_1,v_2)$ in $M$ satisfying $[u_i,v_j]=0$ for
$1\le i,j\le2$. The distances of these representatives from this set,
measured by the sum of four $L^2$ distances, have ultralimit
$\delta_\Phi>0$. Indeed, if this limit were zero, choosing approximate
minimizers would give the forbidden lifts.

Fix $L\ge3$ with $8\pi/L<\delta_\Phi$. Round each of the four unitaries in $M^\omega$ to the nearest
$L$th root of unity by Borel functional calculus, and write
$A,B,C,D$ for their spectral partitions.
Every member of $A$ or $B$ commutes with every member of $C$ or
$D$. Since $\Phi$ is trace-preserving, $U_1,U_2$ are free Haar unitaries, as
are $V_1,V_2$. Thus each spectral projection has
trace $1/L$, and
\[
\tau(A_aB_b)=\tau(C_cD_d)=L^{-2}.
\]
The rounded quadruple is within total operator-norm distance $4\pi/L$
of the original one. Hence these partitions have no coordinate lifts by partitions of $1$, allowing zero projections,
for which every member of the first two commutes with every member of
the last two, even if the coordinate projection traces are allowed to vary.

In $M_1^\omega=M_2(M^\omega)$, set
\[
e=\begin{pmatrix}1&0\\0&0\end{pmatrix},\quad f=1-e,\qquad
g=\frac12\begin{pmatrix}1&1\\1&1\end{pmatrix},\quad h=1-g,
\]
and form
\[
p^0_{ac}=e\otimes A_aC_c,\quad p^1_{bd}=f\otimes B_bD_d,\qquad
q^0_{bc}=g\otimes B_bC_c,\quad q^1_{ad}=h\otimes A_aD_d.
\]
Using $L(\mathbb F_2\times\mathbb F_2)\cong L(\mathbb F_2)\mathbin{\overline\otimes}L(\mathbb F_2)$
and trace preservation of $\Phi$, we obtain $\tau(xy)=\tau(x)\tau(y)$
whenever $x\in W^*(U_1,U_2)$ and $y\in W^*(V_1,V_2)$.
The projections in $A,B$ belong to the first algebra, and those in
$C,D$ to the second. The overlap matrix $T^0$ of this pair of partitions is therefore
\[
\begin{aligned}
\tau_{M_1}(p^0_{ac}q^0_{bc'})&=\frac{\delta_{c,c'}}{4L^3},&
\tau_{M_1}(p^0_{ac}q^1_{a'd})&=\frac{\delta_{a,a'}}{4L^3},\\
\tau_{M_1}(p^1_{bd}q^0_{b'c})&=\frac{\delta_{b,b'}}{4L^3},&
\tau_{M_1}(p^1_{bd}q^1_{ad'})&=\frac{\delta_{d,d'}}{4L^3}.
\end{aligned}
\]
Each row and column has $2L$ nonzero entries and sum $1/n$.
We lift the four partitions
$(A_aC_c)_{a,c}$, $(B_bD_d)_{b,d}$, $(B_bC_c)_{b,c}$, and
$(A_aD_d)_{a,d}$ of $1$ separately by
Lemma~\ref{lem:nested-lifts-general}, and combine these lifts with
$e,f,g,h$ to obtain $(\mathbf p_m,\mathbf q_m)$ representing this
pair. A projection that is a sum of members of two of these partitions
need not receive the same representative from the two separate lifts.

Suppose there were representatives with exact overlap matrix $T^0$.
Zero overlaps give zero products, so
Lemma~\ref{lem:rectangular-partition-decoding} applies.
With $E_m,G_m$ defined as in that lemma, set
\[
\Omega_m=1-(E_m-G_m)^2,\qquad
z_m=1_{[1/4,3/4]}(\Omega_m).
\]
Since $(\Omega_m)_\omega=1/2$,
\[
\tau_{M_1}(1-z_m)\le16\|\Omega_m-\tfrac12\|_2^2\longrightarrow_\omega0.
\]
The function $t\mapsto[t(1-t)]^{-1/2}$ is bounded on $[1/4,3/4]$
and continuous near $1/2$. Applying functional calculus in the corners $E_mz_mM_1E_mz_m$ to the formula for $v_m$, we obtain $(v_m)_\omega=e_{12}\otimes1$.
The four partitions obtained from
Lemma~\ref{lem:rectangular-partition-decoding} have unit $E_mz_m$ and
represent $A,B,C,D$ in the upper-left corner. We add $E_m(1-z_m)$ to the projection with index $0$ in all
four partitions. They retain the four commutation relations and now have unit $E_m$.
Their classes in $M_1^\omega$ are unchanged.
Since the row sums are exact, $\tau_{M_1}(E_m)=1/2=\tau_{M_1}(e)$, so by
Lemma~\ref{lem:equal-trace-unitary-transport} there is a single unitary
$w_m$ satisfying
\[
w_mE_mw_m^*=e,\qquad
\|w_m-1\|_2\le2\|E_m-e\|_2\longrightarrow_\omega0.
\]
Conjugating all four partitions by this same unitary gives
coordinate lifts with all four required commutation relations in
$eM_1e\cong M$, whose normalized trace is $2\tau_{M_1}$. This contradicts
the choice of $A,B,C,D$.

The set $\mathcal Y_{M_1}(T^0)$ is nonempty. Indeed, we may split $1_{M_1}$
into projections of traces $T^0_{ij}$ and take row and column sums.
By the preceding contradiction, we have
\[
\kappa_0=\lim_{m\to\omega}
d_2((\mathbf p_m,\mathbf q_m),\mathcal Y_{M_1}(T^0))>0.
\]
Choose an increasing subsequence with overlap error below $1/k$
and distance at least $\kappa=\kappa_0/2$ at its $k$th term.
Finally, we index the rows by $p_{(s,x,y)}=p^s_{x,y}$ and the
columns by
\[
q_{(0,x,y)}^{\rm new}=q^0_{x,y+1},\qquad
q_{(1,x,y)}^{\rm new}=q^1_{x,y+1}.
\]
This fixed permutation gives the stated $T$, with $T_{ii}=0$ and the
same row and column sums.
\end{proof}

\begin{proposition}
\label{prop:scalar-atomic-table-transfer}
Fix a diffuse tracial von Neumann algebra $(D,\tau_D)$, an integer $n\ge2$, and a matrix
$T\in M_n(\R)$ with nonnegative entries summing to $1$ and $T_{ii}=0$ for every $i$.  Let $0<t<1$ satisfy
\[
E_{ij}:=n^{-2}-tT_{ij}>0\qquad(i,j).
\]
Suppose that a sequence $(\mathbf p^{(m)},\mathbf q^{(m)})\in\mathcal P_{n,n}(D)$ and a constant $\kappa>0$ satisfy
\[
\max_{i,j}\abs{\tau_D(p_i^{(m)}q_j^{(m)})-T_{ij}}\longrightarrow0,
\qquad
d_2((\mathbf p^{(m)},\mathbf q^{(m)}),\mathcal Y_D(T))\ge\kappa
\quad(m\ge1).
\]
Equip the fixed von Neumann algebra
\[
N=\C^{n^2}\oplus D
\]
with the trace
\[
\tau_N((z_{ij}),x)=\sum_{i,j}E_{ij}z_{ij}+t\tau_D(x).
\]
Then there are pairs $(\mathbf P^{(m)},\mathbf Q^{(m)})\in\mathcal P_{n,n}(N)$ such that
\[
\max_{i,j}\abs{\Delta_{ij}(\mathbf P^{(m)},\mathbf Q^{(m)})}\longrightarrow0,
\qquad
d_2((\mathbf P^{(m)},\mathbf Q^{(m)}),Z_{n,n}(N))\ge\rho
\quad(m\ge1),
\]
where
\[
\rho:=\min\left\{\sqrt{\min_{i,j}E_{ij}},\sqrt t\,\kappa\right\}>0.
\]
In particular, $\operatorname{PS}(N;n,n)$ fails.
\end{proposition}

\begin{proof}
On the atom $(k,l)$, set $r_i(k,l)=\delta_{ik}$ and $s_j(k,l)=\delta_{jl}$,
and define
$P_i^{(m)}=r_i\oplus p_i^{(m)}$, $Q_j^{(m)}=s_j\oplus q_j^{(m)}$.
Their overlaps are
\[
\tau_N(P_i^{(m)}Q_j^{(m)})
=E_{ij}+t\tau_D(p_i^{(m)}q_j^{(m)})
=n^{-2}+t\bigl(\tau_D(p_i^{(m)}q_j^{(m)})-T_{ij}\bigr).
\]
Thus each projection trace tends to $1/n$ and the covariance
defects tend to zero.

Take any independent pair $(\mathbf P',\mathbf Q')$.
If any projection changes on an atom, the distance from
$(\mathbf P^{(m)},\mathbf Q^{(m)})$ is at least
$\sqrt{\min_{i,j}E_{ij}}$.
Otherwise write $P'_i=r_i\oplus p'_i$, $Q'_j=s_j\oplus q'_j$ and set
$A_i=\tau_N(P'_i)$, $B_i=\tau_N(Q'_i)$. Positivity and $T_{ii}=0$ give
\[
A_iB_i=E_{ii}+t\tau_D(p'_iq'_i)\ge n^{-2},
\qquad
1\ge\sum_i\sqrt{A_iB_i}\ge1.
\]
Equality in Cauchy--Schwarz and in each lower bound forces
$A_i=B_i=1/n$ for every $i$. Hence, by the covariance equations,
\[
\tau_D(p'_iq'_j)=\frac{A_iB_j-E_{ij}}t=T_{ij}.
\]
Thus $(\mathbf p',\mathbf q')\in\mathcal Y_D(T)$, and since the atomic
parts are unchanged, we have
\[
d_2((\mathbf P^{(m)},\mathbf Q^{(m)}),(\mathbf P',\mathbf Q'))
=\sqrt t\,d_2((\mathbf p^{(m)},\mathbf q^{(m)}),(\mathbf p',\mathbf q'))
\ge\sqrt t\,\kappa.
\]
Taking the minimum of these two bounds, we obtain $\rho$, as desired.
\end{proof}

\begin{corollary}
\label{cor:universal-table-independence-equivalence}
Let $(D,\tau_D)$ be a diffuse tracial von Neumann algebra. If a
fixed matrix $C\in M_{a\times b}(\R)$ with nonnegative entries summing
to $1$ is not stable in $D$, then $\operatorname{PS}(N;a+b,a+b)$ fails
in a von Neumann algebra $N$ consisting of finitely many scalar atoms
and this same $D$, with fixed positive trace weights.

Consequently the following assertions are equivalent.
\begin{enumerate}[label=\textup{(\roman*)}]
\item For every tracial von Neumann algebra $(N,\tau)$ of the form
$N=N_0\oplus D$, where $N_0$ is zero or a direct sum of matrix-factor
central summands and $\tau$ restricts to a positive multiple of $\tau_D$ on $D$, and every
nonprincipal ultrafilter $\cU$ on $\N$, $\operatorname{P30}(N,\cU)$
holds.
\item Every fixed finite matrix with nonnegative entries summing to $1$
is stable in $D$.
\end{enumerate}
For example, both assertions hold when $D$ is of type~I, by
Theorem~\ref{thm:finite-type-I-problem-30}, and both fail for the
factor $M_1$ of Proposition~\ref{prop:fixed-factor-table-obstruction}.
\end{corollary}

\begin{proof}
Let $n=a+b$ and embed $C$ in the off-diagonal rectangle
$T_{i,a+j}=C_{ij}$, setting all other entries to zero.
Then $T_{ii}=0$, and for $t=1/(2n^2)$,
$E_{ij}=n^{-2}-tT_{ij}\ge1/(2n^2)>0$.
Take a sequence of pairs in $\mathcal P_{a,b}(D)$ whose overlap matrices converge
to $C$ and whose distances from the pairs with overlap matrix $C$
are bounded below by a positive constant. We add zero projections in the
rows and columns used to extend $C$ to $T$. In any pair with
overlap matrix $T$, the projections in these rows and columns vanish by faithfulness.
Thus the extended pairs remain separated from all pairs with overlap
matrix $T$. Hence, we may apply Proposition~\ref{prop:scalar-atomic-table-transfer}.
The von Neumann algebra $N=\C^{n^2}\oplus D$ obtained this way has the form
in \textup{(i)}, so by
Proposition~\ref{prop:fixed-ambient-stability-equivalence}, \textup{(i)}
implies \textup{(ii)}.

Conversely, assume \textup{(ii)}, let $N=N_0\oplus D$ be as in
\textup{(i)}, and suppose that $\operatorname{PS}(N;a,b)$ fails.
Choose pairs in $\mathcal P_{a,b}(N)$ whose covariance defects tend to zero and whose
distances from the independent pairs are bounded below by a positive
constant. Let $t>0$ be the weight of $D$, so that $\tau=t\tau_D$ on $D$.
Passing to a subsequence, we may assume that the vectors of traces of the two partitions
converge to $\alpha\in\R^a$ and $\beta\in\R^b$ and, by the compactness in
Fact~\ref{prop:purely-atomic-diagonal-ultrapower}, that the restrictions to
$N_0$ converge in $L^2$. We replace these restrictions by their limits,
at distance tending to zero, and let $H$ be the limiting overlap
matrix on $N_0$, measured with $\tau$. Set $H=0$ when $N_0=0$.
The overlap matrices on $D$, measured with $\tau_D$, converge~to
\[
C'=\frac{\alpha\beta^{\mathsf T}-H}{t},
\]
a matrix with nonnegative entries summing to $1$.
By \textup{(ii)}, we can perturb the restricted pair in $\mathcal P_{a,b}(D)$
so that its overlap matrix is exactly $C'$. Combining the corrected
partitions with the limit on $N_0$ gives the overlap matrix
$H+tC'=\alpha\beta^{\mathsf T}$, hence exact independence, at distance
tending to zero. This contradiction proves $\operatorname{PS}(N;a,b)$, and
\textup{(i)} follows from
Proposition~\ref{prop:fixed-ambient-stability-equivalence}.
\end{proof}

\begin{proof}[Proof of Theorem~\ref{thm:fixed-ambient-problem-30-counterexample}]
Since $n=2L^2$,
\[
E_{ij}\in\left\{\frac1{8L^4},\frac1{4L^4}\right\},\qquad
E_{ii}=n^{-2},\qquad \sum_{i,j}E_{ij}=1-t.
\]
Applying Proposition~\ref{prop:scalar-atomic-table-transfer} to the sequence
in Proposition~\ref{prop:fixed-factor-table-obstruction}, we obtain the
required sequence, with separation constant $\rho$.

Index both partitions by $0,\ldots,n-1$, let $\zeta=e^{2\pi i/n}$, and set
\[
x_m=\sum_i\zeta^iP_i^{(m)},\qquad
y_m=\sum_j\zeta^jQ_j^{(m)}.
\]
For any nonprincipal $\cU$, the spectral projections of
$x=(x_m)_\cU$ and $y=(y_m)_\cU$ each have trace $1/n$, and every
overlap between the two partitions is $1/n^2$. Thus $x,y$ are independent.
If some unitary representatives generated orthogonal unital von Neumann algebras
$C_m,D_m$, their ultraproducts would contain these spectral partitions.
Applying Lemma~\ref{lem:nested-lifts-general} separately inside $\prod_\cU
C_m$ and $\prod_\cU D_m$, we obtain coordinate partitions that are
independent and represent
the same projections, index by index. Their distance $d_2$ from
the constructed pairs therefore tends to zero along $\cU$,
contradicting the uniform separation $\rho$.
\end{proof}

On the atomic summand, the normalized projection traces are $1/n$,
while the diagonal overlaps are $E_{ii}/(1-t)>1/n^2$. These restrictions
are therefore not independent and fail the atomic hypothesis of
Theorem~\ref{thm:independent-atomic-central-assembly}. Adjoining only one
scalar atom to $M_1$, however, would preserve the stability of independence by
Corollary~\ref{cor:diffuse-scalar-problem-30}.

It remains open whether stability holds in every direct sum of a
$\mathrm{II}_1$ factor with two scalar atoms of fixed positive weights, and
what the smallest partition size is for which stability fails.

\part{Exact moment identities for matrix models}
For matrix families that converge to free limits, possibly with amalgamation,
we use the perturbation method of Part~I to make the scalar and amalgamated
free-product identities exact through degrees that grow with the dimension
(Sections~\ref{sec:scalar-models} and~\ref{sec:relative-exactification}).
For GUE families, we show that the degree can grow like $\log N$ and no faster
(Section~\ref{sec:gue-bounds}).

\section{Scalar moment identities and strong matrix models}
\label{sec:scalar-models}

We apply the perturbation theorems of
Sections~\ref{sec:local-correction} and~\ref{sec:higher-order-no-loss}
in matrix algebras. We seek exact moment identities while preserving a
given strong limit, meaning convergence of normalized traces and operator
norms of all fixed noncommutative $^*$-polynomials. To preserve
convergence of these polynomial norms, we require the conjugating
unitaries to tend to $1$ in operator norm.

For subalgebras of fixed finite-dimensional isomorphism classes whose minimal projections
have traces bounded away from zero,
we obtain $n$-independence in every matrix algebra from a small defect
of order $2n$ (Theorem~\ref{thm:fd-noncommutative-exactify}). In all
sufficiently large dimensions, a small defect of order $n$ suffices
(Corollary~\ref{cor:large-matrix-fixed-type-no-loss}).
In Subsection~\ref{sec:stable-strong-models}, we treat strongly convergent matrix families with
free limits, assuming that the polynomial spaces of the individual tuples
eventually have the dimensions of the corresponding spaces in the limiting tracial GNS representation. The exact identities
are those of the free product of the marginal states at each dimension, and the correction preserves the strong limit.

\subsection{Fixed finite-dimensional \texorpdfstring{$C^*$}{C*}-algebras}

\begin{fact}
\label{lem:projection-polynomial-comparison}
Let $(M,\tau)$ be a tracial von Neumann algebra, let $n,a,b\ge1$ be
integers, and let
$1=\sum_{i=1}^a p_i=\sum_{j=1}^b q_j$ be two partitions of $1$,
allowing zero projections. Define the finite-dimensional unital abelian $^*$-subalgebras
\[
A_0:=\operatorname{span}\{p_1,\ldots,p_a\},\qquad
B_0:=\operatorname{span}\{q_1,\ldots,q_b\},
\]
and let $(F,\varphi):=(A_0,\tau|_{A_0})*(B_0,\tau|_{B_0})$ be their
von Neumann algebra free product, in which $\bar p_i,\bar q_j$ denote
the canonical images of $p_i,q_j$.
Then $A_0,B_0$ are $n$-independent if and only if
\[
\tau\bigl(R(p_1,\ldots,p_a,q_1,\ldots,q_b)\bigr)
=
\varphi\bigl(R(\bar p_1,\ldots,\bar p_a,\bar q_1,\ldots,\bar q_b)\bigr)
\]
for every noncommutative $^*$-polynomial $R$ in the two projection
families of total $^*$-degree at most $2n+1$.
\end{fact}

This follows from Lemma~\ref{prop:truncated-free-product} and the
definition of $n$-independence.

\begin{corollary}
\label{cor:matrix-rounding}
Fix $n,a,b\ge1$ and $\delta>0$, and let $\varepsilon_0,C$ be constants from
Corollary~\ref{thm:quantitative-local-exactify}. Let $N\ge1$ and let
\[
1=\sum_{i=1}^a p_i=\sum_{j=1}^b q_j
\]
be two partitions of $1$ in $M_N(\C)$ such that
\[
\operatorname{tr}_N(p_i)\ge\delta,
\qquad
\operatorname{tr}_N(q_j)\ge\delta.
\]
Define the unital abelian $^*$-subalgebras $A:=\operatorname{span}\{p_1,\ldots,p_a\}$ and
$B:=\operatorname{span}\{q_1,\ldots,q_b\}$. If
\[
\varepsilon_{2n}(A,B)\le\varepsilon_0,
\]
then there are $h=h^*\in M_N(\C)$ and $w=e^{ih}$ such that the partitions
\[
p_i':=wp_iw^*,
\qquad
q_j':=q_j
\]
have the same ranks as the original partitions, their abelian spans are $n$-independent, and
\[
\norm{h}_2\le C\varepsilon_{2n}(A,B),
\qquad
\norm{h}_\infty\le C\varepsilon_{2n}(A,B),
\qquad
\norm{w-1}_\infty\le C\varepsilon_{2n}(A,B).
\]
Moreover,
\[
\max_{1\le i\le a}\norm{p_i'-p_i}_\infty
\le2C\varepsilon_{2n}(A,B).
\]
\end{corollary}

\begin{proof}
This follows from Corollary~\ref{thm:quantitative-local-exactify} applied with the normalized trace. For the last bound, note that $\norm{wp_iw^*-p_i}_\infty\le2\norm{w-1}_\infty$.
\end{proof}

Although uniform in $N$, the preceding correction requires a small
defect of order $2n$. We reduce this requirement to order $n$ in large
dimensions by passing to a matrix ultraproduct.

\begin{corollary}
\label{cor:large-matrix-fixed-type-no-loss}
Fix $n\ge1$, finite-dimensional unital $C^*$-algebras $P,Q$, and $\delta,\rho>0$.  There are
\[
N_0=N_0(n,P,Q,\delta,\rho)<\infty,
\qquad
\epsilon=\epsilon(n,P,Q,\delta,\rho)>0
\]
such that the following holds. Suppose $N\ge N_0$ and
$A,B\subseteq M_N(\C)$ are unital $^*$-subalgebras supplied with unital $^*$-isomorphisms $P\cong A$
and $Q\cong B$. If every minimal projection of $A$ and $B$ has
normalized trace at least $\delta$, then
\[
\varepsilon_n(A,B)<\epsilon
\]
implies that there is $u\in U(N)$ with
\[
\norm{u-1}_\infty<\rho
\]
such that $uAu^*,B$ are $n$-independent for $\operatorname{tr}_N$.
\end{corollary}

\begin{proof}
If $P$ or $Q$ is one-dimensional, the identity unitary suffices. Otherwise,
suppose that the corollary fails. Then we can choose $N_k\ge k$ and unital $^*$-embeddings
$\pi_k:P\to M_{N_k}(\C)$, $\sigma_k:Q\to M_{N_k}(\C)$ satisfying the
trace bound, with defects less than $1/k$, but admitting no such
unitary. In a nonprincipal tracial ultraproduct
$M=\prod_\omega(M_{N_k}(\C),\operatorname{tr}_{N_k})$, a
${\rm II}_1$ factor, they induce faithful embeddings $\pi,\sigma$
with $n$-independent images. Indeed, the trace bound gives
faithfulness, and centering fixed elements with respect to the
coordinate traces lets us pass the alternating identities to the limit.

By Proposition~\ref{prop:dense-regular-exact-locus}, we can choose $w=e^{ia}$, with $a=a^*$
and $\norm a_\infty<\rho/3$, such that
the pair $w\pi(P)w^*,\sigma(Q)$ is $n$-independent
and order-$n$ regular.
We lift $a$ to self-adjoint $a_k$ with $\norm{a_k}_\infty\le\norm a_\infty$
and set $w_k=e^{ia_k}$.
Thus $(w_k)_\omega=w$ and $\norm{w_k-1}_\infty<\rho/3$ at every index.
Applying Lemma~\ref{lem:regular-limit-exactification} to
$\operatorname{Ad}(w_k)\circ\pi_k$ and $\sigma_k$, we obtain $h_k=h_k^*$
with $\norm{h_k}_\infty\to_\omega0$ making the coordinate
pairs $n$-independent on an $\omega$-large set.
On a smaller such set, $u_k=e^{ih_k}w_k$ satisfies
\[
\norm{u_k-1}_\infty
\le\norm{h_k}_\infty+\norm{w_k-1}_\infty<\rho,
\]
contradicting the choice of the pairs.
\end{proof}

In the next corollary, $E_A$ and $E_B$ denote the $\operatorname{tr}_N$-preserving conditional expectations of $M_N(\C)$ onto $A$ and $B$.

\begin{corollary}
\label{cor:large-matrix-relative-privacy}
Fix finite-dimensional unital $C^*$-algebras $P,Q$ and $\delta,\rho>0$. There are
$N_0<\infty$ and $\eta>0$ with the following property. Let $N\ge N_0$,
and let $A,B\subseteq M_N(\C)$ be unital $^*$-subalgebras identified with $P,Q$ by fixed
unital $^*$-isomorphisms. If every minimal projection of $A$ and $B$ has
normalized trace at least $\delta$ and
\[
\norm{(E_A-D_N)E_B}_{2\to2}<\eta,
\qquad
D_N(x):=\operatorname{tr}_N(x)1,
\]
then there is $u\in U(N)$ satisfying
\[
\norm{u-1}_\infty<\rho,
\qquad
(E_{uAu^*}-D_N)E_B=0.
\]
Thus $E_{uAu^*}(b)=\operatorname{tr}_N(b)1$ for every $b\in B$.
Here $\norm{\Psi}_{2\to2}$ denotes the operator norm of a linear map $\Psi$ on $M_N(\C)$ with respect to $\norm{\cdot}_2$.
\end{corollary}

\begin{proof}
For centered contractions $a\in A$ and $b\in B$,
\[
\abs{\operatorname{tr}_N(ab)}
=\abs{\operatorname{tr}_N(aE_A(b))}
\le\norm{a}_2\norm{E_A(b)}_2
\le\norm{(E_A-D_N)E_B}_{2\to2}.
\]
Hence $\varepsilon_1(A,B)\le\eta$.  Now choose $\eta$ below the threshold in
Corollary~\ref{cor:large-matrix-fixed-type-no-loss} at $n=1$. By that corollary, there is $u\in U(N)$ with $\norm{u-1}_\infty<\rho$ such that $uAu^*,B$ are $1$-independent, i.e., $(E_{uAu^*}-D_N)E_B=0$.
\end{proof}

Note that the identity $E_{uAu^*}(b)=\operatorname{tr}_N(b)1$, $b\in B$, of
Corollary~\ref{cor:large-matrix-relative-privacy} is the relative privacy condition
for $B$ and $E_{uAu^*}$ studied by Kribs, Levick, Nelson, Pereira, and Rahaman in connection with quantum
privacy~\cite[Definitions~2 and~5]{kribs-levick-nelson-pereira-rahaman}. Moreover, for every $r\ge1$, a two-term
telescoping estimate in trace norm gives
$\norm{(\operatorname{Ad}(u)-\operatorname{id})\otimes\operatorname{id}_{M_r(\C)}}_{1\to1}\le2\norm{u-1}_\infty$,
where $\norm{\cdot}_{1\to1}$ is the operator norm induced by the trace norm. Hence
\[
\norm{\operatorname{Ad}(u)-\operatorname{id}}_\diamond\le2\norm{u-1}_\infty<2\rho
\]
in the diamond norm $\norm{\Psi}_\diamond=\sup_{r\ge1}\norm{\Psi\otimes\operatorname{id}_{M_r(\C)}}_{1\to1}$.

We now construct $n$-independent pairs of partitions of $1$ with prescribed
ranks. First suppose that the normalized ranks converge to positive
limits. By Voiculescu's asymptotic freeness theorem~\cite[Theorems~3.8 and~3.11]{voiculescu-limit-laws},
in the quantitative form of Hiai and Petz~\cite[Corollary~2.2]{hiai-petz}, there are
partitions with small mixed-moment defects, and by Corollary~\ref{cor:matrix-rounding}
we can make them $n$-independent without changing their ranks. Under a uniform lower bound on the normalized ranks, we obtain such partitions in every sufficiently large dimension by compactness.

\begin{corollary}
\label{cor:matrix-sequential}
Fix $n,a,b\ge1$, strictly positive probability vectors
\[
\alpha=(\alpha_1,\ldots,\alpha_a),
\qquad
\beta=(\beta_1,\ldots,\beta_b),
\]
and dimensions $N_k\to\infty$. Suppose positive integer ranks satisfy
\[
\sum_i r_i^{(k)}=\sum_j s_j^{(k)}=N_k,
\qquad
\frac{r_i^{(k)}}{N_k}\longrightarrow\alpha_i,
\qquad
\frac{s_j^{(k)}}{N_k}\longrightarrow\beta_j.
\]
Then, for every sufficiently large $k$, there are partitions
$1=\sum_i p_{i,k}=\sum_j q_{j,k}$ in $M_{N_k}(\C)$ with the prescribed ranks such that their abelian spans are
$n$-independent. The partitions may be chosen to be asymptotically free at all orders:
for every noncommutative $^*$-polynomial $R$,
\[
\operatorname{tr}_{N_k}\bigl(R(p_{1,k},\ldots,p_{a,k},q_{1,k},\ldots,q_{b,k})\bigr)
\longrightarrow
\varphi_{\alpha,\beta}\bigl(R(\bar p_1,\ldots,\bar p_a,\bar q_1,\ldots,\bar q_b)\bigr),
\]
where $\bar p_i,\bar q_j$ are the canonical minimal projections of $\C^a$ and $\C^b$, and $\varphi_{\alpha,\beta}$ is the free-product trace
with atom weights $\alpha$ and $\beta$.
\end{corollary}

\begin{proof}
Choose partitions $(e_{i,k})_i$, $(f_{j,k})_j$ with the prescribed ranks.
For fixed distinct real $\lambda_i,\mu_j$, set
$X_k=\sum_i\lambda_i e_{i,k}$ and $Y_k=\sum_j\mu_j f_{j,k}$.
Their limiting spectral distributions have weights $\alpha,\beta$.
The projections $e_{i,k}$ and $f_{j,k}$ are fixed polynomials in
$X_k$ and $Y_k$, respectively. By the fixed-word estimate~\cite[Corollary~2.2]{hiai-petz},
which is uniform over traceless matrices of bounded norm,
every mixed projection moment converges in probability,
under Haar conjugation, to its free limit.
A diagonal choice over the countably many projection monomials gives
unitaries $u_k$ such that the unital abelian $^*$-subalgebras
\[
A_k=\operatorname{span}\{e_{i,k}\}_i,\qquad
B_k=\operatorname{span}\{u_kf_{j,k}u_k^*\}_j
\]
are asymptotically free with limiting weights $\alpha,\beta$.

By \eqref{eq:labelled-defect-comparison} and the moment convergence, we have
$\varepsilon_{2n}(A_k,B_k)\to0$.
Now choose $0<\delta<\min\{\alpha_i,\beta_j\}$. Then, by
Corollary~\ref{cor:matrix-rounding}, for all sufficiently large $k$ there is
$w_k=e^{ih_k}$ with $\norm{w_k-1}_\infty\to0$ such that the partitions
\[
p_{i,k}=w_ke_{i,k}w_k^*,\qquad q_{j,k}=u_kf_{j,k}u_k^*
\]
have unchanged ranks and $n$-independent spans. If $a=1$ or $b=1$,
take $w_k=1$. For every projection monomial $m$ of degree $d$,
\[
\left|\operatorname{tr}_{N_k}
\bigl(m(p_k,q_k)-m(e_k,u_kf_ku_k^*)\bigr)\right|
\le2d\norm{w_k-1}_\infty\longrightarrow0.
\]
Thus the correction preserves the full limiting joint distribution.
\end{proof}

\begin{theorem}
\label{thm:matrix-truncated-free}
Fix $n,a,b\ge1$ and $\delta>0$. There is an integer
$N_0=N_0(n,a,b,\delta)$ such that the following holds. If $N\ge N_0$ and positive integer ranks satisfy
\[
\sum_{i=1}^a r_i=\sum_{j=1}^b s_j=N,
\qquad
\frac{r_i}{N}\ge\delta,
\qquad
\frac{s_j}{N}\ge\delta,
\]
then there are partitions of $1$ in $M_N(\C)$,
\[
1=\sum_{i=1}^a p_i=\sum_{j=1}^b q_j,
\]
with $\operatorname{rank}(p_i)=r_i$ and $\operatorname{rank}(q_j)=s_j$ such that their abelian spans are $n$-independent for
$\operatorname{tr}_N$.

Equivalently, if $(F_N,\varphi_N)$ is the von Neumann algebra free product of $\C^a$ and $\C^b$ with atom weights
$(r_i/N)_{i=1}^a$ and $(s_j/N)_{j=1}^b$, then the partitions may be chosen so that their joint trace agrees with
$\varphi_N$ on every projection $^*$-polynomial of total degree at most $2n+1$.
\end{theorem}

\begin{proof}
Suppose otherwise. Then we can choose admissible ranks in dimensions $N_k\to\infty$ for which the conclusion fails.
By compactness, we may assume that the normalized rank vectors converge
to probability vectors $\alpha,\beta$ whose entries are all at least $\delta$.
However, by Corollary~\ref{cor:matrix-sequential}, these ranks are realized
with $n$-independent spans for all sufficiently large $k$, a contradiction.
The final equivalence follows from Fact~\ref{lem:projection-polynomial-comparison}.
\end{proof}

Patchell and Shiner~\cite[Corollary~D]{patchell-shiner} independently combined their results with
Voiculescu's asymptotic freeness theorem to show that, for every order, there are tuples of matrices of
some size that are independent to that order.

Chen, Van Voorhis, and Edelman~\cite[Example~9 and Section~3]{chen-van-voorhis-edelman-partial-freeness}
gave an earlier explicit example: trace-zero self-adjoint involutions
$A,B\in M_{2m}(\C)$ with $\operatorname{tr}_{2m}((AB)^k)=0$ for
$1\le k<m$. The rank-$m$ projections $(1+A)/2$ and $(1+B)/2$ therefore
generate unital $^*$-subalgebras that are $n$-independent for every
$1\le n<m$.

Note that in Theorem~\ref{thm:matrix-truncated-free} the partitions
have fixed sizes, and in Corollary~\ref{cor:large-matrix-fixed-type-no-loss}
the algebras lie in fixed isomorphism classes. Hence neither result
covers maximal abelian $^*$-subalgebras of $M_N(\C)$. For rank-one partitions, $1$-independence means that $|\langle e_i,f_j\rangle|^2=1/N$ for all $i,j$, where
$(e_i)$ and $(f_j)$ are the corresponding orthonormal bases of $\C^N$.
Popa~\cite[Theorem~3.2]{popa-orthogonal-pairs} constructed two
orthogonal maximal abelian $^*$-subalgebras of $M_N(\C)$ for every
$N\ge2$.

\subsection{Strong limits and polynomial filtrations}
\label{sec:stable-strong-models}

We now obtain exact moment identities for strongly convergent matrix families.
To do this, we fix a polynomial degree and compare the marginal polynomial spaces in
the matrices with those in the limiting GNS representation, assuming that
their dimensions eventually agree. Polynomials representing a basis at
the limit then give a basis in all sufficiently large matrices. Without
this assumption, a polynomial whose value is nonzero in the matrices but zero in the limiting GNS space may
still have to satisfy an exact mixed-moment equation at each finite stage
(Proposition~\ref{prop:strong-model-rank-obstructions}). If the limiting
families are free, we correct the mixed moments of these bases at each fixed degree, using
Lemma~\ref{lem:finite-dimensional-correction} and the free Gram lower bound of
Proposition~\ref{prop:noncommutative-free-gram}, and choose the degree
by a diagonal argument. In Section~\ref{sec:gue-bounds}, we obtain degrees proportional to $\log N$ for GUE families.

\emph{Strong convergence along arbitrary sizes and couplings.}
Male~\cite[Theorem~1.6]{male-norm-polynomials} proved almost-sure strong
convergence of
independent normalized GUE matrices together with matrices $\mathbf Y_N$ that are independent of them and converge
strongly almost surely to a tuple in a $C^*$-probability space with a faithful trace;
Collins and Male~\cite[Theorem~1.4]{collins-male} proved the same for
independent Haar unitary matrices. We use these
results here and in later sections, always for deterministic $\mathbf Y_N$,
in the following form: the dimensions may run through any strictly increasing
sequence, and the random matrices at different dimensions may be arbitrarily
dependent. These papers do not state this form; it follows from their proofs. First, both proofs, and those of
the asymptotic freeness theorems they quote, apply verbatim along a sequence
of dimensions. In Male's proof, every probabilistic estimate concerns one
dimension at a time, and almost-sure convergence comes from the Borel--Cantelli
lemma applied to summable failure probabilities
\cite[Section~8]{male-norm-polynomials}. Collins and Male reduce their theorem,
dimension by dimension, to Male's theorem and to the convergence of the
spectral distributions of one GUE matrix and of one Haar unitary matrix
\cite[Section~3]{collins-male}. Second, for a fixed polynomial and tolerance,
the failure event at a given dimension involves only the matrices of that
dimension. When the matrices at different dimensions are independent,
almost-sure convergence and the second Borel--Cantelli lemma show that the
failure probabilities are summable. The first Borel--Cantelli lemma then gives
almost-sure convergence for every joint law of the matrices at different
dimensions.

If $\mathbf X$ is a finite tuple in a unital $^*$-algebra, let
\[
\mathcal V_{\mathbf X}(d)
:=
\operatorname{span}_{\C}\{1,w(\mathbf X,\mathbf X^*):w\text{ is a $^*$-monomial of degree at most }d\}.
\]
If $\mathbf x$ lies in a tracial $C^*$-probability space $(\mathcal A,\tau)$, write
\[
\overline{\mathcal V}^{\,2}_{\mathbf x}(d)
\subseteq L^2(\mathcal A,\tau)
\]
for the image of $\mathcal V_{\mathbf x}(d)$ in the tracial GNS space. In the last paragraph of the next theorem, the real moment map is the constraint map of Section~\ref{sec:local-correction} for the bases described in that paragraph: its variable is $h=h^*\in M_{N_j}(\C)$, and the first tuple is conjugated by $e^{ih}$.

\begin{theorem}[Exact moment identities for strong models]
\label{thm:gns-rank-strong-exactification}
Let $(N_j)$ be positive integers.  Let $(\mathbf X_j,\mathbf Y_j)$ be pairs of tuples of fixed lengths in
$M_{N_j}(\C)$ which converge jointly and strongly to $(\mathbf x,\mathbf y)$ in a tracial $C^*$-probability space
$(\mathcal A,\tau)$.  Assume that
$C^*(\mathbf x)$ and $C^*(\mathbf y)$ are freely independent and that, for every fixed $d\ge0$, eventually
\[
\begin{aligned}
\dim_{\C}\mathcal V_{\mathbf X_j}(d)
&=\dim_{\C}\overline{\mathcal V}^{\,2}_{\mathbf x}(d),\\
\dim_{\C}\mathcal V_{\mathbf Y_j}(d)
&=\dim_{\C}\overline{\mathcal V}^{\,2}_{\mathbf y}(d).
\end{aligned}
\eqtag
\label{eq:gns-rank-stability}
\]
Set
\[
A_j:=C^*(\mathbf X_j),
\qquad
B_j:=C^*(\mathbf Y_j),
\]
and let $\varphi_j$ be the canonical state on the reduced $C^*$-free product
\[
(A_j,\operatorname{tr}_{N_j}|_{A_j})*(B_j,\operatorname{tr}_{N_j}|_{B_j}),
\]
with canonical embeddings $\iota_{A,j}$ and $\iota_{B,j}$.  After deleting finitely many indices, there are unitaries
$u_j\in U(N_j)$ and nondecreasing odd integers $D_j\to\infty$ such that
\[
\norm{u_j-1}_\infty\longrightarrow0,
\]
the corrected tuples $(\operatorname{Ad}(u_j)(\mathbf X_j),\mathbf Y_j)$ converge strongly to
$(\mathbf x,\mathbf y)$, and
\[
\begin{split}
&\operatorname{tr}_{N_j}
R(\operatorname{Ad}(u_j)(\mathbf X_j),\mathbf Y_j)\\
&\qquad=
\varphi_j
R(\iota_{A,j}(\mathbf X_j),\iota_{B,j}(\mathbf Y_j))
\end{split}
\eqtag
\label{eq:growing-exact-free-product-law}
\]
for every $^*$-polynomial $R$ of total degree at most $D_j$.

At stage $j$, with $n=(D_j-1)/2$, we may also require that the derivative
of the real moment map be surjective at the corrected tuple. To define
this map, take self-adjoint bases of the centered parts of
$\mathcal V_{\mathbf X_j}(2n+1)$ and $\mathcal V_{\mathbf Y_j}(2n+1)$. Among alternating
formal words in these bases of pair length $k$ with $1\le k\le n$,
choose one representative from each class under cyclic rotation and
adjoint. Conjugate the first tuple and take as coordinates the real
parts of the chosen word traces, together with their imaginary parts
when the adjoint word is not a cyclic rotation of the chosen word.
\end{theorem}

\begin{proof}
Fix $n\ge1$ and set $d=2n+1$. Choose self-adjoint formal polynomials
$f_0=1,f_1,\ldots,f_p$ of degree at most $d$ whose values form an
orthogonal real basis of the self-adjoint part of
$\overline{\mathcal V}^{\,2}_{\mathbf x}(d)$, with $f_i(\mathbf x)$
centered for $i\ge1$. Choose $g_0,\ldots,g_q$ in the same way for
$\mathbf y$. For $1\le i\le p$ and $1\le\ell\le q$, set
\[
a_{i,j}=f_i(\mathbf X_j)-\operatorname{tr}_{N_j}(f_i(\mathbf X_j)),
\qquad
b_{\ell,j}=g_\ell(\mathbf Y_j)-\operatorname{tr}_{N_j}(g_\ell(\mathbf Y_j)).
\eqtag\label{eq:strong-coordinate-bases}
\]
By convergence of the Gram matrices, these lists, together with $1$, are
linearly independent for all large $j$. Hence, by the dimension equalities
in \eqref{eq:gns-rank-stability}, they are bases of
$\mathcal V_{\mathbf X_j}(d)$ and $\mathcal V_{\mathbf Y_j}(d)$,
respectively. After fixed rescalings, we may assume that the centered basis
vectors are contractions.

If either centered list is empty, take $h_j=0$.
Otherwise, let $\Phi_j$ be the moment map through pair length $n$, with
cyclic-gradient Gram matrix $G_j$ at zero.
By moment convergence and Proposition~\ref{prop:noncommutative-free-gram},
we have $\Phi_j(0)\to0$ and $G_j\to G^{\rm fr}>0$.
We then apply Lemma~\ref{lem:finite-dimensional-correction}, using the
inverse of $G_j$ as in the proof of
Theorem~\ref{thm:fd-noncommutative-exactify}, and obtain $h_j$ with $\Phi_j(h_j)=0$ and $\|h_j\|_\infty\to0$.
The gradients change in $L^2$ by $O_n(\|h_j\|_\infty)$, so their
Gram matrix remains positive definite at the corrected point. Hence the derivative of the real moment map is surjective at the corrected tuple. Applying
Corollary~\ref{cor:filtered-truncated-comparison} with $D=d=2n+1$ and $m=n$, we
obtain \eqref{eq:growing-exact-free-product-law} through degree $d$.

Choose $n(j)\uparrow\infty$ sufficiently slowly that $u_j=e^{ih_j}$, with
$h_j$ obtained for $n=n(j)$, satisfies $\|u_j-1\|_\infty<1/n(j)$, and set
$D_j=2n(j)+1$, discarding a finite initial segment.
Since the generators are uniformly bounded and their changes tend to zero in
operator norm, the limiting traces and polynomial norms are unchanged.
\end{proof}

\begin{proposition}
\phantomsection
\label{prop:strong-model-rank-obstructions}
\begin{enumerate}[label=\textup{(\roman*)}]
\item There exist pairs of tuples with a faithful limiting trace
satisfying all hypotheses of
Theorem~\ref{thm:gns-rank-strong-exactification} other than
\eqref{eq:gns-rank-stability}, such that no unitaries $u_N$ with
$\norm{u_N-1}_\infty\to0$ yield the identities
\eqref{eq:growing-exact-free-product-law} through degree two.
\item The same degree-two obstruction occurs for pairs satisfying
every hypothesis of
Theorem~\ref{thm:gns-rank-strong-exactification} except
\eqref{eq:gns-rank-stability}, with a nonfaithful limiting trace and
with algebraic dimensions satisfying
$\dim_\C\mathcal V_{\mathbf X_N}(d)=\dim_\C\mathcal V_{\mathbf x}(d)$ and
$\dim_\C\mathcal V_{\mathbf Y_N}(d)=\dim_\C\mathcal V_{\mathbf y}(d)$
for every $d$.
\end{enumerate}
\end{proposition}

\begin{proof}
For (i), let $N$ run through the even integers and choose centered
self-adjoint unitaries $Y_N$. By~\cite[Theorem~1.4]{collins-male}, in the form recalled in
Subsection~\ref{sec:stable-strong-models}, we can realize Haar unitaries $U_N$
so that, on setting $X_N:=U_NY_NU_N^*$, the pair $(X_N,Y_N)$ converges
strongly to two free centered self-adjoint unitaries $(x,y)$. Let
$\epsilon_N>0$ tend to zero and set
\[
\mathbf A_N:=(X_N,X_N+\epsilon_NY_N),
\qquad
\mathbf B_N:=(Y_N).
\]
Then $(\mathbf A_N,\mathbf B_N)\to((x,x),y)$ strongly, and the limiting
marginal $C^*$-algebras are free. If conjugating $\mathbf A_N$ by $u_N$ gave
the free-product trace identities through degree two, then
\[
0=\operatorname{tr}_N
\bigl(u_N(A_{2,N}-A_{1,N})u_N^*Y_N\bigr)
=\epsilon_N\operatorname{tr}_N(u_NY_Nu_N^*Y_N).
\]
Dividing by $\epsilon_N$ contradicts $u_N\to1$ and
$\operatorname{tr}_N(Y_N^2)=1$. At each finite stage, the difference of
the two generators of $\mathbf A_N$ is $\epsilon_NY_N$, so their
polynomial span contains $Y_N$. The two limiting generators are equal,
and this direction is absent from the limiting polynomial space.

For (ii), take even $N\ge4$ and let $Y_N$ be a diagonal symmetry
with equal $+1$ and $-1$ multiplicities, chosen so that the
range of $P_N=e_{11}$ lies in its $+1$ eigenspace. The joint spectra are the same three points
\[
(1,1),\quad(0,1),\quad(0,-1),
\]
so $(P_N,Y_N)$ converges strongly to a commuting pair $(p,y)$ in $\C^3$ whose limiting trace assigns weight zero to $p$.
The two limiting marginal $C^*$-algebras are free for this nonfaithful trace: every centered element of
$C^*(p)$ is a multiple of $p$, and every alternating word containing such an element has trace zero.  For every $d\ge1$, both algebraic marginal polynomial spaces
have dimension two at the finite stages and at the limit; at $d=0$ both have dimension one. Nevertheless, the degree-two free-product
identity would require
\[
0=\operatorname{tr}_N(u_NP_Nu_N^*Y_N)
=\frac1N\langle u_Ne_1,Y_Nu_Ne_1\rangle.
\]
If $u_N\to1$, the final inner product tends to $1$, a contradiction.  In the limiting GNS space, $p=0$, so
\eqref{eq:gns-rank-stability} excludes the example.
\end{proof}

\begin{corollary}
\label{cor:strong-limit-exactification}
Let $P,Q$ be fixed finite-dimensional unital $C^*$-algebras with self-adjoint generating tuples $\mathbf x,\mathbf y$.  For
each $j$, let
\[
\pi_j:P\longrightarrow M_{N_j}(\C),
\qquad
\rho_j:Q\longrightarrow M_{N_j}(\C)
\]
be unital $^*$-embeddings, where $N_j\to\infty$, and set $A_j:=\pi_j(P)$ and $B_j:=\rho_j(Q)$.  Assume that there is $\delta>0$
such that every minimal projection of every $A_j$ and $B_j$ has normalized trace at least $\delta$, and that
\[
\varepsilon_m(A_j,B_j)\longrightarrow0
\qquad\text{for every fixed }m\ge1.
\]
If the joint tuples $(\pi_j(\mathbf x),\rho_j(\mathbf y))$ converge strongly, then, after discarding finitely many indices, there are unitaries
$u_j\in U(N_j)$ and nondecreasing integers $r_j\to\infty$ such that
\[
\norm{u_j-1}_\infty\longrightarrow0,
\qquad
u_jA_ju_j^*,B_j\text{ are }r_j\text{-independent},
\]
and the corrected tuples
\[
\bigl(\operatorname{Ad}(u_j)(\pi_j(\mathbf x)),\rho_j(\mathbf y)\bigr)
\]
converge strongly to the same limit.
\end{corollary}

\begin{proof}
For each fixed $n$, we have $\varepsilon_{2n}(A_j,B_j)\to0$. Hence, by Theorem~\ref{thm:fd-noncommutative-exactify},
we obtain $n$-independence with conjugating unitaries tending
to $1$ in operator norm. Letting the order grow sufficiently
slowly gives nondecreasing $r_j\to\infty$ and $\|u_j-1\|_\infty<1/r_j$
eventually. The generators
are uniformly bounded and change by amounts tending to zero in operator
norm, so each fixed polynomial retains its limiting trace and norm.
\end{proof}

\begin{theorem}
\label{thm:strong-matrix-models}
Let $(P,\tau_P)$ and $(Q,\tau_Q)$ be finite-dimensional unital $C^*$-algebras with faithful tracial states, and fix
self-adjoint generating tuples $\mathbf x$ of $P$ and $\mathbf y$ of $Q$. Let
\[
(\mathcal F,\varphi):=(P,\tau_P)*(Q,\tau_Q)
\]
be their reduced $C^*$-free product, with canonical embeddings $\iota_P$ and $\iota_Q$. Then there are strictly increasing integers
$N_j\to\infty$, unital $^*$-embeddings
\[
\pi_j:P\to M_{N_j}(\C),
\qquad
\rho_j:Q\to M_{N_j}(\C),
\]
and nondecreasing integers $r_j\to\infty$ such that, upon writing
\[
\tau_{P,j}:=\operatorname{tr}_{N_j}\circ\pi_j,
\qquad
\tau_{Q,j}:=\operatorname{tr}_{N_j}\circ\rho_j,
\]
the following hold.
\begin{enumerate}[label=\textup{(\roman*)}]
\item $\tau_{P,j}\to\tau_P$ and $\tau_{Q,j}\to\tau_Q$ in norm.
\item The unital $^*$-subalgebras $\pi_j(P)$ and $\rho_j(Q)$ are $r_j$-independent with respect to $\operatorname{tr}_{N_j}$.
\item For every fixed noncommutative $^*$-polynomial $R$,
\[
\operatorname{tr}_{N_j}\!\left(R(\pi_j(\mathbf x),\rho_j(\mathbf y))\right)
\longrightarrow
\varphi\!\left(R(\iota_P(\mathbf x),\iota_Q(\mathbf y))\right)
\]
and
\[
\norm{R(\pi_j(\mathbf x),\rho_j(\mathbf y))}_\infty
\longrightarrow
\norm{R(\iota_P(\mathbf x),\iota_Q(\mathbf y))}_{\mathcal F}.
\]
\end{enumerate}
Moreover, let
\[
(\mathcal F_j,\varphi_j):=(P,\tau_{P,j})*(Q,\tau_{Q,j})
\]
be the reduced $C^*$-free product formed from the finite-stage marginal states, with canonical embeddings $\iota_{P,j}$ and
$\iota_{Q,j}$. For every $j$ and every noncommutative $^*$-polynomial $R$ of total $^*$-degree at most $2r_j+1$,
\[
\operatorname{tr}_{N_j}\!\left(R(\pi_j(\mathbf x),\rho_j(\mathbf y))\right)
=
\varphi_j\!\left(R(\iota_{P,j}(\mathbf x),\iota_{Q,j}(\mathbf y))\right).
\]
\end{theorem}

\begin{proof}
Choose faithful representations
$\theta_j:P\to M_{p_j}(\C)$ and $\psi_j:Q\to M_{q_j}(\C)$
with positive integer multiplicities for the matrix summands, such that the
induced normalized traces converge to $\tau_P$ and $\tau_Q$,
respectively. We arrange that $p_j,q_j\to\infty$ and that $N_j=p_jq_j$ is
strictly increasing. The tensor-product representations give
\[
\mathbf X_j=(\theta_j\otimes1)(\mathbf x),\qquad
\mathbf Y_j^0=(1\otimes\psi_j)(\mathbf y)
\]
converging strongly to the canonical tuple in
$(P\otimes Q,\tau_P\otimes\tau_Q)$. Indeed, traces converge and all norms
are already those of this fixed finite-dimensional $C^*$-algebra.

By~\cite[Theorem~1.4]{collins-male}, in the form recalled in
Subsection~\ref{sec:stable-strong-models}, we can choose a realization
$(U_j)_j$ of a sequence of Haar unitary matrices such that
$(\mathbf X_j,\operatorname{Ad}(U_j)(\mathbf Y_j^0))$
converges strongly to free copies of $P,Q$.
We realize this limit using a Haar unitary $u$ free from
$(P\otimes Q,\tau_P\otimes\tau_Q)$ in a tracial $C^*$-algebra
with a faithful trace. The limiting copies of $P,Q$ are
$P\otimes1$ and $u(1\otimes Q)u^*$. Since the trace is faithful, the
$C^*$-algebra they generate is isometrically identified with $\mathcal F$.

Define the unital $^*$-subalgebras $A_j=(\theta_j\otimes1)(P)$ and
$B_j=\operatorname{Ad}(U_j)((1\otimes\psi_j)(Q))$.
The limiting marginal states are faithful, so the traces of minimal projections
of $A_j$ and $B_j$ have a uniform positive lower bound.
Since $P$ and $Q$ are finite-dimensional and the embeddings are
isometric, compactness of their unit balls makes moment convergence
uniform over contractions centered with respect to $\operatorname{tr}_{N_j}$. Thus $\varepsilon_m(A_j,B_j)\to0$ for every fixed $m$.
Now, by Corollary~\ref{cor:strong-limit-exactification}, we obtain unitaries $W_j\to1$
and nondecreasing $r_j\to\infty$. Set
\[
\pi_j=\operatorname{Ad}(W_j)\circ(\theta_j\otimes1),\qquad
\rho_j=\operatorname{Ad}(U_j)\circ(1\otimes\psi_j).
\]
Then (i) follows from the construction and the preservation of marginal
states, and (ii)--(iii) follow from the corollary.
Finally, grouping a monomial of degree at most $2r_j+1$ into marginal
blocks and applying Lemma~\ref{lem:quantitative-truncated-comparison}
at zero defect gives the asserted exact identity for the actual
finite-stage marginal states.
\end{proof}

For finite groups, we choose multiples of the regular representations,
so that the marginal states are the canonical group traces.

\begin{corollary}
\label{cor:exact-character-strong-models}
Let $G,H$ be nontrivial finite groups and set $\Gamma:=G*H$.  Denote its left regular representation by
$\lambda_\Gamma$ and its canonical trace by $\tau_\Gamma$.  There are $N_j\to\infty$ and $r_j\to\infty$, with unitary representations
\[
\Pi_j:\Gamma\longrightarrow U(N_j)
\]
such that, for every $x\in\C[\Gamma]$,
\[
\operatorname{tr}_{N_j}(\Pi_j(x))\longrightarrow\tau_\Gamma(x),
\qquad
\norm{\Pi_j(x)}_\infty\longrightarrow\norm{\lambda_\Gamma(x)}_\infty,
\]
and, for every $g\in\Gamma$ with $\ell_{\rm syl}(g)\le2r_j+1$,
\[
\operatorname{tr}_{N_j}(\Pi_j(g))=\delta_{g,e},
\]
where $\ell_{\rm syl}(g)$ is the number of nonidentity $G$- or $H$-factors
in the reduced expression for $g$, called its syllable length.
\end{corollary}

\begin{proof}
We use the construction of Theorem~\ref{thm:strong-matrix-models} for
$P=L(G)$ and $Q=L(H)$, choosing positive multiples of the regular
representations at every stage. Take the real and imaginary parts
of the canonical group unitaries as self-adjoint generators.
At every stage, the states $\tau_{P,j}$ and $\tau_{Q,j}$ are the canonical traces $\tau_G$ and $\tau_H$,
and the embeddings induce representations $\Pi_j$ of $G*H$. By the identification
$(L(G),\tau_G)*(L(H),\tau_H)\cong C^*_{\rm red}(G*H)$ and part~\textup{(iii)} of the
theorem, we obtain strong convergence. Finally, the identity $\operatorname{tr}_{N_j}(\Pi_j(g))=\delta_{g,e}$
follows from the last identity of the theorem, since a word of
syllable length at most $2r_j+1$
is a polynomial of degree at most $2r_j+1$ in the chosen generators.
\end{proof}

Barbieri and Jezernik~\cite[Theorems~1.1 and~5.7]{barbieri-jezernik}
constructed random-permutation models for free products of finitely many
finite groups. After composition with the standard representation,
the polynomial norms converge in probability to the corresponding norms
in the left regular representation, while the expected normalized traces
tend to the regular character.
In Proposition~\ref{prop:finite-image-log-character-models} below, we obtain representations $\Pi_N$ with $\operatorname{tr}_N(\Pi_N(g))=\delta_{g,e}$ for all $g$ of length at most a constant times $\log N$.

In Corollary~\ref{cor:exact-character-strong-models}, the dimensions may run
through all sufficiently large multiples of
$\ell:=\operatorname{lcm}(\abs{G},\abs{H})$.
To see this, realize multiples of the regular representations of $G,H$ in the same
fixed matrix algebra $M_\ell(\C)$ and amplify both to $M_{m\ell}(\C)$.
Since amplification preserves joint polynomial norms and normalized traces,
the Haar-conjugation argument in the proof of
Theorem~\ref{thm:strong-matrix-models} applies along the full sequence
$N=m\ell$. The marginal states are the canonical group traces, so by the same compactness
argument and Corollary~\ref{cor:strong-limit-exactification}, there are
corrections tending to $1$ in operator norm along this sequence.
The trace identities then follow as in the preceding proof.
The multiples of $\ell$ are precisely the dimensions in which the restrictions to $G$ and $H$ can
both have the normalized regular character, since a representation of a finite group with
this character is a multiple of its regular representation.

In the proof of Proposition~\ref{prop:finite-image-log-character-models}, we use the theorem of Bordenave and Collins~\cite[Theorem~9.1]{bordenave-collins-norm-polynomials} that independent uniform random permutation matrices, restricted to the orthogonal complement of the constant vectors and tensored with deterministic unitaries, converge strongly in probability to free Haar unitaries, provided that the number of letters grows fast enough compared with the dimension of the deterministic unitaries. In the statement of that proposition, $\norm{\cdot}$ is the operator norm.

\begin{proposition}[Finite-image models from tensor convergence]
\label{prop:finite-image-log-character-models}
For each integer $s\ge2$, there are $c_s>0$ and finite-image representations
$\Pi_N:\mathbb F_s\to U(N)$ in every sufficiently large dimension such
that
\[
\operatorname{tr}_N(\Pi_N(g))=\delta_{g,e}
\quad (|g|\le c_s\log N),
\qquad
\norm{\Pi_N(x)}\longrightarrow\norm{\lambda_{\mathbb F_s}(x)}
\quad (x\in\C[\mathbb F_s]).
\]
Here $|g|$ denotes reduced word length in the standard generators.
For nontrivial finite groups $G,H$, the same conclusions hold for
$\Gamma=G*H$ in every sufficiently large dimension divisible by
$\ell=\operatorname{lcm}(\abs{G},\abs{H})$, with syllable length in place
of reduced word length and a constant $c_\Gamma>0$.
\end{proposition}

\begin{proof}
Fix $s\ge2$.

\emph{Finite quotients of $\mathbb F_s$.}
For $p\in\{3,5\}$, set
\[
 A_p=\begin{pmatrix}1&p\\0&1\end{pmatrix},
 \qquad B_p=\begin{pmatrix}1&0\\p&1\end{pmatrix},
 \qquad H_{i,p}=B_p^iA_pB_p^{-i}\quad(0\le i<s).
\]
The projective actions of nonzero powers of $A_p$ send $|t|<1$ into
$|t|>1$, and those of $B_p$ send $|t|>1$ into $|t|<1$. By the ping-pong lemma, $A_p$ and $B_p$
freely generate a subgroup of rank $2$. The $H_{i,p}$ freely generate a subgroup
of rank $s$. Indeed, after adjacent powers with the same index are combined,
a reduced product expands into nonzero powers of $A_p$ separated by
nonzero powers of $B_p$. Choose $C>1$ bounding the maximum absolute
row sum of the $H_{i,p}$ and their inverses, for both values of $p$.
Let $Q_{p,k}$ be the image modulo $p^k$ of the group generated by the matrices $H_{i,p}$, $0\le i<s$. All generators are $1$
modulo $p$, so $Q_{p,k}$ is a $p$-group and
\[
 p^{k-1}\le |Q_{p,k}|\le p^{3(k-1)}.
\]
Indeed, each successive congruence kernel has at most $p^3$ elements,
represented by matrices $1+p^jX$ with $\operatorname{Tr}(X)=0$ modulo
$p$, while $A_p$ has order $p^{k-1}$ modulo $p^k$.
A nonidentity word of length $L$ gives a nonidentity integer matrix
with row norm at most $C^L$. If $C^L+1<p^k$, it therefore remains
nonidentity in $Q_{p,k}$. Margulis~\cite{margulis-girth} used finite quotients of free subgroups of
the modular group to construct Cayley graphs whose girth is logarithmic in the number of vertices.

\emph{Choice of the block dimensions.}
For a large target dimension $N$, set
\[
 k_p=\left\lfloor\frac{\log N}{12\log p}\right\rfloor,
 \quad a=|Q_{3,k_3}|,
 \quad b=|Q_{5,k_5}|,
 \quad
 R_N=\left\lfloor
 \frac{\min_{p\in\{3,5\}}k_p\log p-\log2}{\log C}
 \right\rfloor.
\]
Then $a,b$ are coprime, both tend to infinity, and $a,b\le N^{1/4}$.
The regular quotient representations $\sigma_{3,N},\sigma_{5,N}$
have character zero on every nonidentity word of length at most $R_N$,
since $C^{R_N}+1\le p^{k_p}/2+1<p^{k_p}$. Moreover,
$R_N\ge(\log N)/(12\log C)-O_s(1)$.
Since $a$ and $b$ are coprime, every integer $n\ge(a-1)(b-1)$ has the form $n=ax+by$ with integers $x,y\ge0$.
Applying this to $n=N-(a+b)\lceil N^{1/4}\rceil$, we can write
\[
 N=ad+be,\qquad d,e\ge \lceil N^{1/4}\rceil,
\]
whenever $N\ge(a+b)\lceil N^{1/4}\rceil+(a-1)(b-1)$.
This inequality holds for all large $N$, since its right side is at most $3N^{1/2}+2N^{1/4}$.

\emph{Tensor products with random permutations.}
Let $U_1,\ldots,U_s$ be independent uniform permutations on $d+1$
letters, restricted to the $d$-dimensional orthogonal complement of
the constants, and let $V_1,\ldots,V_s$ be the analogous matrices in
dimension $e$. Send the free generators $g_i$ to
\[
 W_i=(U_i\otimes\sigma_{3,N}(g_i))
       \oplus(V_i\otimes\sigma_{5,N}(g_i)).
\]
This defines a finite-image representation in dimension $N$.
The characters of $\sigma_{3,N}$ and $\sigma_{5,N}$ vanish on all
nonidentity words of length at most $R_N$. Thus each tensor summand,
and hence their direct sum, has the required character on the
radius-$R_N$ ball for every sample.
By the permutation case of the theorem of Bordenave and
Collins~\cite[Theorem~9.1]{bordenave-collins-norm-polynomials}, both tensor
tuples converge strongly in probability to free Haar unitaries, because
\[
 \frac{\log(d+1)}{\sqrt{(\log\log a)\log a}}\longrightarrow\infty,
 \qquad
 \frac{\log(e+1)}{\sqrt{(\log\log b)\log b}}\longrightarrow\infty.
\]
Here $d,e\ge N^{1/4}$ and $a,b\le N^{1/4}$.
The coefficient dimensions $a,b$ may take the same value for several $N$, but the theorem still applies.
Indeed, suppose that for some polynomial and some $\epsilon>0$, the probability that its norm in the
first tensor tuple differs from the limit by more than $\epsilon$ does not tend to zero. Then this
probability is bounded below along a subsequence, and after passing to a further subsequence, $a$ is
strictly increasing. For the integers $n$ that are not values of $a$ along this subsequence, we take
identity matrices of size $n$ and random permutations on $n$ letters. The growth condition above still
holds, which contradicts the theorem. The same argument applies to the second tuple.

The norm of a polynomial in $W_1,\ldots,W_s$ is the maximum of its norms in the two tensor tuples, so it
converges in probability to the same limit. Hence, for every finite set of polynomials with coefficients
in $\mathbb Q+i\mathbb Q$ and every $\epsilon>0$, with positive probability all their norms are within
$\epsilon$ of the limits, in every sufficiently large dimension $N$. Letting the set grow and $\epsilon$
decrease sufficiently slowly, we choose one such sample in each dimension. Since the norm of a polynomial
in unitaries is at most the sum of the absolute values of its coefficients, norm convergence for every
polynomial follows by approximation. Taking, for example, $c_s=1/(24\log C)$ proves the assertion for free groups.

\emph{Free products of finite groups.}
Now let $\Gamma=G*H$ as in the statement. Let $G$ and $H$ each act
freely on the same $\ell$-point set, as disjoint unions of regular
actions. Every $\Gamma$-orbit has size divisible by both group orders,
so this action is transitive. Its stabilizer $F$ has index $\ell$ and
meets every conjugate of either factor trivially. It acts freely on
the Bass--Serre tree of the free product $G*H$, so it is free, of rank
$r=1+\ell-\ell/|G|-\ell/|H|$.
For $r\ge2$, we use the representations $\Pi_n$ just constructed for $F$. If
$r=1$, then $G\cong H\cong C_2$, and we send a generator of $F\cong\mathbb Z$
to the diagonal unitary whose eigenvalues are the $n$th roots of unity,
each occurring once. These representations have finite image, have
zero trace at nonzero exponents of absolute value less than $n$, and
converge strongly because the $n$th roots of unity become dense in
the circle.

Fix coset representatives $t_1,\ldots,t_\ell$ and a free basis of $F$.
For $u\in(G\cup H)\setminus\{e\}$, write
$ut_i=t_j\kappa(u,i)$, where $\kappa(u,i)\in F$, and choose $K\ge1$
bounding the free-basis lengths of these finitely many elements.
A word $g$ of syllable length $L$ then has induced blocks $\Pi_n(h)$
with $|h|_F\le KL$. On a nonzero diagonal block, $h=t_i^{-1}gt_i$, which
is nonidentity when $g\ne e$. Thus
$\operatorname{Ind}_F^\Gamma\Pi_n$ has the exact regular character
on a syllable ball of radius at least a positive constant times
$\log(\ell n)$, for all large $n$.
Strong convergence persists under fixed matrix amplification. Indeed, the map
from $C^*_{\rm red}(F)$ defined by the $\Pi_n$ into any nonprincipal norm
ultraproduct is an isometric $^*$-homomorphism and hence completely
isometric. Applying this to every subsequence gives ordinary
matrix-norm convergence.
Since $\operatorname{Ind}_F^\Gamma\lambda_F\cong\lambda_\Gamma$,
finite-index induction preserves strong convergence by
\cite[Lemma~7.1]{louder-magee-strong-representations}.
The induced representation has finite image and dimension $\ell n$,
giving every sufficiently large dimension divisible by $\ell$.
\end{proof}

Note that the representations $\Pi_N$ of $\mathbb F_s$ constructed in the proof of
Proposition~\ref{prop:finite-image-log-character-models} are reducible: the projection onto the first
direct summand commutes with every generator. Hence the average of the conjugations by the generators and
their inverses has norm $1$ on the traceless matrices with
the $L^2$-norm. In contrast, the clock representations of
Theorem~\ref{thm:irreducible-clock-strong-models} satisfy Hastings's
asymptotically optimal bound for this average~\cite[p.~1 and equations~(10)--(12)]{hastings-quantum-expanders}.

\begin{theorem}
\label{thm:multi-gns-rank-strong-exactification}
Fix $s\ge2$ and positive integers $N_j$.  For $1\le i\le s$, let $\mathbf X_j^{(i)}$ be tuples of fixed length in $M_{N_j}(\C)$, and assume that the joint
tuple converges strongly to mutually free limiting families $\mathbf x^{(1)},\ldots,\mathbf x^{(s)}$.  Assume that, for every $i$ and every fixed degree, the dimension equalities
\eqref{eq:gns-rank-stability} hold eventually.  Set
\[
A_{i,j}:=C^*(\mathbf X_j^{(i)}),
\]
and let $\varphi_j^{(s)}$ be the state of the iterated reduced $C^*$-free product of the finite-stage marginal
$C^*$-probability spaces $(A_{i,j},\operatorname{tr}_{N_j}|_{A_{i,j}})$, with canonical embeddings
$\iota_{i,j}:A_{i,j}\to\mathop{*}_{r=1}^s(A_{r,j},\operatorname{tr}_{N_j}|_{A_{r,j}})$.

After deleting finitely many indices, there are unitaries $u_{i,j}\in U(N_j)$ and nondecreasing odd integers $D_j\to\infty$ such that
$u_{s,j}=1$ for every $j$,
\[
\max_i\norm{u_{i,j}-1}_\infty\longrightarrow0,
\]
the corrected joint tuple has the same strong limit, and
\[
\operatorname{tr}_{N_j}R\bigl((\operatorname{Ad}(u_{i,j})(\mathbf X_j^{(i)}))_{i=1}^s\bigr)
=\varphi_j^{(s)}R\bigl((\iota_{i,j}(\mathbf X_j^{(i)}))_{i=1}^s\bigr)
\]
for every $^*$-polynomial $R$ of total degree at most $D_j$.
\end{theorem}

\begin{proof}
To prove \eqref{eq:gns-rank-stability} for any tuple obtained by grouping some of the marginal families, fix a degree $d$ and such a grouping. In the tracial GNS space
of each limiting marginal $C^*$-algebra, choose an orthogonal real basis of
self-adjoint polynomial vectors compatible with the degree filtration,
and assign to each vector the least degree in which it occurs.
Consider the products of centered basis vectors in which consecutive factors come from different families, with the
sum of their degrees at most $d$. Together with $1$, they span the image of the joint polynomial space through degree $d$
in the limiting tracial GNS space.
By freeness, these products are identified with nonzero orthogonal tensors
in the free-product GNS decomposition.

For all sufficiently large indices, the corresponding marginal
polynomials evaluate to bases in the matrices, by the dimension hypothesis,
as in the proof of Theorem~\ref{thm:gns-rank-strong-exactification}.
Expanding each marginal block in these bases shows that
the products just described span the polynomial space of the grouped
tuple. Their Gram matrices converge to the positive-definite Gram
matrix at the free limit, so they are also linearly independent.
This proves \eqref{eq:gns-rank-stability} for the grouped tuple.
The same argument applies after some families have been conjugated by
unitaries tending to $1$: these conjugations preserve each marginal
polynomial dimension and the joint moment and strong limits.

We first apply Theorem~\ref{thm:gns-rank-strong-exactification} to families
$1$ and $2$. The corrected pair converges strongly and agrees
with its free-product law through an odd degree tending to infinity.
Then we apply the same theorem to this pair, viewed as one tuple, and family
$3$, conjugating the pair by a common unitary close to the identity.
This preserves all trace identities already obtained for the pair.
For a polynomial whose degree is within both bounds, the free-product
formula for the new target trace uses pair moments only within the
earlier bound.
These already agree with the two-family free-product state, so the
new target equals the iterated free-product trace of all three
families. We continue in this way through family $s$.

We take the minimum of the finitely many degree bounds, which still tends
to infinity. Each family is conjugated by a product of at most $s-1$
unitaries tending to $1$, so the product also tends to $1$. The last
family remains fixed. The corrected generators differ from the
original ones by amounts tending to zero in operator norm, which
preserves the joint strong limit.
\end{proof}

In Theorems~\ref{thm:gns-rank-strong-exactification} and
\ref{thm:multi-gns-rank-strong-exactification}, we can also compare the corrected and original tuples on polynomials of increasing degree. Write
\[
\mathbf Z_j:=(\mathbf X_j^{(1)},\ldots,\mathbf X_j^{(s)}),
\qquad
\widetilde{\mathbf Z}_j
:=(\operatorname{Ad}(u_{i,j})(\mathbf X_j^{(i)}))_{i=1}^s,
\]
interpreting Theorem~\ref{thm:gns-rank-strong-exactification} as the case $s=2$, $u_{1,j}=u_j$, and $u_{2,j}=1$.  If $M\ge1$ bounds all original and corrected
generator norms, choose the conjugating unitaries at degree $D_j$ so that
\[
\delta_j:=\max_i\norm{u_{i,j}-1}_\infty,
\qquad
D_jM^{D_j}\delta_j\longrightarrow0.
\]
Telescoping each monomial gives
\[
\sup_{\substack{\deg R\le D_j\\ \norm{R}_{\ell^1}\le1}}
\norm{R(\widetilde{\mathbf Z}_j)-R(\mathbf Z_j)}_\infty
\le 2D_jM^{D_j}\delta_j
\longrightarrow0,
\]
where $\norm{R}_{\ell^1}$ is the sum of the absolute values of its formal monomial coefficients.
However, convergence to the limiting norms is proved only for fixed polynomials.

\begin{corollary}
\label{cor:gue-growing-exactification}
Fix $s\ge2$, and let $G_N^{(1)},\ldots,G_N^{(s)}$ be independent normalized GUE matrices
in $M_N(\C)$.  On a probability-one event,
after deleting finitely many $N$, there are unitaries $u_{i,N}\in U(N)$, with one family left fixed, and nondecreasing odd integers
$D_N\to\infty$ such that
\[
\max_i\norm{u_{i,N}-1}_\infty\longrightarrow0,
\]
the corrected tuple $(\operatorname{Ad}(u_{i,N})(G_N^{(i)}))_{i=1}^s$ converges strongly to a free semicircular family, and its normalized trace agrees through total degree
$D_N$ with the reduced free-product state determined by the normalized
trace states of the $s$ individual matrices.
\end{corollary}

\begin{proof}
By the norm convergence theorem of Haagerup and
Thorbj{\o}rnsen~\cite[Theorem~7.1]{haagerup-thorbjornsen} and Voiculescu's asymptotic freeness theorem~\cite[Theorem~2.2]{voiculescu-limit-laws}, in the almost-sure form proved by
Hiai and Petz~\cite[Corollary~3.3]{hiai-petz} and, by a different method, by Thorbj{\o}rnsen~\cite[Corollary~3.9]{thorbjornsen-mixed-moments}, the
GUE tuple almost surely
converges strongly along
the full sequence to free semicircular variables. Both the norm convergence and the almost-sure convergence of moments need independence only within each dimension.
Almost surely, every $G_N^{(i)}$ has simple spectrum.
For fixed $d<N$, the powers $1,G_N^{(i)},\ldots,(G_N^{(i)})^d$ are therefore linearly independent, while the corresponding
powers of a semicircular variable are linearly independent in $L^2$ because its law is diffuse.  Hence
\eqref{eq:gns-rank-stability} holds for every individual family, and
Theorem~\ref{thm:multi-gns-rank-strong-exactification} applies on
this event. Conjugation leaves each matrix's eigenvalues unchanged,
but can alter the joint probability law.
\end{proof}

In Theorem~\ref{thm:gue-relative-exactification}, we obtain Borel choices of
the conjugating unitaries and exact identities through a degree proportional
to $\log N$. There we also allow several GUE matrices in each
family and amalgamation over a finite-dimensional unital $C^*$-algebra.

\section{Correction of amalgamated moments}
\label{sec:relative-exactification}

We now assume that the marginal $C^*$-algebras contain a common
finite-dimensional unital $^*$-subalgebra $D$. We correct the mixed moments to those of
their reduced amalgamated $C^*$-free product over $D$ by conjugating by unitaries that commute with $D$,
which preserves the $D$-valued law of each family.
Our polynomials have coefficients in $D$: a monomial $d_0Z_1d_1\cdots Z_kd_k$, where $d_j\in D$ and
each $Z_j$ is a tuple entry or its adjoint, has degree $k$.

In Theorem~\ref{thm:relative-quantitative-exactification}, we prove an amalgamated analogue of Theorem~\ref{thm:fd-noncommutative-exactify}, using a free semicircular family over $D$ as the reference point.
When the moments through degree $2L$ are close to those of this family, we obtain the amalgamated
identities through degree $L$ by conjugating by unitaries whose distance
from $1$ is at most $K^L$ times the moment defect, with $K$ independent of $L$.
Combining this with the GUE moment estimates in
Section~\ref{sec:gue-bounds}, we obtain degrees proportional to $\log N$. In
Theorem~\ref{thm:relative-gns-rank-exactification}, we treat arbitrary
limiting families free over $D$, assuming the dimension hypothesis of
Subsection~\ref{sec:stable-strong-models}.

Since the conjugating unitaries must commute with $D$, we
differentiate only in self-adjoint directions in the relative commutant $D'$. To do this, we pass to a corner in which the compression $D_0$ of $D$
is abelian, so that $D_0'$ consists of the elements that are block diagonal with respect to
the minimal projections of $D_0$.

\emph{Standing assumptions.}
Fix
\[
D=\bigoplus_{v=1}^b M_{d_v}(\C),\qquad s\ge2,\qquad m_i\ge1,
\qquad m=\sum_{i=1}^s m_i,\qquad q=\dim_{\C}D.
\]
Choose matrix units $f^v_{jk}$ for $D$.  Throughout the quantitative assertions, every minimal projection of $D$ has
trace at least a fixed $\delta>0$.  Set
\[
\mathcal D=\{1\}\cup\{f^v_{jk}:v,j,k\}.
\eqtag\label{eq:rel-constants}
\]

Let $M$ be a unital $C^*$-algebra containing $D$ as a unital $^*$-subalgebra, equipped with a faithful tracial state $\tau$, and let
$\mathbf H_i=(H_{i,1},\ldots,H_{i,m_i})$ be self-adjoint tuples in $M$ with
$\max_{i,a}\norm{H_{i,a}}_\infty\le3$.  Write $E_D:M\to D$ for the trace-preserving conditional expectation.
Let $\mathbf S=(S_1,\ldots,S_m)$ be a free standard semicircular family free from $(D,\tau|_D)$, with reference
trace $\tau_{\rm ref}$, and write $\mathbf H=(H_\alpha)_{\alpha=1}^m$ for the concatenation of the tuples
$\mathbf H_1,\ldots,\mathbf H_s$.  Define
\[
\begin{split}
\Delta_{2L}(\mathbf H;D)=\max\bigl\{&
\bigl|\tau(d_0H_{\alpha_1}d_1\cdots H_{\alpha_k}d_k)\\
&\quad-\tau_{\rm ref}(d_0S_{\alpha_1}d_1\cdots S_{\alpha_k}d_k)\bigr|:\\[-2pt]
&0\le k\le2L,\quad d_j\in\mathcal D,\quad 1\le\alpha_j\le m\bigr\}.
\end{split}
\eqtag\label{eq:rel-raw-defect}
\]
Set $A_i=C^*(D,\mathbf H_i)$ and let $(\mathcal A,\mathcal E)=*_{D,i}(A_i,E_D|_{A_i})$ be their reduced amalgamated $C^*$-free product, with canonical
embeddings $\jmath_i$. This construction is due to Voiculescu~\cite[Section~5]{voiculescu-symmetries} (see
\cite[Section~2]{ivanov-amalgamated}). The free-product moments in
Theorem~\ref{thm:relative-quantitative-exactification} are those
determined by the actual $D$-valued marginal laws of the given tuples,
not by the semicircular reference family.

\begin{theorem}
\label{thm:relative-quantitative-exactification}
With the notation and standing assumptions above, there is a constant $K\ge1$ depending only on $s$, $m$, $q$ and
$\delta$ for which the following holds. Fix $L\ge2$ and assume that
$K^L\Delta_{2L}(\mathbf H;D)\le1$.
Then there are unitaries $u_i\in D'\cap U(M)$, with $u_s=1$, such that
\[
\max_i\norm{u_i-1}_\infty\le K^L\Delta_{2L}(\mathbf H;D).
\eqtag\label{eq:rel-motion}
\]
Set $\widetilde H_{i,a}=u_iH_{i,a}u_i^*$. Then, for every
$0\le\ell\le L$, every choice of indices, and all
$d_0,\ldots,d_\ell\in D$,
\[
\begin{split}
&E_D(d_0\widetilde H_{i_1,a_1}d_1\cdots\widetilde H_{i_\ell,a_\ell}d_\ell)\\
&\qquad=\mathcal E(d_0\jmath_{i_1}(H_{i_1,a_1})d_1\cdots
\jmath_{i_\ell}(H_{i_\ell,a_\ell})d_\ell).
\end{split}
\eqtag\label{eq:rel-exact-moments}
\]
The complete $D$-valued law of each family is preserved.  The unitaries can be constructed in
$C^*(D,\mathbf H_1,\ldots,\mathbf H_s)$, and in matrix algebras the construction is Borel in the matrix entries.
\end{theorem}

\begin{proof}
Write $e=\Delta_{2L}(\mathbf H;D)$ and set
\[
a=(100+s+m+q+\delta^{-1})^{100},\qquad K=a^{40}.
\]

\emph{Reduction to a corner.}
For a polynomial $P$ with coefficients in $D$ in formal variables $X_1,\ldots,X_m$, let $|P|_3$ be the infimum of
$\sum_w|c_w|3^{\deg w}$ over all ways of writing $P=\sum_wc_ww$, where $c_w\in\C$ and each $w$ is a word
$d_0X_{\alpha_1}d_1\cdots X_{\alpha_k}d_k$ with $d_j\in\mathcal D$ and $k\ge0$. We call $|P|_3$ the coefficient norm of $P$.
Since a product of two elements of $\mathcal D$ is zero or an element of $\mathcal D$, we have $|PQ|_3\le|P|_3|Q|_3$ and $|P^*|_3=|P|_3$.
Since the $H_\alpha$ and the $S_\alpha$ have norm at most $3$, we have $\norm{P(\mathbf H)}_\infty\le|P|_3$ and $\norm{P(\mathbf S)}_\infty\le|P|_3$, and the traces satisfy
\[
\abs{\tau(P(\mathbf H))-\tau_{\rm ref}(P(\mathbf S))}\le e|P|_3
\qquad(\deg P\le2L).
\eqtag\label{eq:rel-coefficient-control}
\]

Set
\[
\begin{gathered}
p=\sum_v f^v_{11},\qquad M_0=pMp,\qquad t=\tau(p),\qquad
\tau_0=t^{-1}\tau|_{M_0},\\
D_0=pDp=\bigoplus_v\C e_v,\qquad e_v=f^v_{11},\qquad\beta_v=\tau_0(e_v).
\end{gathered}
\eqtag\label{eq:rel-full-corner}
\]
Both $t$ and every $\beta_v$ are at least $\delta$.  A unitary $v\in D_0'\cap U(M_0)$ lifts by
\[
\mathcal L(v)=\sum_{w=1}^b\sum_{j=1}^{d_w}f^w_{j1}v f^w_{1j}
\in D'\cap U(M),\qquad
\norm{\mathcal L(v)-1}_\infty=\norm{v-p}_\infty.
\eqtag\label{eq:rel-unitary-lift}
\]
The formula places a copy of $e_wve_w$ in each diagonal corner of
the $w$th matrix summand. It gives
$f^w_{1j}\mathcal L(v)=vf^w_{1j}$ and
$\mathcal L(v)^*f^w_{k1}=f^w_{k1}v^*$.
Inserting $1=\sum_{w,j}f^w_{j1}pf^w_{1j}$ between successive letters
expresses every moment in terms of corner polynomials of no greater
degree. In these polynomials, the $i$th family is conjugated by its
chosen unitary in $M_0$.
For $x\in M_0$, bimodularity gives
$E_{D_0}(x)=pE_D(x)p=E_D(x)$.  In particular, a $D_0$-centered corner element is $D$-centered, so the unital $C^*$-subalgebras
$p\jmath_i(A_i)p$ are free over $D_0$ in $p\mathcal A p$.  For arbitrary $x\in M$, writing $\theta_v=\tau(f^v_{11})$, we have
\[
E_D(x)=\sum_{v,j,k}\theta_v^{-1}\tau(f^v_{kj}x)f^v_{jk},\qquad
\tau(f^v_{kj}x)=\tau(f^v_{1j}xf^v_{k1}).
\eqtag\label{eq:rel-expectation-reconstruction}
\]
Thus the corner trace identities determine all $D$-valued expectations in the original $C^*$-algebra $M$.

\emph{Reference bases.}
We realize the reference family on a Fock space, as is standard for operator-valued semicircular families
(see Speicher~\cite{speicher-operator-valued} and Shlyakhtenko~\cite{shlyakhtenko-a-valued}).
For family $i$, we use the correspondence
\[
\mathcal H_i=(D\otimes_{\C}D)^{\oplus m_i},\qquad
\langle a\otimes b,c\otimes d\rangle_D=\tau(a^*c)b^*d.
\]
The full Fock construction on $\mathcal H_i$ yields a Hilbert $C^*$-module
over $D$. Creation by $1\otimes1$ in the $a$th summand of $\mathcal H_i$,
together with its adjoint, gives $S_a$.
After tensoring over $D$ with $L^2(D,\tau)$, the degree-$k$ subspace is the
direct sum of copies of $L^2(D)^{\otimes(k+1)}$ indexed by strings
of $k$ generator indices. The noncrossing-pairing expansion shows that the
vacuum expectation, followed by $\tau|_D$, is the tracial semicircular
reference state. For $k\ge2$, the Wick
polynomials satisfy
\[
\begin{split}
W_{a_1\cdots a_k}(d_0,\ldots,d_k)
={}&d_0S_{a_1}W_{a_2\cdots a_k}(d_1,\ldots,d_k)\\
&-\mathbf1_{a_1=a_2}d_0\tau(d_1)
 W_{a_3\cdots a_k}(d_2,\ldots,d_k),
\end{split}
\eqtag\label{eq:rel-wick-recurrence}
\]
with $W_a(d_0,d_1)=d_0S_ad_1$ and $W_\varnothing(d_0)=d_0$.
Creation adds the first tensor coordinate; the second term cancels the annihilation contribution.  Hence
\[
\begin{split}
&\langle W_{\mathbf a}(d_0,\ldots,d_k),W_{\mathbf b}(d'_0,\ldots,d'_\ell)\rangle\\
&\qquad=\mathbf1_{k=\ell}\mathbf1_{\mathbf a=\mathbf b}
\prod_{j=0}^k\tau(d_j^*d'_j).
\end{split}
\eqtag\label{eq:rel-wick-orthogonality}
\]
The recurrence is triangular in degree, so these polynomials span
the entire polynomial space at each degree. Expanding the recurrence by contracting
disjoint adjacent pairs of generator occurrences shows that taking
adjoints reverses the generator indices and replaces the coefficient
list $d_0,\ldots,d_k$ by $d_k^*,\ldots,d_0^*$.
In the Fock space built in the same way from $\bigoplus_i\mathcal H_i$, a product of Wick polynomials of
positive degree in which consecutive factors come from different families sends the vacuum to the concatenation
of their tensors, since the annihilation terms vanish at each change
of family. This tensor is orthogonal to the vacuum sector, so
the $D$-valued vacuum expectation of the product is zero. Expanding in Wick polynomials
and taking norm limits, we obtain that the marginal $C^*$-algebras are free
over $D$.

We now choose bases of the corner polynomial spaces in the reference model. In the Wick polynomials
$W_{a_1\cdots a_k}(d_0,\ldots,d_k)$ with $k\ge1$, we let $d_0$ and $d_k$ range over the normalized matrix units
in $pD$ and $Dp$, respectively, and $d_1,\ldots,d_{k-1}$ over a self-adjoint basis of $D$ that is orthonormal
for $\tau$, consisting of normalized diagonal matrix units and normalized real and imaginary
parts of off-diagonal matrix units. We multiply the resulting Wick polynomials
by $\sqrt t$ to normalize their reference values for the corner trace
$t^{-1}\tau_{\rm ref}$.

We sort these polynomials by their left and right block supports
$e_v,e_w$. In each diagonal block, we replace every pair $R\ne R^*$
by $(R+R^*)/\sqrt2$ and $(R-R^*)/(i\sqrt2)$, and keep
the self-adjoint polynomials. For each pair of opposite off-diagonal
blocks, we keep the basis in one and take its adjoint in the other.
At degree zero, we include $e_v/\sqrt{\beta_v}$. Thus the bases of the diagonal
blocks are self-adjoint, and the bases with supports $(e_v,e_w)$ and
$(e_w,e_v)$ are adjoints of one another. A product of these basis elements
vanishes unless successive block indices match.

Fix a family and a block. For $k\ge1$, every degree-$k$ basis
polynomial $R_\lambda$ and the number of basis vectors through degree
$k$ satisfy
\[
|R_\lambda|_3\le a^k,\qquad \#\{\text{basis vectors through degree }k\}\le a^k.
\eqtag\label{eq:rel-reference-cost}
\]
At degree $k$, the number of basis vectors is bounded by a fixed
constant times $m_i^kq^{k+1}$. In each normalized basis vector of
$pD$, $D$, or $Dp$ used above, the sum of the absolute values of the
matrix-unit coefficients is at most $2\delta^{-1/2}$.
By the recurrence \eqref{eq:rel-wick-recurrence}, the coefficient sum of a degree-$k$ Wick polynomial is at most $4^k$ times the product of these sums.
Passing to the self-adjoint combinations adds a factor of at most $\sqrt2$. The choice of $a$
also accounts for the degree-zero normalization and summation over lower degrees.

\emph{Bases adapted to degree.}
Write $\|\cdot\|_{\rm F}$ for the Frobenius norm, i.e., the
unnormalized Hilbert--Schmidt norm on scalar matrices. For each
family and block, we order the basis by degree and let $G_k$ be the
Gram matrix of the evaluated basis polynomials of degree at most
$k$. Then
\[
\|G_k-I\|_{\rm F}\le a^{4k}e\qquad(k\le L)
\eqtag\label{eq:rel-prefix-gram}
\]
by \eqref{eq:rel-coefficient-control}, \eqref{eq:rel-reference-cost}, and $t\ge\delta$.
For a Gram matrix $G=I+E$ with $\|E\|_{\rm F}\le1/8$, write
$G=C^*C$, where $C=I+V$ is upper triangular with positive diagonal. Then
$V=\mathcal T(E-V^*V)$, where $\mathcal T$ retains the strict upper
triangle and half the diagonal. The map on the right is a contraction
on $\|V\|_{\rm F}\le2\|E\|_{\rm F}$ with constant at most
$4\|E\|_{\rm F}$. This and the Neumann series give
\[
\|C-I\|_{\rm F}\le2\|E\|_{\rm F},\qquad
\|C^{-1}-I\|_{\rm F}\le4\|E\|_{\rm F}.
\eqtag\label{eq:rel-cholesky}
\]
For each $k$, we multiply the row of reference polynomials of degree at
most $k$ by $C^{-1}$ and obtain polynomials $U_\lambda$. Since $C^{-1}$ is upper triangular, this
preserves degrees and the spans through every lower degree. Since the triangular factor with positive diagonal is unique, the factor
of $G_k$ is the leading block of the factor of $G_L$, and real or complex-conjugate
Gram matrices have real or complex-conjugate factors. Real triangular
factors on diagonal blocks and conjugate factors on opposite blocks
preserve the adjoint symmetries. At each degree, the evaluated $U_\lambda$
form a basis of the entire marginal polynomial space in that block.
Since the trace is faithful, $L^2$ expansions in these bases are
identities in $M_0$, and the positive-degree basis vectors are
$D_0$-centered.
Each column of $C^{-1}-I$ has at most $a^k$ entries. Thus, by
Cauchy--Schwarz, \eqref{eq:rel-reference-cost},
\eqref{eq:rel-prefix-gram}, and \eqref{eq:rel-cholesky}, we have
$|U_\lambda-R_\lambda|_3\le4a^{k/2}a^{4k}e\,a^k\le a^{6k}e$.
Together with $|R_\lambda|_3\le a^k$ and $e\le a^{-40L}$, this gives
\[
|U_\lambda|_3\le a^{2k},\qquad
|U_\lambda-R_\lambda|_3\le a^{6k}e
\quad(\deg\lambda=k\ge1).
\eqtag\label{eq:rel-basis-cost}
\]
Set
\[
x_\lambda=a^{-2k}U_\lambda(\mathbf H),\qquad
y_\lambda=a^{-2k}R_\lambda(\mathbf S),\qquad \kappa=a^{-2}.
\eqtag\label{eq:rel-scaled-letters}
\]
The $x_\lambda$ and $y_\lambda$ are $D_0$-centered contractions
of $L^2$-norm $\kappa^k$ for $\tau_0$ and $t^{-1}\tau_{\rm ref}$,
respectively.
The formal polynomials $a^{-2k}U_\lambda$ and $a^{-2k}R_\lambda$
have coefficient norms at most one, and their difference has
coefficient norm at most $a^{4k}e$. We keep these polynomials fixed
when the families are conjugated.

\emph{The moment equations and their gradients.}
We consider the products $w=x_1\cdots x_\ell$, $\ell\ge2$, of
positive-degree scaled basis elements $x_\lambda$ from
\eqref{eq:rel-scaled-letters} such that
$x_j\in e_{v_{j-1}}M_0e_{v_j}$ with $v_\ell=v_0$,
the total degree is at most $L$, and consecutive letters,
including the last and first, come from different families.
We choose one representative from each class under cyclic rotation and
adjoint. If a cyclic class differs from its adjoint class, our equations are
$\Ree\tau_0(w)=0$ and $\Imm\tau_0(w)=0$; otherwise, the only equation is $\Ree\tau_0(w)=0$.
There are at most $a^{2k}$ letters of degree $k$ and at most
$2^{n-1}$ ways to express a total degree $n$ as a sum of positive
degrees. Summing over $n\le L$ and counting real and imaginary parts
shows that the number $d$ of real equations is at most $a^{3L}$.

Let $c_j$ be the index of the family containing $x_j$, set
$c_0=c_\ell$, and let $R_j(w)$ be the rotation beginning at $x_j$.
For $i<s$, set
\[
A_i(w)=\sum_{j=1}^{\ell}
(\mathbf1_{c_j=i}-\mathbf1_{c_{j-1}=i})R_j(w),\qquad
\Xi_i(w)=-\frac{i}{2}A_i(w).
\eqtag\label{eq:rel-cyclic-gradients}
\]
Conjugating family $i$ by $e^{ih_i}$, for $h_i\in(D_0'\cap M_0)_{\rm sa}$ and $h_s=0$, gives
\[
D\tau_0(w)[h]=-2\sum_{i<s}\tau_0(h_i\Xi_i(w)).
\eqtag\label{eq:rel-cyclic-derivative}
\]
To obtain this formula, we expand each commutator and use cyclicity to
place $h_i$ first in each trace. Every rotated word belongs to
a diagonal block $e_vM_0e_v$, so each $\Xi_i(w)$ commutes with $D_0$.
For each chosen real equation $\alpha$, let $g_{\alpha,i}$ be the corresponding
real or imaginary part of $\Xi_i(w)$. Each $g_{\alpha,i}$ is self-adjoint
and has coefficient norm and operator norm at most $L/2$.

In the reference model, distinct reduced block words are orthogonal.
Successive conditional expectations give, for a word with block
indices $v_0,\ldots,v_\ell$,
\[
\norm{y_1\cdots y_\ell}_2^2
=\frac{\prod_{j=1}^{\ell}\norm{y_j}_2^2}
{\prod_{j=1}^{\ell-1}\beta_{v_j}}
\ge\kappa^{2L},
\eqtag\label{eq:rel-path-norm}
\]
since $E_{D_0}(y_j^*y_j)=\norm{y_j}_2^2\beta_{v_j}^{-1}e_{v_j}$.
Set $\Xi(w)=(\Xi_i(w))_{i<s}$, with componentwise adjoint and the
direct-sum $L^2$-norm. Equal rotations of a periodic word have the
same first and last family indices, so their coefficients in
\eqref{eq:rel-cyclic-gradients} add without cancellation.
At every change of family, the vector
$(\mathbf1_{c_j=i}-\mathbf1_{c_{j-1}=i})_{i<s}$ is nonzero,
including at a change involving the fixed family $s$. Hence
$\norm{\Xi(w)}_2^2\ge\kappa^{2L}/4$.
Gradients from distinct classes under rotation and adjoint use distinct
reduced block words.

If the cyclic class differs from its adjoint, $\Xi(w)$ and
$\Xi(w)^*$ are orthogonal and have equal norms. The real and
imaginary parts of $\Xi(w)$ are therefore orthogonal, each with
half its squared norm. Otherwise $\Xi(w)$ is self-adjoint, since
$\Xi_i(w^*)=\Xi_i(w)^*$ for each $i<s$.
Consequently the real Gram matrix in the reference model satisfies
\[
G_{\rm ref}\ge\frac18a^{-4L}I_d.
\eqtag\label{eq:rel-reference-gradient-gram}
\]

Let $\Phi(h)$ be the real vector of moments after conjugation, and let
$G_{\alpha\beta}=\sum_{i<s}\tau_0(g_{\alpha,i}g_{\beta,i})$ be the Gram
matrix of the tuples $(g_{\alpha,i})_{i<s}$ at zero.
The formal polynomials representing $x_\lambda$ and $y_\lambda$
differ in coefficient norm by at most $a^{4L}e$.
Telescoping the rotated words shows that the corresponding gradient
polynomials in the two models differ in coefficient norm by at most
$L^2a^{4L}e$.
Every polynomial used here has degree at most $2L$, including the products of gradients.
For such a product, the additional error from evaluating the trace is at most
$\delta^{-1}sL^2e$ by \eqref{eq:rel-coefficient-control}.
These bounds give an error of at most $a^{8L}e$ in each Gram entry.
The reference values of the chosen moments are zero. Applying the same
estimates to these moments and using the bound $d\le a^{3L}$ on the number
of coordinates,
we obtain
\[
\norm{G-G_{\rm ref}}\le a^{11L}e,\qquad
\eta:=\norm{\Phi(0)}_{\ell_2}\le a^{8L}e.
\eqtag\label{eq:rel-gram-residual-transfer}
\]
Since $e\le a^{-40L}$,
\[
G\ge\lambda I_d,\qquad \lambda=\frac1{16}a^{-4L}.
\eqtag\label{eq:rel-gradient-gap}
\]

\emph{Solving the equations.}
If $d=0$ or $\eta=0$, we take all the unitaries to be identities.
Otherwise, for $h=(h_i)_{i<s}$ with
$h_i\in(D_0'\cap M_0)_{\rm sa}$, we use the Hilbert norm
$\norm h_2^2=\sum_{i<s}\norm{h_i}_2^2$.
Each conjugated word has at most $2L$ unitary factors. Differentiating twice and
applying H\"older's inequality for $\tau_0$ gives, for every $h$,
\[
\norm{D^2\Phi(h)[k,l]}_{\ell_2}
\le B_L\norm k_2\norm l_2,\qquad
B_L=4L^2\sqrt d\le a^{4L}.
\eqtag\label{eq:rel-hessian-bound}
\]
Define a right inverse for $-D\Phi(0)$ by
\[
J(z)_i=\frac12\sum_\alpha(G^{-1}z)_\alpha g_{\alpha,i},
\qquad D\Phi(0)J=-I,
\]
with
\[
\norm{Jz}_2\le\frac{\norm z_{\ell_2}}{2\sqrt\lambda},\qquad
\max_i\norm{J(z)_i}_\infty
\le\frac{L\sqrt d}{2\lambda}\norm z_{\ell_2}.
\eqtag\label{eq:rel-inverse-bounds}
\]
For $R(h)=\Phi(h)-\Phi(0)-D\Phi(0)h$, the remainder constant is
$C=B_L/2$. Since $\eta\le\lambda/(8B_L)$, we may apply
Lemma~\ref{lem:finite-dimensional-correction} on the real span of the
gradient tuples $(g_{\alpha,i})_{i<s}$, with
$A=1/(2\sqrt\lambda)$ and $B=L\sqrt d/(2\lambda)$.
The iteration from zero converges to a solution of
\[
h=J(\Phi(0)+R(h))
\eqtag\label{eq:rel-fixed-point}
\]
satisfying $\Phi(h)=0$ and
\[
\max_i\norm{h_i}_\infty\le\frac{3L\sqrt d\,\eta}{4\lambda}
\le\frac{L\sqrt d\,\eta}{\lambda}
\le a^{20L}e\le K^Le.
\eqtag\label{eq:rel-final-generator-bound}
\]
Set $v_i=e^{ih_i}$, $v_s=p$, and $u_i=\mathcal L(v_i)$. By
\eqref{eq:rel-unitary-lift}, the same bound then holds for
$\norm{u_i-1}_\infty$.
The iterates lie in this fixed finite-dimensional span, so their
$L^2$ convergence also gives convergence in operator norm in the
generated $C^*$-algebra. When $M$ is a matrix algebra, triangular factorization and Gram inversion in
the chosen ordered bases are Borel, so the iterates and resulting unitaries
are Borel functions of the matrix entries.
The set of inputs satisfying the hypotheses is closed, because $\Delta_{2L}(\mathbf H;D)$ depends continuously on the
matrix units of $D$ and the entries of $\mathbf H$. Indeed, the reference moments are polynomials in the
values of $\tau$ on the minimal projections of $D$.

\emph{The amalgamated moment identities.}
We expand each same-family block of a corner word in the basis of its
marginal polynomial space, separating the $D_0$ term from the centered
terms. We absorb each $D_0$ term into the neighboring blocks and merge
neighbors from the same family. If the first and last block
indices coincide, we also use cyclicity of the trace to merge the first and last
blocks when they come from the same family.
Repeating these steps does not increase the degree and reduces the trace to
marginal moments and the centered cyclically alternating moments
that we have set equal to zero. A product in $e_vM_0e_w$ with
$v\ne w$ has trace and $D_0$-expectation zero.

The same reductions hold in the reduced amalgamated $C^*$-free product formed from
the actual marginal laws, whose canonical state $\tau|_D\circ\mathcal E$
is tracial (compare~\cite[Section~2]{dykema-amalgamated-fd}). Indeed, if $x\in\jmath_i(A_i)$ and $y$ is an alternating
product of $n\ge1$ centered elements, then, by freeness, $\mathcal E(xy)=\mathcal E(yx)=0$ unless $n=1$ and
$y\in\jmath_i(A_i)$. In that case, and for $y\in D$, the trace property of $\tau$ on $A_i$ applies. Such $x$ generate
$\mathcal A$, and such $y$ span a dense subspace. Testing with $e_v$ and using
\eqref{eq:rel-expectation-reconstruction} therefore gives
\eqref{eq:rel-exact-moments}. Finally, conjugation by a unitary in
$D'$ preserves the $E_D$-value of every word $x$ in one family with coefficients from $D$, since for such a unitary $u$
and $d\in D$ we have $\tau(uxu^*d)=\tau(xd)$, and $\tau$ is faithful on $D$. Hence every marginal law
remains
unchanged.
\end{proof}

We can replace the semicircular reference family by arbitrary families
that are free over $D$, provided the following polynomial dimensions
eventually agree. This
gives an amalgamated version of
Theorem~\ref{thm:multi-gns-rank-strong-exactification}.

\begin{theorem}[Exact amalgamated moments under a dimension hypothesis]
\label{thm:relative-gns-rank-exactification}
Let $M_j$ be unital $C^*$-algebras with faithful tracial states $\tau_j$,
and let $M$ be a unital $C^*$-algebra with a tracial state $\tau$.
Assume that these algebras contain specified unital $^*$-isomorphic copies of a fixed
finite-dimensional unital $C^*$-algebra $D$.

Assume that $\tau_j|_D\to\tau|_D$ and that $\tau|_D$ is faithful.
Fix $s\ge2$. For $1\le i\le s$, let $\mathbf X_j^{(i)}$ be tuples in $M_j$ of fixed
length, uniformly bounded in operator norm, and let $\mathbf x^{(i)}$ be
tuples in $M$, such that
\[
(D,\mathbf X_j^{(1)},\ldots,\mathbf X_j^{(s)})
\longrightarrow(D,\mathbf x^{(1)},\ldots,\mathbf x^{(s)})
\]
in all joint tracial $^*$-moments. Assume that the $C^*$-algebras
$C^*(D,\mathbf x^{(i)})$ are free over $D$ with respect to the
trace-preserving conditional expectation.
For each family $i$, let $\mathcal V_{i,j}(d)$ and
$\mathcal V_i(d)$ be the spans in the respective $L^2$ spaces of all
$D$-coefficient $^*$-polynomials of tuple degree at most $d$ in that family. Suppose that, for every fixed $d$ and $i$,
\[
\dim_{\C}\mathcal V_{i,j}(d)=\dim_{\C}\mathcal V_i(d)
\quad\text{for all sufficiently large }j.
\eqtag\label{eq:rel-gns-rank-hypothesis}
\]
After deleting finitely many indices, there are nondecreasing integers $L_j\to\infty$ and unitaries
\[
\begin{gathered}
u_{i,j}\in D'\cap U\bigl(C^*(D,\mathbf X_j^{(1)},\ldots,\mathbf X_j^{(s)})\bigr),
\qquad u_{s,j}=1,\\
\max_i\norm{u_{i,j}-1}_\infty\longrightarrow0,
\end{gathered}
\eqtag\label{eq:rel-gns-corrections}
\]
such that the $D$-valued moments of the corrected families $\operatorname{Ad}(u_{i,j})(\mathbf X_j^{(i)})$ through degree
$L_j$ agree with those of the reduced amalgamated $C^*$-free product formed from their
actual finite-stage $D$-valued marginal laws.  Every marginal law is preserved at every degree.
If the original tuples converge strongly together with $D$, the corrected tuples have the same strong limit.
\end{theorem}

\begin{proof}
After discarding finitely many indices, the traces of all minimal
projections of $D$ are uniformly bounded away from zero. We pass to the full corners
in \eqref{eq:rel-full-corner}. Identifying equivalent corners by the
matrix units gives the following dimension formula, both at the
finite stages and in the limiting GNS spaces:
\[
\dim\mathcal V_{i,j}(d)
=\sum_{v,w}d_vd_w\dim(e_v\mathcal V_{i,j}(d)e_w).
\eqtag\label{eq:rel-edge-rank-decomposition}
\]
Convergence of the moment Gram matrices implies that, for fixed $d$,
the dimension of each finite-stage block is eventually at least its
limiting GNS dimension. Since all weights $d_vd_w$ in
\eqref{eq:rel-edge-rank-decomposition} are positive, equality of the
total dimensions in \eqref{eq:rel-gns-rank-hypothesis} forces eventual
equality in every block.

Fix $d$. In each block, choose an orthonormal basis of the limiting
GNS polynomial space compatible with the degree filtration, taking
self-adjoint vectors on diagonal blocks and adjoint bases on opposite
blocks, with degree zero in $D_0$. Represent these vectors by
polynomials of the corresponding degrees and with the same adjoint
symmetries. Their evaluated Gram matrices converge to the identity.
Since the dimensions agree in each block, triangular factorization in
degree order gives bases of the marginal polynomial spaces at
the finite stages. These bases are centered in positive degree, and
their polynomial coefficients converge. Because the generators are uniformly bounded, we can choose
one $0<\kappa_d<1$ so that multiplying each degree-$k$ letter by
$\kappa_d^k$ makes it a contraction for all sufficiently large $j$.

For products whose first and last block indices coincide, the calculation in
\eqref{eq:rel-cyclic-gradients}--\eqref{eq:rel-reference-gradient-gram}
bounds the reference Gram matrix below by $\kappa_d^{2d}I/8$. By moment
convergence and the contraction argument
\eqref{eq:rel-hessian-bound}--\eqref{eq:rel-fixed-point}, we obtain
conjugating unitaries in the generated $C^*$-algebras, commuting
with $D$ and tending to $1$. Expanding and reducing same-family blocks as above
gives all the required $D$-valued moments without increasing degree.
The finite-stage traces are faithful, so these polynomial expansions are
algebraic identities.

Choose the degree by a diagonal argument to obtain nondecreasing
$L_j\to\infty$ and $\max_i\|u_{i,j}-1\|_\infty\le1/L_j$.
Since these unitaries tend to $1$ in norm, the corrected tuples have the same strong limit as the original ones, as desired.
\end{proof}

\section{Degree and size of corrections for GUE families}
\label{sec:gue-bounds}

We apply Theorem~\ref{thm:relative-quantitative-exactification} to
independent GUE families and ask how
many exact moment identities can be obtained by unitaries at a given
distance from $1$. Using estimates for GUE moments, we impose identities through a degree proportional to $\log N$, with a deterministic bound on the distance of the unitaries from $1$ and summable failure probabilities (Theorem~\ref{thm:gue-relative-exactification}). For degrees $o(\log N)$, the deterministic bound on the distance can be taken to be $N^{-1+o(1)}$ (Proposition~\ref{prop:gue-degree-motion-profile}), and for a fixed degree, $N$ times the distance is bounded in probability (Corollary~\ref{cor:gue-relative-polynomial-images}).

Conversely, almost surely and for all sufficiently large $N$, the identities can hold only through a degree that is at most a constant times $\log N$ (Proposition~\ref{prop:relative-degree-ceiling}). No deterministic bound $C/N$ on the distance holds with probability tending to one, already for the identities of degree two (Proposition~\ref{prop:gue-motion-obstruction}). If $D\ne\C$ or one of the families has at least two generators, we prove a lower bound on the distance that grows exponentially with the degree, and hence a deterministic bound $N^{-1+o(1)}$ can hold with probability tending to one only at degrees $o(\log N)$ (Theorem~\ref{thm:relative-degree-motion-lower-bound}).

In Subsection~\ref{subsec:optimal-correction-costs}, we consider covariance equations for a fixed finite list of real
polynomials in one GUE matrix, with a second independent GUE matrix
left fixed. Assuming the polynomials are linearly independent modulo
constants, we determine the limiting distribution of $N$ times the
minimum of $\|u-1\|_\infty$ over the unitaries $u$ such that all these equations hold after the first matrix is conjugated by $u$ (Corollary~\ref{cor:optcost-polynomial-cost}). For the two polynomials $x$ and $x^2$, this limiting
distribution differs from that of $N$ times the larger of the two separate minima
(Proposition~\ref{prop:optcost-genuine-joint-cost}).

\subsection{Bounds on degree and distance}
\label{subsec:gue-degree-bounds}

At each dimension, we equip $D$ with the trace induced by its
embedding in $M_N(\C)$. As in
Section~\ref{sec:relative-exactification}, degree counts the GUE
generators, and coefficients in $D$ have degree zero.

\begin{theorem}[Exact amalgamated moments for GUE families]
\label{thm:gue-relative-exactification}
Fix $D=\bigoplus_{v=1}^bM_{d_v}(\C)$, $s\ge2$, $m_1,\ldots,m_s\ge1$, and $\delta>0$.  Let $N$ range over an unbounded
set of dimensions, and let $\iota_N:D\to M_N(\C)$ be deterministic unital $^*$-embeddings such that
$\operatorname{tr}_N(\iota_N(e))\ge\delta$ for every minimal projection $e\in D$.  Set $D_N=\iota_N(D)$, and let
$E_N:M_N(\C)\to D_N$ be the trace-preserving conditional expectation.  For each $N$, let
$\mathbf G_{i,N}=(G_{i,a,N})_{a=1}^{m_i}$, $1\le i\le s$, be tuples of normalized GUE matrices, with all generators mutually independent at each dimension.  Define their marginal $C^*$-algebras and reduced amalgamated $C^*$-free product by
\[
A_{i,N}=C^*(D_N,\mathbf G_{i,N}),\qquad
(\mathcal A_N,\mathcal E_N)
=\mathop{*}_{D_N,\,i=1}^s(A_{i,N},E_N|_{A_{i,N}}),
\]
with canonical embeddings $\jmath_{i,N}$.

For every $0<\gamma<1$ and $B>1$, there are $c>0$ and $N_0$, depending only on $D$, $s$, $m_1,\ldots,m_s$, $\delta$, $\gamma$, and $B$, and
Borel functions of the matrix entries
\[
u_{i,N}\in D_N'\cap U(N),\qquad u_{s,N}=1,
\]
defined for every tuple of self-adjoint matrices and satisfying, for $N\ge N_0$,
\[
\max_i\norm{u_{i,N}-1}_\infty\le N^{-\gamma}.
\eqtag\label{eq:relgauss-global-cap}
\]
Set $L_N=\lfloor c\log N\rfloor$ and $\widetilde G_{i,a,N}=u_{i,N}G_{i,a,N}u_{i,N}^*$.  For each $N\ge N_0$, with probability at least
$1-N^{-B}$, simultaneously for every $0\le\ell\le L_N$, all indices $i_1,\ldots,i_\ell$ and $a_1,\ldots,a_\ell$, and all
$d_0,\ldots,d_\ell\in D_N$,
\[
\begin{split}
&E_N\bigl(d_0\widetilde G_{i_1,a_1,N}d_1\cdots
\widetilde G_{i_\ell,a_\ell,N}d_\ell\bigr)\\
&\quad=\mathcal E_N\bigl(d_0\jmath_{i_1,N}(G_{i_1,a_1,N})d_1\cdots
\jmath_{i_\ell,N}(G_{i_\ell,a_\ell,N})d_\ell\bigr).
\end{split}
\eqtag\label{eq:relgauss-exact-identities}
\]
Each family's complete $D_N$-valued law is preserved for every realization.  On a probability-one event,
\eqref{eq:relgauss-exact-identities} holds at all sufficiently large dimensions in the chosen set.

If the induced traces converge, $\operatorname{tr}_N\circ\iota_N\to\tau_D$, then the corrected tuple together with the matrix units
of $D$ has the same almost-sure strong limit as the original tuple: a free semicircular family free from
$(D,\tau_D)$. The norm convergence also holds for polynomials with
matrix coefficients of any fixed size.
\end{theorem}

Applying \eqref{eq:relgauss-exact-identities} entry by entry gives
the same identities for polynomials with matrix coefficients. Because
the unitaries depend on the matrices, the corrected matrices need not
have the GUE distribution.

To prove Theorem~\ref{thm:gue-relative-exactification}, we bound the probability that the moment defect of the GUE tuple is large or that one of the GUE matrices has norm larger than $3$ (Lemma~\ref{lem:relgauss-raw-moments}), and we apply Theorem~\ref{thm:relative-quantitative-exactification} on the complementary event.

\begin{lemma}
\label{lem:relgauss-raw-moments}
Set $q=\dim_{\C}D$, $m=\sum_i m_i$, and
\[
D_*=[2m(q+1)]^4,\qquad
b_{L,N}=\frac{4^LL^3}{N^2},\qquad J_p=D_*(128p)^p.
\]
Let $e>0$, let $\mathbf G=(G_\alpha)_{\alpha=1}^m$ be the concatenation of the tuples $\mathbf G_{1,N},\ldots,\mathbf G_{s,N}$ of Theorem~\ref{thm:gue-relative-exactification}, and let $\Delta_{2L}(\mathbf G;D_N)$ be the moment defect used in
Theorem~\ref{thm:relative-quantitative-exactification}.
If $L\ge2$, $(2L)^3\le N^2$, and $b_{L,N}\le e/2$, then
\[
\mathbb P\{\Delta_{2L}(\mathbf G;D_N)>e\}
\le J_p^L(Ne)^{-p}
\eqtag\label{eq:relgauss-dictionary-tail}
\]
for every even integer $p\ge2$.  For one normalized GUE matrix and every integer
$r\ge1$,
\[
\mathbb P\{\norm G_\infty>3\}
\le N(4/9)^r\exp\left(\frac{r^3-r}{12N^2}\right).
\eqtag\label{eq:relgauss-norm-tail}
\]
The bounds are uniform over the embeddings of $D$.
\end{lemma}

\begin{proof}
\emph{A genus bound.}
For $r\ge1$, let $\varepsilon_g(r)$ count pairings $\pi$ of $2r$ cyclic
positions of genus $g$, where $2g=r+1-\#(\zeta\pi)$,
$\zeta=(1\ \cdots\ 2r)$, and $\#$ denotes the number of cycles.
Let $C_r$ be the $r$th Catalan number. We use the generating-function identity
of Harer and Zagier~\cite[Theorem~3]{harer-zagier}:
\[
1+2\sum_{r\ge0}\frac{x^{r+1}}{(2r-1)!!}
\sum_g\varepsilon_g(r)y^{r+1-2g}
=\left(\frac{1+x}{1-x}\right)^y,
\eqtag\label{eq:relgauss-harer-zagier}
\]
with $(-1)!!=1$ and $\varepsilon_0(0)=1$.
Then we differentiate in $x$, multiply by $1-x^2$, and compare coefficients.
Solving the resulting recurrence for $\varepsilon_g(r)$ and dividing by $C_r$, we obtain
\[
a_g(r)=a_g(r-1)+\frac{r(r-1)}4a_{g-1}(r-2),
\qquad a_g(r)=\frac{\varepsilon_g(r)}{C_r},\quad r\ge2.
\]
For $r=0,1$, only genus zero occurs.
Iterating the recurrence expresses $a_g(r)$ as a sum of products of $g$
nonnegative weights $j(j-1)/4$ over indices separated by at least two.
If we drop that restriction and bound the elementary symmetric sum of degree $g$ in these weights by the $g$th power of their sum divided by $g!$, we obtain
\[
a_g(r)\le\frac1{g!}\left(\sum_{j=1}^r\frac{j(j-1)}4\right)^g
=\frac1{g!}\left(\frac{r^3-r}{12}\right)^g.
\]
Consequently, for every $r,N$,
\[
\sum_g\varepsilon_g(r)N^{-2g}
\le C_r\exp\left(\frac{r^3-r}{12N^2}\right).
\eqtag\label{eq:relgauss-genus-majorant}
\]

\emph{Expectations.}
The independent GUE matrices have covariance
\[
\mathbb E(G_\alpha)_{ij}(G_\beta)_{kl}
=N^{-1}\delta_{\alpha\beta}\delta_{il}\delta_{jk}.
\]
For deterministic contractions $d_j\in D_N$, we consider
$F=\operatorname{tr}_N(d_0G_{\alpha_1}d_1\cdots G_{\alpha_k}d_k)$.
When $k=2r$, we merge its endpoint coefficients cyclically and expand by Wick's formula.
If $\pi$ pairs positions carrying the same GUE matrix, then after we sum over the matrix indices, we obtain one trace of a product of coefficients for each cycle of $\zeta\pi$. If we write these traces with $\operatorname{tr}_N$, the remaining power of $N$ is $N^{\#(\zeta\pi)-r-1}=N^{-2g}$.  Each coefficient trace is that of a product of contractions and hence has absolute value at most one.  The noncrossing pairings give the moment $\tau_{\rm ref}(d_0S_{\alpha_1}d_1\cdots S_{\alpha_k}d_k)$, computed with the normalized trace of $M_N(\C)$ restricted to $D_N$.  By \eqref{eq:relgauss-genus-majorant}, we have
\[
\begin{split}
&\abs{\mathbb EF-\tau_{\rm ref}
 (d_0S_{\alpha_1}d_1\cdots S_{\alpha_k}d_k)}\\
&\qquad\le C_r\left[\exp\left(\frac{r^3-r}{12N^2}\right)-1\right].
\end{split}
\eqtag\label{eq:relgauss-bias}
\]
For odd $k$, both $\mathbb EF$ and $\tau_{\rm ref}(d_0S_{\alpha_1}d_1\cdots S_{\alpha_k}d_k)$ are zero.

\emph{Variances.}
We express each GUE matrix in an orthonormal basis $(E_a)$ of the real
space of self-adjoint matrices, for the inner product
$(X,Y)\mapsto\operatorname{Tr}(XY)$, with coordinates $N^{-1/2}\xi_a$ for independent
standard real Gaussian variables $\xi_a$.
For $k\ge1$, each differentiated occurrence in $F$ contributes $N^{-3/2}\operatorname{Tr}(BE_a)$ to the derivative of $F$ with respect to $\xi_a$, where $B$ is the cyclic product
of the remaining $k-1$ random letters and coefficients. Parseval's identity gives
$\sum_a|\operatorname{Tr}(BE_a)|^2=\operatorname{Tr}(B^*B)$.
By Cauchy--Schwarz, applied to the sum over differentiated
occurrences, and the Gaussian Poincar\'e
inequality, we have
\[
\begin{aligned}
\mathbb E\abs{F-\mathbb EF}^2
&\le\mathbb E\norm{\nabla F}^2
\le\frac{k^2}{N^2}\max_B\mathbb E\operatorname{tr}_N(B^*B)\\
&\le\frac{k^2}{N^2}C_{k-1}
\exp\left(\frac{(k-1)^3}{12N^2}\right).
\end{aligned}
\eqtag\label{eq:relgauss-variance}
\]
For the last inequality, we apply the Wick expansion above to $\operatorname{tr}_N(B^*B)$, with $r=k-1$, and use \eqref{eq:relgauss-genus-majorant}.
For complex $F$, we apply the Gaussian Poincar\'e inequality to its real and
imaginary parts and sum.

\emph{Tail bounds.}
For a real Gaussian polynomial $F$ of degree at most $k$, write
$F=\sum_{j=0}^kF_j$ for its orthogonal expansion in products of
Hermite polynomials, grouped by total degree. By Nelson's hypercontractivity
theorem~\cite[Theorem~2]{nelson-free-markoff}, we have
\[
\norm F_{L^p}
\le\left(\sum_{j=0}^k(p-1)^j\norm{F_j}_{L^2}^2\right)^{1/2}
\le(p-1)^{k/2}\norm F_{L^2},\qquad p\ge2.
\eqtag\label{eq:relgauss-hypercontractivity}
\]
For complex polynomials, an extra factor $2$ suffices.  For $1\le k\le2L$ and $(2L)^3\le N^2$, equations
\eqref{eq:relgauss-variance}--\eqref{eq:relgauss-hypercontractivity} imply
\[
\norm{F-\mathbb EF}_{L^p}
\le\frac{2k\,2^{k-1}\exp(1/12)\,p^{k/2}}N
\le\frac{(64p)^L}{N}.
\]
Words of degree zero lie in $D_N$ and contribute no defect. The list of words defining the defect contains at most
\[
(q+1)\sum_{k=0}^{2L}[m(q+1)]^k\le D_*^L
\]
words.  By \eqref{eq:relgauss-bias} and the bounds $C_r\le4^r$ and $r\le L$, the left side of \eqref{eq:relgauss-bias} is at most $b_{L,N}$ for every word in the list.  By Markov's inequality at threshold $e/2$, the probability on the left of \eqref{eq:relgauss-dictionary-tail} is at most
\[
D_*^L\left(\frac{2(64p)^L}{Ne}\right)^p
\le J_p^L(Ne)^{-p},
\]
proving \eqref{eq:relgauss-dictionary-tail}.  Finally, \eqref{eq:relgauss-norm-tail} follows from
$\norm G_\infty^{2r}\le\operatorname{Tr}(G^{2r})$
and \eqref{eq:relgauss-genus-majorant}.  In particular, for a fixed number of GUE matrices and any prescribed power of $N^{-1}$, if we choose $r$ to be a sufficiently large constant times $\log N$, then, for all sufficiently large $N$, the probability that any of their norms exceeds $3$ is smaller than this power.
\end{proof}

\begin{proof}[Proof of Theorem~\ref{thm:gue-relative-exactification}]
Let $K$ be the constant of Theorem~\ref{thm:relative-quantitative-exactification}; enlarging it if necessary, we may assume that $K>1$. We set
$\epsilon=(1-\gamma)/3$, choose an even integer $p$ with $p\epsilon>2(B+3)$, and take
\[
\begin{split}
c&=\min\left\{\frac{\epsilon}{2\log K},\,
\frac{p\epsilon}{4\log J_p},\,\frac1{16\log4}\right\},\\
L_N&=\lfloor c\log N\rfloor,\qquad e_N=N^{-1+\epsilon}.
\end{split}
\eqtag\label{eq:relgauss-parameters}
\]
For all sufficiently large $N$, we have $L_N\ge2$ and $(2L_N)^3\le N^2$.
Moreover, $b_{L_N,N}=o(e_N)$. By Lemma~\ref{lem:relgauss-raw-moments}, we
have
\[
\mathbb P\{\Delta_{2L_N}>e_N\}
\le N^{-3p\epsilon/4}\le N^{-B-3}.
\]
We consider the event on which $\Delta_{2L_N}\le e_N$ and $\max_{i,a}\norm{G_{i,a,N}}_\infty\le3$. We choose the moment order $r$ in \eqref{eq:relgauss-norm-tail} to be a sufficiently large constant times $\log N$, so that the probability that either bound fails is at most $N^{-B}$. Since $c\le\epsilon/(2\log K)$, we have
\[
K^{L_N}e_N\le N^{-1+3\epsilon/2}
=N^{-(1+\gamma)/2}\le N^{-\gamma}.
\]
We apply Theorem~\ref{thm:relative-quantitative-exactification} on this event, and set $u_{i,N}=1$ outside it and for the finitely many $N<N_0$. These
choices are Borel, by the measurability assertion of that theorem and because the event is closed
($\Delta_{2L_N}$ and the norms are continuous functions of the matrices). For every realization, \eqref{eq:relgauss-global-cap} holds.

Since $B>1$, the probabilities $N^{-B}$ are summable over the chosen set of dimensions, so by the first Borel--Cantelli lemma, \eqref{eq:relgauss-exact-identities} holds almost surely at all sufficiently large dimensions in this set.

For the strong-convergence assertion, assume additionally that
$\operatorname{tr}_N\circ\iota_N\to\tau_D$. Every fixed polynomial in
the matrix units of $D$ has the same norm under the embeddings
$\iota_N$, and its trace converges. The trace $\tau_D$ is faithful
since every minimal projection has trace at least $\delta$.
By~\cite[Theorem~1.6 and Corollary~2.3]{male-norm-polynomials}, in the form recalled in
Subsection~\ref{sec:stable-strong-models},
these matrix units and the independent GUE tuple converge
jointly and strongly almost surely, also for fixed matrix coefficients.
Almost surely, the generators are eventually norm bounded, and
\[
\norm{u_{i,N}G_{i,a,N}u_{i,N}^*-G_{i,a,N}}_\infty
\le2\norm{G_{i,a,N}}_\infty\norm{u_{i,N}-1}_\infty\longrightarrow0.
\]
By telescoping the monomials, we see that, for every polynomial, the trace and the norm at the corrected tuples have the same limits as at the original tuples.  In the same way, we obtain the norm convergence for polynomials with matrix coefficients of each fixed size.
\end{proof}

If we conjugate families of matrices by independent Haar unitaries in the relative commutant $D_N'\cap U(N)$, then they become approximately free over the unital $^*$-subalgebra $D_N$, and their marginal laws are preserved. Brown, Dykema, and Jung~\cite[Theorem~3.8 and
Corollary~3.9]{brown-dykema-jung} proved this for moments of a fixed degree, up to a fixed positive error, uniformly
for
tuples with a fixed number of entries and a common norm bound. Their
setup allows the induced traces on $D$ to converge to a faithful
trace~\cite[Subsection~3.2]{brown-dykema-jung}.
Their conclusion also holds when the induced traces on $D$ vary with $N$ and give every minimal projection of $D$ trace at least $\delta$. Indeed, suppose that it fails. Since the set of traces on $D$ is compact, we can then choose a counterexample along a subsequence on which the induced traces converge to a trace on $D$, and this trace is faithful. This contradicts the uniform statement of Brown, Dykema, and Jung.

Hayes, Jekel, Nelson, and Sinclair~\cite[Theorem~6.10 and
Remarks~6.11--6.12]{hayes-jekel-nelson-sinclair} allowed atomic amalgamation and countably many families
in their theorem on
asymptotic freeness under conjugation. In their Remark~6.12, they noted that one family
may be kept fixed. This follows by conjugating all families by one further unitary and using the invariance of the Haar measure.

\begin{proposition}
\label{prop:gue-degree-motion-profile}
Under the assumptions of Theorem~\ref{thm:gue-relative-exactification},
for every $B>1$ there are $C_B>1$, $c_B>0$ and $N_0$, depending only on $D$, $s$, $m_1,\ldots,m_s$, $\delta$, and $B$, with the following property. For every
integer sequence $(L_N)$ satisfying
$2\le L_N\le c_B\log N$ for $N\ge N_0$, there are Borel choices
of unitaries $u_{i,N}\in D_N'\cap U(N)$, with $u_{s,N}=1$, such
that, for each $N\ge N_0$,
\eqref{eq:relgauss-exact-identities} holds through degree $L_N$ with
probability at least $1-N^{-B}$. For every realization,
\[
\max_i\norm{u_{i,N}-1}_\infty
\le\frac1N\left(\frac{C_B\log N}{L_N}\right)^{L_N}.
\eqtag\label{eq:relgauss-degree-motion}
\]
In particular, for every prescribed $L_N=o(\log N)$ with $L_N\ge2$, the correcting unitaries can be chosen within distance
$N^{-1+o(1)}$ from the identity, with all marginal laws preserved. The identities hold almost surely at all sufficiently large dimensions in the chosen set. If the induced traces on $D$ converge, the corrected tuples have the strong limit stated in
Theorem~\ref{thm:gue-relative-exactification}.
\end{proposition}

\begin{proof}
In \eqref{eq:relgauss-dictionary-tail}, we take
\[
\begin{split}
p_N&=2\left\lceil\frac{(B+3)\log N}{2L_N}\right\rceil,\\
e_N&=\frac{(128p_N)^{L_N}}N
D_*^{L_N/p_N}N^{(B+2)/p_N}.
\end{split}
\]
With these choices,
$\mathbb P\{\Delta_{2L_N}>e_N\}\le N^{-B-2}$.
For $L_N\le\log N$,
\[
\begin{gathered}
p_N\le(B+5)\frac{\log N}{L_N},\qquad
N^{(B+2)/p_N}\le\exp(L_N),\\
D_*^{L_N/p_N}\le D_*^{L_N/2}.
\end{gathered}
\]
Consequently $K^{L_N}e_N$ is bounded by the right side of \eqref{eq:relgauss-degree-motion}, for example with
\[
C_B=256\exp(1)K(B+5)\sqrt{D_*},\qquad
c_B\le\min\{(4C_B)^{-1},(16\log4)^{-1}\}.
\]
Since $x\log(C_B/x)$ is increasing on $(0,(4C_B)^{-1}]$,
for $0<x\le c_B$ we have
$x\log(C_B/x)\le\log(4C_B^2)/(4C_B)\le1/2$.  Thus, with $x=L_N/\log N$, the right side of \eqref{eq:relgauss-degree-motion} is at most $N^{-1/2}$.  The conditions $(2L_N)^3\le N^2$ and $b_{L_N,N}\le e_N/2$
of Lemma~\ref{lem:relgauss-raw-moments}
also hold uniformly for large $N$: the first follows from
$L_N\le c_B\log N$, and the second from $b_{L_N,N}=o(N^{-1})$
and $e_N\ge N^{-1}$.  We apply Theorem~\ref{thm:relative-quantitative-exactification} when
$\Delta_{2L_N}\le e_N$ and $\max_{i,a}\norm{G_{i,a,N}}_\infty\le3$,
and set $u_{i,N}=1$ otherwise and for $N<N_0$.  By \eqref{eq:relgauss-norm-tail} and \eqref{eq:relgauss-dictionary-tail}, \eqref{eq:relgauss-exact-identities} holds through degree $L_N$ with probability at least $1-N^{-B}$ for all $N\ge N_0$, where $N_0$ does not depend on the sequence $(L_N)$.

If $x_N=L_N/\log N\to0$, then
\[
\log\left[\left(\frac{C_B\log N}{L_N}\right)^{L_N}\right]
=\log N\,x_N\log(C_B/x_N)=o(\log N).
\]
Since $B>1$, the probabilities $N^{-B}$ are summable, so by the first Borel--Cantelli lemma, the identities hold almost surely for all large dimensions. When the induced traces converge, we see from the telescoping estimate in the proof of Theorem~\ref{thm:gue-relative-exactification} that the joint strong limit is preserved.
\end{proof}

\begin{corollary}
\label{cor:gue-relative-polynomial-images}
Under the assumptions of Theorem~\ref{thm:gue-relative-exactification},
for every fixed integer $L\ge0$ there are Borel choices of unitaries
$u_{i,N}\in D_N'\cap U(N)$, with $u_{s,N}=1$, such that \eqref{eq:relgauss-exact-identities} holds through degree $L$ with probability tending to one and
\[
\max_i\norm{u_{i,N}-1}_\infty=O_{\mathbb P}(N^{-1}).
\eqtag\label{eq:relgauss-fixed-degree-cost}
\]
Here $Z_N=O_{\mathbb P}(N^{-1})$ means that $NZ_N$ is bounded in probability.
Moreover, suppose that we replace each family in either Theorem~\ref{thm:gue-relative-exactification} or Proposition~\ref{prop:gue-degree-motion-profile} by a fixed finite list of $D_N$-coefficient $^*$-polynomials in that same family, of degree at most a fixed positive integer $f$. Then the same unitaries preserve each new family's $D_N$-valued law for every realization. Whenever the original identities hold, the $D_N$-valued joint moments of the new families agree with those of their amalgamated free product through degree $\lfloor L_N/f\rfloor$.
For each polynomial value $Y$ in a new family, every relation
$pYq=0$ with $p,q\in D_N$ is preserved exactly.
\end{corollary}

\begin{proof}
For fixed $L\ge2$, only finitely many words occur in the defect, and
\eqref{eq:relgauss-bias}--\eqref{eq:relgauss-variance} give
$\Delta_{2L}=O_{\mathbb P}(N^{-1})$, uniformly over embeddings with $\operatorname{tr}_N(\iota_N(e))\ge\delta$ for every minimal projection $e\in D$.
We apply Theorem~\ref{thm:relative-quantitative-exactification} on the event
$\Delta_{2L}\le K^{-L}$ and $\max_{i,a}\norm{G_{i,a,N}}_\infty\le3$,
and set $u_{i,N}=1$ elsewhere. This event is closed, hence Borel, and its probability tends to
one by the bound just stated and \eqref{eq:relgauss-norm-tail}. On it, $\max_i\norm{u_{i,N}-1}_\infty\le K^L\Delta_{2L}$,
which proves \eqref{eq:relgauss-fixed-degree-cost}.
For $L\le1$, we take $u_{i,N}=1$, since the moments of degrees zero and one
already agree with the free-product moments.

Since the unitaries commute with every $D_N$ coefficient, the value of a polynomial at the conjugated variables is the conjugate of its original value by the same unitary.  If we substitute the defining polynomials into a word of degree at most $\lfloor L_N/f\rfloor$ in the new families, we obtain a polynomial of degree at most $L_N$ in the original GUE matrices.  Inside $\mathcal A_N$, the unital $C^*$-subalgebras generated by $D_N$ and the
values of the defining polynomials at $\jmath_{i,N}(\mathbf G_{i,N})$ are
free over $D_N$. The restrictions of $\mathcal E_N$ give the
$D_N$-valued laws of the new families. Thus their joint moments are
those of the required amalgamated free product.  Finally, for such a polynomial value $Y$ and
the unitary $u$ assigned to its family, we have $p(uYu^*)q=u(pYq)u^*$, so every relation $pYq=0$ with $p,q\in D_N$ is preserved.
\end{proof}

In the rest of this section, we write $\norm{\cdot}$ for the operator norm $\norm{\cdot}_\infty$. In the rest of this subsection, all matrix $L^2$-norms use $\operatorname{tr}_N$.

\begin{proposition}
\label{prop:relative-degree-ceiling}
Let $D\subset M_N(\C)$ be a unital $^*$-subalgebra, let $s\ge2$, let
$\mathbf X_1,\ldots,\mathbf X_s$ be tuples of self-adjoint matrices in
$M_N(\C)$, and let $X_i$ be an entry of $\mathbf X_i$. Suppose $e$ is a
minimal projection of $D$ of rank $n$ and every $eX_ie$ has simple
spectrum in $eM_Ne$. If there are unitaries $u_1,\ldots,u_s\in D'\cap U(N)$
for which the tuples $u_i\mathbf X_iu_i^*$ satisfy all amalgamated
free-product identities for the original marginal laws through degree
$L$, then
\[
 \sum_{j=0}^{r}s^j\le n^2,
 \qquad r=\min\{\lfloor L/2\rfloor,n-1\}.
 \eqtag\label{eq:relcost-degree-count}
\]
Consequently, for sufficiently large $n$,
\[
 L\le\frac{4\log n}{\log s}+1
 \le\frac{4\log N}{\log s}+1.
 \eqtag\label{eq:relcost-degree-ceiling}
\]
For deterministic $D$ and GUE matrices, the simple-spectrum hypothesis holds almost surely.
\end{proposition}

\begin{proof}
We equip $eM_Ne$ with its normalized trace.  Since $e$ is a minimal projection of $D$, we have $eDe=\C e$. Hence, in the amalgamated free product, the compressed families are free for the normalized corner trace. Indeed, an element of the corner is centered for the normalized corner trace exactly when its conditional expectation onto $D$ is zero, so alternating products of centered elements of the compressed families have zero expectation. Since a unitary in $D'$ commutes with $e$, conjugating $X_i$ by it does not change the spectrum of $eX_ie$. For each $i$, since $eX_ie$ has simple spectrum,
there are polynomials $P_{i,j}$, $0\le j<n$, with
$\deg P_{i,j}=j$ and $P_{i,0}=1$, whose values $P_{i,j}(eX_ie)$ are
orthonormal for the normalized trace on $eM_Ne$.

For every word of length $j$ with $1\le j\le r$ in the alphabet $\{1,\ldots,s\}$, we write its maximal constant runs as $i_1^{j_1}\cdots i_t^{j_t}$ and consider the product $P_{i_1,j_1}(eX_{i_1}e)\cdots P_{i_t,j_t}(eX_{i_t}e)$ in the amalgamated free product.  Along with $1$,
these form an orthonormal family in the free-product $L^2$ space.
The products used to compute
the Gram matrix have degree at most $2r\le L$. By the amalgamated free-product identities through degree $L$ for the tuples $u_i\mathbf X_iu_i^*$, the corresponding Gram matrix in $L^2(eM_Ne)$ is the identity. Since $L^2(eM_Ne)$ has complex dimension $n^2$, this proves
\eqref{eq:relcost-degree-count}.  For large $n$, we have $s^{n-1}>n^2$, which
excludes $r=n-1$. Then $s^{\lfloor L/2\rfloor}\le n^2$, and
\eqref{eq:relcost-degree-ceiling} follows.
Finally, a deterministic rank-$n$ compression of a normalized GUE matrix has a nondegenerate Gaussian density
on the self-adjoint $n\times n$ matrices.  The discriminant of its characteristic polynomial, a nonzero polynomial in the entries, therefore
vanishes with probability zero.
\end{proof}

For $n=1$, \eqref{eq:relcost-degree-count} imposes no restriction on
the degree, but under the hypotheses of
Theorem~\ref{thm:gue-relative-exactification} we have $n\ge\delta N$.
Almost surely, for all sufficiently large $N$, the largest $L$ for which unitaries as in Proposition~\ref{prop:relative-degree-ceiling} exist is bounded above and below by positive
multiples of $\log N$. We do not know the optimal constants in these bounds.

\begin{proposition}
\label{prop:gue-motion-obstruction}
Let $G_N,H_N$ be independent normalized GUE matrices and let $D_N\subset M_N(\C)$ be
deterministic unital $^*$-subalgebras isomorphic to a fixed finite-dimensional $C^*$-algebra $D$.  For every $a_N\to0$,
the probability that there are unitaries $u,v\in D_N'\cap U(N)$
giving the degree-two identities for the amalgamated free product of
the marginal laws and satisfying
$\max\{\norm{u-1},\norm{v-1}\}\le a_N/N$ tends to zero.  More precisely, let
$\mathcal R_N(r)$ be the event that there are such unitaries with $\max\{\norm{u-1},\norm{v-1}\}\le r$.
For every fixed $C\ge0$,
\[
 \limsup_{N\to\infty}\mathbb P(\mathcal R_N(C/N))
 \le\mathbb P\{\abs{Z}\le4C\}<1,
 \quad Z\sim\mathcal N(0,1).
 \eqtag\label{eq:relcost-fixed-cap}
\]
Thus, if $\mathbb P(\mathcal R_N(r_N))\to1$ for deterministic $r_N$, then $Nr_N\to\infty$.
\end{proposition}

\begin{proof}
Write $E_N=E_{D_N}$, $A=G_N-E_N(G_N)$, and $C_1=H_N-E_N(H_N)$.
Every pair of unitaries $u,v\in D_N'\cap U(N)$ giving the required
degree-two identities satisfies $\operatorname{tr}_N(uAu^*vC_1v^*)=0$.
Let $r=\max\{\norm{u-1},\norm{v-1}\}$.
Conjugating both factors by $v^*$ and applying Cauchy--Schwarz gives
\[
 \abs{\operatorname{tr}_N(AC_1)}\le4r\norm{A}_2\norm{C_1}_2.
 \eqtag\label{eq:relcost-covariance-bound}
\]
Since $A$ is orthogonal to $D_N$ and $H_N$ is a normalized GUE matrix independent of $G_N$, we have, conditionally on $G_N$,
\[
 N\operatorname{tr}_N(AC_1)=\operatorname{Tr}(AH_N)
 \sim\mathcal N(0,\operatorname{tr}_N(A^2)).
 \eqtag\label{eq:relcost-conditional-normal}
\]
In the real Hilbert space of self-adjoint matrices, with inner product
$(X,Y)\mapsto\operatorname{tr}_N(XY)$, a normalized GUE matrix has covariance
$N^{-2}I$.  Since $(D_N)_{\rm sa}$ has real dimension
$q=\dim_\C D$,
\[
 \mathbb E\norm{E_N(G_N)}_2^2=q/N^2,
 \qquad \mathbb E\norm{G_N}_2^2=1,
 \qquad \operatorname{Var}(\norm{G_N}_2^2)=2/N^2.
\]
Hence $\norm{A}_2,\norm{C_1}_2\to1$ in probability, and the random variable in
\eqref{eq:relcost-conditional-normal} converges in distribution to $Z$.
Since every pair of unitaries satisfying the required moment identities obeys
\eqref{eq:relcost-covariance-bound}, continuity of the Gaussian distribution gives
\eqref{eq:relcost-fixed-cap} and the assertion for $a_N\to0$.  If $Nr_N$ did not tend
to infinity, a bounded subsequence would contradict \eqref{eq:relcost-fixed-cap}.
\end{proof}

We strengthen the preceding lower bound by projecting an independent
GUE matrix onto the space of all self-adjoint $D$-coefficient polynomials in the other family, through
a given degree, with zero conditional expectation onto $D$. To count the independent Gaussian coordinates of this
projection, we construct a tuple on which evaluation of all
$D$-coefficient polynomials through the given degree is injective.
Klep and \v Spenko~\cite[Theorem~2.2 and Remark~3.4(a)]{klep-spenko} gave a related construction. For $g\ge2$, they obtained
$n^2$ words of degree
$2\lceil\log_g n\rceil$ whose evaluations at a single $g$-tuple of $n\times n$ matrices are linearly
independent, and the matrices may be taken symmetric.
Here we also take into account the multiplicities of $D$.

\begin{lemma}
\label{lem:full-filtered-relative-rank}
Let $D=\bigoplus_v M_{d_v}(\C)\subset M_N(\C)$ be a unital $^*$-subalgebra, let
$q=\dim_\C D$, and, for matrix units $f^v_{jk}$ of $D$, write $n_v=\operatorname{rank}(f^v_{11})$ for its multiplicities.
For $m_0\ge1$ and $k\ge0$, set
\[
 \rho=m_0q,\qquad h_+=\lceil k/2\rceil,\qquad h_-=\lfloor k/2\rfloor,\qquad
 T_k=1+\sum_{j=1}^{h_+}\rho^j+\sum_{j=1}^{h_-}\rho^j.
 \eqtag\label{eq:relcost-rank-multiplicity}
\]
Consider the formal $D$-coefficient polynomials of degree at most $k$ in $m_0$ variables: the linear combinations of monomials $d_0x_{a_1}d_1\cdots x_{a_j}d_j$, with $j\le k$ and $d_0,\ldots,d_j\in D$, which depend multilinearly on $(d_0,\ldots,d_j)$ and satisfy no other relations. They form a space of complex dimension
\[
 q\sum_{j=0}^{k}\rho^j.
 \eqtag\label{eq:relcost-filtered-dimension}
\]
If $n_v\ge d_vT_k$ for every $v$, then, for every self-adjoint $m_0$-tuple $\mathbf X$ in $M_N(\C)$ outside a set of Lebesgue measure zero, a formal polynomial of degree at most $k$ that vanishes at $\mathbf X$ is zero.
\end{lemma}

\begin{proof}
We equip $D$ with a faithful tracial state $\tau$. For each $j\ge0$, let
$\mathcal K_j$ be the direct sum of $m_0^j$ copies of
$L^2(D)^{\otimes(j+1)}$, indexed by strings of generator indices.
We take two copies of these spaces, the plus and minus branches, up to levels $h_+$ and $h_-$, and identify their degree-zero spaces:
\[
 \mathcal K=\mathcal K_0\oplus\bigoplus_{j=1}^{h_+}\mathcal K_j^+
                    \oplus\bigoplus_{j=1}^{h_-}\mathcal K_j^-.
\]
The left action of $D$ is on the first tensor factor.  On $\mathcal K_j$ its multiplicity
in block $v$ is $d_v\rho^j$, so its multiplicity on $\mathcal K$ is $d_vT_k$.
Since $n_v\ge d_vT_k$ for every $v$, we can realize this representation as a
reducing subrepresentation of the given action of $D$ on $\C^N$.

At a level of a branch below its top level, let $\ell_a$ send $\xi$ to
$1\otimes\xi$ in the summand of the next level of that branch whose
string of generator indices is $a$ followed by that of $\xi$. On
$\mathcal K_0$, let $\ell_a\xi$ be the sum of the vectors
$1\otimes\xi$ in the summands indexed by $a$ in the first levels of the branches of positive height.
On the top level of each branch, $\ell_a=0$.
Let $H_a=\ell_a+\ell_a^*$, extended by zero on the complementary
$D$-subrepresentation.
Suppose a formal polynomial $P$ of maximal degree $r\le k$ evaluates to zero.  The case
$r=0$ follows from faithfulness of the action of $D$.  Otherwise set
$u=\min\{r,h_+\}$ and $v_0=r-u\le h_-$.
We consider the block from minus level $v_0$ to plus level $u$, using $\mathcal K_0$ as the domain
when $v_0=0$.  Multiplication by $D$ preserves level and branch.  A path between these
spaces requires at least $u+v_0=r$ occurrences of the operators
$H_a$. At that length, the path must first move toward the common
degree-zero space and then outward along the plus branch. Thus all
lower-degree terms vanish on this block.

For the monomial $b_0H_{a_1}b_1\cdots H_{a_r}b_r$, the only
surviving path goes from the summand indexed by
$(a_r,\ldots,a_{u+1})$ to the summand indexed by
$(a_1,\ldots,a_u)$. Write $T_{b_0,\ldots,b_r}$ for this block operator.
Its value on $\zeta_0\otimes\cdots\otimes\zeta_{v_0}$ is
\[
 \left(\prod_{j=u+1}^{r}\tau(b_j\zeta_{r-j})\right)
 \bigl(b_0\otimes\cdots\otimes b_{u-1}\otimes b_u\zeta_{v_0}\bigr).
 \eqtag\label{eq:relcost-shortest-path}
\]
For fixed $a_1,\ldots,a_r$, the linear map
$b_0\otimes\cdots\otimes b_r\mapsto T_{b_0,\ldots,b_r}$ on
$D^{\otimes(r+1)}$ is injective.
Indeed, the first $u$ coefficients remain as tensor factors. Each of the last
$v_0$ coefficients determines a linear functional by the nondegenerate
pairing $(b,\zeta)\mapsto\tau(b\zeta)$, and the middle coefficient acts
by the faithful left regular representation on $L^2(D)$.  The tensor product of these injective maps is
injective, and the input and output summands distinguish all strings
of generator indices. Thus the degree-$r$ coefficients vanish. Descending in $r$, we obtain injectivity
on the whole polynomial space through degree $k$.

The space of formal homogeneous polynomials of degree $j$ with coefficients
in $D$ is the direct sum of $m_0^j$ copies of $D^{\otimes(j+1)}$,
giving \eqref{eq:relcost-filtered-dimension}.
Fix a basis $(P_\beta)$ of the formal polynomials of degree at most $k$, and, for self-adjoint $m_0$-tuples $\mathbf X$ in $M_N(\C)$, set
\[
\phi(\mathbf X)=\det\bigl(\operatorname{Tr}(P_\beta(\mathbf X)^*P_{\beta'}(\mathbf X))\bigr)_{\beta,\beta'}.
\]
Then $\phi$ is a real polynomial in the real coordinates of $\mathbf X$, and $\phi(H_1,\ldots,H_{m_0})\ne0$ for the tuple constructed above. Hence the set where $\phi$ vanishes has Lebesgue measure zero, and the matrices $P_\beta(\mathbf X)$ are linearly independent whenever $\phi(\mathbf X)\ne0$.
\end{proof}

In the next theorem, $d_v$ is again the size of the $v$th matrix summand of $D$; the numbers $d_k$ and $\mathsf d_k$, indexed by a degree $k$, are defined in the statement.

\begin{theorem}
\label{thm:relative-degree-motion-lower-bound}
Let $D\subset M_N(\C)$ be a deterministic unital $^*$-subalgebra, let $\mathbf X$ be any self-adjoint
tuple, and let $H$ be a normalized GUE matrix independent of $\mathbf X$.  In the real Hilbert space of
self-adjoint matrices with the $L^2(\operatorname{tr}_N)$ inner product, set
\[
 \begin{split}
 V_k&=\{P(\mathbf X):\deg P\le k,\ P(\mathbf X)^*=P(\mathbf X),\
                         E_D(P(\mathbf X))=0\},\\
 d_k&=\dim_\R V_k,
 \end{split}
 \eqtag\label{eq:relcost-centered-polynomial-space}
\]
where $P$ has coefficients in $D$.  Write $P_{V_k}$ for the real orthogonal projection onto
$V_k$. Suppose $u,v\in D'\cap U(N)$ and the conjugated families
$u\mathbf Xu^*$ and $vHv^*$ satisfy the amalgamated free-product moment
identities through degree $k+1$ for the original marginal laws. Then
\[
 4Nr\norm{H}_2\ge N\norm{P_{V_k}H}_2,
 \qquad r=\max\{\norm{u-1},\norm{v-1}\}.
 \eqtag\label{eq:relcost-conditional-chi}
\]
Conditionally on $\mathbf X$, the right side has the chi distribution with
$d_k$ degrees of freedom, with value zero when $d_k=0$. If one family stays fixed, $4$ may be replaced by $2$.

Now let $N\to\infty$, with deterministic unital $^*$-embeddings of a fixed
$C^*$-algebra $D=\bigoplus_v M_{d_v}(\C)$ into $M_N(\C)$. Suppose $\mathbf X_N$
has $m_0$ generators and a joint law absolutely continuous with respect to Lebesgue measure on the full space of self-adjoint $m_0$-tuples at every $N$.
Let $q$, $n_v$, and $T_k$ be as in Lemma~\ref{lem:full-filtered-relative-rank}.
Assume $n_v\ge\delta N$ for every $v$ and some fixed $\delta>0$.
Set $\rho=m_0q$ and assume $\rho>1$.  For every $k_N\to\infty$ satisfying
$n_v\ge d_vT_{k_N}$ for all $v$, and every $\eta>0$, with probability tending to one,
every pair $u,v\in D'\cap U(N)$ giving these amalgamated identities
through degree $k_N+1$ obeys
\[
 r\ge\frac{(1-\eta)\sqrt{\mathsf d_{k_N}}}{4N},
 \qquad \mathsf d_k=q\sum_{j=1}^k\rho^j.
 \eqtag\label{eq:relcost-unsaturated-bound}
\]
The condition $n_v\ge d_vT_{k_N}$ for all $v$ holds eventually for $k_N=\lfloor c\log N\rfloor$ whenever $0<c<2/\log\rho$.

Moreover, there is $c_0>0$, depending only on $D,\delta,m_0$, such
that, for every integer sequence $L_N\to\infty$, with probability tending
to one, every pair $u,v\in D'\cap U(N)$ giving the identities
through degree $L_N$ satisfies
\[
 r\ge c_0\min\left\{\frac{\rho^{(L_N-1)/2}}{N},1\right\}.
 \eqtag\label{eq:relcost-saturated-bound}
\]

Let $L_N\to\infty$ be integers, and let $r_N>0$ be deterministic
with $r_N\to0$. If, with
probability tending to one, there are $u,v\in D'\cap U(N)$ giving
the identities through degree $L_N$,
with both unitaries at operator-norm distance at most $r_N$ from
$1$, then
\[
 L_N\le\frac{2\log(Nr_N)}{\log\rho}+O(1).
 \eqtag\label{eq:relcost-degree-motion-bound}
\]
For $r_N=N^{-\gamma}$, this is
$L_N\le2(1-\gamma)\log N/\log\rho+O(1)$. The additive constant is uniform over
embeddings satisfying the multiplicity bound $n_v\ge\delta N$. Each probability bound holds on an event
where the stated inequality is valid for every correcting pair of
unitaries, including pairs chosen as arbitrary functions of the given matrices.
\end{theorem}

\begin{proof}
If $P\in V_k$, then $E_D(P)=0$, so the amalgamated free-product expectation of $PH$ vanishes.
Hence, by the moment identities for $u,v$, we have
$\operatorname{tr}_N(PwHw^*)=0$, where $w=u^*v$.  Thus
\[
 \norm{P_{V_k}H}_2\le\norm{H-wHw^*}_2
 \le2\norm{w-1}\norm{H}_2\le4r\norm{H}_2.
\]
If one unitary is $1$, the final constant is $2$.
Conditionally on $\mathbf X$, we choose an orthonormal basis
$P_1,\ldots,P_{d_k}$ of $V_k$ for the real inner product induced by
$\operatorname{tr}_N$. The variables $N\operatorname{tr}_N(P_jH)$ are
independent standard real Gaussians.  Hence, conditionally on $\mathbf X$, $N\norm{P_{V_k}H}_2$ has the chi distribution with $d_k$ degrees of freedom, and \eqref{eq:relcost-conditional-chi} holds simultaneously for all pairs $u,v\in D'\cap U(N)$ for which $u\mathbf Xu^*$ and $vHv^*$ satisfy the identities through degree $k+1$.

If $n_v\ge d_vT_k$ for every $v$, then by Lemma~\ref{lem:full-filtered-relative-rank}, since the law of $\mathbf X$ is absolutely continuous, almost surely no nonzero formal $D$-coefficient polynomial of degree at most $k$ vanishes at $\mathbf X$. By \eqref{eq:relcost-filtered-dimension}, the space of the values of these polynomials at $\mathbf X$ then has complex dimension $q\sum_{j=0}^k\rho^j$. It contains $D$ and is $^*$-closed, so its complex dimension equals the real dimension of its self-adjoint part. The restriction of $E_D$ to this self-adjoint part is onto $D_{\rm sa}$, of dimension $q$, so
$d_k=\mathsf d_k$ almost surely. A chi variable with $d$ degrees of
freedom, divided by $\sqrt d$, tends to $1$ in probability as
$d\to\infty$, and $\norm H_2\to1$. Together, these give
\eqref{eq:relcost-unsaturated-bound}.

For $\rho>1$, the geometric sums in
\eqref{eq:relcost-rank-multiplicity} give $T_k\asymp\rho^{k/2}$,
with constants depending only on $\rho$. Let $k_N^*$ be the largest integer $k$ such that $n_v\ge d_vT_k$ for every $v$. Since
$\delta N\le n_v\le N$,
\[
k_N^*=\frac{2\log N}{\log\rho}+O(1),
\]
uniformly over the embeddings, with the error depending only on
$D,\delta,m_0$. This also proves that $n_v\ge d_vT_{k_N}$ for all $v$ holds eventually for $k_N=\lfloor c\log N\rfloor$ when $c<2/\log\rho$.
Set $k_N=\min\{L_N-1,k_N^*\}$. Then $k_N\to\infty$ and
$\mathsf d_{k_N}\ge c_2\min\{\rho^{L_N-1},N^2\}$ for a fixed
$c_2>0$. Applying \eqref{eq:relcost-unsaturated-bound} with $\eta=1/2$
gives \eqref{eq:relcost-saturated-bound}.

Suppose that, with probability tending to one, there are unitaries
$u,v\in D'\cap U(N)$ satisfying the identities through degree $L_N$
and $\max\{\norm{u-1},\norm{v-1}\}\le r_N$.
For all sufficiently large $N$, this event intersects the event in
\eqref{eq:relcost-saturated-bound}, since both probabilities tend to one.
We choose a realization $(\mathbf X_N,H_N)$ in their intersection and unitaries $u,v$ with these properties. The lower bound for these unitaries gives
\[
r_N\ge c_0\min\left\{\frac{\rho^{(L_N-1)/2}}{N},1\right\}.
\]
Since $r_N\to0$, eventually $r_N<c_0$, so
$\rho^{(L_N-1)/2}\le Nr_N/c_0$.
Taking logarithms, we obtain \eqref{eq:relcost-degree-motion-bound}, as desired.
\end{proof}

If $\rho=1$, then necessarily $D=\C$ and $m_0=1$.  Then $T_k=k+1$, and by
Lemma~\ref{lem:full-filtered-relative-rank}, $d_k=k$ for $k\le N-1$ almost surely.
For $L_N\to\infty$, the conditional Gaussian estimate
\eqref{eq:relcost-conditional-chi} gives, with probability tending
to one, a lower bound of order
$\sqrt{\min\{L_N-1,N-1\}}/N$ on $\max\{\norm{u-1},\norm{v-1}\}$ for every pair $u,v$ giving the identities through degree $L_N$.
If both matrices have simple spectrum, then by
Proposition~\ref{prop:relative-degree-ceiling}, the joint identities can hold only through a degree that is at most a constant times $\log N$.

\subsection{Optimal correction costs}
\label{subsec:optimal-correction-costs}

By Corollary~\ref{cor:gue-relative-polynomial-images},
at every fixed degree there are unitaries at distance $O_{\mathbb P}(N^{-1})$ from $1$ that give the identities through that degree with probability tending to one.
However, by Proposition~\ref{prop:gue-motion-obstruction}, no
deterministic bound $C/N$ on the distance holds with probability tending
to one, already for the degree-two covariance equation. We now determine the limiting law of $N$ times the least distance for this equation and
for finite lists of polynomial covariance equations.

\begin{theorem}[Optimal degree-two correction cost]
\label{thm:gue-scalar-exact-cost}
For two independent normalized GUE matrices $G_N,H_N$, define
\[
 \mathfrak c_N=\min\bigl\{\norm{u-1}:u\in U(N),\quad
 \operatorname{tr}_N(uG_Nu^*H_N)
 =\operatorname{tr}_N(G_N)\operatorname{tr}_N(H_N)\bigr\}.
 \eqtag\label{eq:relcost-one-sided-minimum}
\]
The minimum exists for every pair of self-adjoint matrices $(G_N,H_N)$
and is a Borel function of this pair. Let $S,T$ be
free standard semicircular variables in a tracial von Neumann algebra with trace $\tau$, and set
\[
 \chi=\tau\abs{i[S,T]}=\frac12+\frac2\pi.
 \eqtag\label{eq:relcost-chi-value}
\]
Then
\[
 N\mathfrak c_N\ \xrightarrow{\ d\ }\ \frac{\abs{Z}}{\chi},
 \qquad Z\sim\mathcal N(0,1).
 \eqtag\label{eq:relcost-one-sided-law}
\]
If $\mathfrak b_N$ is instead the minimum of $\max\{\norm{u-1},\norm{v-1}\}$
over $u,v\in U(N)$ subject to the same identity for $uG_Nu^*,vH_Nv^*$, then
\[
 N\mathfrak b_N\ \xrightarrow{\ d\ }\ \frac{\abs{Z}}{2\chi}.
 \eqtag\label{eq:relcost-two-sided-law}
\]
\end{theorem}

\begin{proof}
We subtract the scalar traces, setting $A=G_N-\operatorname{tr}_N(G_N)$,
$C=H_N-\operatorname{tr}_N(H_N)$, and
\[
 \delta_N=\operatorname{tr}_N(AC),\qquad Q_N=i[A,C],\qquad
 \chi_N=\operatorname{tr}_N\abs{Q_N}.
\]
By Lemma~\ref{lem:optcost-feasibility-selection} below, with $P_1(t)=t$, the
one-sided minimum is attained and is Borel measurable.
Taking $v=1$, we see that the two-sided identity also has a solution. By the proof of Lemma~\ref{lem:optcost-feasibility-selection}, with $U(N)^2$ in place of $U(N)$, the two-sided minimum is attained and is Borel measurable.

The conditional Gaussian calculation \eqref{eq:relcost-conditional-normal}, now with
$D=\C$, gives $N\delta_N\xrightarrow{d}Z$.  Since the joint trace moments of $(G_N,H_N)$ converge to those of $(S,T)$, the moments of $Q_N$ converge to those of $i[S,T]$. By \eqref{eq:relgauss-norm-tail}, the norm $\norm{Q_N}$ is bounded by a fixed constant with probability tending to one.
Uniform polynomial approximation of $\abs{x}$ on a fixed interval then
yields $\chi_N\to\chi$ in probability.
Since $\tau((i[S,T])^2)=2$ and $\norm{i[S,T]}\le8$,
we have $\chi\ge1/4$.

For fixed $(G_N,H_N)$, let $B_N=\operatorname{sgn}(Q_N)$ and
$f_N(t)=\operatorname{tr}_N(e^{itB_N}Ae^{-itB_N}C)$.  Then
\[
 f_N(0)=\delta_N,\qquad f_N'(0)=\chi_N,\qquad
 \abs{f_N''(t)}\le4\norm{A}\norm{C}.
 \eqtag\label{eq:relcost-sign-derivative}
\]
On an event where $\|A\|,\|C\|$ are bounded by fixed constants and
$\chi_N\ge\chi/2$, the intermediate value theorem and the second-derivative
bound give a root
$t_N=-\delta_N/\chi_N+O(\delta_N^2)$ for sufficiently small
$\abs{\delta_N}$, with a uniform remainder constant.
Therefore
\[
 \mathfrak c_N\le\frac{\abs{\delta_N}}{\chi_N}+O(\delta_N^2)
 =O_{\mathbb P}(N^{-1}).
 \eqtag\label{eq:relcost-sign-upper}
\]
Conversely, let $u\in U(N)$ satisfy $\operatorname{tr}_N(uAu^*C)=0$,
with $r=\norm{u-1}<1$. Write $u=e^{ih}$ using its principal
self-adjoint logarithm. Then $\norm{h}=2\arcsin(r/2)=r+O(r^3)$,
and Taylor expansion gives
\[
 0=\delta_N+\operatorname{tr}_N(hQ_N)+O(\norm{h}^2),
 \qquad \abs{\operatorname{tr}_N(hQ_N)}\le\norm{h}\chi_N,
\]
with a uniform remainder constant on the same event. Applying this to a minimizing unitary and
using \eqref{eq:relcost-sign-upper}, we obtain
\[
 N\mathfrak c_N-\frac{\abs{N\delta_N}}{\chi_N}\longrightarrow0
 \quad\text{in probability}.
 \eqtag\label{eq:relcost-asymptotic-linearization}
\]
By Slutsky's theorem, we obtain \eqref{eq:relcost-one-sided-law}.

To compute $\chi$, we use the formula $\kappa_n(Q)=\operatorname{Tr}(M^n)$ for the free cumulants $\kappa_n(Q)$ of a quadratic form $Q=\sum_{a,b}M_{ab}S_aS_b$ in free standard semicircular variables $S_a$
\cite[Remark~2.6]{ejsmont-lehner}.
Nica and Speicher~\cite[Example~1.5(2)]{nica-speicher-commutators}
computed the law of the semicircular commutator. For $Q=i[S,T]$, the
coefficient matrix
$\left(\begin{smallmatrix}0&i\\-i&0\end{smallmatrix}\right)$
has eigenvalues $1,-1$.
Thus the odd cumulants of $Q$ vanish and its even cumulants equal $2$.  With the convention
$R_Q(w)=\sum_{n\ge1}\kappa_n(Q)w^{n-1}$ and Cauchy transform
$G_Q(z)=\tau((z-Q)^{-1})$,
\[
 R_Q(w)=\frac{2w}{1-w^2},\qquad
 z=\frac1{G_Q(z)}+\frac{2G_Q(z)}{1-G_Q(z)^2}.
 \eqtag\label{eq:relcost-cauchy-equation}
\]
The second equation is the relation $z=1/G_Q(z)+R_Q(G_Q(z))$ between the Cauchy transform and the $R$-transform, which holds near infinity and extends, after
clearing denominators, throughout the upper half-plane by analyticity.  Since the law of $Q$ is symmetric, we have $G_Q(it)=-is(t)$ with $s(t)>0$, hence
\[
 t=\frac{1-s(t)^2}{s(t)(1+s(t)^2)}.
\]
For $t>0$ this forces $0<s(t)<1$.  The function $s\mapsto(1-s^2)/(s(1+s^2))$ decreases strictly from infinity to zero on $(0,1)$.  Applying Tonelli's theorem to $\abs{x}=(2/\pi)\int_0^\infty x^2/(x^2+t^2)\,dt$, followed by
the change of variables given by the preceding display, we obtain
\[
 \begin{split}
 \tau\abs{Q}
 &=\frac2\pi\int_0^\infty(1-ts(t))\,dt\\
 &=\frac2\pi\int_0^1\left(\frac2{1+s^2}
       +\frac{4s^2(1-s^2)}{(1+s^2)^3}\right)\,ds
 =\frac2\pi\left(1+\frac\pi4\right).
 \end{split}
 \eqtag\label{eq:relcost-absolute-moment}
\]
For the last equality, the substitution $s=\tan\theta$, $0\le\theta\le\pi/4$, makes the two summands
integrate to $\pi/2$ and $1-\pi/4$.  This proves \eqref{eq:relcost-chi-value}.

Finally, if $u,v$ satisfy the identity for $uG_Nu^*,vH_Nv^*$, then $v^*u$ satisfies the identity in \eqref{eq:relcost-one-sided-minimum} and $\norm{v^*u-1}\le2\max\{\norm{u-1},\norm{v-1}\}$. Hence
$\mathfrak b_N\ge\mathfrak c_N/2$.  On the event $\mathfrak c_N<1$, we write a minimizing one-sided unitary as $w=e^{ih}$ and set $u=e^{ih/2}$, $v=e^{-ih/2}$.
This pair satisfies the two-sided identity, and $\max\{\norm{u-1},\norm{v-1}\}=\mathfrak c_N/2+O(\mathfrak c_N^3)$.
By \eqref{eq:relcost-one-sided-law} and \eqref{eq:relcost-sign-upper}, we then obtain \eqref{eq:relcost-two-sided-law}.
\end{proof}

We now fix real polynomials $P_1,\ldots,P_k$, write $P=(P_1,\ldots,P_k)$, and, for self-adjoint $G,H\in M_N(\C)$, impose the equations
\[
 \operatorname{tr}_N(uP_j(G)u^*H)
 =\operatorname{tr}_N(P_j(G))\operatorname{tr}_N(H)
\]
using one unitary $u$ and leaving $H$ fixed. We call such a unitary feasible.

By a \emph{symmetric convex body} we mean a compact convex set
$K=-K\subseteq\R^k$ containing a neighborhood of zero. Its
support function and gauge are defined by
\[
 s_K(z)=\max_{x\in K}z\cdot x,\qquad
 p_K(y)=\inf\{t\ge0:y\in tK\},
\]
and the gauge is a norm. We write $|\cdot|$ for the Euclidean norm on
$\R^k$ and $\|\cdot\|$ for the operator norm.

\begin{lemma}
\label{lem:optcost-feasibility-selection}
Fix integers $k,N\ge1$ and real polynomials $P_1,\ldots,P_k$. For self-adjoint
$G,H\in M_N(\C)$ set
\[
 A_j=P_j(G)-\operatorname{tr}_N(P_j(G)).
\]
There is a unitary $u$ for which
$\operatorname{tr}_N(uA_ju^*H)=0$ for every $j$. Consequently the cost
\[
 c_P(G,H)=\min\{\|u-1\|:u\in U(N),\
             \operatorname{tr}_N(uA_ju^*H)=0\ (1\le j\le k)\}
\]
is attained for every pair of self-adjoint matrices and lies in $[0,2]$.
As a function of $(G,H)$, $c_P$ is lower semicontinuous, hence Borel,
and a minimizing unitary can be chosen Borel measurably.
\end{lemma}

\begin{proof}
We diagonalize $G=V\operatorname{diag}(\lambda_b)V^*$ and $H=W\operatorname{diag}(\mu_a)W^*$. Set $F=(N^{-1/2}e^{2\pi iab/N})_{a,b=1}^N$. Since
$F$ is unitary and $|F_{ab}|^2=1/N$ for all $a,b$, the unitary $u=WFV^*$ satisfies, for every polynomial $P$,
\[
\operatorname{tr}_N(uP(G)u^*H)
=\frac1{N^2}\sum_{a,b}\mu_aP(\lambda_b)
=\operatorname{tr}_N(H)\operatorname{tr}_N(P(G)).
\]
Thus the set of feasible unitaries is nonempty and compact.

Write $x=(G,H)$ and $f(x,u)=(\operatorname{tr}_N(uA_ju^*H))_{j=1}^k$, which is continuous in $(x,u)$.
For each compact $D\subseteq U(N)$, the extended minimum
\[
v_D(x)=\min\{\|u-1\|:u\in D,\ f(x,u)=0\},
\qquad \min\varnothing=+\infty,
\]
has closed sublevel sets, by compactness of $D$.
In particular, $c_P=v_{U(N)}$ is lower semicontinuous.
Starting with $U(N)$, we recursively cover each compact set by a finite ordered list of nonempty compact subsets, with diameters tending to zero as the stage increases. A set $D$ in this construction meets the minimizing
set exactly when $v_D(x)=c_P(x)$, a Borel condition. Successively
choosing the first such subset gives a nested sequence whose
intersection is a single minimizing unitary. Fix a representative
point in each compact set. The representatives of the selected sets are Borel functions of $x$ and converge to that unitary. Hence a minimizing unitary can be chosen Borel measurably.
\end{proof}

\begin{lemma}
\label{lem:optcost-tangent-recovery}
Let $A_1,\ldots,A_k,H$ be self-adjoint matrices in $M_N(\C)$. Set
\[
 \begin{split}
 \Phi(h)&=\bigl(\operatorname{tr}_N(e^{ih}A_je^{-ih}H)\bigr)_{j=1}^k,
       \qquad h=h^*,\\
 v&=\Phi(0),\qquad B_j=i[A_j,H],\\
 Lh&=\bigl(\operatorname{tr}_N(hB_j)\bigr)_{j=1}^k,\qquad
 Q=\bigl(\operatorname{tr}_N(B_jB_\ell)\bigr)_{j,\ell=1}^k,\\
 K_N&=L\{h=h^*:\|h\|\le1\}.
 \end{split}
\]
Assume $Q>0$. The set $K_N$ is a symmetric convex body,
with support function and gauge
\begin{equation}
 \begin{split}
 s_N(z)&=\operatorname{tr}_N\left|\sum_jz_jB_j\right|,\\
 p_N(y)&=\min_{Lh=y}\|h\|
       =\sup_{z\ne0}\frac{|z\cdot y|}{s_N(z)}.
 \end{split}
 \label{eq:optcost-finite-tangent-duality}
\end{equation}
Define
\[
 \begin{split}
 C&=\max\left\{1,\ 2\|H\|
                     \left(\sum_j\|A_j\|^2\right)^{1/2}\right\},\\
 M&=\max\left\{1,\ \|Q^{-1}\|
                     \left(\sum_j\|B_j\|^2\right)^{1/2}\right\},
 \qquad r=p_N(v).
 \end{split}
\]
If $r\le\min\{1/4,(8MC)^{-1}\}$ and $h_0$ solves
$Lh_0=-v$ with the smallest possible operator norm
$\|h_0\|=r$, then there is $w$ such that $h=h_0+w$ satisfies
\begin{equation}
 \Phi(h)=0,\qquad \|w\|\le4MC r^2,
 \qquad \|e^{ih}-1\|\le r+4MC r^2.
 \label{eq:optcost-exact-tangent-recovery}
\end{equation}
Conversely, every unitary $u$ with $\operatorname{tr}_N(uA_ju^*H)=0$ for all $j$ and $a=\|u-1\|<1$ satisfies
\begin{equation}
 p_N(v)\le 2\arcsin(a/2)+MC\bigl(2\arcsin(a/2)\bigr)^2.
 \label{eq:optcost-arbitrary-unitary-lower}
\end{equation}
It follows that the minimum $c$ of $\|u-1\|$ over the unitaries $u$ with $\operatorname{tr}_N(uA_ju^*H)=0$ for all $j$ satisfies
$c=p_N(v)+O_{M,C}(p_N(v)^2)$ as $p_N(v)\to0$. At each fixed
$N$, the gauge is continuous in $(B_1,\ldots,B_k,y)$ on $Q>0$. The
equation $Lh=y$ has a solution of minimal operator norm that depends Borel measurably on $(B_1,\ldots,B_k,y)$. If $r\le\min\{1/4,(8MC)^{-1}\}$,
the exact solution $h=h_0+w$ can also be chosen as a Borel function of $(A_1,\ldots,A_k,H)$.
\end{lemma}

\begin{proof}
By cyclicity, $D\Phi(0)=L$. Trace duality gives
$s_N(z)=\operatorname{tr}_N|\sum_jz_jB_j|$. Define $J$ by $Jy=\sum_j(Q^{-1}y)_jB_j$. Then
\[
LJ=I,\qquad \|Jy\|\le M|y|.
\]
Thus $K_N$ is compact, convex, symmetric, and contains the Euclidean
ball of radius $1/M$. Compactness and separation give the minimum and the
dual formula in \eqref{eq:optcost-finite-tangent-duality}. In particular,
$p_N(y)\le M|y|$.

Differentiating the exponentials twice gives $|D^2\Phi(h)[b,b']|\le2C\|b\|\|b'\|$ globally.
For $R(h)=\Phi(h)-v-Lh$, it follows that
\begin{equation}
 |R(h)|\le C\|h\|^2,\qquad
 |R(h)-R(h')|\le C(\|h\|+\|h'\|)\|h-h'\|.
 \label{eq:optcost-simultaneous-remainder}
\end{equation}
We choose an optimal $h_0$ with $Lh_0=-v$ and $\|h_0\|=r$. If $r=0$, we take
$w=0$. Otherwise, we set $d:=4MC r^2\le r/2$.
If $\|w\|\le d$, then
$\|JR(h_0+w)\|\le MC(r+d)^2\le(9/4)MC r^2\le d$.
Thus $w\mapsto-JR(h_0+w)$ maps this ball into itself and has
Lipschitz constant at most $2MC(r+d)\le3/8$. For its fixed point $w$, the matrix $h=h_0+w$ satisfies \eqref{eq:optcost-exact-tangent-recovery}.

Conversely, let $u$ be a unitary with $\operatorname{tr}_N(uA_ju^*H)=0$ for all $j$ and $a=\|u-1\|<1$. We write $u=e^{ih}$, using the principal branch. Then $h=h^*$ has norm
$\ell=2\arcsin(a/2)$. Since $v=-Lh-R(h)$,
\[
p_N(v)\le p_N(Lh)+p_N(R(h))\le\ell+MC\ell^2.
\]
This proves \eqref{eq:optcost-arbitrary-unitary-lower}. Applying it to a
minimizer, using the upper bound from \eqref{eq:optcost-exact-tangent-recovery}
and $2\arcsin(c/2)=c+O(c^3)$,
we obtain $c=r+O_{M,C}(r^2)$.

On $Q>0$, the dual formula is a maximum over the unit sphere with a
locally uniformly positive denominator, so the gauge is continuous.
For $y\ne0$, we choose $g$ in the self-adjoint unit ball with $Lg=y/p_N(y)$ by the selection argument in the proof of Lemma~\ref{lem:optcost-feasibility-selection}, applied to the self-adjoint unit ball in place of $U(N)$: a compact set $D$ in the construction meets this solution set exactly when $\min_{g\in D}|Lg-y/p_N(y)|=0$, a Borel condition. Then $p_N(y)g$ solves $Lh=y$ with smallest operator norm. When $y=0$, we use
$h=0$. Finally, we apply this selection to $y=-v$ and iterate
$w_{n+1}=-JR(h_0+w_n)$ from zero. The limit gives a Borel choice
of an exact solution.
\end{proof}

\begin{lemma}
\label{lem:optcost-support-convergence}
For fixed $k$, let $B_{1,N},\ldots,B_{k,N}$ be self-adjoint matrices that
are uniformly bounded in operator norm for all sufficiently large $N$.
Suppose their joint normalized trace moments converge to those of
self-adjoint $B_1,\ldots,B_k$ in a tracial von Neumann algebra $(\mathcal M,\tau)$,
and that $Q=(\tau(B_jB_\ell))_{j,\ell}>0$. Define
\[
 \begin{split}
 K_N&=\{(\operatorname{tr}_N(hB_{j,N}))_{j=1}^k:
                         h=h^*,\ \|h\|\le1\},\\
 K&=\{(\tau(hB_j))_{j=1}^k:h=h^*\in\mathcal M,\ \|h\|\le1\}.
 \end{split}
\]
Then $K$ is a symmetric convex body with support function
$s(z)=\tau|\sum_jz_jB_j|$. Moreover,
\begin{equation}
 \sup_{|z|=1}\left|
     \operatorname{tr}_N\left|\sum_jz_jB_{j,N}\right|-s(z)\right|
 \longrightarrow0.
 \label{eq:optcost-uniform-support-limit}
\end{equation}
The bodies $K_N$ are full-dimensional for all sufficiently large $N$,
converge to $K$ in Hausdorff distance, and their gauges converge to that
of $K$ uniformly on every bounded subset of $\R^k$.
\end{lemma}

\begin{proof}
Since the maps $h\mapsto\tau(hB_j)$ are weak-$*$ continuous and the self-adjoint unit ball of $\mathcal M$ is weak-$*$ compact, $K$ is compact. Trace duality gives its support
function. With
$b=(\sum_j\|B_j\|^2)^{1/2}$ and
$\lambda=\lambda_{\min}(Q)>0$, for $z\ne0$ we have
\[
s(z)\ge\frac{\tau((\sum_jz_jB_j)^2)}{\|\sum_jz_jB_j\|}
\ge\rho|z|,\qquad \rho=\lambda/b>0.
\]
Thus $K$ contains the Euclidean ball of radius $\rho$.

Uniform polynomial approximation of $|t|$ on a common spectral interval
gives $s_N(z)=\operatorname{tr}_N|\sum_jz_jB_{j,N}|\to s(z)$ for each fixed
$z$. By the trace-norm triangle
inequality, these support functions are uniformly Lipschitz on the unit
sphere. For each $\theta>0$, we choose a finite $\theta$-net of the unit sphere. Since $s_N\to s$ at each point of this net and $\theta$ is arbitrary, we obtain \eqref{eq:optcost-uniform-support-limit}.
Let $\epsilon_N$ be the supremum in \eqref{eq:optcost-uniform-support-limit} and $\eta_N=\epsilon_N/\rho$.
For $\eta_N<1$, homogeneity and separation give
\begin{equation}
 (1-\eta_N)K\subseteq K_N\subseteq(1+\eta_N)K,\qquad
 \frac{p_K(y)}{1+\eta_N}\le p_{K_N}(y)
 \le\frac{p_K(y)}{1-\eta_N}.
 \label{eq:optcost-homothetic-body-comparison}
\end{equation}
Since $K$ is bounded and $p_K(y)\le |y|/\rho$, these imply the asserted
Hausdorff and locally uniform gauge convergence.
Finally, the matrices $(\operatorname{tr}_N(B_{j,N}B_{\ell,N}))_{j,\ell}$ converge to $Q$, and their smallest eigenvalues are eventually at least $\lambda/2$.
\end{proof}

\begin{theorem}[Optimal cost of two Gaussian moment identities]
\label{thm:optcost-two-moment-cost}
Let $G_N,H_N$ be independent normalized GUE matrices in
$M_N(\C)$, with
$\mathbb E(H_N)_{ab}(H_N)_{cd}=N^{-1}\delta_{ad}\delta_{bc}$,
and set
\[
 A_{1,N}=G_N-\operatorname{tr}_N(G_N),
 \qquad A_{2,N}=G_N^2-\operatorname{tr}_N(G_N^2).
\]
Define
\begin{equation}
 \mathfrak c_N^{(2)}=
 \min\{\|u-1\|:u\in U(N),\
       \operatorname{tr}_N(uA_{j,N}u^*H_N)=0\ (j=1,2)\}.
 \label{eq:optcost-two-moment-minimum}
\end{equation}
The minimum is attained for every pair of self-adjoint matrices
$(G_N,H_N)$. Its value is a lower semicontinuous, hence Borel,
function of this pair, and a minimizing unitary can be chosen Borel
measurably.

Let $S,T$ be free standard semicircular variables in a tracial von Neumann algebra with trace $\tau$, and set
\begin{equation}
 \begin{split}
 B_1&=i[S,T],\qquad B_2=i[S^2,T],\\
 \mathcal K&=\{(\tau(hB_1),\tau(hB_2)):
                  h=h^*\in W^*(S,T),\ \|h\|\le1\},\\
 s(z)&=\tau\left|i[z_1S+z_2S^2,T]\right|.
 \end{split}
 \label{eq:optcost-limiting-body}
\end{equation}
The set $\mathcal K$ is a symmetric convex body with support
function $s$. Writing $\mathbb B_2$ for the Euclidean
unit ball in $\R^2$, we have
\begin{equation}
 \frac1{2\sqrt{13}}\mathbb B_2\subseteq\mathcal K
                         \subseteq\sqrt2\,\mathbb B_2.
 \label{eq:optcost-body-nondegeneracy}
\end{equation}
Its gauge is the norm
\[
 p_{\mathcal K}(x)=\inf\{t\ge0:x\in t\mathcal K\}
                 =\sup_{z\ne0}\frac{|z\cdot x|}{s(z)}.
\]
For $Z=(Z_1,Z_2)$ with independent standard real Gaussian coordinates,
\begin{equation}
 N\mathfrak c_N^{(2)}\ \xrightarrow{\ d\ }\ p_{\mathcal K}(Z).
 \label{eq:optcost-cost-limit}
\end{equation}
More precisely, if
\[
 v_N=(\operatorname{tr}_N(A_{j,N}H_N))_{j=1}^2,
 \quad B_{j,N}=i[A_{j,N},H_N],
 \quad\mathcal K_N=\{(\operatorname{tr}_N(hB_{j,N}))_{j=1}^2:
                                  h=h^*,\ \|h\|\le1\},
\]
then $\mathcal K_N$ is full-dimensional almost surely for $N\ge3$ and
\begin{equation}
 \mathfrak c_N^{(2)}=p_{\mathcal K_N}(v_N)+O_{\mathbb P}(N^{-2}),
 \qquad
 N\mathfrak c_N^{(2)}-p_{\mathcal K}(Nv_N)
                    \longrightarrow0\quad\text{in probability}.
 \label{eq:optcost-finite-tangent-cost}
\end{equation}
On the null set where $\mathcal K_N$ is not full-dimensional, any
Borel real-valued extension of the gauge gives the same probability
statements.
The limiting distribution is specified, for every $t\ge0$, by
\begin{equation}
 \mathbb P\{p_{\mathcal K}(Z)\le t\}
 =\frac1{2\pi}\int_{t\mathcal K}e^{-|x|^2/2}\,dx.
 \label{eq:optcost-gaussian-body-distribution}
\end{equation}
\end{theorem}

\begin{proof}
By Lemma~\ref{lem:optcost-feasibility-selection}, with $P_1(t)=t$ and
$P_2(t)=t^2$, the minimum is attained for every pair of self-adjoint
matrices, its value is lower semicontinuous, and a minimizing unitary can be
chosen Borel measurably.
Since the conclusions are in distribution and in probability, we may realize the matrices of all
dimensions on one probability space, for instance independently.
By Lemma~\ref{lem:relgauss-raw-moments}, Chebyshev's inequality and the
Borel--Cantelli lemma, there is one event of probability one on which all joint trace moments of $(G_N,H_N)$ converge and $\|G_N\|,\|H_N\|$ are eventually at most $3$. For these norm bounds, we use $r=\lceil6\log N\rceil$ in \eqref{eq:relgauss-norm-tail}.

Conditionally on $G_N$, the vector $Nv_N$ is centered Gaussian with
covariance
\begin{equation}
 \Gamma_N=
 \begin{pmatrix}
 m_{2,N}-m_{1,N}^2&m_{3,N}-m_{1,N}m_{2,N}\\
 m_{3,N}-m_{1,N}m_{2,N}&m_{4,N}-m_{2,N}^2
 \end{pmatrix}
 \longrightarrow I_2,\qquad m_{\ell,N}=\operatorname{tr}_N(G_N^\ell).
 \label{eq:optcost-exact-conditional-covariance}
\end{equation}
Conditional characteristic functions therefore give $Nv_N\xrightarrow{d}Z$.

For centered $a,b\in W^*(S)$, freeness gives
\begin{equation}
 \tau\bigl(i[a,T]\,i[b,T]\bigr)=2\tau(ab).
 \label{eq:optcost-free-commutator-gram}
\end{equation}
Since $S,S^2-1$ have Gram matrix $I_2$, we have $(\tau(B_jB_\ell))_{j,\ell=1}^2=2I_2$. Moreover,
\[
\|z_1B_1+z_2B_2\|\le8|z_1|+12|z_2|\le4\sqrt{13}|z|,
\qquad \tau((z_1B_1+z_2B_2)^2)=2|z|^2.
\]
With $B=z_1B_1+z_2B_2$, the inequalities $\tau|B|\ge\tau(B^2)/\|B\|$ and $\tau|B|\le\tau(B^2)^{1/2}$ give $s(z)\ge|z|/(2\sqrt{13})$ and $s(z)\le\sqrt2|z|$, and hence
\eqref{eq:optcost-body-nondegeneracy}.
We apply Lemma~\ref{lem:optcost-support-convergence} to the commutators $B_{1,N},B_{2,N}$ on the preceding probability-one event. Then the gauges of
$\mathcal K_N$ converge to the gauge of $\mathcal K$ uniformly
on bounded sets. The constants $M,C$ in
Lemma~\ref{lem:optcost-tangent-recovery} are eventually bounded on this
event.

For each $N\ge3$, almost surely $G_N$ has simple spectrum and all
off-diagonal entries of $H_N$ in an eigenbasis of $G_N$ are nonzero. If a real
linear combination $z_1B_{1,N}+z_2B_{2,N}$ of the two commutators vanished, the polynomial
$z_1t+z_2t^2$ would be constant on at least three distinct eigenvalues.
This forces $z_1=z_2=0$, so $\mathcal K_N$ is full-dimensional.

Now $r_N=p_{\mathcal K_N}(v_N)\le M|v_N|=O_{\mathbb P}(N^{-1})$.
By Lemma~\ref{lem:optcost-tangent-recovery}, we have $\mathfrak
c_N^{(2)}=r_N+O_{\mathbb P}(N^{-2})$.
Since $Nv_N$ is bounded in probability and the gauges converge uniformly on
bounded sets, we have
$p_{\mathcal K_N}(Nv_N)-p_{\mathcal K}(Nv_N)\to0$ in probability.
By homogeneity and the continuous mapping theorem, we obtain \eqref{eq:optcost-finite-tangent-cost} and \eqref{eq:optcost-cost-limit}.
The Gaussian density gives \eqref{eq:optcost-gaussian-body-distribution}.
\end{proof}

\begin{corollary}
\label{cor:optcost-polynomial-cost}
Let $P_1,\ldots,P_k\in\R[t]$ be fixed, with $1,P_1,\ldots,P_k$
linearly independent over $\R$, and let $G_N,H_N$ have the law
in Theorem~\ref{thm:optcost-two-moment-cost}. Set
\[
 A_{j,N}=P_j(G_N)-\operatorname{tr}_N(P_j(G_N)),
 \qquad v_{P,N}=(\operatorname{tr}_N(A_{j,N}H_N))_{j=1}^k,
\]
and let $\mathfrak c_{P,N}$ be the minimum of $\|u-1\|$ over the unitaries $u\in U(N)$ with $\operatorname{tr}_N(uA_{j,N}u^*H_N)=0$ for $1\le j\le k$. It is attained for every
pair of self-adjoint matrices $(G_N,H_N)$. Its value is a lower semicontinuous, hence Borel, function of this pair, and a minimizing unitary can be chosen Borel measurably.

For free standard semicircular variables $S,T$ in a tracial von Neumann algebra with trace $\tau$, set
\[
 \begin{split}
 a_j&=P_j(S)-\tau(P_j(S)),\qquad
 \Gamma_P=(\tau(a_ja_\ell))_{j,\ell=1}^k,\\
 \mathcal K_P&=\{(\tau(h\,i[P_j(S),T]))_{j=1}^k:
                    h=h^*\in W^*(S,T),\ \|h\|\le1\}.
 \end{split}
\]
Then $\Gamma_P>0$, $\mathcal K_P$ is a symmetric convex body, and
\begin{equation}
 s_{\mathcal K_P}(z)=\tau\left|i\left[\sum_jz_jP_j(S),T\right]\right|,
 \qquad
 N\mathfrak c_{P,N}\ \xrightarrow{\ d\ }\ p_{\mathcal K_P}(Z_P),
 \quad Z_P\sim\mathcal N_k(0,\Gamma_P).
 \label{eq:optcost-polynomial-cost-law}
\end{equation}
Writing
\[
 B_{j,N}=i[P_j(G_N),H_N],\qquad
 \mathcal K_{P,N}
 =\{(\operatorname{tr}_N(hB_{j,N}))_{j=1}^k:
             h=h^*\in M_N(\C),\ \|h\|\le1\},
\]
we have
\[
 \mathfrak c_{P,N}
 =p_{\mathcal K_{P,N}}(v_{P,N})+O_{\mathbb P}(N^{-2}),\qquad
 N\mathfrak c_{P,N}-p_{\mathcal K_P}(Nv_{P,N})\longrightarrow0
 \quad\text{in probability}.
\]
The body $\mathcal K_{P,N}$ is full-dimensional almost surely whenever
$N>\max_j\deg P_j$. As in Theorem~\ref{thm:optcost-two-moment-cost}, its gauge
may be extended by a Borel real-valued function on the null event where
$\mathcal K_{P,N}$ is not full-dimensional.
\end{corollary}

\begin{proof}
Feasibility, attainment and Borel selection follow from
Lemma~\ref{lem:optcost-feasibility-selection}. A polynomial that is constant
almost everywhere
for the semicircular measure is identically constant. Thus linear
independence of the $P_j$ modulo constants is equivalent to linear
independence of the centered vectors $a_j$ in $L^2$, and
$\Gamma_P>0$.

Conditionally on $G_N$, the vector $Nv_{P,N}$ is centered Gaussian with covariance
\[
(\Gamma_{P,N})_{j\ell}
=\operatorname{tr}_N(P_j(G_N)P_\ell(G_N))
-\operatorname{tr}_N(P_j(G_N))\operatorname{tr}_N(P_\ell(G_N))
\longrightarrow(\Gamma_P)_{j\ell}.
\]
By \eqref{eq:optcost-free-commutator-gram}, we have $(\tau(i[P_j(S),T]\,i[P_\ell(S),T]))_{j,\ell=1}^k=2\Gamma_P$.
We now argue as in the proof of Theorem~\ref{thm:optcost-two-moment-cost}. Since the $P_j$ are fixed polynomials, on the event of probability one in that proof, the matrices $A_{j,N}$ and $B_{j,N}$ are eventually bounded in norm and their joint trace moments converge. By Lemmas~\ref{lem:optcost-support-convergence} and \ref{lem:optcost-tangent-recovery}, we obtain the two relations for $\mathfrak c_{P,N}$ in the statement and \eqref{eq:optcost-polynomial-cost-law}.

For $N>\max_j\deg P_j$, we argue in an eigenbasis of $G_N$ as in the proof of Theorem~\ref{thm:optcost-two-moment-cost}. If a real linear combination $\sum_jz_jB_{j,N}$ of the commutators vanished, then the polynomial $\sum_jz_jP_j$, of degree less than $N$, would be constant on the $N$ distinct eigenvalues of $G_N$. Since $1,P_1,\ldots,P_k$ are linearly independent, this forces $z=0$, so $\mathcal K_{P,N}$ is full-dimensional.
\end{proof}

For $k=1$ and $P_1(t)=t$, we have $\Gamma_P=1$ and $\mathcal K_P=[-\chi,\chi]$, with $\chi$ as in \eqref{eq:relcost-chi-value}, so \eqref{eq:optcost-polynomial-cost-law} gives the limit \eqref{eq:relcost-one-sided-law} of Theorem~\ref{thm:gue-scalar-exact-cost}. For $k=2$ and $(P_1,P_2)=(t,t^2)$, we have $\Gamma_P=I_2$ and $\mathcal K_P=\mathcal K$, so \eqref{eq:optcost-polynomial-cost-law} gives the limit \eqref{eq:optcost-cost-limit} of Theorem~\ref{thm:optcost-two-moment-cost}.

Replacing the $P_j$ by invertible real linear combinations of them,
or adding constants, gives the same constraint set and the same
minimum cost. Explicitly, for $\widetilde P=CP+d$ with
$C\in GL_k(\R)$ and $d\in\R^k$,
\[
 \Gamma_{\widetilde P}=C\Gamma_PC^{\mathsf T},\qquad
 \mathcal K_{\widetilde P}=C\mathcal K_P,\qquad
 p_{\mathcal K_{\widetilde P}}(Cx)=p_{\mathcal K_P}(x).
\]

In the next proposition and its proof, $G_N$, $H_N$, $A_{1,N}$, $A_{2,N}$, $\mathfrak c_N^{(2)}$, $v_N$, $B_{j,N}$, $S$, $T$, $\tau$, and $Z$ are as in Theorem~\ref{thm:optcost-two-moment-cost}.

\begin{proposition}
\label{prop:optcost-genuine-joint-cost}
For the body $\mathcal K$ in
Theorem~\ref{thm:optcost-two-moment-cost}, let
\[
 \chi_1=\tau|i[S,T]|=\frac12+\frac2\pi,
 \qquad \chi_2=\tau|i[S^2,T]|>0,
 \qquad R_\chi=[-\chi_1,\chi_1]\times[-\chi_2,\chi_2].
\]
The only point of $\mathcal K$ with first coordinate $\chi_1$ is
$(\chi_1,0)$. In particular, $\mathcal K$ is a proper subset of
$R_\chi$, and
\[
 p_{\mathcal K}(x)\ge
     \max\{|x_1|/\chi_1,|x_2|/\chi_2\},
\]
with strict inequality on a nonempty open subset of $\R^2$.
If $\Phi_{\rm G}$ denotes the standard real Gaussian distribution
function, then for every $t>0$,
\begin{equation}
 \mathbb P\{p_{\mathcal K}(Z)\le t\}
 <\bigl(2\Phi_{\rm G}(t\chi_1)-1\bigr)
  \bigl(2\Phi_{\rm G}(t\chi_2)-1\bigr).
 \label{eq:optcost-strict-joint-distribution}
\end{equation}
Moreover, let $\mathfrak c_{j,N}^{\rm sep}$ be the minimum of $\|u-1\|$ over the unitaries $u\in U(N)$ with $\operatorname{tr}_N(uA_{j,N}u^*H_N)=0$. Then jointly
\begin{equation}
 (N\mathfrak c_{1,N}^{\rm sep},N\mathfrak c_{2,N}^{\rm sep},
       N\mathfrak c_N^{(2)})
 \ \xrightarrow{\ d\ }\
 (|Z_1|/\chi_1,|Z_2|/\chi_2,p_{\mathcal K}(Z)).
 \label{eq:optcost-joint-versus-separate}
\end{equation}
Thus the limiting law of $N\mathfrak c_N^{(2)}$ differs from that of $N\max\{\mathfrak c_{1,N}^{\rm sep},\mathfrak c_{2,N}^{\rm sep}\}$.
\end{proposition}

\begin{proof}
Let $B_1=i[S,T]$ and $B_2=i[S^2,T]$. By the Cauchy-transform calculation
\eqref{eq:relcost-cauchy-equation}, with $s(t)$ as in the proof of
Theorem~\ref{thm:gue-scalar-exact-cost}, we have
\[
\tau(1_{\{0\}}(B_1))=\lim_{t\downarrow0}t s(t)=0,\qquad 0<s(t)<1.
\]
Thus $B_1$ has zero kernel (see also~\cite[Section~2.5.3]{ejsmont-lehner}).
If $h$ is a self-adjoint contraction with $\tau(hB_1)=\tau|B_1|$, then $h=\operatorname{sgn}(B_1)$, since $\tau$ is faithful.
Since $(-S,T)$ is again a pair of free standard semicircular variables, there is a trace-preserving
automorphism of $W^*(S,T)$ sending $(S,T)$ to $(-S,T)$. It negates $B_1$, hence $\operatorname{sgn}(B_1)$, and fixes $B_2$, so
$\tau(\operatorname{sgn}(B_1)B_2)=0$.
Thus $(\chi_1,0)$ is the only point of $\mathcal K$ whose first
coordinate is $\chi_1$.

Since $\abs{\tau(hB_j)}\le\tau\abs{B_j}=\chi_j$ for every self-adjoint contraction $h$, we have $\mathcal K\subseteq R_\chi$. However,
$(\chi_1,\chi_2)\notin\mathcal K$. Since $\mathcal K$ is compact and $\chi_2>0$, the difference
$R_\chi\setminus\mathcal K$ has positive area, and the inequality $p_{\mathcal K}(x)\ge\max\{|x_1|/\chi_1,|x_2|/\chi_2\}$ is strict on a nonempty open set.
Every positive dilation of this difference has positive Gaussian
measure, proving \eqref{eq:optcost-strict-joint-distribution}.

Applying Lemma~\ref{lem:optcost-tangent-recovery} to each separate constraint and writing $v_N=(v_{1,N},v_{2,N})$, we have
\[
\mathfrak c_{j,N}^{\rm sep}
=\frac{|v_{j,N}|}{\operatorname{tr}_N|B_{j,N}|}
+O_{\mathbb P}(N^{-2}),\qquad
\operatorname{tr}_N|B_{j,N}|\longrightarrow\chi_j>0.
\]
Together with \eqref{eq:optcost-finite-tangent-cost} and $Nv_N\xrightarrow{d}Z$,
these simultaneous approximations prove
\eqref{eq:optcost-joint-versus-separate}.
\end{proof}

\part{Free groups and \texorpdfstring{$q$}{q}-Gaussian families}
For free groups, we round independent Haar unitaries to clocks and use small
conjugations to make the traces of all short nontrivial words vanish
(Section~\ref{sec:clock-representations}). We then study the traces of words in the Haar unitaries
themselves (Section~\ref{subsec:haar-word-fluctuations}) and construct strongly convergent Gaussian matrix
models for mixed $q$-Gaussian families (Section~\ref{sec:mixed-q-strong-models}).

\section{Clock representations}
\label{sec:clock-representations}

In Corollary~\ref{cor:exact-character-strong-models}, we saw that, for the free product $G*H$ of two nontrivial finite groups, small conjugations make the normalized traces of short nontrivial words vanish while retaining strong convergence. In Proposition~\ref{prop:finite-image-log-character-models}, we obtained strongly convergent finite-image representations of $\mathbb F_s$ whose normalized traces vanish on the nontrivial words of length at most a constant times $\log N$. The representations constructed there are reducible. Here, for $\mathbb F_s$, we also prescribe the spectra of the generators.
We show that, for each fixed $s\ge2$, small perturbations of independent
Haar unitary matrices give strongly convergent irreducible representations of $\mathbb F_s$ in
every
sufficiently large dimension, with each generator sent to a conjugate of the cyclic
shift (Theorem~\ref{thm:irreducible-clock-strong-models}). Their normalized traces vanish on all nontrivial words of length at most a constant times $\log N$, and this order of growth is optimal (Proposition~\ref{prop:exact-character-radius-ceiling}). They satisfy Hastings's asymptotically optimal estimate for the norm, on the traceless matrices, of the average of the conjugations by the generators and their inverses, and from it we identify the relative commutants of the induced embeddings of $L(\mathbb F_s)$ into ultraproducts (Corollary~\ref{cor:ramanujan-clock-ultraproduct-commutants}).

A \emph{clock} in $U(N)$ is a unitary whose eigenvalues are the $N$th roots
of unity, each with multiplicity one; a \emph{clock representation} of $\mathbb F_s$ is a
homomorphism $\mathbb F_s\to U(N)$ that sends each of the standard generators to a clock. Write $\zeta_N=e^{2\pi i/N}$.
A unitary is a clock if and only if it is conjugate to the clock matrix
$\operatorname{diag}(1,\zeta_N,\ldots,\zeta_N^{N-1})$, or equivalently to the cyclic shift. A unitary
$U\in U(N)$ satisfies $\operatorname{tr}_N(U^m)=0$ for $1\le m<N$
precisely when $U$ is a scalar multiple of a clock. No $U\in U(N)$
can satisfy these identities for all $1\le m\le N$. Thus a clock has the same moments as a Haar unitary in a tracial von Neumann algebra up to order $N-1$, the largest order possible in $U(N)$.
For the free group $\mathbb F_s$ with standard generators $a_1,\ldots,a_s$,
we write $\abs g$ for reduced word length and $B_{\mathbb F_s}(r)$ for the
corresponding closed ball. For unitaries $V_1,\ldots,V_s$, we write $\mathbf V=(V_1,\ldots,V_s)$ and $g(\mathbf V)$ for the image of $g\in\mathbb F_s$ under the homomorphism sending $a_i$ to $V_i$, and similarly for other letters. As in Corollary~\ref{cor:large-matrix-relative-privacy}, $\norm{\Psi}_{2\to2}$ denotes the operator norm, with respect to $\norm{\cdot}_2$, of a linear map $\Psi$ on $M_N(\C)$ or on a subspace of it.

\begin{theorem}[Strongly convergent clock representations]
\label{thm:irreducible-clock-strong-models}
Fix $s\ge2$, let $a_1,\ldots,a_s$ be the standard generators of $\mathbb F_s$, and set
\[
\lambda_s:=\frac{\sqrt{2s-1}}s,
\qquad
R_N:=\left\lfloor\frac{\log N}{16\log(2s)}\right\rfloor.
\]
For some $N_s^0$ depending only on $s$, there are representations $\pi_N:\mathbb F_s\to U(N)$, $N\ge N_s^0$, such that, with $V_{i,N}:=\pi_N(a_i)$:
\begin{enumerate}[label=\textup{(\roman*)}]
\item each $V_{i,N}$ is a clock;
\item $\operatorname{tr}_N(\pi_N(g))=\delta_{g,e}$ whenever $\abs g\le R_N$;
\item $\pi_N$ converges strongly to the left regular representation: for every $f\in\C[\mathbb F_s]$,
\[
\operatorname{tr}_N(\pi_N(f))\longrightarrow\tau_{\mathbb F_s}(f),
\qquad
\norm{\pi_N(f)}_\infty\longrightarrow\norm{\lambda_{\mathbb F_s}(f)}_\infty;
\]
\item the averaging operators
\[
T_N(x):=\frac1{2s}\sum_{i=1}^s
\bigl(V_{i,N}xV_{i,N}^*+V_{i,N}^*xV_{i,N}\bigr)
\]
on $L^2(M_N(\C),\operatorname{tr}_N)$ are self-adjoint, leave the subspace $L^2_0(M_N)$ of traceless matrices invariant, and satisfy
\[
\rho_N:=\left\lVert T_N\big|_{L^2_0(M_N)}\right\rVert_{2\to2}
\longrightarrow\lambda_s.
\]
\end{enumerate}

For every traceless $x\in M_N(\C)$,
\[
\frac1s\sum_{i=1}^s\norm{[V_{i,N},x]}_2^2
\ge2(1-\rho_N)\norm{x}_2^2.
\eqtag
\label{eq:clock-poincare}
\]
Consequently, for all sufficiently large $N$, we have $C^*(\pi_N(\mathbb F_s))=M_N(\C)$, and $\operatorname{tr}_N\circ\pi_N$ is an extreme point of the convex set of normalized characters of $\mathbb F_s$. Let $\operatorname{Red}_{s,N}$ be the set of tuples $(W_1,\ldots,W_s)\in U(N)^s$ with $C^*(W_1,\ldots,W_s)\ne M_N(\C)$. Then, for the maximum operator-norm metric on $U(N)^s$, the tuple $\mathbf V_N:=(V_{1,N},\ldots,V_{s,N})$ satisfies
\[
\operatorname{dist}_\infty(\mathbf V_N,\operatorname{Red}_{s,N})
\ge\sqrt{\frac{1-\rho_N}{2}}.
\eqtag
\label{eq:clock-distance-reducible}
\]
\end{theorem}

To prove this theorem, we start with independent Haar unitary matrices. Since conjugation does not change the traces of the powers of a unitary, we first replace the eigenvalues of each, in their circular order, by the $N$th roots of unity, each occurring once. We then conjugate the first $s-1$ of the resulting clocks by unitaries close to $1$, leaving the last fixed, so that the traces of mixed words of length at most $R_N$ vanish as well (Proposition~\ref{prop:quantitative-haar-clock-exactification}). We prove the theorem after that proposition.

\subsection{Group-word correction}

Fix integers $s,D\ge2$.
We consider the mixed words
\[
w=b_1\cdots b_k,
\qquad
b_j=a_{c_j}^{m_j},
\qquad
k\ge2,
\qquad
\sum_{j=1}^k\abs{m_j}\le D,
\]
where $m_j\in\mathbb Z\setminus\{0\}$ and $c_j\ne c_{j+1}$ cyclically.
We call these maximal blocks $b_j$ the syllables of $w$.
Let $\mathcal W_D$ contain exactly one representative from each class modulo cyclic syllable rotation and inversion.

For $w=b_1\cdots b_k\in\mathcal W_D$, set
\[
R_j(w):=b_j\cdots b_kb_1\cdots b_{j-1},
\]
and, for $1\le i<s$, define
\[
A_i(w):=\sum_{j=1}^k
\bigl(\mathbf1_{c_j=i}-\mathbf1_{c_{j-1}=i}\bigr)R_j(w).
\]
Here $c_0=c_k$. By cyclicity of the trace, for unitaries $\mathbf C=(C_1,\ldots,C_s)$ in a tracial von Neumann algebra $(M,\tau)$ and self-adjoint $h_1,\ldots,h_{s-1}\in M$,
with $h_s=0$, we have
\[
\frac d{dt}\bigg|_{t=0}\tau\bigl(w\bigl((e^{ith_i}C_ie^{-ith_i})_{i=1}^s\bigr)\bigr)
=i\sum_{i<s}\tau\bigl(h_iA_i(w)(\mathbf C)\bigr).
\]
Set $A(w):=(A_i(w))_{i<s}\in\C[\mathbb F_s]^{s-1}$.
Define the self-adjoint gradient tuples
\[
\gamma_w^{\Re}:=\frac i2\bigl(A(w)-A(w)^*\bigr),
\qquad
\gamma_w^{\Im}:=\frac12\bigl(A(w)+A(w)^*\bigr).
\]
Since the $h_i$ are self-adjoint, the inner products of $(h_i)_{i<s}$ with $\gamma_w^{\Re}(\mathbf C)$ and $\gamma_w^{\Im}(\mathbf C)$ in $\bigoplus_{i<s}L^2(M,\tau)$ are the real and imaginary parts, respectively, of the derivative above.
Let $d_D:=2\abs{\mathcal W_D}$ be the number of the gradient tuples $\gamma_w^{\Re}$ and $\gamma_w^{\Im}$.
For an element of the group $^*$-algebra $\C[\mathbb F_s]$, its coefficient $\ell^1$-norm is the sum of the absolute values of its coefficients, and its degree is the largest reduced word length in its support.
For tuples, we use the direct-sum $L^2$-norm.

\begin{lemma}
\label{lem:group-word-gradient-gram}
For integers $s,D\ge2$, the number $d_D=2\abs{\mathcal W_D}$ satisfies
\[
1\le d_D\le4(2s)^D.
\]
In $\bigoplus_{i<s}L^2(L(\mathbb F_s),\tau_{\mathbb F_s})$, the real Gram matrix $G_D^{\mathrm{fr}}$ of the $d_D$ self-adjoint gradient tuples is
diagonal and satisfies
\[
G_D^{\mathrm{fr}}\ge I_{d_D}.
\]
Every gradient has degree at most $D$ and total coefficient $\ell^1$-norm, summed over its components, at most $2D$.
Consequently, if the gradients are ordered as $(\gamma_\alpha)_{\alpha=1}^{d_D}$, every polynomial
\[
\sum_{i<s}\gamma_{\alpha,i}\gamma_{\beta,i}
\]
has degree at most $2D$ and coefficient $\ell^1$-norm at most $4D^2$.
\end{lemma}

\begin{proof}
Let $q=2s$. There are
$\sum_{\ell=1}^Dq(q-1)^{\ell-1}\le2q^D$ nontrivial reduced words
of length at most $D$. Since $\mathcal W_D\ne\varnothing$, this gives
$1\le d_D\le4q^D$.

No nontrivial cyclically reduced word is a cyclic rotation of its
inverse. Indeed, such a rotation would give a reflection of the cyclic order of the letters that sends each letter to its inverse. This reflection would either fix a letter, which would then equal its inverse, or exchange two adjacent letters, which would then cancel.
Thus $A(w)$ and $A(w)^*$ have disjoint supports.

Set $v_j=(\mathbf1_{c_j=i}-\mathbf1_{c_{j-1}=i})_{i<s}$, so
$1\le\norm{v_j}_{\ell_2}^2\le2$ and $\norm{v_j}_{\ell^1}\le2$.
Equal rotations have the same first and last syllables,
hence the same $v_j$. There are at least two distinct rotations.
Writing $r_R$ for the number of $j$ with $R_j(w)=R$ and $v_R$ for the corresponding $v_j$, we obtain
\[
\norm{A(w)}_2^2=\sum_R r_R^2\norm{v_R}_{\ell_2}^2\ge2.
\]
Since $A(w)$ and $A(w)^*$ have disjoint supports, the tuples $\gamma_w^{\Re}$ and $\gamma_w^{\Im}$ are orthogonal and each has squared norm
$\norm{A(w)}_2^2/2\ge1$. Different representatives have disjoint
rotation-and-inversion orbits, so their gradients are also orthogonal.
Hence $G_D^{\mathrm{fr}}$ is diagonal and $G_D^{\mathrm{fr}}\ge I_{d_D}$.

Each gradient has degree at most $D$ and total coefficient
$\ell^1$-norm at most $\sum_j\norm{v_j}_{\ell^1}\le2k\le2D$.
The coefficient $\ell^1$-norm is submultiplicative, and the reduced length of a product of two words is at most the sum of their lengths. Hence each $\sum_{i<s}\gamma_{\alpha,i}\gamma_{\beta,i}$ has coefficient $\ell^1$-norm at most $(2D)^2$ and degree at most $2D$.
\end{proof}

For a tuple $\mathbf h=(h_1,\ldots,h_{s-1})$ of self-adjoint elements of a tracial von Neumann algebra, we use
\[
\norm{\mathbf h}_2:=\left(\sum_{i<s}\norm{h_i}_2^2\right)^{1/2},
\qquad
\norm{\mathbf h}_{\max,\infty}:=\max_{i<s}\norm{h_i}_\infty,
\qquad
h_s:=0.
\]

\begin{proposition}
\label{prop:group-word-quantitative-exactification}
Fix integers $s,D\ge2$, and let $d=d_D=2\abs{\mathcal W_D}$.  Let $(M,\tau)$ be a tracial von Neumann
algebra, and let $\mathbf C=(C_1,\ldots,C_s)\in U(M)^s$ satisfy
\[
\tau(C_i^m)=0
\qquad(1\le\abs m\le D,\ 1\le i\le s).
\]
Suppose that, for some $E\ge0$,
\[
\max_{\substack{g\in\mathbb F_s\setminus\{e\}\\\abs g\le2D}}
\abs{\tau(g(\mathbf C))}
\le E\le\frac1{64D^2d}.
\]
Then there are self-adjoint $h_1,\ldots,h_s\in M$, with $h_s=0$, such that, for
\[
u_i:=e^{ih_i},
\qquad
V_i:=u_iC_iu_i^*, \] we have
\[
\tau(g(\mathbf V))=\delta_{g,e}
\qquad(g\in\mathbb F_s,\ \abs g\le D)
\]
and
\[
\norm{\mathbf h}_2\le2\sqrt2\sqrt d\,E,
\qquad
\norm{\mathbf h}_{\max,\infty}\le6DdE.
\]
The same respective bounds hold for
\[
\left(\sum_{i<s}\norm{u_i-1}_2^2\right)^{1/2},
\qquad
\max_{i<s}\norm{u_i-1}_\infty.
\]
Moreover, $V_s=C_s$, and each $V_i$ has the same spectral
distribution as $C_i$.
\end{proposition}

\begin{proof}
Set $C(\mathbf h)_i=e^{ih_i}C_i e^{-ih_i}$, with $h_s=0$, and let $\Phi(\mathbf h)$ consist of the real and imaginary parts
of $\tau(w(C(\mathbf h)))$ for $w\in\mathcal W_D$.
The derivative of $e^{ith_i}C_i^me^{-ith_i}$ at $t=0$ is $i[h_i,C_i^m]$, with no factor of $m$.
Hence, by cyclicity of the trace, we have
\[
\frac d{dt}\bigg|_{t=0}\tau(w(C(t\mathbf h)))
=i\sum_{i<s}\tau(h_iA_i(w)(\mathbf C)),\qquad
D\Phi(0)\mathbf h=\left(\sum_{i<s}\tau(h_i\gamma_{\alpha,i})\right)_\alpha.
\]
Here the gradients of Lemma~\ref{lem:group-word-gradient-gram} are
evaluated at $\mathbf C$. Let $G$ be their real Gram matrix. By the last assertion of Lemma~\ref{lem:group-word-gradient-gram} and the hypothesis on $E$, every entry differs from the corresponding entry of $G_D^{\mathrm{fr}}$ by at most $4D^2E$. Thus $\norm{G-G_D^{\mathrm{fr}}}\le4D^2dE\le1/16$ and $G\ge I/2$.
Let $H$ be the real span of the gradients and define $L:\R^d\to H$ by
\[
L(y)=\sum_\alpha(G^{-1}y)_\alpha\gamma_\alpha.
\]
Then $D\Phi(0)L=I$, $LD\Phi(0)|_H=I$, and
\[
\norm{Ly}_2\le\sqrt2\norm y_{\ell_2},\qquad
\norm{Ly}_{\max,\infty}\le4D\sqrt d\,\norm y_{\ell_2}.
\]
For the second bound, we use that $\norm{G^{-1}y}_{\ell^1}\le\sqrt d\,\norm{G^{-1}y}_{\ell_2}\le2\sqrt d\,\norm y_{\ell_2}$ and that each gradient has total coefficient $\ell^1$-norm at most $2D$ and is evaluated at unitaries.

Arguing as for the remainder in the proof of Proposition~\ref{prop:pair-linearization}, now with the direct-sum $L^2$-norm, we obtain
\[
\norm{R(\mathbf h)-R(\mathbf h')}_{\ell_2}
\le C_{\rm nl}(\norm{\mathbf h}_2+\norm{\mathbf h'}_2)\norm{\mathbf h-\mathbf h'}_2,
\qquad C_{\rm nl}=4D^2\sqrt d,
\]
where $R(\mathbf h)=\Phi(\mathbf h)-\Phi(0)-D\Phi(0)\mathbf h$. Indeed, each of the at most $D$ rotations of $w$ has at most $D$ conjugated syllables, and its coefficient vector $(\mathbf1_{c_j=i}-\mathbf1_{c_{j-1}=i})_{i<s}$ has $\ell_2$-norm at most two.
This bound holds for all tuples $\mathbf h,\mathbf h'$ of self-adjoint elements of $M$. Since
$\eta:=\norm{\Phi(0)}_{\ell_2}\le\sqrt dE\le1/(16C_{\rm nl})$,
we may apply Lemma~\ref{lem:finite-dimensional-correction} to the
restriction of $\Phi$ to the finite-dimensional space $H$, with
$T=D\Phi(0)|_H$, $J=-L$, and the constants $\sqrt2$, $4D\sqrt d$ and $C_{\rm nl}$ in the roles of $A$, $B$ and $C$. We obtain $\mathbf h\in H$
with $\Phi(\mathbf h)=0$ and
\[
\norm{\mathbf h}_2\le2\sqrt2\sqrt dE,\qquad
\norm{\mathbf h}_{\max,\infty}\le6DdE.
\]

After cyclic reduction, a nontrivial reduced word of length at most $D$ either is a nonzero
power of one generator, whose trace remains zero under conjugation, or
belongs up to rotation and inversion to $\mathcal W_D$.
Since $\mathbf V=C(\mathbf h)$ and $\Phi(\mathbf h)=0$, we have $\tau(g(\mathbf V))=0$ for all $g\in\mathbb F_s\setminus\{e\}$ with $\abs g\le D$. Finally, the unitaries
$u_i=e^{ih_i}$, with $u_s=1$, satisfy the asserted bounds because
$\norm{e^{ih_i}-1}_2\le\norm{h_i}_2$ and
$\norm{e^{ih_i}-1}_\infty\le\norm{h_i}_\infty$.
\end{proof}

\subsection{The clock construction}

Collins and Male~\cite[Lemma~3.3 and Section~3.4]{collins-male}
approximated Haar unitaries in operator norm by clocks in the same
eigenbases and obtained strongly convergent clock models. Because the approximation error tends
to zero in operator norm, the average of the conjugations by the clocks and their adjoints has the same limiting norm on the traceless matrices as Hastings found for the Haar unitaries~\cite[p.~1]{hastings-quantum-expanders}.
To apply Proposition~\ref{prop:group-word-quantitative-exactification} with $D$ of order $\log N$, we need quantitative bounds for the operator-norm distance between each Haar unitary and its clock and for the traces of the mixed words in the clocks. We obtain them from
a moment estimate for products of traces of words in independent Haar
unitaries, with explicit dependence on their total length (Lemma~\ref{lem:haar-growing-word-moments}). We also use this estimate in Section~\ref{subsec:haar-word-fluctuations}.

\begin{lemma}
\label{lem:haar-growing-word-moments}
Let $s\ge1$ be an integer and let $U_1,\ldots,U_s$ be independent Haar unitaries in $U(N)$.
If $K\ge1$ is an integer and $N\ge K^{3/2}$, then every ordered collection $(v_1,\ldots,v_j)$ of nontrivial cyclically
reduced words of total length at most $K$ satisfies
\[
\abs{\mathbb E\prod_{\ell=1}^j\operatorname{Tr}v_\ell(\mathbf U)}\le K^{j/2},
\]
with the empty product equal to $1$.  In particular, for every nontrivial $w\in\mathbb F_s$ and every even integer $p\ge2$,
\[
\norm{\operatorname{Tr}w(\mathbf U)}_{L^p(\mathbb P)}\le\sqrt{p\abs w}
\qquad\text{if }N\ge(p\abs w)^{3/2}.
\]
Consequently, if $L\ge1$ and $N\ge16L^2$, every nontrivial $w\in\mathbb F_s$ with $\abs w\le L$ satisfies
\[
\mathbb E\abs{\operatorname{Tr}w(\mathbf U)}^2\le4L,
\qquad
\mathbb E\abs{\operatorname{tr}_Nw(\mathbf U)}^2\le\frac{4L}{N^2}.
\]
Then $\abs{\mathbb E\operatorname{tr}_Nw(\mathbf U)}\le2\sqrt L/N$, and the variance of $\operatorname{tr}_Nw(\mathbf U)$ is at most $4L/N^2$.
\end{lemma}

\begin{proof}
Let $\mathcal C_K$ be the finite set of ordered collections
$\mathcal V=(v_1,\ldots,v_j)$ of nontrivial cyclically reduced words of total length at most $K$, including the
empty collection.  Set
\[
M(\mathcal V):=\mathbb E\prod_{\ell=1}^j\operatorname{Tr}v_\ell(\mathbf U),
\qquad
\beta:=\max_{\mathcal V\in\mathcal C_K}K^{-j/2}\abs{M(\mathcal V)}.
\]
The empty collection gives $\beta\ge1$.

We use the identity obtained by differentiating under Haar
measure~\cite[equations~(15)--(20)]{hastings-quantum-expanders}.
For a nonempty collection, write $v_1=xX$. If $x$ is an inverse generator,
we replace the corresponding Haar matrix by its adjoint and
interchange that generator with its inverse in every word.
This preserves the product Haar law and
permutes $\mathcal C_K$, so we may assume that $x$ is a generator
without changing $\beta$.
For a real orthonormal basis $(T_\alpha)$ of self-adjoint matrices in the unnormalized Hilbert--Schmidt inner product,
we have
\[
\begin{gathered}
\sum_\alpha T_\alpha^2=N,\\
\sum_\alpha\operatorname{Tr}(T_\alpha A T_\alpha B)
=\operatorname{Tr}(A)\operatorname{Tr}(B),\\
\sum_\alpha\operatorname{Tr}(T_\alpha A)\operatorname{Tr}(T_\alpha B)=\operatorname{Tr}(AB)
\end{gathered}
\]
for arbitrary complex matrices $A,B$. Now we differentiate
\[
\mathbb E\left[\operatorname{Tr}(T_\alpha v_1(\mathbf U))
\prod_{\ell>1}\operatorname{Tr}v_\ell(\mathbf U)\right]
\]
under $x\mapsto e^{itT_\alpha}x$. The expectation does not depend on $t$, by Haar invariance, so its derivative at $t=0$ is zero. At $t=0$, each occurrence of $x$ has derivative $iT_\alpha x$, and each occurrence of $x^{-1}$ has derivative $-ix^{-1}T_\alpha$. After summing over $\alpha$, the first letter contributes $iN M(\mathcal V)$.
Solving for $M(\mathcal V)$, we obtain an identity expressing $M(\mathcal V)$ as the sum of the following terms, one for each remaining occurrence of $x$ or $x^{-1}$. In the last two rows, we rotate $v_\ell$ to start at that occurrence. All traces not displayed
remain as factors:
\[
\begin{array}{c|c|c}
\text{participating words}&\text{replacement traces}&\text{coefficient}\\ \hline
v_1=xA\,xB&\operatorname{Tr}(xA)\operatorname{Tr}(xB)&-N^{-1}\\
v_1=xA\,x^{-1}B&\operatorname{Tr}(A)\operatorname{Tr}(B)&N^{-1}\\
v_1=xX,\ v_\ell=xY&\operatorname{Tr}(xXxY)&-N^{-1}\\
v_1=xX,\ v_\ell=x^{-1}Y&\operatorname{Tr}(XY)&N^{-1}
\end{array}
\]
There are at most $K-1$ terms, counted with multiplicity. The first two rows are the splits of $v_1$, at an occurrence of $x$ or of $x^{-1}$. The last two are the joins of $v_1$ with another word $v_\ell$, at an occurrence of $x$ (a same-sign join) or of $x^{-1}$ (an opposite-sign join).

We cyclically reduce every resulting word. Both words in a split are nontrivial. For a split at $x^{-1}$, the words $A$ and $B$ are
nonempty reduced segments, since an empty segment would force an
ordinary or cyclic cancellation in $v_1$. For a split at $x$, the words $xA$
and $xB$ are nonempty reduced segments.  A same-sign join cannot be the identity either, since
otherwise $xY=(xX)^{-1}$ would force the last letter of the cyclically reduced word $xX$ to be
$x^{-1}$.
Only a join of opposite signs can produce the identity. When this occurs, the trace is $\operatorname{Tr}(1)=N$,
and both participating traces disappear from the remaining collection. Every resulting collection lies in $\mathcal C_K$:
the total length is unchanged or decreases by two before reduction, and cannot increase afterwards.

After multiplying by $K^{-j/2}$, a split, a nontrivial join, and an identity join have coefficients of absolute value
\[
\frac{\sqrt K}N,\qquad \frac1{\sqrt K\,N},\qquad \frac1K,
\]
respectively, when each term is written as a multiple of $K^{-j'/2}M(\mathcal V')$, where $\mathcal V'$ is the resulting collection and $j'$ is the number of its words. This includes an identity join with $j=2$,
which leaves the empty collection. Since $N\ge K^{3/2}$, each coefficient has absolute value at most $1/K$.
If $\beta>1$, a nonempty collection attains the maximum, and the identity above for this collection gives
\[
\beta\le\frac{K-1}{K}\beta<\beta,
\]
a contradiction.  This also covers $K=1$, when the identity above has no remaining terms.  Hence $\beta=1$,
which proves the first estimate.

We cyclically reduce $w$. This leaves a nontrivial word and does not change $\operatorname{Tr}w(\mathbf U)$.  For even $p$, the product of
$p/2$ copies of this trace and $p/2$ copies of its conjugate has $p$ nontrivial traces and total length at most
$K=p\abs w$.  By the first estimate,
$\mathbb E\abs{\operatorname{Tr}w(\mathbf U)}^p\le(p\abs w)^{p/2}$.
For $\abs w\le L$ and $N\ge16L^2$, we have $N\ge(2\abs w)^{3/2}$, so the case $p=2$ gives
$\mathbb E\abs{\operatorname{Tr}w(\mathbf U)}^2\le2\abs w\le4L$.
Dividing by $N^2$ gives the bound for $\operatorname{tr}_N$. The remaining assertions follow from Cauchy--Schwarz and the fact that the variance is at most the second moment.
\end{proof}

\begin{lemma}
\label{lem:fourier-clock-rounding}
Let $U\in U(N)$ and let $H\ge1$ be an integer. There is a clock $C$,
diagonal in an eigenbasis of $U$, such that
\[
\norm{C-U}_\infty
\le2\pi\left(\frac6{H+1}
+\frac4\pi\sum_{m=1}^H\frac{\abs{\operatorname{tr}_N(U^m)}}m+\frac1N\right).
\]
\end{lemma}

\begin{proof}
We write the eigenvalues of $U$, with multiplicity, as $e^{2\pi i\theta_j}$,
where $0\le\theta_1\le\cdots\le\theta_N<1$. Since
$N^{-1}\sum_je^{2\pi im\theta_j}=\operatorname{tr}_N(U^m)$, the
Erd\H{o}s--Tur\'an
inequality~\cite[Chapter~2, Theorem~2.5]{kuipers-niederreiter} bounds
the discrepancy
\[
\Delta=\sup_{0\le\alpha<\beta\le1}
\abs{N^{-1}\#\{j:\alpha\le\theta_j<\beta\}-(\beta-\alpha)}
\]
by $6/(H+1)+(4/\pi)\sum_{m=1}^H\abs{\operatorname{tr}_N(U^m)}/m$. From
the definition of $\Delta$, we obtain $j/N-\Delta\le\theta_j\le(j-1)/N+\Delta$.
In a corresponding eigenbasis, replacing $e^{2\pi i\theta_j}$ by $e^{2\pi
ij/N}$ produces a clock $C$ satisfying
$\norm{C-U}_\infty\le2\pi(\Delta+1/N)$.
\end{proof}

Note that if unitary matrices $U_N\in U(N)$ satisfy $\operatorname{tr}_N(U_N^m)\to0$ for every fixed integer $m\ne0$, then the clocks $C_N$ constructed in the proof of Lemma~\ref{lem:fourier-clock-rounding} satisfy $\norm{C_N-U_N}_\infty\to0$. Indeed, these clocks do not depend on $H$, so we may let $N\to\infty$ in the bound of the lemma for each fixed $H$ and then let $H\to\infty$.

\begin{corollary}
\label{cor:quantitative-haar-clock-rounding}
For a Haar unitary $U\in U(N)$ with $N\ge256$, set
\[
H_N:=\left\lfloor\frac{\sqrt N}{4}\right\rfloor,\qquad
\mathcal A_N^{\mathrm{round}}
:=\left\{\sum_{m=1}^{H_N}\frac{\abs{\operatorname{tr}_N(U^m)}}m\le N^{-1/4}\right\}.
\]
Then $\mathbb P((\mathcal A_N^{\mathrm{round}})^c)\le16N^{-2}$.
On $\mathcal A_N^{\mathrm{round}}$, there is a clock $C$, diagonal in an eigenbasis of $U$, such that
\[
\norm{C-U}_\infty\le100N^{-1/4}.
\]
\end{corollary}

\begin{proof}
For $m\le H_N$, we have $N\ge(4m)^{3/2}$. Hence, by
Lemma~\ref{lem:haar-growing-word-moments} with $p=4$ and Minkowski's inequality,
\[
\left\|\sum_{m=1}^{H_N}\frac{\abs{\operatorname{tr}_N(U^m)}}m\right\|_{L^4(\mathbb P)}
\le\frac2N\sum_{m=1}^{H_N}m^{-1/2}
\le2N^{-3/4}.
\]
Markov's inequality gives
$\mathbb P((\mathcal A_N^{\mathrm{round}})^c)\le16N^{-2}$.
On $\mathcal A_N^{\mathrm{round}}$, we apply Lemma~\ref{lem:fourier-clock-rounding} with $H=H_N$. Since $H_N+1>\sqrt N/4$ and $N\ge256$, we obtain
\[
\norm{C-U}_\infty
\le2\pi\left(24N^{-1/2}+\frac4\pi N^{-1/4}+\frac1N\right)
\le2\pi\left(6+\frac4\pi+\frac1{64}\right)N^{-1/4}<100N^{-1/4}.
\qedhere
\]
\end{proof}

\begin{proposition}
\label{prop:quantitative-haar-clock-exactification}
Fix an integer $s\ge2$, and set
\[
D_N:=\left\lfloor\frac{\log N}{16\log(2s)}\right\rfloor,
\qquad
N_s^0:=\left\lceil\exp\bigl(\max\{256,32\log(2s)\}\bigr)\right\rceil.
\]
For every $N\ge N_s^0$ and independent Haar unitaries $U_1,\ldots,U_s\in U(N)$, there is an event
$\mathcal A_N$ such that
\[
\mathbb P(\mathcal A_N^c)
\le128D_N^2(2s)^{2D_N}N^{-2}+16sN^{-2}
\le128D_N^2N^{-15/8}+16sN^{-2}<\frac14.
\]
For every tuple in $\mathcal A_N$, there are clocks $C_1,\ldots,C_s$ and self-adjoint $h_1,\ldots,h_s\in M_N(\C)$, with $h_s=0$,
such that, for $u_i=e^{ih_i}$ and $V_i=u_iC_iu_i^*$,
\[
\operatorname{tr}_N(g(\mathbf V))=\delta_{g,e}
\qquad(g\in\mathbb F_s,\ \abs g\le D_N),
\]
\[
\max_i\norm{C_i-U_i}_\infty\le100N^{-1/4},
\qquad
\max_i\norm{h_i}_\infty\le5000D_N^2N^{-3/16},
\]
and
\[
\left(\sum_{i=1}^s\norm{h_i}_2^2\right)^{1/2}
\le804\sqrt2\,D_NN^{-7/32}.
\]
The same two bounds hold with $h_i$ replaced by $u_i-1$.
\end{proposition}

\begin{proof}
Write $q=2s$ and $D=D_N$, which is the radius $R_N$ of Theorem~\ref{thm:irreducible-clock-strong-models}. For $N\ge N_s^0$, we have
\[
\begin{gathered}
\log N\ge\max\{256,32\log q\},\qquad 2\le D\le\log N,\\
q^D\le N^{1/16},\qquad N\ge64D^2,\qquad N>2D.
\end{gathered}
\]
There are at most $2q^{2D}$ nonidentity reduced words of length at most
$2D$. Since $N\ge64D^2\ge(8D)^{3/2}$, by Lemma~\ref{lem:haar-growing-word-moments} with $p=4$, we have $\mathbb E\abs{\operatorname{tr}_N(g(\mathbf U))}^4\le64D^2/N^4$ for each such word. Let $\mathcal A_N$ be the intersection of the event
that all these traces have modulus at most $N^{-1/2}$ with the $s$
events defined in Corollary~\ref{cor:quantitative-haar-clock-rounding},
one for each $U_i$. By Markov's inequality and the bound in that corollary, we have
\[
\mathbb P(\mathcal A_N^c)
\le128D^2q^{2D}N^{-2}+16sN^{-2}
\le128D^2N^{-15/8}+16sN^{-2}<\frac14.
\]
On $\mathcal A_N$, let $C_i$ be the clock given by Corollary~\ref{cor:quantitative-haar-clock-rounding} for $U_i$. By telescoping each word, we see that the normalized traces of the nontrivial words of length at most $2D$ in the $C_i$ have modulus at most
\[
N^{-1/2}+2D\max_i\norm{C_i-U_i}_\infty\le201DN^{-1/4}=:\kappa_N.
\]
Since $D<N$, the clocks satisfy
$\operatorname{tr}_N(C_i^m)=0$ for $1\le |m|\le D$ and $1\le i\le s$.
Since $d=d_D\le4q^D$, the condition $E\le1/(64D^2d)$ of Proposition~\ref{prop:group-word-quantitative-exactification} with $E=\kappa_N$ follows from
\[
64D^2d\kappa_N\le51456D^3N^{-3/16}<1
\]
for $N\ge N_s^0$. Hence, by
Proposition~\ref{prop:group-word-quantitative-exactification}, we obtain $\operatorname{tr}_N(g(\mathbf V))=\delta_{g,e}$ for $\abs g\le D$, and
\[
\max_i\norm{h_i}_\infty\le6Dd\kappa_N\le5000D^2N^{-3/16},
\qquad
\left(\sum_i\norm{h_i}_2^2\right)^{1/2}
\le2\sqrt2\sqrt d\,\kappa_N\le804\sqrt2\,DN^{-7/32}.
\]
Since $\norm{e^{ih_i}-1}_2\le\norm{h_i}_2$ and $\norm{e^{ih_i}-1}_\infty\le\norm{h_i}_\infty$, the same bounds hold for $u_i-1$.
\end{proof}

\begin{proof}[Proof of Theorem~\ref{thm:irreducible-clock-strong-models}]
\emph{The almost-sure event.}
The bounds on $\mathbb P(\mathcal A_N^c)$ in Proposition~\ref{prop:quantitative-haar-clock-exactification} are summable in $N$. By the Borel--Cantelli lemma, almost surely, the Haar tuple of dimension $N$ lies in $\mathcal A_N$ for every sufficiently large $N$, for any coupling of the Haar tuples of different dimensions.
We intersect this event with the almost-sure strong-convergence event of
Collins and Male~\cite[Theorem~1.4]{collins-male}, in the form recalled in Subsection~\ref{sec:stable-strong-models}, and the almost-sure
event on which the norm on $L^2_0(M_N)$ of the average of the conjugations by the Haar unitaries and their adjoints converges to $\lambda_s$. Hastings's lower bound for this norm~\cite[equations~(10)--(12)]{hastings-quantum-expanders} holds for every tuple of unitaries. For the upper bound, we use the form of Hastings's theorem given by Hastings and Harrow~\cite[Theorem~6, with $k=1$]{hastings-harrow}: for every $c>1$ and all sufficiently large $N$, this norm exceeds $c\bigl(1+O(N^{-1/6}\log N)\bigr)\lambda_s$ with probability at most $c^{-N^{1/6}/4}$. For each $c$, these probabilities are summable in $N$, so, by the Borel--Cantelli lemma and the lower bound, this norm converges to $\lambda_s$ almost surely, as Farah, Jekel, and Pi observed for the second largest eigenvalue~\cite[Theorem~4.6]{farah-jekel-pi}.

\emph{The corrected clocks.}
Choose a realization $(\mathbf U_N)_N$ of the Haar tuples in the
probability-one event just established, and set
\[
S_N=\frac1{2s}\sum_{i=1}^s
(\operatorname{Ad}U_{i,N}+\operatorname{Ad}U_{i,N}^*).
\]
Then the tuples converge strongly to free Haar unitaries and
$\|S_N|_{L^2_0(M_N)}\|_{2\to2}\to\lambda_s$.
Since $\mathcal A_N$ is nonempty for every $N\ge N_s^0$, we may replace the
finitely many tuples outside $\mathcal A_N$ by arbitrary elements of
$\mathcal A_N$ without changing either limit.

We now apply Proposition~\ref{prop:quantitative-haar-clock-exactification}
to these tuples, obtaining clocks $C_{i,N}$ and correcting unitaries
$u_{i,N}$ with $u_{s,N}=1$. Set
\[
\delta_N=\max_i\|C_{i,N}-U_{i,N}\|_\infty\le100N^{-1/4},
\qquad
\epsilon_N=\max_i\|u_{i,N}-1\|_\infty
\le5000D_N^2N^{-3/16}.
\]
For $V_{i,N}=u_{i,N}C_{i,N}u_{i,N}^*$, we have
$\max_i\|V_{i,N}-U_{i,N}\|_\infty\le\delta_N+2\epsilon_N\to0$, and
\[
\|f(\mathbf V_N)-f(\mathbf U_N)\|_\infty
\le(\delta_N+2\epsilon_N)\sum_g|c_g||g|
\quad\left(f=\sum_gc_gg\right).
\]
Define $\pi_N:\mathbb F_s\to U(N)$ by $\pi_N(a_i)=V_{i,N}$.
By the preceding estimate, $\mathbf V_N$ has the same strong limit as
$\mathbf U_N$.
Each $V_{i,N}$ is a clock, and by that proposition, we have $\operatorname{tr}_N(\pi_N(g))=\delta_{g,e}$ for
$\abs g\le R_N=D_N$.

\emph{The spectral gap.}
Let $S_N^{\rm clk}$ be the averaging operator defined as $S_N$ above,
with $U_{i,N}$ replaced by $C_{i,N}$. The inequality
$\|\operatorname{Ad}U-\operatorname{Ad}V\|_{2\to2}\le2\|U-V\|_\infty$
gives
\[
\|S_N^{\rm clk}-S_N\|_{2\to2}\le2\delta_N\to0,
\eqtag\label{eq:clock-channel-rounding}
\]
\[
\|T_N-S_N^{\rm clk}\|_{2\to2}\le4\epsilon_N\to0.
\eqtag\label{eq:clock-channel-exactification}
\]
Since all three operators preserve $L^2_0(M_N)$, it follows that $\rho_N\to\lambda_s$.

\emph{Irreducibility.}
For traceless $x$,
\[
\frac1s\sum_i\|[V_{i,N},x]\|_2^2
=2\langle(1-T_N)x,x\rangle
\ge2(1-\rho_N)\|x\|_2^2.
\]
Eventually $\rho_N<1$, so the commutant of $\pi_N(\mathbb F_s)$ is scalar. By Burnside's theorem, $C^*(\pi_N(\mathbb F_s))=M_N(\C)$. Hence the GNS von Neumann algebra $\pi_N(\mathbb F_s)''=M_N(\C)$ of the character $\operatorname{tr}_N\circ\pi_N$ is a factor, so the character is extreme.
Finally, let $\mathbf W_N=(W_{1,N},\ldots,W_{s,N})\in\operatorname{Red}_{s,N}$, choose a common reducing
projection $0<p<1$, and set
$t=\operatorname{tr}_N(p)$ and $\eta_N=\max_i\|V_{i,N}-W_{i,N}\|_\infty$.
Each off-diagonal corner of $[V_{i,N},p]=[V_{i,N}-W_{i,N},p]$
has operator norm at most $\eta_N$ and rank at most
$N\min\{t,1-t\}$, hence squared $L^2$-norm at most
$\eta_N^2\min\{t,1-t\}$. Applying \eqref{eq:clock-poincare} to $p-t$
gives
\[
2\eta_N^2\min\{t,1-t\}\ge2(1-\rho_N)t(1-t).
\]
Thus $\eta_N^2\ge(1-\rho_N)\max\{t,1-t\}\ge(1-\rho_N)/2$,
which proves \eqref{eq:clock-distance-reducible}.
\end{proof}

Since $V_{i,N}^N=1$ for every $i$, the representations factor through
the free product of $s$ cyclic groups of order $N$. Because $R_N<N$, every syllable of a nontrivial word of length at most $R_N$ has a nonzero exponent modulo $N$, so such a word has nontrivial image in this free product. Thus, by Theorem~\ref{thm:irreducible-clock-strong-models}\textup{(ii)}, on the images of the words of length at most $R_N$, the normalized character of the induced representation of this free product equals the free product of the canonical traces of the $s$ cyclic groups.

By Hastings's lower bound~\cite[equations~(10)--(12)]{hastings-quantum-expanders}, $\lambda_s$ is asymptotically optimal among the norms on $L^2_0(M_N)$ of the uniform averages of the conjugations by $s$ unitaries and their inverses.

We now use the limit $\rho_N\to\lambda_s<1$ to identify the relative commutants of the embeddings of $L(\mathbb F_s)$ into ultraproducts induced by these representations. For a separable ${\rm II}_1$ factor, Brown~\cite{brown-topological-dynamics} defined a convex-like structure on the set of its unital $^*$-homomorphisms into $R^\cU$ modulo unitary conjugacy; we call this set Brown's convex space.

\begin{corollary}
\label{cor:ramanujan-clock-ultraproduct-commutants}
Let $\cU$ be any nonprincipal ultrafilter on $\N$.  The representations in
Theorem~\ref{thm:irreducible-clock-strong-models} induce a normal unital trace-preserving embedding of von Neumann algebras
\[
\Theta:L(\mathbb F_s)\longrightarrow
\prod_{\cU}\bigl(M_N(\C),\operatorname{tr}_N\bigr)
\]
such that
\[
\Theta(L(\mathbb F_s))'
\cap\prod_{\cU}\bigl(M_N(\C),\operatorname{tr}_N\bigr)
=\C1.
\]

Let $R$ be the hyperfinite ${\rm II}_1$ factor. For each $N$, choose a trace-preserving isomorphism
$\Phi_N:M_N(\C)\mathbin{\overline\otimes} R\to R$ and set $B_N:=\Phi_N(1\otimes R)$, a unital von Neumann subalgebra of $R$.  The assignments
\[
\lambda_{\mathbb F_s}(a_i)\longmapsto
\bigl(\Phi_N(V_{i,N}\otimes1)\bigr)_{\cU}
\]
induce a normal unital trace-preserving embedding of von Neumann algebras $\widetilde\Theta:L(\mathbb F_s)\to R^\cU$ satisfying
\[
\widetilde\Theta(L(\mathbb F_s))'\cap R^\cU
=\prod_{\cU}B_N
\cong R^\cU.
\]
In particular, the relative commutant is a ${\rm II}_1$ factor and
the class of $\widetilde\Theta$ is an extreme point of Brown's convex space
of trace-preserving embeddings $L(\mathbb F_s)\to R^\cU$
modulo unitary conjugacy.
\end{corollary}

\begin{proof}
By the trace limits in Theorem~\ref{thm:irreducible-clock-strong-models}\textup{(iii)}, we obtain embeddings of the reduced group $C^*$-algebra that send $\lambda_{\mathbb F_s}(a_i)$ to $(V_{i,N})_\cU$ and to $\bigl(\Phi_N(V_{i,N}\otimes1)\bigr)_\cU$, respectively. These extend to normal trace-preserving
embeddings of $L(\mathbb F_s)$. Let $(x_N)_\cU$ belong to the first relative
commutant, with $(x_N)$ uniformly bounded. Since
$\rho_N\to\lambda_s<1$, it follows from \eqref{eq:clock-poincare} that
$\|x_N-\operatorname{tr}_N(x_N)\|_2\to_\cU0$.
Thus $(x_N)_\cU=\lim_\cU\operatorname{tr}_N(x_N)$, proving
that this commutant is scalar.

For the second relative commutant, set
\[
E_N=\Phi_N\circ(\operatorname{tr}_N\otimes\operatorname{id})
\circ\Phi_N^{-1}:R\to B_N.
\]
Tensoring \eqref{eq:clock-poincare} with the identity on $L^2(R)$ gives, for $x_N\in R$,
\[
2(1-\rho_N)\|x_N-E_N(x_N)\|_2^2
\le\frac1s\sum_i
\|[\Phi_N(V_{i,N}\otimes1),x_N]\|_2^2.
\]
Thus, if $(x_N)_\cU$ belongs to the commutant and $(x_N)$ is a
bounded representative, then
$(x_N)_\cU=(E_N(x_N))_\cU$. The elements $E_N(x_N)$ are uniformly
bounded and lie in $B_N$. Conversely, $B_N$ commutes with $\Phi_N(V_{i,N}\otimes1)$ for every $i$. This proves that the commutant is
$\prod_\cU B_N\cong R^\cU$.
Extremality follows from the implication (2)~$\Rightarrow$~(1)
of~\cite[Proposition~5.2]{brown-topological-dynamics}, applied to the separable ${\rm II}_1$ factor $L(\mathbb F_s)$.
\end{proof}

To identify the relative commutants, we use only that $\operatorname{tr}_N(\pi_N(g))\to\delta_{g,e}$ for every $g\in\mathbb F_s$ and that $\rho_N\le c$ for some $c<1$ and all sufficiently large $N$. The tensor amplification argument is standard (see
Brown~\cite[Lemma~6.5]{brown-topological-dynamics} and the proof of Farah, Jekel, and Pi~\cite[Lemma~4.8]{farah-jekel-pi}).
For Haar tuples, which are asymptotically free by Voiculescu's theorem~\cite[Theorem~3.8]{voiculescu-limit-laws}, and almost surely so by a theorem of Hiai and Petz~\cite[Theorem~2.1]{hiai-petz},
the construction of Farah, Jekel, and Pi with ambient factor $R$ already gives a normal unital
trace-preserving embedding of von Neumann algebras $L(\mathbb F_s)\to R^\cU$ with factorial
relative commutant. Indeed, for three independent Haar unitaries $U_{j,N}\in U(N)$, they proved that, almost surely, the relative commutant of the elements $(U_{j,N}\otimes1)_\cU$ in $\prod_\cU(M_N(\C)\mathbin{\overline\otimes}R)$ is $\prod_\cU(1\otimes R)$~\cite[Lemma~4.8]{farah-jekel-pi}, and their proof applies to $s$ unitaries for every $s\ge2$, by their Corollary~4.7. Thus Popa's
question~\cite[Question~1.1]{goldbring-factorial-commutant}, recorded as
\cite[Problem~37]{kunnawalkam-elayavalli-open-problems}, whether every separable ${\rm II}_1$ factor that embeds into $R^\cU$ admits an embedding with factorial relative commutant, has a positive answer for
$L(\mathbb F_s)$, although Farah, Jekel, and Pi~\cite{farah-jekel-pi} did not state this consequence.

For $\Gamma=\mathrm{SL}(n,\mathbb Z)$ with odd $n\ge3$, Brown combined the tensor amplification argument with a theorem of
Bekka~\cite{bekka-superrigidity} to embed $L(\Gamma)$ into $R^\cU$ with factorial relative
commutant~\cite[Lemma~6.5 and the paragraph before Corollary~6.11]{brown-topological-dynamics}.
Popa observed, using results of Bekka~\cite{bekka-superrigidity}, that, for $\Gamma=\mathrm{PSL}(n,\mathbb
Z)$ with $n\ge3$, finite-dimensional irreducible
representations of $\Gamma$ with dimensions tending to infinity generate a copy of $L(\Gamma)$ with trivial
relative commutant in the ultraproduct of the corresponding matrix algebras~\cite[Section~1.7]{popa-independence}.
By tensor amplification, this observation also gives embeddings of $L(\Gamma)$ into $R^\cU$ with factorial relative commutant;
Goldbring~\cite[p.~5007]{goldbring-factorial-commutant} listed the case $n=3$ of this observation as a positive answer
to Popa's question.
Jekel embedded $L(\mathbb F_s)$, $2\le s<\infty$, into
$\prod_{\cU}(M_N(\C),\operatorname{tr}_N)$ with trivial relative commutant, for every
nonprincipal ultrafilter $\cU$ on $\N$~\cite[Theorem~1.2 and the discussion after it]{jekel-free-probability-model-theory}.
Assuming the Continuum Hypothesis, Chifan, Drimbe, and Ioana embedded every separable non-Gamma
${\rm II}_1$ factor, in particular $L(\mathbb F_s)$, into the ultrapower of any infinitely generic
${\rm II}_1$ factor, with factorial relative commutant~\cite[Theorem~6.2]{chifan-drimbe-ioana}.

We now show that the radius $R_N$ in Theorem~\ref{thm:irreducible-clock-strong-models} is optimal up to a constant factor.

\begin{proposition}
\label{prop:exact-character-radius-ceiling}
Fix integers $s\ge2$ and $r\ge0$, and let
$\pi:\mathbb F_s\to U(N)$ be a representation. If
\[
\operatorname{tr}_N(\pi(g))=\delta_{g,e}
\qquad(\abs g\le2r),
\]
then
\[
\abs{B_{\mathbb F_s}(r)}\le N^2
\]
and
\[
(2s-1)^r
\le1+\frac{s-1}{s}(N^2-1).
\]
Consequently, if $\operatorname{tr}_N(\pi(g))=\delta_{g,e}$ whenever $\abs g\le R$, then
\[
R\le\frac{4\log N}{\log(2s-1)}+O_s(1).
\]
\end{proposition}

\begin{proof}
For $g,h\in B_{\mathbb F_s}(r)$,
\[
\langle\pi(g),\pi(h)\rangle_2
=\operatorname{tr}_N(\pi(g^{-1}h))
=\delta_{g,h}.
\]
These matrices form an orthonormal subset of $L^2(M_N(\C),\operatorname{tr}_N)$, which has dimension $N^2$.  Since
\[
\abs{B_{\mathbb F_s}(r)}
=1+\frac{2s((2s-1)^r-1)}{2s-2},
\]
the first two displayed bounds follow.  Finally, if $\operatorname{tr}_N(\pi(g))=\delta_{g,e}$ whenever $\abs g\le R$, we apply these bounds with $r=\lfloor R/2\rfloor$ and take logarithms to obtain the last estimate.
\end{proof}

By Theorem~\ref{thm:irreducible-clock-strong-models} and this proposition, the largest $R$ for which some representation $\pi:\mathbb F_s\to U(N)$ satisfies $\operatorname{tr}_N(\pi(g))=\delta_{g,e}$ whenever $\abs g\le R$ is of order $\log N$. We do not know the optimal constant.

We now bound the effect of the perturbations in Proposition~\ref{prop:quantitative-haar-clock-exactification} on polynomials whose degree grows with $N$.

\begin{corollary}
\label{cor:clock-growing-polynomial-perturbation}
Let $\mathbf U_N$ and $\mathbf V_N$ be any original and corrected tuples supplied by
Proposition~\ref{prop:quantitative-haar-clock-exactification}, including those used for
Theorem~\ref{thm:irreducible-clock-strong-models}.  For every pair of integers $r,L\ge1$ and every
\[
P=\sum_{\abs g\le L}A_g\otimes g\in M_r(\C)\otimes\C[\mathbb F_s],
\qquad
\sum_g\norm{A_g}_\infty\le1, \] we have
\[
\norm{P(\mathbf V_N)-P(\mathbf U_N)}_\infty
\le L\bigl(100N^{-1/4}+10000D_N^2N^{-3/16}\bigr).
\]
Consequently, if $L_ND_N^2N^{-3/16}\to0$, the difference tends to
zero uniformly over all matrix sizes $r$ and all polynomials of degree
at most $L_N$ satisfying the coefficient bound above. This
applies to $L_N=D_N$ and $L_N=\lfloor N^{1/8}\rfloor$.
\end{corollary}

\begin{proof}
We apply the estimate for $\norm{f(\mathbf V_N)-f(\mathbf U_N)}_\infty$ in the proof of
Theorem~\ref{thm:irreducible-clock-strong-models} to each term, using
$\norm{A_g\otimes X}_\infty=\norm{A_g}_\infty\norm X_\infty$ and
$\sum_g\norm{A_g}_\infty|g|\le L$. This gives the stated bound independently
of the matrix size $r$. Since
$N^{-1/4}=o(D_N^2N^{-3/16})$, the convergence assertion follows from the condition $L_ND_N^2N^{-3/16}\to0$.
\end{proof}

\section{Traces of words in Haar unitaries}
\label{subsec:haar-word-fluctuations}

We now return to the original independent Haar matrices and study the size and joint distribution of the traces of short nontrivial group words. For each fixed
nontrivial cyclically reduced word $w$, the trace
$\operatorname{Tr}w(\mathbf U_N)$ converges in distribution to a nondegenerate
complex Gaussian variable (Theorem~\ref{thm:haar-word-gaussian}),
while the trace of $w$ on the corrected clocks is eventually zero
(Theorem~\ref{thm:irreducible-clock-strong-models}). Thus the correction
retains the strong limit but changes these fluctuations. For unitary tuples $\mathbf U$ and $\mathbf V$ with $\max_i\norm{U_i-V_i}_\infty\le\eta$, by telescoping, we have $\abs{\operatorname{Tr}w(\mathbf U)-\operatorname{Tr}w(\mathbf V)}\le N\abs w\,\eta$.
Because of the factor $N$, this bound need not tend to zero even when
$\eta\to0$.

For fixed words, joint Gaussian convergence, with covariance determined by
cyclic rotations and inversion, was obtained by
R\u{a}dulescu~\cite[Theorem~4.1]{radulescu-haar-traces}, for words in two independent Haar matrices,
and independently by Mingo, \'Sniady, and
Speicher~\cite[Remark~3.11, Theorem~3.13, and Corollary~3.14]{mingo-sniady-speicher}, for any number of matrices.
For powers of a single Haar matrix this goes back to Diaconis and Shahshahani~\cite{diaconis-shahshahani}.
We give bounds for collections whose size and word lengths increase with
the matrix dimension. Johansson~\cite[Theorem~2.6]{johansson-compact-groups} proved a superexponential rate of Gaussian
approximation for fixed real linear combinations of the real and imaginary parts of the power traces of a single Haar matrix, and Johansson and
Lambert~\cite{johansson-lambert} extended such bounds to the vector of these traces;
see also~\cite{dobler-stolz,fulman-stein-traces}.

Throughout this section, length counts individual generators and their inverses in $\mathbb F_s$, and $\operatorname{Tr}$ is unnormalized. As in Section~\ref{sec:clock-representations}, $a_1,\ldots,a_s$ are the standard generators of $\mathbb F_s$. For a tuple $\mathbf U=(U_1,\ldots,U_s)$ of independent Haar unitary matrices in $U(N)$, we write $T_w:=\operatorname{Tr}w(\mathbf U)$.
For a nontrivial cyclically reduced word $w$, let
\[
r(w):=\#\{0\le j<\abs w:\rho^j(w)=w\},
\]
where $\rho$ rotates the letters of $w$ cyclically by one place.  Equivalently, $r(w)$ is the largest $r$ such that $w=v^r$ for some word $v$.
For probability measures with finite $p$th moments on $\R^d$ or $\C^k$, $1\le p<\infty$, write
\[
W_p(\mu,\nu):=\inf\left\{\bigl(\mathbb E\norm{X-Y}_2^p\bigr)^{1/p}:
\mathcal L(X)=\mu,\ \mathcal L(Y)=\nu\right\},
\]
where the infimum is over couplings.

For $k$ nontrivial cyclically reduced words of length at most $L$, no two related by cyclic rotation or by cyclic rotation followed by inversion, the law of the vector of the traces $\operatorname{Tr}w(\mathbf U_N)/\sqrt{r(w)}$ is within $W_p$ distance $3pkL^2/N$ of the law of a vector of independent circular complex Gaussians of variance one, for every even integer $p\ge2$ with $N\ge(2pL)^{3/2}$ (Theorem~\ref{thm:haar-word-gaussian}). We also bound the expectation of the trace of a word (Corollary~\ref{cor:haar-word-bias}).

We study the largest absolute values of word traces, with and without division by
$\sqrt{r(w)}$. With this division, the complex trace of every fixed word has limiting
variance one. Let $s\ge2$ be fixed and let
$L=L_N\to\infty$ be integers. We take both maxima over the classes, modulo cyclic rotation and inversion, of
nontrivial cyclically reduced words of length at most $L_N$. We show that if $L_N=o(N^{2/7})$, then, with
probability tending to one, every class
attaining the normalized maximum is represented by a word that is not a proper power (Corollary~\ref{cor:haar-word-full-ball-extremes}). If
$L_N^3/N\to0$, then, with probability tending to one, every class attaining the unnormalized maximum is
represented by a power of a generator (Theorem~\ref{thm:haar-word-raw-extremes}).
In each case, after centering and scaling, the squared absolute values form a point process converging to a Poisson point process with intensity $e^{-x}\,dx$, and their maximum has limiting distribution function $x\mapsto\exp(-e^{-x})$. For the normalized values, we prove this for every collection of $k_N\to\infty$ such words of length at most $L_N$, no two related by cyclic rotation or by cyclic rotation followed by inversion, with $L_N^2(\log k_N)^{3/2}/N\to0$ (Theorem~\ref{thm:haar-word-gumbel}). If $L_N=o(N^{2/7})$, the two maxima are asymptotically independent
(Corollary~\ref{cor:haar-word-joint-extremes}).

\subsection{Gaussian approximation}

\begin{theorem}[Gaussian approximation for word traces]
\label{thm:haar-word-gaussian}
Let $s,k,L\ge1$ be integers, let $U_1,\ldots,U_s$ be independent Haar matrices in $U(N)$, and let
$w_1,\ldots,w_k$ be nontrivial cyclically reduced words of lengths at most $L$.
Assume no two words are related by cyclic rotation or by cyclic rotation followed by inversion.  Set
\[
r_i:=r(w_i),\qquad Z_i:=\frac{\operatorname{Tr}w_i(\mathbf U)}{\sqrt{r_i}},
\qquad Z=(Z_1,\ldots,Z_k),
\]
and let $G=(G_1,\ldots,G_k)$ be a tuple of independent circular complex Gaussians with
$\mathbb E\abs{G_i}^2=1$. Thus the real and imaginary parts of each
$G_i$ are independent centered real Gaussians of variance $1/2$.

Fix a cyclic rotation of each word. Let $d_i$ be the number of letters of $w_i$, other than the first, that are equal to the first letter of $w_i$ or its inverse, and let $n_{ij}$ be the number of letters of $w_j$ that are equal to the first letter of $w_i$ or its inverse. Set
\[
S_1:=\left(\sum_i\frac{d_i^2}{r_i}\right)^{1/2},
\qquad
S_2:=\left(\sum_{i,j}\frac{n_{ij}^2}{r_ir_j}\right)^{1/2}.
\]

For every even integer $p\ge2$ with $N\ge(2pL)^{3/2}$,
\[
\begin{aligned}
W_p(\mathcal L(Z),\mathcal L(G))
&\le\frac{pL}{N}S_1+\frac{2c_p\sqrt{pL}}N S_2\\
&\le\frac pN(LS_1+2\sqrt L\,S_2)\\
&\le\frac pN(\sqrt k\,L^2+2kL^{3/2})
\le\frac{3pkL^2}{N},
\end{aligned}
\]
where $c_p:=(\mathbb E\abs g^p)^{1/p}$ for a standard real Gaussian $g$.  In particular, when $N\ge(4L)^{3/2}$,
\[
W_2(\mathcal L(Z),\mathcal L(G))
\le\frac{2L}{N}S_1+\frac{2\sqrt{2L}}N S_2.
\]

If $k=k_N$, $L=L_N$ and the words depend on $N$, with $k_NL_N^2/N\to0$, then $W_p(\mathcal L(Z),\mathcal L(G))\to0$ for each fixed even $p\ge2$.
\end{theorem}

We deduce Theorem~\ref{thm:haar-word-gaussian} from
Lemma~\ref{lem:gaussian-stein-hs}, which we prove using the
Ornstein--Uhlenbeck flow, as in
\cite[proof of Proposition~3.4]{ledoux-nourdin-peccati}
and~\cite[Section~6.2]{fang-koike}. Since Gaussian convolution gives a smooth density at every positive time, we need not assume that the bounded random vector $V$ of the lemma has a density; compare \cite[Appendix~A.2]{fang-koike}.

The $W_1$ bound in the lemma also follows from the dimension-free
Hilbert--Schmidt Hessian bound of Chatterjee and
Meckes~\cite[Lemma~2.2(i)]{chatterjee-meckes}.

\begin{lemma}
\label{lem:gaussian-stein-hs}
Let $V$ be a bounded random vector in $\R^d$.  Suppose integrable real random coefficients
$\beta\in\R^d$ and $K\in M_d(\R)$ satisfy
\[
\mathbb E[V_a\phi(V)]
=\mathbb E[\beta_a\phi(V)]+\sum_b\mathbb E[K_{ab}\partial_b\phi(V)]
\qquad(1\le a\le d)
\]
for every smooth bounded real-valued function $\phi$ with bounded derivatives.  If $\gamma_d$ is the standard Gaussian law, then
\[
W_1(\mathcal L(V),\gamma_d)
\le\mathbb E\norm{\beta}_2+\mathbb E\norm{K-I}_{\rm HS}.
\]
The matrix $K$ need be neither symmetric nor positive, and the random coefficients need not be functions of $V$.
If $2\le p<\infty$ and
$\mathbb E\norm{\beta}_2^p+\mathbb E\norm{K-I}_{\rm HS}^p<\infty$, then, writing
$c_p:=(\mathbb E\abs{g}^p)^{1/p}$ for $g\sim\mathcal N(0,1)$,
\[
W_p(\mathcal L(V),\gamma_d)
\le\bigl(\mathbb E\norm{\beta}_2^p\bigr)^{1/p}
+c_p\bigl(\mathbb E\norm{K-I}_{\rm HS}^p\bigr)^{1/p}.
\]
Moreover, $c_2=1$, and $c_p\le\sqrt p$ for even integers $p$.
\end{lemma}

\begin{proof}
Fix $p=1$ or $2\le p<\infty$, and set $\kappa_1=1$ and $\kappa_p=c_p$ for $p\ge2$.
Let $G$ be a standard Gaussian vector independent of $(V,\beta,K)$, and set
$a_t=e^{-t}$, $q_t=1-e^{-2t}$ and $X_t=a_tV+\sqrt{q_t}G$.
The density $\rho_t$ of $X_t$ is positive and smooth for $t>0$. Applying the assumed identity to $y\mapsto\mathbb E_G f(a_ty+\sqrt{q_t}G)$ for every smooth bounded real-valued function $f$ with bounded derivatives and then integrating by parts in $G$ gives
\[
g_t(X_t):=(\nabla\log \rho_t)(X_t)+X_t
=a_t\mathbb E[\beta\mid X_t]
+\frac{a_t^2}{\sqrt{q_t}}\mathbb E[(K-I)G\mid X_t].
\]
For every fixed matrix $A\in M_d(\R)$, we have
$\norm{AG}_{L^p(\ell_2)}\le \kappa_p\norm A_{\rm HS}$.
Indeed, for $p=1$, this follows from Cauchy--Schwarz. For $p\ge2$, we diagonalize $A^*A$ and reduce the estimate to Minkowski's inequality in $L^{p/2}$.
Hence, by Jensen's inequality applied conditionally, we have
\[
\norm{g_t(X_t)}_{L^p}
\le a_t\norm{\beta}_{L^p(\ell_2)}+\kappa_p\frac{a_t^2}{\sqrt{q_t}}\norm{K-I}_{L^p({\rm HS})}.
\]

Choose $R$ so that $\norm V_2\le R$ almost surely. Since $\rho_t$ is a mixture of Gaussian densities, differentiating $\log \rho_t$ gives
\[
\nabla\log \rho_t(x)=-\frac{x-a_t\mathbb E[V\mid X_t=x]}{q_t},\qquad
D^2\log \rho_t(x)=-\frac I{q_t}
+\frac{a_t^2}{q_t^2}\operatorname{Cov}(V\mid X_t=x).
\]
The conditional mean and covariance are bounded by $R$ and $R^2$, respectively.
Thus the vector field $v_t=-g_t$ is globally Lipschitz with linear
growth on every $[\epsilon,T]\subset(0,\infty)$. Since $\partial_t\rho_t=-\operatorname{div}(\rho_tv_t)$, the flow of $v_t$ carries the law with density $\rho_\epsilon$ to the law with density $\rho_t$.
Coupling the endpoints along this flow and applying Minkowski's
inequality gives
\[
W_p(\mathcal L(X_\epsilon),\mathcal L(X_T))
\le\int_\epsilon^T\norm{g_t(X_t)}_{L^p}\,dt.
\]
Moreover, $X_\epsilon\to V$ and $X_T\to G$ in $L^p$, since
$X_t=a_tV+\sqrt{q_t}G$.
Both $\int_0^\infty a_t\,dt$ and
$\int_0^\infty a_t^2/\sqrt{q_t}\,dt$ equal one, so letting
$\epsilon\downarrow0$ and $T\to\infty$, we obtain $W_p(\mathcal L(V),\gamma_d)\le\norm{\beta}_{L^p(\ell_2)}+\kappa_p\norm{K-I}_{L^p({\rm HS})}$.
Finally, $c_2=1$ and $c_p^p=(p-1)!!\le p^{p/2}$ for even $p$.
\end{proof}

\begin{proof}[Proof of Theorem~\ref{thm:haar-word-gaussian}]
\emph{Integration by parts.}
Write $T_v=\operatorname{Tr}v(\mathbf U)$, with $T_1=\operatorname{Tr}(1)=N$.
For the chosen rotation $w=xW$ of one of the words $w_i$, with first letter $x$, we may assume that $x$ is a generator. Otherwise, $x^{-1}$ is a generator, and we replace the corresponding Haar matrix by its adjoint and interchange $x$ and $x^{-1}$ in every word. Set
\[
R_w=\frac1N\left(
\sum_{w=xAx^{-1}A'}T_AT_{A'}-\sum_{w=xAxA'}T_{xA}T_{xA'}\right),
\]
\[
J_{w,v}=\frac1N\left(
\sum_{v\ {\rm cut\ as}\ x^{-1}Y}T_{WY}
-\sum_{v\ {\rm cut\ as}\ xY}T_{xWxY}\right),
\]
where the sums in $R_w$ and in $J_{w,v}$ run over the occurrences of $x$ and $x^{-1}$ in $W$ and in $v$, respectively. A cut of $v$ means the cyclic rotation of $v$ beginning at such an occurrence.
We choose a self-adjoint basis $(H_\alpha)$ orthonormal for the
unnormalized Hilbert--Schmidt inner product. Then we apply the Haar
differentiation and basis identities from the proof of Lemma~\ref{lem:haar-growing-word-moments}
to
$\operatorname{Tr}(H_\alpha w)F$, where $F$ is a smooth bounded real-valued function of $(T_{w_1},\ldots,T_{w_k})$ with bounded derivatives. After summation over $\alpha$, differentiating the first letter of $w$ and its
remaining letters gives $iNT_wF$ and $-iNR_wF$, respectively. Differentiating $F$ gives, for each $j$, $-iNJ_{w,w_j}$ times its derivative with respect to $T_{w_j}$ and $-iNJ_{w,w_j^{-1}}$ times its derivative with respect to $\overline{T_{w_j}}=T_{w_j^{-1}}$. By Haar invariance, the expectations of these terms sum to zero.
For $Z_i=T_{w_i}/\sqrt{r_i}$, we write $Z_j=\xi_j+i\eta_j$ and use
\[
\partial_{Z_j}=\tfrac12(\partial_{\xi_j}-i\partial_{\eta_j}),\qquad
\partial_{\bar Z_j}=\tfrac12(\partial_{\xi_j}+i\partial_{\eta_j}).
\]
By the chain rule, we have
\[
\mathbb E[Z_iF(Z)]
=\mathbb E[b_iF(Z)]
+\sum_j\mathbb E[B_{ij}\partial_{Z_j}F(Z)
+C_{ij}\partial_{\bar Z_j}F(Z)],
\]
where
\[
b_i=\frac{R_{w_i}}{\sqrt{r_i}},\qquad
B_{ij}=\frac{J_{w_i,w_j}}{\sqrt{r_ir_j}},\qquad
C_{ij}=\frac{J_{w_i,w_j^{-1}}}{\sqrt{r_ir_j}}.
\]
If the chosen rotation of $w_i$ begins with an inverse generator, we now undo the change of variables for that $i$. All these identities then hold for the original
trace vector on the original product probability space.

\emph{Joins that reduce to the identity.}
As in the proof of Lemma~\ref{lem:haar-growing-word-moments},
splits and same-sign joins are nontrivial.
An opposite-sign join $WY$ reduces to the identity if and only if
$Y=W^{-1}$, or equivalently, the cut $x^{-1}Y$ of $v$ equals
$x^{-1}W^{-1}$. Such a cut exists exactly when $v$ is a cyclic rotation
of $w^{-1}$. Once one such cut is chosen, the rotations fixing $w$ yield the
$r(w)$ possible cuts.
As we saw in the proof of Lemma~\ref{lem:group-word-gradient-gram}, no
nontrivial cyclically reduced word is a rotation of its inverse.
By this and the assumption that no two of the words are related by cyclic rotation followed by inversion, no join contributing to any $B_{ij}$ reduces to the group identity.
In $C_{ij}$, the joins that reduce to the identity occur only for $j=i$, since no two of the words are related by cyclic rotation. There are $r_i$ of them, each contributing $1/r_i$, so their total contribution is $\delta_{ij}$.

\emph{Moment bounds.}
For even $p$, the $p$th absolute moment of a split product is the
expectation of a product of $2p$ nontrivial traces, with total word
length at most $pL$. A joined word has length at most $2L$.
Hence, by Lemma~\ref{lem:haar-growing-word-moments} and the assumption $N\ge(2pL)^{3/2}$, we have
\[
\|T_AT_{A'}\|_{L^p}\le pL,\qquad
\|T_v\|_{L^p}\le\sqrt{2pL}\quad(1\ne v,\ |v|\le2L).
\]
Counting the $d_i$ splits and at most $n_{ij}$ joins yields
\[
\|b_i\|_{L^p}\le\frac{pL\,d_i}{N\sqrt{r_i}},\qquad
\|B_{ij}\|_{L^p},\ \|C_{ij}-\delta_{ij}\|_{L^p}
\le\frac{\sqrt{2pL}\,n_{ij}}{N\sqrt{r_ir_j}}.
\]

\emph{Conclusion.}
Set $V=\sqrt2(\Ree Z_1,\Imm Z_1,\ldots,\Ree Z_k,\Imm Z_k)$
and $D_{ij}=C_{ij}-\delta_{ij}$. Taking the real and imaginary parts of the identity for
$\mathbb E[Z_iF(Z)]$ gives the formula assumed in
Lemma~\ref{lem:gaussian-stein-hs}, with
$\beta_{2i-1}=\sqrt2\Ree b_i$ and $\beta_{2i}=\sqrt2\Imm b_i$.
The $(i,j)$-block of $K-I$ is
\[
\begin{pmatrix}
\Ree(B_{ij}+D_{ij})&\Imm(B_{ij}-D_{ij})\\
\Imm(B_{ij}+D_{ij})&\Ree(D_{ij}-B_{ij})
\end{pmatrix},
\]
and its squared Hilbert--Schmidt norm is $2(|B_{ij}|^2+|D_{ij}|^2)$.
By Minkowski's inequality in $L^{p/2}$, we have
\[
\|\beta\|_{L^p(\ell_2)}\le\frac{\sqrt2\,pL}{N}S_1,\qquad
\|K-I\|_{L^p({\rm HS})}\le\frac{2\sqrt{2pL}}N S_2.
\]
Since $V$, $\beta$ and $K$ are bounded at each $N$, we may apply Lemma~\ref{lem:gaussian-stein-hs}. Scaling distances back by $1/\sqrt2$, we obtain the first bound.
The others follow from $c_p\le\sqrt p$, $d_i<L$, $n_{ij}\le L$,
and $r_i\ge1$, using $c_2=1$ when $p=2$.
For fixed $p$, the condition $k_NL_N^2/N\to0$ also ensures that
$N\ge(2pL_N)^{3/2}$ for all sufficiently large $N$.
\end{proof}

\begin{corollary}
\label{cor:haar-word-bias}
For independent Haar matrices $U_1,\ldots,U_s$ in $U(N)$ and every nontrivial $w\in\mathbb F_s$,
\[
\abs{\mathbb E\operatorname{Tr}w(\mathbf U)}
\le\frac{\abs w(\abs w-1)}N
\qquad\text{if }N\ge\abs w^{3/2}.
\]
\end{corollary}

\begin{proof}
We cyclically reduce $w$ to a nontrivial word of length $m\le\abs w$.
Taking $F=1$ in the identity for $\mathbb E[Z_iF(Z)]$ from the proof of Theorem~\ref{thm:haar-word-gaussian}, we obtain $\mathbb ET_w=\mathbb ER_w$.
Each split produces two nontrivial words of total length at most $m$.
By Lemma~\ref{lem:haar-growing-word-moments} and the assumption $N\ge\abs
w^{3/2}$, the absolute value of the
expected product of their traces is at most $m$. There are at most $m-1$
terms, each with coefficient $\pm1/N$, so
$\abs{\mathbb ET_w}\le m(m-1)/N\le\abs w(\abs w-1)/N$.
\end{proof}

\subsection{Extreme values}

Fang and Koike~\cite[Theorems~2.2 and~4.2]{fang-koike} derived moderate
deviations from $W_p$ bounds by choosing a moment order of logarithmic
size and applying Markov's inequality to the coupling error. We use this device in the proof of Theorem~\ref{thm:haar-word-gumbel} to obtain the estimate \eqref{eq:haar-word-relative-tails} for joint tail
probabilities of word traces.

\begin{theorem}[Poisson limits for word traces]
\label{thm:haar-word-gumbel}
For each $N$, let $s_N\ge1$ be an integer and let $U_1,\ldots,U_{s_N}$ be independent Haar unitary matrices in $U(N)$.
Let $L_N\ge1$ be integers and let $\mathcal C_N\subseteq\mathbb F_{s_N}$ be a collection of $k_N\ge1$
nontrivial cyclically reduced words of lengths at most $L_N$, with no two related by cyclic rotation or by cyclic
rotation followed by inversion.  Suppose
\[
k_N\longrightarrow\infty,
\qquad
\frac{L_N^2(\log k_N)^{3/2}}N\longrightarrow0.
\]
Set $T_w:=\operatorname{Tr}w(\mathbf U)$ and
\[
\Xi_N:=\sum_{w\in\mathcal C_N}\delta_{\abs{T_w}^2/r(w)-\log k_N}.
\]
Then, in the vague topology, $\Xi_N$ converges in distribution as
$N\to\infty$ to a Poisson point process on $\R$ with intensity $e^{-x}\,dx$. Moreover, for any finite family of pairwise
disjoint intervals $I_j=(\alpha_j,\beta_j]$, where $-\infty<\alpha_j<\beta_j\le\infty$
(with $(\alpha_j,\infty]$ read as $(\alpha_j,\infty)$), the counts $\Xi_N(I_j)$
converge jointly to independent Poisson variables with means
$e^{-\alpha_j}-e^{-\beta_j}$, where $e^{-\infty}=0$. In particular, for every $x\in\R$,
\[
\mathbb P\left\{\max_{w\in\mathcal C_N}\frac{\abs{T_w}^2}{r(w)}-\log k_N\le x\right\}
\longrightarrow\exp(-e^{-x}).
\]
\end{theorem}

\begin{proof}
Write $k=k_N$, $L=L_N$, and $Z_w=T_w/\sqrt{r(w)}$.  Fix $h\ge1$ and let $w_1,\ldots,w_h$ be distinct words in $\mathcal C_N$.  Set
\[
p:=2\left\lceil\frac{(h+2)\log k}{2}\right\rceil,
\qquad
D:=\frac pN\bigl(\sqrt h\,L^2+2hL^{3/2}\bigr),
\qquad
\epsilon:=2eD.
\]
For all sufficiently large $N$, we have $p\ge2$ and $N\ge(2pL)^{3/2}$, since $k\to\infty$, $L\ge1$ and $L^2(\log k)^{3/2}/N\to0$.  By Theorem~\ref{thm:haar-word-gaussian}, there is a coupling of $Z=(Z_{w_1},\ldots,Z_{w_h})$ with a vector $G=(G_1,\ldots,G_h)$ of independent circular Gaussians such that $\norm{Z-G}_{L^p(\ell_2)}\le2D$. Hence, by  Markov's inequality,
\[
\mathbb P\{\max_{j\le h}\abs{Z_{w_j}-G_j}>\epsilon\}\le e^{-p}\le k^{-h-2}.
\]
For $t_j>\epsilon$, the triangle inequality therefore gives
\[
e^{-\sum_j(t_j+\epsilon)^2}-e^{-p}
\le\mathbb P\{\abs{Z_{w_j}}>t_j\text{ for all }j\}
\le e^{-\sum_j(t_j-\epsilon)^2}+e^{-p}.
\]
We take $t_j=\sqrt{\log k+x_j}$ for fixed real levels $x_j$.  Since
\[
\epsilon\sqrt{\log k}=O_h(L^2(\log k)^{3/2}/N)=o(1),
\]
these inequalities give, uniformly over all $h$-tuples of distinct words,
\[
\mathbb P\{\abs{Z_{w_j}}^2-\log k>x_j\text{ for all }j\}
=k^{-h}e^{-\sum_jx_j}(1+o(1)).
\eqtag\label{eq:haar-word-relative-tails}
\]

Let $C_{j,N}=\Xi_N(I_j)$ and $\lambda_j=e^{-\alpha_j}-e^{-\beta_j}$.
Inclusion--exclusion in \eqref{eq:haar-word-relative-tails} and summation over tuples of distinct words in $\mathcal C_N$ give, for every fixed tuple of nonnegative integers
$m_j$,
\[
\mathbb E\prod_j(C_{j,N})_{m_j}\longrightarrow\prod_j\lambda_j^{m_j},
\]
where $(n)_m$ denotes the falling factorial.
Since Poisson laws are determined by their moments, it follows that the $C_{j,N}$ converge jointly to independent
Poisson variables with means $\lambda_j$. The bounded half-open intervals form a dissecting semiring of sets whose
boundaries are null for the intensity $e^{-x}\,dx$, so convergence of the joint counts implies vague
convergence in distribution~\cite[Theorem~4.11]{kallenberg-random-measures}.
Taking the probability that $\Xi_N((x,\infty))=0$, we obtain the last assertion.
\end{proof}

\begin{corollary}
\label{cor:haar-word-full-ball-extremes}
Fix $s\ge2$ and set $q=2s-1$. For every integer $L\ge1$, let $\mathcal C_L$ contain exactly one representative
from each class of nontrivial cyclically reduced words of length at most $L$, modulo cyclic rotation and inversion, and set
$k_L=\abs{\mathcal C_L}$.  If $L=L_N\to\infty$ and $L_N=o(N^{2/7})$, the conclusions of
Theorem~\ref{thm:haar-word-gumbel} hold for this collection.  Moreover,
\[
k_L\sim\frac{q^{L+1}}{2(q-1)L},
\qquad
\log k_L=L\log q-\log L+\log\frac q{2(q-1)}+o(1).
\]
With probability tending to one, every word maximizing $\abs{T_w}/\sqrt{r(w)}$ over $\mathcal C_L$ has $r(w)=1$.
\end{corollary}

\begin{proof}
The matrix of allowed transitions between the $2s$ signed letters in a reduced word is $J-P$, where $J$ has all entries one and $P$ exchanges each
letter with its inverse.  Its eigenvalues are $q$ once, $1$ with multiplicity $s$, and $-1$ with multiplicity $s-1$.
Thus the number of cyclically reduced strings of length $m$ is
\[
c_m=q^m+s+(s-1)(-1)^m.
\]
A string fixed by a nontrivial rotation has the form $v^{m/d}$ for a string $v$ of length $d$, where $d\mid m$ and $d\le m/2$, so it is determined by its initial segment of length $d$.  Since there are at most $2sq^{d-1}$ reduced strings of length $d$, there are $O_s(q^{m/2})$ such strings.  Every remaining rotation orbit has size $m$,
and inversion pairs rotation orbits, fixing none of them, as shown in the
proof of Lemma~\ref{lem:group-word-gradient-gram}. Consequently the number of classes of
exact length $m$ is $q^m/(2m)+O_s(q^{m/2})$, and
\[
k_L=\frac12\sum_{m=1}^L\frac{q^m}{m}+O_s(q^{L/2})
\sim\frac{q^{L+1}}{2(q-1)L}.
\]
Taking logarithms gives the stated expansion of $\log k_L$. By this expansion, the condition $L^2(\log k_L)^{3/2}/N\to0$ of Theorem~\ref{thm:haar-word-gumbel} becomes $L^{7/2}/N\to0$.

There are $O_s(q^{L/2})=o(k_L)$ proper-power classes. By \eqref{eq:haar-word-relative-tails} with $h=1$,
for each fixed $x$, the probability that any of them satisfies
$|T_w|^2/r(w)-\log k_L>x$ tends to zero. If a proper power attains the
maximum, either this event occurs or
$\max_{w\in\mathcal C_L}|T_w|^2/r(w)\le\log k_L+x$. Hence the upper limit of the probability that a proper power attains the
maximum is at most $\exp(-e^{-x})$. Letting $x\to-\infty$ proves the last assertion.
\end{proof}

Without division by $r(w)$, both the centering and the maximizing
words change. In the Gaussian approximation, the trace of a generator power $a_i^m$ has variance $m$, while the trace of every other class $w=v^{r(w)}$ in $\mathcal C_L$, with $\abs v=d\ge2$, has variance $r(w)\le L/d$. The largest of the $sL$ values $\abs{T_{a_i^m}}^2$ is of order $L\log L$. For each $d\ge2$, there are at most of order $L(2s-1)^d$ such classes $w$, and $\abs{T_w}^2$ exceeds this level with probability at most of order $L^{-d}$, up to logarithmic factors. Hence the expected number of such classes above this level tends to zero.
This suggests that generator powers
determine the maximum.

\begin{theorem}
\label{thm:haar-word-raw-extremes}
Fix $s\ge2$ and let $\mathcal C_L$ be the collections from Corollary~\ref{cor:haar-word-full-ball-extremes}.  Suppose the integers $L=L_N\to\infty$ satisfy $L_N^3/N\to0$, and set
\[
b_L:=\log(sL)-\log\log(sL),
\qquad
\Psi_N:=\sum_{w\in\mathcal C_L}\delta_{\abs{T_w}^2/L-b_L}.
\]
Then $\Psi_N$ converges vaguely in distribution to a Poisson point
process of intensity $e^{-x}\,dx$. Moreover, for any finite family of pairwise disjoint intervals $I_j$ as in Theorem~\ref{thm:haar-word-gumbel}, the counts $\Psi_N(I_j)$ converge jointly to independent Poisson variables with the same means as in that theorem. In particular,
\[
\mathbb P\left\{\frac1L\max_{w\in\mathcal C_L}\abs{T_w}^2-b_L\le x\right\}
\longrightarrow\exp(-e^{-x})
\qquad(x\in\R).
\]
With probability tending to one, every word maximizing $\abs{T_w}$ over $\mathcal C_L$ belongs to the class of $a_i^m$ for some
$1\le i\le s$ and $1\le m\le L$.
\end{theorem}

\begin{proof}
First let $G_{i,m}$, $1\le i\le s$, $1\le m\le L$, be independent
circular complex Gaussians with $\mathbb E|G_{i,m}|^2=1$, and set
\[
Q_{i,m}:=\frac mL\abs{G_{i,m}}^2-b_L.
\]
For fixed $x$, let $\theta=b_L+x\sim\log L$.  The function $f_\theta(t)=e^{-\theta/t}$ on $(0,1]$, extended by $f_\theta(0)=0$, is
increasing for large $L$.  Hence
\[
\begin{aligned}
\left|\sum_{m=1}^Le^{-\theta L/m}-L\int_0^1e^{-\theta/t}\,dt\right|&\le e^{-\theta},\\
\int_0^1e^{-\theta/t}\,dt
&=\frac{e^{-\theta}}\theta\int_0^\infty e^{-u}(1+u/\theta)^{-2}\,du
\sim\frac{e^{-\theta}}\theta.
\end{aligned}
\]
Since $\theta/L\to0$, the error $e^{-\theta}$ is negligible compared with $Le^{-\theta}/\theta$.  It follows that
\[
\begin{gathered}
\sum_{i,m}\mathbb P\{Q_{i,m}>x\}
\sim\frac{sLe^{-b_L-x}}{b_L+x}\longrightarrow e^{-x},\\
\max_{i,m}\mathbb P\{Q_{i,m}>x\}=e^{-b_L-x}\longrightarrow0.
\end{gathered}
\]
For any fixed family of disjoint intervals $I_j$, the sums $\sum_{i,m}\mathbb P\{Q_{i,m}\in I_j\}$ converge to the corresponding Poisson means, and $\max_{i,m}\mathbb P\{Q_{i,m}\in I_j\}\to0$. Since the $Q_{i,m}$ are independent, the joint probability generating function of the counts is a product over the pairs $(i,m)$. Its logarithm differs from its first-order approximation by at most a constant times the sum of the squared probabilities, which is at most their maximum times their sum and hence tends to zero. Thus the counts have
independent Poisson limits. The vague convergence of the point process and the limit law of the maximum follow as in the proof of Theorem~\ref{thm:haar-word-gumbel}.

We now transfer this limit to the traces of the generator powers. Set
$Z_{i,m}=\operatorname{Tr}(U_i^m)/\sqrt m$.
For the $sL$ words $a_i^m$, $1\le i\le s$ and $1\le m\le L$,
the quantities in Theorem~\ref{thm:haar-word-gaussian} satisfy
$S_1^2\le sL^2$ and $S_2^2\le sL^3$.
Hence, by the $p=2$ estimate of that theorem, we have
\[
W_2(\mathcal L(Z),\mathcal L(G))\le C_s L^2/N
\qquad(N\ge(4L)^{3/2}).
\]
We choose a coupling of $Z$ and $G$ with $\norm{Z-G}_{L^2(\ell_2)}\le2C_sL^2/N$ and set $\Delta_L=\norm{Z-G}_{\ell_2}$.
Since $L^3/N\to0$, we have
$\Delta_L\sqrt{\log L}\to0$ in probability, while
$\max_{i,m}\abs{G_{i,m}}=O_{\mathbb P}(\sqrt{\log L})$.
Thus
\[
\max_{i,m}\left|
\frac mL\abs{Z_{i,m}}^2-\frac mL\abs{G_{i,m}}^2
\right|
\le2\Delta_L\max_{i,m}\abs{G_{i,m}}+\Delta_L^2
\longrightarrow0\quad\text{in probability}.
\eqtag\label{eq:haar-power-score-coupling}
\]
The limiting Poisson point process has no atom at any fixed endpoint,
and the number of $Q_{i,m}$ above each fixed level is bounded in
probability. Together with \eqref{eq:haar-power-score-coupling}, these
facts give convergence of the interval counts for the generator powers,
including intervals unbounded above. Vague convergence then follows from~\cite[Theorem~4.11]{kallenberg-random-measures} applied to the interval counts.

It remains to exclude the words that are not generator powers.
Write a cyclically reduced word as $w=v^r$, with $r=r(w)$, so that
$v$ is not a proper power. Let $d=\abs v$. Then $rd\le L$, and a
class outside the generator powers has $d\ge2$. For $d\ge2$,
set
\[
p_d:=2\lceil d\log q+2\log L\rceil,
\qquad q=2s-1.
\]
Uniformly for $2\le d\le L$, we have $p_d=O_s(d+\log L)=O_s(L)$ and $(2p_dL)^{3/2}=O_s(L^3)=o(N)$.
By Theorem~\ref{thm:haar-word-gaussian} with $k=1$,
there is a circular Gaussian coupled to $Z_w=T_w/\sqrt r$ so that
their $L^{p_d}$ distance is at most $6p_dL^2/N$. With
$\epsilon_d=6ep_dL^2/N$, the probability that the coupling error exceeds
$\epsilon_d$ is at most
\[
e^{-p_d}\le q^{-2d}L^{-4}.
\]
Fix $x\in\R$ and set $\theta=b_L+x$. For all sufficiently large $L$,
the inequality $\abs{T_w}^2/L-b_L>x$ is equivalent to $\abs{Z_w}>t_w$, where
\[
t_w:=\sqrt{L\theta/r}\ge\sqrt{\theta d}.
\]
Uniformly in $d$ and $r$, we have
$\epsilon_d=O_s(L^3/N)=o(1)$ and $t_w\ge\sqrt{2\theta}\to\infty$.
Hence $\epsilon_d/t_w\to0$ uniformly.
Thus, eventually, $(t_w-\epsilon_d)^2\ge\tfrac34\theta d$, and the triangle inequality gives
\[
\mathbb P\{\abs{T_w}^2/L-b_L>x\}
\le e^{-3d(b_L+x)/4}+q^{-2d}L^{-4}.
\]
There are at most $C_sq^d$ classes of words $v$ of length $d$ that
are not proper powers, and at most $L$ powers of each such word to
consider. Summing over all classes other than generator powers, we see that the expected number of such classes with $\abs{T_w}^2/L-b_L>x$ is at most
\[
\begin{aligned}
C_sL\sum_{d=2}^Lq^d\bigl(e^{-3d(b_L+x)/4}+q^{-2d}L^{-4}\bigr)
&=O_{s,x}\bigl(L^{-1/2}(\log L)^{3/2}\bigr)+O_s(L^{-3})\\
&\longrightarrow0.
\end{aligned}
\]
Above every fixed level, $\Psi_N$ agrees with its restriction to
generator powers with probability tending to one. The two processes
therefore have the same limits for interval counts and maxima, and
the same vague limit. If a class outside the generator powers attains the maximum, then either $\abs{T_w}^2/L-b_L>x$ for some class $w$ outside the generator powers, or $\abs{T_w}^2/L-b_L\le x$ for every generator power $w$. Letting first $N\to\infty$ and then $x\to-\infty$ proves the
last assertion.
\end{proof}

\begin{corollary}
\label{cor:haar-word-joint-extremes}
Fix $s\ge2$, and let the integers $L=L_N\to\infty$ satisfy $L_N=o(N^{2/7})$.
For the collections $\mathcal C_L$ from Corollary~\ref{cor:haar-word-full-ball-extremes}, let $\Xi_N$ and $\Psi_N$ be the point processes in Theorems~\ref{thm:haar-word-gumbel} and \ref{thm:haar-word-raw-extremes}, respectively.
Then $(\Xi_N,\Psi_N)$ converges jointly in distribution, for the product vague topology, to two independent Poisson point processes
of intensity $e^{-x}\,dx$.  Moreover, for any two finite families $(I_j)$ and $(I'_l)$ of intervals $(\alpha,\beta]$ with $-\infty<\alpha<\beta\le\infty$, each family pairwise disjoint, the counts $\Xi_N(I_j)$ and $\Psi_N(I'_l)$ converge jointly to independent Poisson variables, with mean $e^{-\alpha}-e^{-\beta}$ for the interval $(\alpha,\beta]$, where $e^{-\infty}=0$.  In particular, for all
$x,y\in\R$,
\[
\mathbb P\left\{\begin{aligned}
\max_{w\in\mathcal C_L}\frac{\abs{T_w}^2}{r(w)}-\log k_L&\le x,\\
\frac1L\max_{w\in\mathcal C_L}\abs{T_w}^2-b_L&\le y
\end{aligned}
\right\}
\longrightarrow\exp(-e^{-x}-e^{-y}).
\]
\end{corollary}

\begin{proof}
Write $k=k_L$, $q=2s-1$, and let $\mathcal P_L$ be the $sL$ generator-power classes.
By \eqref{eq:haar-word-relative-tails} and $sL/k\to0$, the probability
that any class in $\mathcal P_L$ contributes a point to $\Xi_N$ above
a fixed level tends to zero.
As shown in the proof of Theorem~\ref{thm:haar-word-raw-extremes}, we may
restrict $\Psi_N$ to $\mathcal P_L$. We may further keep only the classes
$a_i^m$ with $m\ge\lceil L/2\rceil$.
Indeed, for the variables $Q_{i,m}$ defined in the proof of
Theorem~\ref{thm:haar-word-raw-extremes}, the expected number with
$m<\lceil L/2\rceil$ and $Q_{i,m}>y$ is at most
\[
sL e^{-2(b_L+y)}=O_{s,y}((\log L)^2/L)\longrightarrow0,
\]
and the probability that any of the generator powers $a_i^m$ with $m<\lceil L/2\rceil$ contributes a point above $y$ also tends to zero, since
\eqref{eq:haar-power-score-coupling} reduces this assertion to the
preceding Gaussian bound at $y-\epsilon$, for any fixed
$\epsilon>0$. We call the classes $a_i^m$ with $m\ge\lceil L/2\rceil$ the retained powers.

Fix $h,\ell\ge0$ with $h+\ell\ge1$, let $w_1,\ldots,w_h$ be distinct words outside $\mathcal P_L$ and $a_{i_1}^{m_1},\ldots,a_{i_\ell}^{m_\ell}$ distinct retained powers, and set
\[
p:=2\left\lceil\frac{(h+2)\log k+(\ell+2)\log L}{2}\right\rceil.
\]
Since $\log k\sim L\log q$, we have $p=O_{s,h,\ell}(L)$ and
$N\ge(2pL)^{3/2}$ eventually. By Theorem~\ref{thm:haar-word-gaussian} and
Markov's inequality, we can couple the selected
variables $T_w/\sqrt{r(w)}$ with independent circular Gaussians of
variance one so that the maximum coordinate error exceeds
$\eta=O_{s,h,\ell}(pL^2/N)$ with probability at most
$k^{-h-2}L^{-\ell-2}$.

Fix real levels $x_1,\ldots,x_h$ and $y_1,\ldots,y_\ell$. For a selected word $w_j$, the point $\abs{T_{w_j}}^2/r(w_j)-\log k$ of $\Xi_N$ lies above $x_j$ if and only if $\abs{T_{w_j}}/\sqrt{r(w_j)}$ exceeds the threshold $\sqrt{\log k+x_j}=O(\sqrt L)$. For a retained power $a_{i_j}^{m_j}$, the corresponding threshold at level $y_j$ in $\Psi_N$ is $\sqrt{(L/m_j)(b_L+y_j)}=O(\sqrt{\log L})$, where $1\le L/m_j\le2$.
Since $L=o(N^{2/7})$, we have $\eta\sqrt L\to0$. Thus shifting
any of these thresholds by $\eta$ changes its square by $o(1)$. Comparing the Gaussian
probabilities as in the proof of Theorem~\ref{thm:haar-word-gumbel} gives, uniformly over
the selected coordinates,
\[
\begin{aligned}
&\mathbb P\left\{\begin{aligned}
\abs{T_{w_j}}^2/r(w_j)-\log k&>x_j &&(1\le j\le h),\\
\abs{T_{a_{i_j}^{m_j}}}^2/L-b_L&>y_j &&(1\le j\le\ell)
\end{aligned}\right\}\\
&\qquad
=k^{-h}e^{-\sum_{j=1}^hx_j}
\prod_{j=1}^{\ell}e^{-(L/m_j)(b_L+y_j)}(1+o(1))
+O(k^{-h-2}L^{-\ell-2}).
\end{aligned}
\]
Summing the additive error over at most
$k^h(sL)^\ell$ tuples gives $O_{s,h,\ell}(k^{-2}L^{-2})=o(1)$.

For a retained power $a_i^m$, the probability
that $Q_{i,m}$ lies in a bounded interval $(\alpha,\beta]$, divided by the
probability that it exceeds $\alpha$, is
\[
1-e^{-(L/m)(\beta-\alpha)}\ge1-e^{-(\beta-\alpha)}>0.
\]
For the selected words $w_j$, the corresponding ratio is the same expression with $L/m$ replaced by one; for an interval unbounded above, it is one.
Since these ratios are bounded below by a positive constant, inclusion--exclusion preserves relative errors. Hence the estimate above holds, with a factor $1+o(1)$ uniform over the selected coordinates, for the probabilities that the quantities $\abs{T_{w_j}}^2/r(w_j)-\log k$ and $\abs{T_{a_{i_j}^{m_j}}}^2/L-b_L$ lie in fixed intervals.
For $\Xi_N$ with the generator powers removed, the number of ordered
$h$-tuples of distinct indices is $(k-sL)_h\sim k^h$.
For the retained powers, the sums of the individual interval
probabilities converge to the corresponding Poisson means, and the largest
individual probability tends to zero.
In the product of these sums, the terms with a repeated power index contribute $o(1)$, so the sum over tuples of distinct retained powers has the same limit. Indeed, for each pair
of positions, we may bound one probability by the largest individual
probability
and use the bounded sums for the other positions. Thus the joint
factorial moments of the interval counts for the two restricted
processes converge to those of independent Poisson variables.

Arguing as in the proof of Theorem~\ref{thm:haar-word-gumbel}, we now obtain the joint convergence of the counts and, from the probability that $\Xi_N((x,\infty))=\Psi_N((y,\infty))=0$, the joint limit law of the two maxima. Joint vague convergence in distribution follows
from~\cite[Theorem~4.11]{kallenberg-random-measures}, applied to the pair viewed as one point process
on
$\R\times\{1,2\}$, with the dissecting semiring of the sets $(\alpha,\beta]\times\{i\}$.
Restoring the omitted points does not change these limits,
since their counts above each fixed level vanish in probability.
\end{proof}

\section{Strong matrix models for mixed \texorpdfstring{$q$}{q}-Gaussian families}
\label{sec:mixed-q-strong-models}

The strong models in
Sections~\ref{sec:scalar-models}--\ref{sec:clock-representations} have
limits that are free, possibly with amalgamation. We now construct strong matrix approximations to mixed $q$-Gaussian families, which need not be
free (Theorem~\ref{thm:mixed-q-strong-gaussian-models}). This construction does not use moment correction. For the
approximations, we use independent centered self-adjoint Gaussian matrices,
choosing the correlations among the entries of each matrix so that, in the limit, the
joint moments and polynomial norms are those of the prescribed mixed $q$-Gaussian
family. Hence the $C^*$-algebra generated by such a family is MF (Corollary~\ref{cor:mixed-q-spectra-mf}). For $-1<q<1$, the Brown--Douglas--Fillmore extension semigroup of the $C^*$-algebra generated by a finite $q$-Gaussian family with at least two generators is not a group (Corollary~\ref{cor:scalar-q-mf-ext}).

Let $r\ge1$ and let
$Q=(q_{ab})_{a,b=1}^r$ be real symmetric, with
\[
\rho:=\max_{a,b}\abs{q_{ab}}<1.
\]
Write $s_a^Q=\ell_a^Q+(\ell_a^Q)^*$ for the mixed $q$-Gaussian
variables of Speicher~\cite{speicher-generalized-statistics} with parameter matrix $Q$, in the deformed Fock
construction of Bo\.{z}ejko and
Speicher~\cite[Corollary~3.2(ii) and Theorems~4.3--4.4]{bozejko-speicher-coxeter}; we use the formulation in \cite[Section~2]{nou-noninjectivity}, with vacuum vector $\Omega$ and
\[
(\ell_a^Q)^*\ell_b^Q
=\delta_{ab}+q_{ab}\ell_b^Q(\ell_a^Q)^*,
\qquad
(\ell_a^Q)^*\Omega=0.
\eqtag\label{eq:mixed-q-ccr}
\]
The vacuum state $\tau_Q$ is faithful on the unital $C^*$-algebra $A_Q:=C^*(s_1^Q,\ldots,s_r^Q)$ and satisfies
$\tau_Q(s_a^Qs_b^Q)=\delta_{ab}$. Since $Q$ is real symmetric, $\tau_Q$ is tracial. Note that $Q$ need not be positive semidefinite.
Throughout this section, the degree of a polynomial is its total degree in the
self-adjoint variables, and $\operatorname{tr}_D=D^{-1}\operatorname{Tr}$. We write $\mathbf s^Q=(s_1^Q,\ldots,s_r^Q)$, and likewise for other parameter matrices and for other tuples, such as $\mathbf X^{(m)}=(X_1^{(m)},\ldots,X_r^{(m)})$. A word is a monomial $w=x_{a_1}\cdots x_{a_n}$ in noncommuting variables, $|w|=n$ is its length, and $w(\mathbf X)$ is its value at a tuple $\mathbf X$.

\begin{theorem}[Strong Gaussian models of mixed $q$-Gaussian families]
\label{thm:mixed-q-strong-gaussian-models}
Let $r\ge1$ and let $Q=(q_{ab})_{a,b=1}^r$ be a real symmetric matrix with $\max_{a,b}\abs{q_{ab}}<1$.
There are integers $D_m\uparrow\infty$ and, for each $m$, independent centered Gaussian
self-adjoint matrices $X_1^{(m)},\ldots,X_r^{(m)}\in M_{D_m}(\C)$, all defined on one probability space, such that
\[
\mathbb E[X_a^{(m)}X_b^{(m)}]=\delta_{ab}.
\eqtag\label{eq:mixed-q-matrix-variance}
\]
On one event of probability one, for every word $w$, every fixed $h\ge1$, and every
$P\in M_h(\C)\otimes\C\langle x_1,\ldots,x_r\rangle$,
\[
\begin{aligned}
\operatorname{tr}_{D_m}w(\mathbf X^{(m)})&\longrightarrow\tau_Q(w(\mathbf s^Q)),\\
\norm{P(\mathbf X^{(m)})}&\longrightarrow\norm{P(\mathbf s^Q)}.
\end{aligned}
\eqtag\label{eq:mixed-q-strong-convergence}
\]
For fixed $h\ge1$ and $d\ge0$, write $P_A=\sum_{|w|\le d}A_w\otimes w$, with $A_w\in M_h(\C)$.
On the same event,
\[
\sup_{\sum_w\norm{A_w}\le1}
\abs{\norm{P_A(\mathbf X^{(m)})}-\norm{P_A(\mathbf s^Q)}}\longrightarrow0.
\eqtag\label{eq:mixed-q-fixed-coefficient-ball}
\]
\end{theorem}

The real and imaginary parts of the entries of each matrix $X_a^{(m)}$ are jointly
Gaussian, with covariance that may be degenerate.

When $q_{ab}=q$ for all $a,b$, we obtain strong Gaussian models of the $q$-Gaussian family for every
$-1<q<1$.
Caspers and Chen~\cite[Theorems~6.11 and~7.1]{caspers-chen} obtained strong $q$-Gaussian models and the MF property for
$\abs q<\sqrt2-1$ (see Corollary~\ref{cor:mixed-q-spectra-mf}
for the definition of MF). Their theorem on the models also applies beyond this range whenever the Fock
representation of the universal $C^*$-algebra of the $q$-commutation relations is
faithful~\cite[Remark~5.6 and Theorem~6.11]{caspers-chen}.
For each fixed degree, their results also allow matrix coefficients whose
dimension exceeds that of the random matrices. We follow their graph-product argument, which builds on
\cite[Theorem~9.8]{chen-garza-vargas-van-handel}, but only for
deformations with zero diagonal.

For fixed $r\ge2$ and $\rho$ below a bound depending on $r$, free transport identifies $A_Q$ with the $C^*$-algebra generated
by $r$ free semicircular variables. This was proved for scalar $q$ by
Guionnet and Shlyakhtenko~\cite{guionnet-shlyakhtenko-free-transport}
and for mixed $Q$ by Nelson and Zeng~\cite{nelson-zeng-mixed-q}.
Combining these trace-preserving isomorphisms with norm approximation by
polynomials and the strong semicircular matrix
models~\cite{haagerup-thorbjornsen}, we obtain MF and strong matrix
approximations in this range; for scalar $q$, Caspers and Chen noted the strong approximation~\cite[Section~1.1]{caspers-chen}.
However, the resulting matrices need not be Gaussian or independent.

Moment approximations are available from the Gaussian matrix models of
\'Sniady~\cite[Theorems~1--3]{sniady-gaussian-models}, for scalar $0\le q\le1$, and from
Speicher's central limit construction with random commutation signs
\cite[Theorem~2]{speicher-generalized-statistics}, for mixed parameters.
Junge and Zeng~\cite[Section~2B and Theorem~3.1]{junge-zeng} gave a
presentation using spin matrices and their trace.
Liu and Shen~\cite[Theorems~1.1--1.2]{liu-shen-syk} obtained mixed $q$-Gaussian
limits in expected joint moments from SYK Hamiltonians with partial
interactions. Their parameters are constrained by interaction lengths,
support overlaps, and parities. With Gaussian coefficients, their
construction gives independent Gaussian matrix families, but without a
norm-convergence conclusion.

To prove Theorem~\ref{thm:mixed-q-strong-gaussian-models}, we first replace each variable index $a$ by $m$ indices $(a,i)$,
$1\le i\le m$, and form a larger mixed $q$-Gaussian family with zero
diagonal parameters. For these parameters, the Fock representation of
the universal $C^*$-algebra generated by elements satisfying the first relation in \eqref{eq:mixed-q-ccr} is faithful
\cite[Theorem~6.2]{kuzmin-pochekai}, which we use below to compare norms in
graph products. We approximate the $a$th original variable by the
normalized sum over the copy index $i$ (Theorem~\ref{thm:mixed-q-diagonal-deletion}), using a rapid-decay bound whose constant does not depend on the number of variables (Lemma~\ref{lem:mixed-q-uniform-rapid-decay}).

We then approximate each fixed family with zero diagonal parameters by normalized
sums of graph-product semicircular variables tensored with sign matrices,
i.e., self-adjoint unitaries with prescribed commutation signs.
In the reduced graph product, the unital $C^*$-algebras at adjacent vertices commute
(see the Fock construction at \eqref{eq:mixed-q-graph-creators}).
M\l{}otkowski~\cite[Theorem~4 and Corollary~3]{mlotkowski-lambda-free}
obtained scalar-$q$ moment limits from normalized sums in graph products
for $0\le q\le1$. The joint moment limit for nonnegative mixed parameters
with zero diagonal is the special case of~\cite[Theorem~5.1]{cebron-santos-youssef-graphon} in which the limit of the graphs is constant between each pair of parts.
Tensoring those sums with the fixed sign matrices gives the signed
moment limit in Theorem~\ref{thm:mixed-q-graph-limit}, where we also
prove norm convergence.

The graph-product variables are in turn approximated by GUE matrices that act on some factors of a tensor product and as the identity on the other factors. By a diagonal argument through these approximations, we obtain the matrices $X_a^{(m)}$ and the dimensions $D_m$ in Theorem~\ref{thm:mixed-q-strong-gaussian-models}.

\subsection{Reduction to zero diagonal parameters}

Set
\[
C_\rho:=\prod_{j=1}^\infty(1-\rho^j)^{-1},
\qquad K_\rho:=C_\rho^{3/2}.
\]
We use the Wick-product bound in~\cite[equation~(3) and Lemma~2]{nou-noninjectivity},
whose constant depends only on $\rho$ and is independent of the number of variables.
Bounds of this type go back to Bo\.{z}ejko~\cite{bozejko-ultracontractivity} and, for Yang--Baxter deformations,
to Kr\'olak~\cite{krolak-wick-product}. As Amrutam, Jekel, and Wasilewski noted~\cite[Remark~1.4]{amrutam-jekel-wasilewski},
Kr\'olak~\cite{krolak-contractivity} also extended Bo\.{z}ejko's rapid-decay estimate to mixed $q$-Gaussian variables.

\begin{lemma}
\label{lem:mixed-q-uniform-rapid-decay}
Let $t,h\ge1$ and $d\ge0$ be integers.  Let $R=(r_{ij})\in M_t(\R)$ be symmetric with $\max_{i,j}\abs{r_{ij}}\le\rho<1$.
For every scalar polynomial $Z$ of degree at most $d$ in the mixed $q$-Gaussian variables with parameter matrix $R$ and vacuum state $\tau_R$,
\[
\norm{Z}\le K_\rho(d+1)^{3/2}\norm{Z}_{2,\tau_R}.
\eqtag\label{eq:mixed-q-scalar-rapid-decay}
\]
For an $M_h(\C)$-valued polynomial $Z$ of degree at most $d$,
\[
\norm{Z}\le h^{3/2}K_\rho(d+1)^{3/2}
\norm{Z}_{2,\operatorname{tr}_h\otimes\tau_R}.
\eqtag\label{eq:mixed-q-matrix-rapid-decay}
\]
\end{lemma}

\begin{proof}
The operator $T_R$ defined by $T_R(e_i\otimes e_j)=r_{ij}e_j\otimes e_i$ is
self-adjoint, has norm at most $\rho<1$, and satisfies the Yang--Baxter
relation on $(\C^t)^{\otimes3}$:
\[
(T_R\otimes 1)(1\otimes T_R)(T_R\otimes 1)
=(1\otimes T_R)(T_R\otimes 1)(1\otimes T_R).
\]
Since $R$ is real symmetric, the vacuum state $\tau_R$ is tracial.
Let $H_n$ be the degree-$n$ tensor summand of the deformed Fock space $\mathcal F_R$ with parameter matrix $R$,
with $H_0=\C\Omega$, and let
$J_{a,b}:H_a\otimes H_b\to H_{a+b}$ be concatenation. Write
$\Theta_\alpha\xi=\alpha\otimes\xi$.
With $S$ the conjugate-linear tensor reversal, set
\mbox{$U_{a,b}(\alpha\otimes\beta)=\Theta_\alpha\Theta_{S\beta}^*$}.
For $\xi\in H_n$, the Wick product $W_R(\xi)$ is the polynomial in
the mixed $q$-Gaussian variables characterized by $W_R(\xi)\Omega=\xi$.
Nou~\cite[Lemma~2, equations~(8)--(9)]{nou-noninjectivity} gave the
following Wick decomposition and stated the norm estimate below:
\[
W_R(\xi)=\sum_{b=0}^nU_{n-b,b}J_{n-b,b}^*\xi,
\qquad \xi\in H_n,
\eqtag\label{eq:mixed-q-wick-decomposition}
\]
and $\norm{W_R(\xi)}\le C_\rho^{3/2}(n+1)\norm\xi$.
Here $J^*$ denotes the adjoint with respect to the deformed scalar products.

The argument for this norm estimate in~\cite[p.~23]{nou-noninjectivity}
bounds $U_{a,b}$ using the multiplication map
$B(\mathcal F_R)\otimes_{\min}B(\mathcal F_R)\to B(\mathcal F_R)$,
which is not bounded in general. Indeed, for matrix units $e_{ij}$, $1\le i,j\le k$, in $B(\mathcal F_R)$, the element $\sum_{i=1}^ke_{1i}\otimes e_{i1}$ has norm one, but the multiplication map sends it to $ke_{11}$. However, Nou also bounded $U_{a,b}$
on the Haagerup tensor product of the column and row Hilbert
spaces~\cite[Corollary~1]{nou-noninjectivity}, using the multiplication map on
the Haagerup tensor product of $B(\mathcal F_R)$ with itself, which is contractive. On $H_a\otimes H_b$, the Haagerup
norm is at most the Hilbert space norm, so Nou's Corollary~1 also gives the
estimate. We include a direct verification.
By~\cite[equation~(3)]{nou-noninjectivity}, concatenation
$J_a:H_a\otimes\mathcal F_R\to\mathcal F_R$
has norm at most $C_\rho^{1/2}$. Using $S$ in the second factor, we identify
$\eta\in H_a\otimes H_b$ with the Hilbert--Schmidt operator $T_\eta:H_b\to
H_a$. Then
$U_{a,b}(\eta)=J_a(T_\eta\otimes 1)J_b^*$, and hence
$\|U_{a,b}\|\le C_\rho$.
Since $\|J_a\|\le C_\rho^{1/2}$, we also have $\|J_{a,b}\|\le C_\rho^{1/2}$,
and with the bounds on $U_{a,b}$ and $J_{a,b}$ the estimate follows from \eqref{eq:mixed-q-wick-decomposition}.

For $Z$ of degree at most $d$, write $Z\Omega=\sum_{n=0}^d\xi_n$
in the orthogonal summands $H_n$. Since $\Omega$ is separating,
$Z=\sum_{n=0}^dW_R(\xi_n)$, so, by Cauchy--Schwarz, we have
\[
\norm Z\le K_\rho\left(\sum_{n=0}^d(n+1)^2\right)^{1/2}
\norm{Z}_{2,\tau_R}
\le K_\rho(d+1)^{3/2}\norm{Z}_{2,\tau_R}.
\]
Applying this entrywise and using
\[
\norm Z\le\sum_{i,j}\norm{Z_{ij}},\qquad
\sum_{i,j}\norm{Z_{ij}}_{2,\tau_R}
\le h^{3/2}\norm{Z}_{2,\operatorname{tr}_h\otimes\tau_R}
\]
gives \eqref{eq:mixed-q-matrix-rapid-decay}.
\end{proof}

For a single variable, the moment limit in Theorem~\ref{thm:mixed-q-diagonal-deletion} also follows from~\cite[Theorem~6.1]{cebron-santos-youssef-graphon}, applied to the parameter matrices with zero diagonal and all other entries equal to $q_{aa}$, which may be negative. However, we also need joint strong convergence.

\begin{theorem}
\label{thm:mixed-q-diagonal-deletion}
For each integer $m\ge1$, consider the mixed $q$-Gaussian
variables $s_{a,i}^{Q^{[m]}}$ indexed by
$(a,i)\in\{1,\ldots,r\}\times\{1,\ldots,m\}$, with the parameter
matrix and the averages defined by
\[
Q^{[m]}_{(a,i),(b,j)}:=
\begin{cases}
0,&(a,i)=(b,j),\\
q_{ab},&(a,i)\ne(b,j),
\end{cases}
\qquad
Y_a^{[m]}:=\frac1{\sqrt m}\sum_{i=1}^m s_{a,i}^{Q^{[m]}}.
\eqtag\label{eq:mixed-q-replicated-system}
\]
Then $\mathbf Y^{[m]}=(Y_a^{[m]})_{a=1}^r$ converges jointly in moments and strongly to $\mathbf s^Q$: for every fixed $h$ and every
$P\in M_h(\C)\otimes\C\langle x_1,\ldots,x_r\rangle$,
\[
\norm{P(\mathbf Y^{[m]})}\longrightarrow\norm{P(\mathbf s^Q)}.
\]
For every word $w$ of length $2p$, $p\ge1$,
\[
\abs{\tau_{Q^{[m]}}(w(\mathbf Y^{[m]}))-\tau_Q(w(\mathbf s^Q))}
\le\frac{(2p-1)!!\binom p2}{m}.
\eqtag\label{eq:mixed-q-replication-moments}
\]
\end{theorem}

\begin{proof}
Let $\mathcal P_2(2p)$ be the set of pair partitions of $\{1,\ldots,2p\}$,
and let $\operatorname{Cr}(\pi)$ be the set of unordered pairs of crossing blocks
of $\pi$.  The mixed Wick moment formula is due to
Speicher, who obtained it from a central limit construction~\cite[Theorems~1--3]{speicher-generalized-statistics};
see also the proof of~\cite[Theorem~4.4]{bozejko-speicher-coxeter}.
In the formulation of
\cite[Theorem~3.1 and Remark~3.17]{junge-zeng}, it reads
\[
\tau_Q(s_{a_1}^Q\cdots s_{a_{2p}}^Q)
=\sum_{\substack{\pi\in\mathcal P_2(2p)\\a_u=a_v\text{ for }\{u,v\}\in\pi}}
\prod_{\{B,C\}\in\operatorname{Cr}(\pi)}q_{a_B,a_C},
\eqtag\label{eq:mixed-q-pairing-moments}
\]
where $a_B$ is the common variable index on $B$.
Expanding $w(\mathbf Y^{[m]})$ and applying \eqref{eq:mixed-q-pairing-moments} to $Q^{[m]}$, we see that the contribution of each pairing $\pi$ that joins equal variable indices is an average over copy indices $i_B\in\{1,\ldots,m\}$, chosen independently and uniformly for the blocks $B$ of $\pi$. For each variable index
$a$,
let $G_{\pi,a}$ be the graph on the $p_a$ blocks with index $a$, with
edges joining crossing blocks. The copy indices color its vertices; the coloring is proper
when adjacent vertices have different copy indices. The product of crossing
weights for $Q^{[m]}$ equals that for $Q$ when all these colorings are proper, and is zero
otherwise. Thus $\pi$ contributes
\[
\left(\prod_{\{B,C\}\in\operatorname{Cr}(\pi)}q_{a_B,a_C}\right)
\prod_{a=1}^r\frac{\chi_{G_{\pi,a}}(m)}{m^{p_a}},
\eqtag\label{eq:mixed-q-chromatic-moments}
\]
where $\chi_G(m)$ counts proper colorings of $G$ with $m$ colors. By a union bound, the probability that some coloring is
improper is at most $\binom p2/m$. Since the product of crossing weights for $Q$ has absolute value at most one, summing over at most
$(2p-1)!!$ pairings proves \eqref{eq:mixed-q-replication-moments}.

To prove norm convergence, we apply rapid decay to powers, as in~\cite[Proposition~3.4]{miyagawa-strong-convergence}.
Fix $P$ of degree at most $d$ with coefficients in $M_h(\C)$, and write $B_m=P(\mathbf Y^{[m]})$.
By Lemma~\ref{lem:mixed-q-uniform-rapid-decay}, applied to $(B_m^*B_m)^n$ for each integer $n\ge1$, we have
\[
\norm{B_m}
\le\left[h^{3/2}K_\rho(2dn+1)^{3/2}\right]^{1/(2n)}
\left[(\operatorname{tr}_h\otimes\tau_{Q^{[m]}})((B_m^*B_m)^{2n})\right]^{1/(4n)}.
\eqtag\label{eq:mixed-q-replication-power}
\]
For fixed $n$, expanding into words and using \eqref{eq:mixed-q-replication-moments} and the vanishing of odd moments on both sides, we see that the last trace converges to $(\operatorname{tr}_h\otimes\tau_Q)((P(\mathbf s^Q)^*P(\mathbf s^Q))^{2n})$ as $m\to\infty$. We take the upper limit as $m\to\infty$ and then let $n\to\infty$.
The prefactor tends to one, and the remaining trace factor tends to
$\norm{P(\mathbf s^Q)}$ because $\operatorname{tr}_h\otimes\tau_Q$ is
faithful. Hence
$\limsup_m\norm{P(\mathbf Y^{[m]})}\le\norm{P(\mathbf s^Q)}$.
Conversely, for every fixed positive integer $n$,
\[
\left[(\operatorname{tr}_h\otimes\tau_{Q^{[m]}})((B_m^*B_m)^n)\right]^{1/(2n)}\le\norm{B_m}.
\]
Letting $m\to\infty$ and then $n\to\infty$, we obtain $\liminf_m\norm{P(\mathbf Y^{[m]})}\ge\norm{P(\mathbf s^Q)}$.
\end{proof}

From the proof, we obtain an explicit comparison for finite $m$.
For integers $d,n,m\ge1$ and a scalar polynomial $P=\sum_wc_ww$ of degree at most $d$, set
\[
\Lambda(P):=\sum_w\abs{c_w},\qquad
F_{\rho,d,n}:=\left[K_\rho(2dn+1)^{3/2}\right]^{1/(2n)},
\]
\[
E_{d,n}:=(4dn-1)!!\binom{2dn}{2}.
\]
Then
\[
\norm{P(\mathbf Y^{[m]})}
\le F_{\rho,d,n}
\left(\norm{P(\mathbf s^Q)}^{4n}+\frac{\Lambda(P)^{4n}E_{d,n}}m\right)^{1/(4n)},
\eqtag\label{eq:mixed-q-finite-replication}
\]
and the same inequality holds with the two norms interchanged.
Indeed, $(P^*P)^{2n}$ has degree at most $4dn$, and $\Lambda((P^*P)^{2n})\le\Lambda(P)^{4n}$.
By \eqref{eq:mixed-q-replication-moments} and the vanishing of odd moments, the trace of each word in the expansion of $(P^*P)^{2n}$ changes by at most $E_{d,n}/m$ when $\mathbf Y^{[m]}$ is replaced by $\mathbf s^Q$.
The two inequalities then follow by applying~\eqref{eq:mixed-q-replication-power} in the unital
$C^*$-algebra generated by
$\mathbf Y^{[m]}$ or by $\mathbf s^Q$, respectively.

\begin{corollary}
\label{cor:mixed-q-norm-continuity}
For every fixed $r,h$, every $\rho_0\in[0,1)$, and every $P\in M_h(\C)\otimes\C\langle x_1,\ldots,x_r\rangle$, the function
$R\mapsto\norm{P(\mathbf s^R)}$ is continuous on the symmetric matrices $R\in M_r(\R)$ satisfying
$\max_{a,b}\abs{r_{ab}}\le\rho_0$.
\end{corollary}

When all $r_{ab}$ equal a common $q$ and $h=1$, the conclusion is Miyagawa's continuity
theorem~\cite[Theorem~1.1]{miyagawa-strong-convergence}. We follow his argument.

\begin{proof}
For two parameter matrices $R,R'$ in this set, telescoping shows that a product of crossing weights changes by at most
the number of crossings times $\max_{a,b}\abs{r_{ab}-r'_{ab}}$.  Thus, for every word $w$ of length $2p$, $p\ge1$,
\[
\abs{\tau_R(w(\mathbf s^R))-\tau_{R'}(w(\mathbf s^{R'}))}
\le(2p-1)!!\binom p2\max_{a,b}\abs{r_{ab}-r'_{ab}}.
\]
Moreover, odd moments vanish and $\tau_R(1)=1$. Since the constant in Lemma~\ref{lem:mixed-q-uniform-rapid-decay}, with $\rho_0$ in place of $\rho$, is uniform on this parameter set, we obtain $\norm{P(\mathbf s^{R'})}\to\norm{P(\mathbf s^R)}$ as $R'\to R$ by the argument with \eqref{eq:mixed-q-replication-power} in the proof of Theorem~\ref{thm:mixed-q-diagonal-deletion}.
\end{proof}

\subsection{Graph products and Gaussian matrices}

It remains to approximate a mixed $q$-Gaussian family with zero diagonal
parameters.
Fix an integer $t\ge1$ and a real symmetric matrix
$R=(r_{ab})_{a,b=1}^t$ with $r_{aa}=0$ and
$\max_{a,b}\abs{r_{ab}}<1$. Set
\[
p_{ab}:=\abs{r_{ab}},\qquad \sigma_{ab}:=\operatorname{sgn}(r_{ab}),
\]
where $\operatorname{sgn}(0)=1$.  Let $\Gamma_k$ be a finite simple graph with parts
\[
V_a=\{(a,1),\ldots,(a,k)\}.
\]
There are no edges within a part. The adjacency matrices between parts are
$A_{ab}\in M_k(\R)$, so $A_{ba}=A_{ab}^{\mathsf T}$.  Write $J_k$ for the all-ones matrix and set
\[
\delta_k:=\max_{a<b}\frac{\norm{A_{ab}-p_{ab}J_k}}k,
\eqtag\label{eq:mixed-q-graph-discrepancy}
\]
with $\delta_k=0$ if $t=1$.

Let $\Gamma_k^+$ be the monoid of words in the vertices, with words identified under interchanges of consecutive letters joined by an edge of $\Gamma_k$.
On $\mathcal F_{\Gamma_k}:=\ell^2(\Gamma_k^+)$, with vacuum $\Omega_{\Gamma_k}=\delta_{\varnothing}$,
left multiplication defines isometries $\ell_v\delta_u=\delta_{vu}$ for vertices $v$ and words $u\in\Gamma_k^+$.  They satisfy
\[
\ell_v^*\ell_v=1,\qquad
\ell_v^*\ell_{v'}=
\begin{cases}
\ell_{v'}\ell_v^*,&\{v,v'\}\in E(\Gamma_k),\\
0,&v\ne v',\ \{v,v'\}\notin E(\Gamma_k).
\end{cases}
\eqtag\label{eq:mixed-q-graph-creators}
\]
Indeed, two distinct letters can both occur first in a word of this monoid precisely when they commute. In that case,
either first letter can be removed.

The operators $\ell_v+\ell_v^*$ are the standard semicircular generators of
the reduced graph product of copies of the unital $C^*$-algebra
$C([-2,2])$, with states given by the standard semicircular measure
$\mu_{\rm sc}$. The orthonormal polynomials are $p_n(x)=U_n(x/2)$,
where $U_n$ is the Chebyshev polynomial of the second kind. They satisfy
$p_0=1$, $p_1=x$, and $xp_n=p_{n+1}+p_{n-1}$ for $n\ge1$.
Using these bases in the vertex spaces $L^2(\mu_{\rm sc})$, we identify
the graph-product Hilbert space with $\ell^2(\Gamma_k^+)$. By the
recurrence, multiplication by $x$ at vertex $v$ corresponds to
$\ell_v+\ell_v^*$. Here $\ell_v$ adds an initial $v$, while $\ell_v^*$
removes one if the word has a representative beginning with $v$, and
is zero otherwise~\cite[Sections~2.2 and~3.1]{caspers-chen}.
The vacuum state $\psi_k$ on the unital $C^*$-algebra generated by $(\ell_v)_v$ restricts to the
graph-product trace on the unital $C^*$-subalgebra generated by $(\ell_v+\ell_v^*)_v$.

To include negative parameters, we tensor the isometries $\ell_v$ with self-adjoint unitaries that commute or anticommute according to
the signs $\sigma_{ab}$. Choose $\epsilon_{ab}\in\{0,1\}$ with
$\sigma_{ab}=(-1)^{\epsilon_{ab}}$, and set
\[
X_0:=\begin{pmatrix}0&1\\1&0\end{pmatrix},\qquad
Z_0:=\begin{pmatrix}1&0\\0&-1\end{pmatrix},
\]
\[
c_a:=\bigotimes_{b<a}Z_0^{\epsilon_{ba}}\otimes X_0\otimes\bigotimes_{b>a}1
\in M_{2^t}(\C).
\]
These self-adjoint unitaries satisfy
$c_ac_b=\sigma_{ab}c_bc_a$ because, for $a<b$, only the $a$th tensor
coordinate can contribute a minus sign.  On $\mathcal H_k:=\C^{2^t}\otimes\mathcal F_{\Gamma_k}$ define
\[
\tilde\ell_{a,i}:=c_a\otimes\ell_{a,i},\qquad
L_{a,k}:=\frac1{\sqrt k}\sum_{i=1}^k \tilde\ell_{a,i},\qquad
S_{a,k}:=L_{a,k}+L_{a,k}^*.
\eqtag\label{eq:mixed-q-graph-averages}
\]
The state $\varphi_k:=\operatorname{tr}_{2^t}\otimes\psi_k$ is tracial on the unital $C^*$-algebra generated by
$\mathbf S_k=(S_{1,k},\ldots,S_{t,k})$.

\begin{lemma}
\label{lem:mixed-q-graph-error}
For the operators in \eqref{eq:mixed-q-graph-averages}, $L_{a,k}^*L_{a,k}=1$, $\norm{S_{a,k}}\le2$, and, for $a\ne b$,
\[
\norm{L_{a,k}^*L_{b,k}-r_{ab}L_{b,k}L_{a,k}^*}
=\frac{\norm{A_{ab}-p_{ab}J_k}}k.
\eqtag\label{eq:mixed-q-graph-error}
\]
\end{lemma}

\begin{proof}
Since there are no edges within $V_a$, we have $\tilde\ell_{a,i}^*\tilde\ell_{a,j}=\delta_{ij}$.  Thus
\[
W_a:=[\tilde\ell_{a,1},\ldots,\tilde\ell_{a,k}]:\mathcal H_k^k\to\mathcal H_k
\]
is an isometry, and
$L_{a,k}^*L_{a,k}=1$.  For different parts,
\[
\tilde\ell_{a,i}^*\tilde\ell_{b,j}=\sigma_{ab}(A_{ab})_{ij}\,\tilde\ell_{b,j}\tilde\ell_{a,i}^*.
\]
Consequently,
\[
L_{a,k}^*L_{b,k}-r_{ab}L_{b,k}L_{a,k}^*
=\frac{\sigma_{ab}}k W_b\bigl[(A_{ab}-p_{ab}J_k)^{\mathsf T}\otimes 1_{\mathcal H_k}\bigr]W_a^*.
\eqtag\label{eq:mixed-q-graph-factorization}
\]
Since $W_a,W_b$ are isometries, the left side of \eqref{eq:mixed-q-graph-error} is at most the right side. Multiplying \eqref{eq:mixed-q-graph-factorization} by
$W_b^*$ on the left and by $W_a$ on the right gives the reverse inequality.
\end{proof}

\begin{theorem}
\label{thm:mixed-q-graph-limit}
Let $t$ and $R$ be as above.  If $\delta_k\to0$, then $\mathbf S_k$ converges jointly and strongly to the mixed $q$-Gaussian
family $\mathbf s^R$ with parameter matrix $R$.  More precisely, for every word $w$, every fixed $h\ge1$, and every
$P\in M_h(\C)\otimes\C\langle x_1,\ldots,x_t\rangle$,
\[
\varphi_k(w(\mathbf S_k))\longrightarrow\tau_R(w(\mathbf s^R)),\qquad
\norm{P(\mathbf S_k)}\longrightarrow\norm{P(\mathbf s^R)}.
\]
There are graphs satisfying $\delta_k\le4/\sqrt k$ for every sufficiently large $k$.
\end{theorem}

\begin{proof}
We pass to a norm ultraproduct of the unital $C^*$-algebras $B(\mathcal H_k)$ along any nonprincipal ultrafilter on
$\N$, and let $L_a$ be the image of $(L_{a,k})_k$. By
Lemma~\ref{lem:mixed-q-graph-error} and $\delta_k\to0$, these operators
satisfy
\[
L_a^*L_a=1,\qquad L_a^*L_b=r_{ab}L_bL_a^*\quad(a\ne b).
\]
For zero-diagonal $R$, the Fock creation operators $\ell_a^R$
faithfully represent the universal unital $C^*$-algebra of these isometry
relations~\cite[Theorem~6.2]{kuzmin-pochekai}. Thus the assignment
$\ell_a^R\mapsto L_a$ extends to a unital $^*$-homomorphism from
$C^*(\ell_1^R,\ldots,\ell_t^R)$ into the norm ultraproduct. Since this
$^*$-homomorphism is contractive at each matrix level and the ultrafilter is
arbitrary, we obtain
\[
\limsup_k\|P(\mathbf S_k)\|\le\|P(\mathbf s^R)\|.
\]

In the GNS representation of the ultralimit of the states $\varphi_k$, the cyclic vector
is annihilated by every $L_a^*$. Indeed,
$L_{a,k}^*(\xi\otimes\Omega_{\Gamma_k})=0$ for every $\xi\in\C^{2^t}$,
so $\varphi_k(L_{a,k}L_{a,k}^*)=0$ for every $k$.
By applying the relations $L_a^*L_a=1$ and $L_a^*L_b=r_{ab}L_bL_a^*$ ($a\ne b$) repeatedly to adjacent factors $L_a^*L_b$, we write every word in the $L_a$ and $L_a^*$ as a scalar multiple of a word in which every $L_a$ stands to the left of every $L_b^*$. Among these words, only the empty word has a nonzero value under the ultralimit of the states $\varphi_k$. Since the operators $\ell_a^R$ satisfy the same relations and $(\ell_a^R)^*\Omega=0$, the ultralimit of $\varphi_k(w(\mathbf S_k))$ is $\tau_R(w(\mathbf s^R))$ for every word $w$. These limits are the same for every
nonprincipal ultrafilter, so the moments converge along the ordinary
sequence. As at the end of the proof of Theorem~\ref{thm:mixed-q-diagonal-deletion}, since the moments converge and $\tau_R$ is faithful, we have $\liminf_k\norm{P(\mathbf S_k)}\ge\norm{P(\mathbf s^R)}$.

To choose the graphs, we include each edge between $V_a$ and $V_b$, $a<b$,
independently with probability $p_{ab}$. For $B=A_{ab}-p_{ab}J_k$,
Hoeffding's inequality and two $1/4$-nets of the Euclidean unit sphere
in $\R^k$, each of size at most $9^k$, give
\[
\mathbb P(|x^{\mathsf T}By|\ge u)\le2e^{-2u^2}
\quad(\|x\|=\|y\|=1),\qquad
\mathbb P(\|B\|>4\sqrt k)
\le2e^{-(8-2\log9)k}\le2e^{-3k}.
\]
The probability that $\|A_{ab}-p_{ab}J_k\|>4\sqrt k$ for at least one pair $a<b$ is at most $t(t-1)e^{-3k}<1$ for all sufficiently large $k$, with $t$ fixed. Hence there is a graph with $\delta_k\le4/\sqrt k$.
\end{proof}

\begin{proof}[Proof of Theorem~\ref{thm:mixed-q-strong-gaussian-models}]
\emph{Tensor GUE models of graph products.}
We approximate the semicircular variables in each fixed graph product by matrices, using strong convergence for tensor GUE
matrices~\cite[Section~9.4 and Theorem~9.8]{chen-garza-vargas-van-handel}.
Fix $L\ge1$, a finite set of vertices and, for each vertex $v$, a nonempty subset $K_v\subseteq\{1,\ldots,L\}$ of the
tensor factors of $(\C^n)^{\otimes L}$.
Let $H_v^{(n)}$ be independent normalized GUE matrices of
size $n^{|K_v|}$, and let $\widehat G_v^{(n)}$ act as $H_v^{(n)}$ on
the factors in $K_v$ and as the identity on the others. With the supports $K_v$ fixed, the norm of every fixed polynomial with matrix coefficients in this tuple converges almost surely to its norm at the graph-product semicircular family of the graph in which $v$ and $v'$ are adjacent exactly when $K_v\cap K_{v'}=\varnothing$.

Morampudi and Laumann derived the limiting mixed moments for local
Haar-conjugation models, including semicircular
marginals~\cite[Sections~III.B and~IV]{morampudi-laumann}.
For convergence of expected normalized word traces in the tensor GUE
model used here, see
Charlesworth and Collins~\cite[Theorem~4]{charlesworth-collins}, stated there for Haar-conjugated
families (the GUE form follows by unitary invariance),
and~\cite[Section~9.4]{chen-garza-vargas-van-handel}.

Fix a word $w$ and set
$F=\operatorname{tr}_{n^L}w(\widehat{\mathbf G}^{(n)})$ and
$d_v=n^{|K_v|}$. For each occurrence of $v$ in $w$, we form the product of
the remaining matrices in cyclic order, starting immediately after that
occurrence, and take its normalized partial trace over the tensor factors
outside $K_v$. Let $C_v\in M_{d_v}(\C)$ be the sum of
these partial traces. Then the derivative of $F$ with respect to $H_v^{(n)}$ in a direction $B$ is $\operatorname{tr}_{d_v}(C_vB)$.
The coordinates of $H_v^{(n)}$ in a self-adjoint basis of $M_{d_v}(\C)$ orthonormal for the unnormalized Hilbert--Schmidt inner product are independent centered Gaussian variables with variance $d_v^{-1}$.
By the Gaussian Poincar\'e inequality and Parseval's identity, we have, for a constant $C(w)$ independent of $n$,
\[
\operatorname{Var}F
\le\sum_v d_v^{-2}\mathbb E\norm{C_v}_{2,\operatorname{tr}_{d_v}}^2
\le C(w)n^{-2}.
\]
Indeed, using the $L^2$ contraction of the normalized partial trace and then
H\"older's inequality, we reduce the bound to moments
$\mathbb E\operatorname{tr}_{d_v}((H_v^{(n)})^{2p})$,
where $p$ depends only on $w$. These moments are bounded uniformly in
$n$ by \eqref{eq:relgauss-genus-majorant}, established in the proof of
Lemma~\ref{lem:relgauss-raw-moments}; together with $d_v\ge n$,
this gives the estimate. For complex $F$, we apply the Gaussian Poincar\'e
inequality to its real and imaginary parts. By Chebyshev's inequality and the Borel--Cantelli lemma, the traces converge almost surely to the limits of their expectations. Together with the convergence of the norms, this shows that the tuple converges strongly almost surely.

Every finite simple graph with vertex set $V$ arises in this way from subsets $K_v$ of $L=|V|+\binom{|V|}2$ coordinates, as in
Caspers and Chen~\cite[Section~6.3]{caspers-chen} (compare Model~A of Charlesworth and
Collins~\cite[p.~587]{charlesworth-collins}): one coordinate
for each vertex and one for each unordered pair of vertices.  The support $K_v$ contains the coordinate of $v$ and the coordinate of each pair
$\{v,v'\}$ that is not an edge.  Two supports are disjoint exactly when their vertices are adjacent.

\emph{The Gaussian matrices.}
At stage $m$, we use the construction of \eqref{eq:mixed-q-graph-averages} with the $rm$ variables $(a,i)$ and the parameter matrix $R=Q^{[m]}$, which has zero diagonal, so that a graph with $k_m$ vertices per part has the $rmk_m$ vertices $(a,i,j)$, $1\le j\le k_m$.
For such a graph, set
\[
L_m=rmk_m+\binom{rmk_m}{2}.
\]
Using the sign matrices $c_{a,i}\in M_{2^{rm}}(\C)$ for $Q^{[m]}$ and independent GUE matrices acting on the factors of $(\C^{n_m})^{\otimes L_m}$ in the supports $K_{(a,i,j)}$ just described for this graph, with the tensor GUE size $n_m$ chosen below, set
\[
X_a^{(m)}=\frac1{\sqrt{mk_m}}\sum_{i=1}^m\sum_{j=1}^{k_m}
c_{a,i}\otimes\widehat G_{a,i,j}^{(n_m)},
\qquad D_m=2^{rm}n_m^{L_m}.
\eqtag\label{eq:mixed-q-finite-ensemble}
\] For distinct $a$, the
sums defining $X_a^{(m)}$ involve disjoint collections of independent
GUE matrices, so the $X_a^{(m)}$ are independent. Since $c_{a,i}^2=1$, each $\widehat G_{a,i,j}^{(n_m)}$ satisfies $\mathbb E[(\widehat G_{a,i,j}^{(n_m)})^2]=1$, and the GUE matrices are independent and centered, we have $\mathbb E[X_a^{(m)}X_b^{(m)}]=\delta_{ab}$.

\emph{The diagonal argument.}
We enumerate the traces of words and the norms of matrix polynomials in the
original $r$ variables, with coefficient entries in $\mathbb Q+i\mathbb Q$,
allowing every finite matrix size and including the generator norms. At stage $m$, we evaluate these traces and norms at the $r$ normalized sums of the $rm$ variables, as in \eqref{eq:mixed-q-replicated-system}.
By Theorem~\ref{thm:mixed-q-graph-limit}, there is a graph with $k_m\ge m$
for which the first $m$ of these traces and norms, evaluated at the graph-product sums, differ from their values
at $\mathbf Y^{[m]}$ by at most $1/m$.

We then keep the graph, the matrices $c_{a,i}\in M_{2^{rm}}(\C)$, and the sets
$K_v$ fixed while choosing the tensor GUE size $n_m$, and evaluate the same traces and norms at the Gaussian sums in
\eqref{eq:mixed-q-finite-ensemble}.
By~\cite[Theorem~9.8]{chen-garza-vargas-van-handel} and the trace
convergence proved above, for every sufficiently large $n$, these traces
and norms differ from their values at the corresponding graph-product
sums by at most $1/m$, outside an event of probability at most $2^{-m}$. We choose such an $n_m$, large enough that $D_m$
is strictly increasing.

On a product probability space, by the Borel--Cantelli lemma and Theorem~\ref{thm:mixed-q-diagonal-deletion}, all the enumerated traces and norms converge on one event of probability one.
Since the generator norms converge, there is $M$ with $\norm{X_a^{(m)}}\le M$ for all $a$ and all large $m$, and $\norm{s_a^Q}\le M$ for all $a$. On this event,
\[
\|P_A(\mathbf X^{(m)})-P_B(\mathbf X^{(m)})\|
\le\max(1,M)^d\sum_{|w|\le d}\|A_w-B_w\|,
\]
and likewise in the limit. Approximating the coefficient entries by elements
of $\mathbb Q+i\mathbb Q$ therefore gives norm convergence for every
fixed matrix polynomial. Finally, by the compactness of $\{(A_w):\sum_{|w|\le
d}\norm{A_w}\le1\}$ and the same Lipschitz bound, we obtain
\eqref{eq:mixed-q-fixed-coefficient-ball} for every fixed $h,d$, on this same event.
\end{proof}

\subsection{MF approximation and covariance maps}

\begin{corollary}[Spectra and the MF property]
\label{cor:mixed-q-spectra-mf}
Let $Q$, $\mathbf s^Q$, and $\mathbf X^{(m)}$ be as in Theorem~\ref{thm:mixed-q-strong-gaussian-models}.
On the same probability-one event, for every fixed $h\ge1$ and every $P\in M_h(\C)\otimes\C\langle x_1,\ldots,x_r\rangle$ with $P^*=P$,
\[
d_{\rm H}\bigl(\operatorname{spec}P(\mathbf X^{(m)}),\operatorname{spec}P(\mathbf s^Q)\bigr)\longrightarrow0,
\]
where $d_{\rm H}$ denotes Hausdorff distance.  For an arbitrary $P\in M_h(\C)\otimes\C\langle x_1,\ldots,x_r\rangle$,
the spectra of $\abs{P(\mathbf X^{(m)})}:=(P(\mathbf X^{(m)})^*P(\mathbf X^{(m)}))^{1/2}$ converge in Hausdorff distance
to the spectrum of $\abs{P(\mathbf s^Q)}$.  If $P(\mathbf s^Q)$ is invertible, then $P(\mathbf X^{(m)})$ is eventually
invertible and
\[
\norm{P(\mathbf X^{(m)})^{-1}}\longrightarrow\norm{P(\mathbf s^Q)^{-1}}.
\]

Moreover, the unital $C^*$-algebra $A_Q$ is MF: it admits a unital injective $^*$-homomorphism
\[
A_Q\longrightarrow
\prod_m M_{D_m}(\C)\Big/\bigoplus_m M_{D_m}(\C),
\]
where the product consists of bounded sequences and the direct sum consists of sequences tending to zero in operator norm.
\end{corollary}

\begin{proof}
We work on the probability-one event of
Theorem~\ref{thm:mixed-q-strong-gaussian-models}. The assertions about
spectra and inverse norms follow from strong convergence and
continuous functional calculus, applied to $P$ when it is self-adjoint
and to $P^*P$ in general.

By strong convergence, the map
\[
P(\mathbf s^Q)\longmapsto[(P(\mathbf X^{(m)}))_m]
\]
is an isometric unital $^*$-homomorphism on the norm-dense polynomial $^*$-subalgebra of $A_Q$.
Its extension to $A_Q$ is the required MF embedding.
\end{proof}

When $q_{ab}=0$ for all $a$ and $b$ in different groups of variables, the corresponding subfamilies of $\mathbf s^Q$ are free, and we apply Theorem~\ref{thm:multi-gns-rank-strong-exactification} to the models of Theorem~\ref{thm:mixed-q-strong-gaussian-models}.

\begin{corollary}
\label{cor:mixed-q-exact-free-blocks}
Let $Q$ and $\mathbf X^{(m)}$ be as in
Theorem~\ref{thm:mixed-q-strong-gaussian-models}. Suppose that
$\{1,\ldots,r\}=I_1\sqcup\cdots\sqcup I_t$, where $t\ge2$ and every $I_k$
is nonempty, and that $q_{ab}=0$ whenever $a$ and $b$ lie in different
sets $I_k$. On the probability-one event of that theorem, after deleting
finitely many $m$, there are unitaries $u_{k,m}\in U(D_m)$, with
$u_{t,m}=1$ and $\max_k\norm{u_{k,m}-1}_\infty\to0$, and nondecreasing
odd integers $d_m\to\infty$ with the following properties.

The tuples
$\bigl(u_{k,m}X_a^{(m)}u_{k,m}^*\bigr)_{a\in I_k}$, $1\le k\le t$,
converge jointly and strongly to $\mathbf s^Q$. On every polynomial of
degree at most $d_m$, their joint normalized trace agrees with the state
of the reduced free product of their $t$ marginal $C^*$-probability
spaces.
\end{corollary}

\begin{proof}
For $a\in I_k$ and $b\in I_l$ with $k\ne l$,
\eqref{eq:mixed-q-ccr} gives $(\ell_a^Q)^*\ell_b^Q=0$. In the deformed
inner product on the $n$-particle space, a permutation of the tensor
positions contributes zero if it reverses the order of two positions
whose indices lie in different sets. Hence concatenation identifies the deformed Fock space with the
free-product Hilbert space of the Fock spaces of the subfamilies. Indeed, tensors split
into maximal runs from the sets $I_k$, and the inner product factors over
these runs. Creation and annihilation act as the corresponding
free-product operators, since $q_{ab}=0$ for $a$ and $b$ in different sets. Thus the families
$(s_a^Q)_{a\in I_k}$, $1\le k\le t$, are free.

Since $\max_{a,b}|q_{ab}|<1$, the deformed inner product is strictly
positive on every algebraic tensor power~\cite[Theorem~2.3]{bozejko-speicher-coxeter}.
A word in the variables sends the vacuum vector to the
corresponding elementary tensor plus vectors in tensor powers of lower
degree. Hence, in each family, the words of length at most $d$ are
linearly independent in $L^2(A_Q,\tau_Q)$, and their Gram matrix is
positive definite. By the trace convergence in \eqref{eq:mixed-q-strong-convergence}, the Gram matrices of the corresponding words in the Gaussian matrices converge to these Gram matrices, so these matrix words are linearly independent for all sufficiently large $m$.
Both sides of \eqref{eq:gns-rank-stability} then equal the number of
words of length at most $d$, for every family and every $d$.
The result now follows from Theorem~\ref{thm:multi-gns-rank-strong-exactification}.
\end{proof}

Jekel, Lee, Nelson, and Pi established convergence theorems for Gaussian matrix
models with covariance maps whose Choi matrices $T_m$ have small
norm~\cite[Proposition~1.4 and Definition~1.5]{jekel-lee-nelson-pi}.
The definition of $T_m$ is recalled in Proposition~\ref{prop:mixed-q-choi-obstruction}, where $E_{ab}\in M_r(\C)$ and $E_{uv}\in M_{D_m}(\C)$ are the standard matrix units.
Their weak-convergence theorem~\cite[Theorem~3.15]{jekel-lee-nelson-pi}
assumes $\norm{T_m}\to0$, and their strong-convergence
theorem~\cite[Theorem~4.12]{jekel-lee-nelson-pi} assumes
$(\log D_m)^3\norm{T_m}\to0$. By a fourth-moment calculation in Proposition~\ref{prop:mixed-q-choi-obstruction}, when $Q\ne0$, the models of Theorem~\ref{thm:mixed-q-strong-gaussian-models} satisfy neither hypothesis. The estimate for the crossing term in the proof of Proposition~\ref{prop:mixed-q-choi-obstruction} is a special case of the bound of
Bandeira, Boedihardjo, and van Handel on the matrix alignment parameter~\cite[Proposition~4.6]{bandeira-boedihardjo-van-handel};
we include a direct proof for the present normalization.

\begin{proposition}
\label{prop:mixed-q-choi-obstruction}
For the Gaussian models in Theorem~\ref{thm:mixed-q-strong-gaussian-models}, set
\[
\eta_{ab}^{(m)}(B):=\mathbb E[X_a^{(m)}BX_b^{(m)}],
\qquad
T_m:=\sum_{a,b=1}^r\sum_{u,v=1}^{D_m}
E_{ab}\otimes\eta_{ab}^{(m)}(E_{uv})\otimes E_{uv}.
\]
Then
\[
\liminf_m\norm{T_m}\ge\rho:=\max_{a,b}\abs{q_{ab}}.
\]
\end{proposition}

\begin{proof}
Fix $m$, write $D=D_m$ and $T=T_m$, and suppress the superscript. We express each matrix as
\[
X_a=\sum_\alpha\xi_{a\alpha}A_{a\alpha},
\qquad A_{a\alpha}=A_{a\alpha}^*,
\qquad \sum_\alpha A_{a\alpha}^2=1,
\]
with independent standard real Gaussian coefficients. The last identity follows from \eqref{eq:mixed-q-matrix-variance}.
With $\operatorname{vec}(A):=\sum_{i,j}A_{ij}e_i\otimes e_j$, define
\[
C_a:=\sum_\alpha\operatorname{vec}(A_{a\alpha})\operatorname{vec}(A_{a\alpha})^*.
\]
By independence, $T=\operatorname{diag}(C_1,\ldots,C_r)$. Moreover,
$C_a\ge0$ and $\operatorname{Tr}C_a=D$.  Hence
\[
\norm{C_a}_{\rm HS}^2\le D\norm{C_a},
\qquad \norm{C_a}_{\rm HS}^2:=\operatorname{Tr}(C_a^*C_a).
\]
The crossing pairing in
$\mathbb E\operatorname{tr}_D(X_aX_bX_aX_b)$ contributes
\[
S_{ab}:=\frac1D\sum_{\alpha,\beta}
\operatorname{Tr}(A_{a\alpha}A_{b\beta}A_{a\alpha}A_{b\beta}),
\]
which satisfies
\[
DS_{ab}=\sum_{i,j,k,l}
(C_a)_{(i,j),(l,k)}(C_b)_{(j,k),(i,l)}.
\]
Indeed, $(C_a)_{(i,j),(l,k)}=\sum_\alpha(A_{a\alpha})_{ij}\overline{(A_{a\alpha})_{lk}}$, and $\overline{(A_{a\alpha})_{lk}}=(A_{a\alpha})_{kl}$ since $A_{a\alpha}$ is self-adjoint.
By Cauchy--Schwarz for the sum over these four indices, we have
\[
\abs{S_{ab}}
\le D^{-1}\norm{C_a}_{\rm HS}\norm{C_b}_{\rm HS}
\le\sqrt{\norm{C_a}\norm{C_b}}
\le\norm{T}.
\]
For $a\ne b$, independence yields
$\mathbb E\operatorname{tr}_D(X_aX_bX_aX_b)=S_{ab}$.
For $a=b$, each of the two noncrossing Wick pairings contributes one, so
$\mathbb E\operatorname{tr}_D(X_a^4)=2+S_{aa}$.

To pass from almost-sure convergence of these word traces to
convergence of their expectations, we prove uniform integrability. Let $X=\sum_\alpha\xi_\alpha A_\alpha$ with
$\sum_\alpha A_\alpha^2=1$, and set
$\eta(B)=\sum_\alpha A_\alpha BA_\alpha$.
By Gaussian integration by parts, H\"older's inequality, and the fact that the unital trace-preserving completely positive map $\eta$ is a contraction for every Schatten norm, we have, for integers $p\ge2$,
\[
\begin{aligned}
\mathbb E\operatorname{tr}_D(X^{2p})
&=\sum_{j=0}^{2p-2}\mathbb E\operatorname{tr}_D
\bigl(\eta(X^j)X^{2p-2-j}\bigr)\\
&\le(2p-1)\mathbb E\operatorname{tr}_D(X^{2p-2}).
\end{aligned}
\]
Since $\mathbb E\operatorname{tr}_D(X^2)=1$, induction yields
$\mathbb E\operatorname{tr}_D(X^{2p})\le(2p-1)!!$. Thus
\[
\mathbb E(\operatorname{tr}_D(X_a^4))^2\le105,\qquad
\mathbb E\abs{\operatorname{tr}_D(X_aX_bX_aX_b)}^2\le9\quad(a\ne b),
\]
by H\"older's inequality for the normalized trace and, for the second
bound, independence.
These second-moment bounds give uniform integrability as $m$ varies.
The almost-sure limits supplied by
Theorem~\ref{thm:mixed-q-strong-gaussian-models}, together with
\[
\tau_Q((s_a^Q)^4)=2+q_{aa},
\qquad
\tau_Q(s_a^Qs_b^Qs_a^Qs_b^Q)=q_{ab}\quad(a\ne b),
\]
imply $S_{ab}^{(m)}\to q_{ab}$. Since $\abs{S_{ab}^{(m)}}\le\norm{T_m}$, the result follows.
\end{proof}

The nonquasidiagonality argument of Caspers and Chen
\cite[Section~1.3 and proof of Theorem~7.2]{caspers-chen} applies to the
scalar $q$-Gaussian $C^*$-algebra on $r\ge2$ generators for every
$-1<q<1$. With their MF models, this argument shows that
$\operatorname{Ext}$ is not a group when $\abs q<\sqrt2-1$. By combining
their argument with Corollary~\ref{cor:mixed-q-spectra-mf}, we obtain that
conclusion throughout $-1<q<1$.

\begin{corollary}[MF and $\operatorname{Ext}$ for $q$-Gaussian algebras]
\label{cor:scalar-q-mf-ext}
Let $r\ge2$ be an integer and $-1<q<1$.  For the $q$-Gaussian family $s_1^q,\ldots,s_r^q$, i.e., the family $\mathbf s^Q$ with $q_{ab}=q$ for all $a,b$,
the unital $C^*$-algebra
\[
A_q:=C^*(s_1^q,\ldots,s_r^q)
\]
is MF and is not quasidiagonal.  Its Brown--Douglas--Fillmore extension semigroup $\operatorname{Ext}(A_q)$ is not a group.
\end{corollary}

\begin{proof}
By Corollary~\ref{cor:mixed-q-spectra-mf}, $A_q$ is MF.
Kennedy and Nica~\cite[Corollary~1.4]{kennedy-nica-exactness} proved that
the $C^*$-algebra generated by the creation operators $\ell_a^Q$, where $q_{ab}=q$ for all $a,b$, is exact for every $-1<q<1$.
Its $C^*$-subalgebra $A_q$ is therefore exact.
For $r\ge2$, the vacuum trace $\tau_q$ is the unique trace on $A_q$
\cite[Theorem~C and its proof]{amrutam-jekel-wasilewski}.
Note that the bound \eqref{eq:mixed-q-scalar-rapid-decay} satisfies the growth
condition used in~\cite[proof of Theorem~A and Remark~2.2]{amrutam-jekel-wasilewski}, that the bound at degree $n^2$, raised to the power $1/n$, tends to $1$. Indeed,
$[K_{|q|}(n^2+1)^{3/2}]^{1/n}\to1$. Hence their argument applies with this bound.

If $A_q$ were quasidiagonal, then a weak-$*$ limit of the normalized matrix traces composed with its quasidiagonal approximations would be an amenable trace, so the unique trace $\tau_q$ would be amenable
\cite[Theorem~4.2 and Proposition~6.1]{brown-quasidiagonal},
\cite[Theorem~3.1.6]{brown-invariant-means}.
Since $A_q$ is exact, it is locally reflexive, by a theorem of
Kirchberg~\cite[p.~42 and the remarks after Proposition~7.1, pp.~71--72]{kirchberg-exactness}.
The GNS von Neumann
algebra of $\tau_q$ would therefore be hyperfinite
\cite[Theorem~3.1.6 and Corollary~4.3.4]{brown-invariant-means}, hence injective.
But this algebra is the $q$-Gaussian von Neumann algebra, which is not injective
\cite[Theorem~2]{nou-noninjectivity}.
Thus $A_q$ is not quasidiagonal. Since an MF $C^*$-algebra whose extension
semigroup is a group is quasidiagonal
\cite[Corollary~13.5]{brown-quasidiagonal},
$\operatorname{Ext}(A_q)$ is not a group.
\end{proof}

\appendix
\section{A second-order proof of partition stability}
\label{app:second-order-geometry}

We give another proof of
Corollary~\ref{thm:faithful-relative-partition-stability}. In the main text, we first conjugate so that the derivative of the moment map is surjective and the defect is still small. In this step, we use Theorem~\ref{thm:finite-order-free-positioning}, which relies on
Popa's asymptotic freeness theorem~\cite[Theorem~2.1]{popa-free-independent-1995},
in the formulation of~\cite[Corollary~0.2]{popa-independence}. Here we do not use that theorem. We allow
that derivative not to be surjective and use the second derivative of the moment map to correct the components outside the range of the
first derivative. We use the identities $\tau(p_iq_j)=\tau(p_i)\tau(q_j)$ in \eqref{eq:hessian-common-corner-mean}, so this proof does not apply to other prescribed overlap matrices.

We argue by contradiction. In a tracial ultraproduct, the limiting pair is independent, and we correct the approximating pairs by Lemma~\ref{lem:second-order-inverse}, with vectors built from the infinite-dimensional spaces $p_rMp_s$ in a ${\rm II}_1$ factor (Lemmas~\ref{lem:partition-hessian-indefiniteness} and~\ref{lem:partition-quadratic-loop}). Together with the
induction in the proof of Theorem~\ref{thm:finite-partition-stability}, which removes the lower bound on projection traces,
this also gives another proof of that theorem.

\begin{lemma}\label{lem:second-order-inverse}
Let $X,Y$ be finite-dimensional real normed spaces, let $V\subseteq X$ be an open neighborhood of $0$, and let
$f_n,f:V\to Y$ be $C^3$ maps. Assume that $f(0)=0$, that the values and first two derivatives at zero converge,
\begin{equation}
\norm{f_n(0)-f(0)}
+\norm{Df_n(0)-Df(0)}
+\norm{D^2f_n(0)-D^2f(0)}
\longrightarrow0,
\label{eq:second-order-jet-convergence}
\end{equation}
and that
\begin{equation}
\sup_n\sup_{x\in V}\norm{D^3f_n(x)}<\infty.
\label{eq:second-order-third-derivative-bound}
\end{equation}
Set $A:=Df(0)$ and $Y_0:=\operatorname{ran}A$, choose a fixed decomposition
$Y=Y_0\oplus R$, and denote its projections by $P_0,P_R$.

If $R=0$, then there are $x_n\to0$ with $f_n(x_n)=0$ for all sufficiently large $n$. If $R\ne0$, the same conclusion holds
provided there are vectors $k_1,\ldots,k_N\in\ker A$ and numbers $\lambda_1,\ldots,\lambda_N>0$ such that, with
\[
z_m:=P_RD^2f(0)[k_m,k_m],
\]
we have \begin{equation} \begin{aligned} P_RD^2f(0)[k_m,k_\ell]&=0 &&(m\ne\ell),\\
\operatorname{span}\{z_1,\ldots,z_N\}&=R,
&
\sum_m\lambda_mz_m&=0.
\end{aligned}
\label{eq:second-order-positive-loop-hypotheses}
\end{equation}
\end{lemma}

When $f_n=f$ for all $n$ and $R\ne0$, the vector $k_0=\sum_m\sqrt{\lambda_m}\,k_m$ satisfies the hypotheses of a theorem of Avakov~\cite[Theorem~1]{avakov-extremum-conditions}, who proved that $f(tk_0+r(t))=f(0)$ for some $r(t)=o(t)$. Tret'yakov~\cite{tretyakov-p-th-order} had proved the corresponding description of the tangent cone when the first derivative vanishes and every nonzero direction of the cone is regular; see also Izmailov and Solodov~\cite[Theorem~2.1]{izmailov-solodov-2-regularity}. In Lemma~\ref{lem:second-order-inverse}, the maps $f_n$ vary, so we give a proof.

\begin{proof}
We work on a fixed smaller ball with closure in $V$. If $R=0$, we choose a
complement $E$ of $\ker A$. The maps $Df_n(0)|_E$ have uniformly bounded
inverses, and by \eqref{eq:second-order-jet-convergence} and \eqref{eq:second-order-third-derivative-bound}, the second derivatives $D^2f_n$ are uniformly bounded on this ball.
Hence, by Lemma~\ref{lem:finite-dimensional-correction} applied to the restriction of $f_n$ to $E$, with $T=Df_n(0)|_E$ and $J=-(Df_n(0)|_E)^{-1}$, there are $x_n$ with $f_n(x_n)=0$ and $\norm{x_n}=O(\norm{f_n(0)})$.

Assume $R\ne0$, set $H=D^2f(0)$ and $k_0=\sum_m\sqrt{\lambda_m}\,k_m$, and choose $E$ so that
$A|_E:E\to Y_0$ is an isomorphism. By the definitions of $k_0$ and $z_m$ and the first identity in \eqref{eq:second-order-positive-loop-hypotheses}, we have
\begin{equation}
P_RH[k_0,k_m]=\sqrt{\lambda_m}\,z_m.
\label{eq:second-order-loop-cross-terms}
\end{equation}
Define
\begin{equation}
e_0=-\tfrac12(A|_E)^{-1}P_0H[k_0,k_0],\qquad
x(t)=tk_0+t^2e_0.
\label{eq:second-order-trial-curve}
\end{equation}
By the choice of $e_0$, the $Y_0$ component of the coefficient of $t^2$ in the Taylor expansion of $f(x(t))$ vanishes. The $R$ component of this coefficient is $\frac12\sum_m\lambda_mz_m=0$. Since $f(0)=0$ and $k_0\in\ker A$, we obtain
\begin{equation}
f(x(t))=O(t^3).
\label{eq:second-order-cubic-residual}
\end{equation}
Let $\epsilon_n\downarrow0$ be such that the left side of \eqref{eq:second-order-jet-convergence} is at most $\epsilon_n$, and set
\begin{equation}
t_n=\max\{\epsilon_n^{1/3},n^{-1}\}.
\label{eq:second-order-scale-choice}
\end{equation}
Then $\epsilon_n=o(t_n^2)$, and by the uniform third-derivative bound, we have
\begin{equation}
\eta_n:=\norm{f_n(x(t_n))}=o(t_n^2).
\label{eq:second-order-perturbed-residual}
\end{equation}
Set
\begin{equation}
\Sigma_t=P_0+t^{-1}P_R.
\label{eq:second-order-scaled-derivative}
\end{equation}
On the subspace $X_0=E\oplus\operatorname{span}_{\R}\{k_m\}$ of $X$,
Taylor expansion gives
\begin{equation}
\Sigma_tDf(x(t))[e+k]\longrightarrow
\Gamma_0(e+k):=Ae+P_RH[k_0,e+k].
\label{eq:second-order-surjective-limit}
\end{equation}
The map $\Gamma_0$ is onto. Indeed, $Ae$ ranges over $Y_0$ as $e$ ranges over $E$, and, by \eqref{eq:second-order-loop-cross-terms}, the vectors $P_RH[k_0,k_m]$ span $R$. The rescaled derivatives
$T_n=\Sigma_{t_n}Df_n(x(t_n))|_{X_0}$ have the same limit, with error
$O(t_n+\epsilon_n/t_n)$. We choose uniformly bounded right inverses $J_n$ of $T_n$ and set
\begin{equation}
S_n=J_n\Sigma_{t_n}.
\label{eq:second-order-right-inverse}
\end{equation}
Then
\begin{equation}
Df_n(x(t_n))S_n=I_Y,\qquad \norm{S_n}=O(t_n^{-1}).
\label{eq:second-order-inverse-bound}
\end{equation}
The second derivatives $D^2f_n$ are uniformly bounded on the chosen ball, and
$\norm{S_n}^2\eta_n=o(1)$. Hence, by Lemma~\ref{lem:finite-dimensional-correction} applied to $h\mapsto f_n(x(t_n)+h)$ on $X_0$, with $T=Df_n(x(t_n))|_{X_0}$ and $J=-S_n$,
there are $h_n\in X_0$ such that
\[
f_n(x(t_n)+h_n)=0,\qquad
\norm{h_n}=O(\norm{S_n}\eta_n)=o(t_n).
\]
Finally, the ball centered at $x(t_n)$ on which we apply that lemma lies in the chosen ball for large $n$, and
$x_n=x(t_n)+h_n\to0$, as desired.
\end{proof}

To obtain the vectors required in Lemma~\ref{lem:second-order-inverse}, we compute the first two derivatives at zero of the moment map of an independent pair.
Fix partitions
$\mathbf p=(p_i)_{i=1}^a$ and $\mathbf q=(q_j)_{j=1}^b$ of $1$ in a
${\rm II}_1$ factor $(M,\tau)$, with
\[
\alpha_i=\tau(p_i)>0,
\qquad
\beta_j=\tau(q_j)>0.
\]
Assume that $\tau(p_iq_j)=\alpha_i\beta_j$ for all $i,j$.
Set $Y:=\R^{(a-1)(b-1)}$ and define the real analytic map
\begin{equation}
F(h)_{ij}
:=\tau(e^{ih}p_i e^{-ih}q_j)-\alpha_i\beta_j,
\qquad i<a,\ j<b,
\label{eq:faithful-covariance-map}
\end{equation}
for $h=h^*$ near $0$. Differentiating conjugation gives
\begin{equation}
DF(0)[h]_{ij}
=\tau(i[h,p_i]q_j)
=\tau\bigl(h\,i[p_i,q_j]\bigr),
\label{eq:faithful-covariance-derivative}
\end{equation}
and the symmetric bilinear Hessian is
\begin{equation}
D^2F(0)[h,k]_{ij}
=-\frac12\tau\!\left(
\bigl([h,[k,p_i]]+[k,[h,p_i]]\bigr)q_j
\right).
\label{eq:faithful-covariance-bilinear-hessian}
\end{equation}
In particular,
\begin{equation}
D^2F(0)[h,h]_{ij}
=-\tau([h,[h,p_i]]q_j).
\label{eq:faithful-covariance-quadratic-hessian}
\end{equation}
Let $A=DF(0)$, $K=\ker A$, $Y_0=\operatorname{ran}A$, $R=Y_0^\perp$, and let $P_R$ be the orthogonal projection of $Y$ onto $R$.

We identify $C=(c_{ij})\in R$ with its coefficient matrix in $M_{(a-1)\times(b-1)}(\R)$ and set
\begin{equation}
B_i:=\sum_{j<b}c_{ij}q_j\quad(i<a),
\qquad
B_a:=0.
\label{eq:hessian-coefficient-representatives}
\end{equation}
Since $C$ is orthogonal to $Y_0$, pairing \eqref{eq:faithful-covariance-derivative} with $C$ and using the nondegeneracy of the real $L^2$ pairing, we obtain
\begin{equation}
\sum_{i=1}^a[p_i,B_i]=0,
\qquad
p_r(B_r-B_s)p_s=0\quad(r\ne s).
\label{eq:hessian-cokernel-relations}
\end{equation}

For $C\in R$, define the scalar Hessian
\begin{equation}
Q_C(h):=\langle C,D^2F(0)[h,h]\rangle
=-\sum_{i=1}^a\tau([h,[h,p_i]]B_i).
\label{eq:hessian-scalar-form}
\end{equation}

\begin{lemma}
\label{lem:partition-hessian-indefiniteness}
For every nonzero $C\in R$, there are $\gamma>0$ and
infinite-dimensional real subspaces $K_+,K_-\subseteq K$ such that
\[
Q_C(h)\ge\gamma\norm h_2^2\quad(h\in K_+),
\qquad
Q_C(h)\le-\gamma\norm h_2^2\quad(h\in K_-).
\]
\end{lemma}

\begin{proof}
For $r\ne s$, set $D_{rs}=B_r-B_s$,
$X_r=p_rD_{rs}p_r$ and $X_s=p_sD_{rs}p_s$.
Block multiplication gives, for $x\in p_rMp_s$,
\begin{equation}
 Q_C(x+x^*)=2\bigl(\tau(x^*xX_s)-\tau(xx^*X_r)\bigr).
 \label{eq:hessian-corner-quadratic-form}
\end{equation}
Since $\tau(p_iq_j)=\alpha_i\beta_j$, the two compressions have
the same normalized trace:
\begin{equation}
 \frac{\tau(X_r)}{\alpha_r}
 =\frac{\tau(X_s)}{\alpha_s}
 =\tau(D_{rs})=:\mu_{rs}.
 \label{eq:hessian-common-corner-mean}
\end{equation}
If $X_r\notin\C p_r$ or $X_s\notin\C p_s$ for some $r\ne s$, we choose nonzero spectral
projections $e_+\le p_r$ and $f_+\le p_s$ for $X_r$ and $X_s$,
respectively, so that $X_re_+\le\mu_{rs}e_+$ and
$X_sf_+\ge\mu_{rs}f_+$, with a positive spectral gap in at least one
inequality. Making the same choices with the inequalities reversed, we
obtain projections $e_-\le p_r$ and $f_-\le p_s$.
By \eqref{eq:hessian-corner-quadratic-form}, $Q_C$ is positive
definite on the real space $\{x+x^*:x\in e_+Mf_+\}$ and negative
definite on $\{x+x^*:x\in e_-Mf_-\}$. The spectral gaps give a
positive lower bound for $Q_C(h)/\norm h_2^2$ on the first space and
for $-Q_C(h)/\norm h_2^2$ on the second, for nonzero $h$.
Both spaces are infinite dimensional because $M$ is a ${\rm II}_1$ factor.

It remains to consider the case
\begin{equation}
 p_r(B_r-B_s)p_r=\mu_{rs}p_r,\qquad
 p_s(B_r-B_s)p_s=\mu_{rs}p_s\quad(r\ne s).
 \label{eq:hessian-diagonal-degeneracy}
\end{equation}
Since $B_a=0$, applying these identities to pairs involving $a$ and
then to the pair $(r,i)$ gives
\begin{equation}
 p_rB_ip_r=\tau(B_i)p_r\quad(1\le i,r\le a).
 \label{eq:hessian-diagonal-block-identity}
\end{equation}
For distinct $r,s,t$, set $W_{rst}=p_s(2B_t-B_r-B_s)p_r$.
For $x\in p_rMp_t$, $y\in p_sMp_t$, the cross term is
\begin{equation}
\begin{aligned}
&Q_C(x+x^*+y+y^*)-Q_C(x+x^*)-Q_C(y+y^*)\\
&\hspace{1cm}=2\Ree\tau(xy^*W_{rst}).
\end{aligned}
\label{eq:hessian-two-corner-cross-term}
\end{equation}
Indeed, for $h=x+x^*$ and $k=y+y^*$, we have $hk=xy^*$ and $kh=yx^*$, while $hp_ik$ and $kp_ih$ equal $xy^*$ and $yx^*$ for $i=t$ and vanish for $i\ne t$. Hence, inserting $[h,[k,p_i]]+[k,[h,p_i]]=(hk+kh)p_i+p_i(hk+kh)-2hp_ik-2kp_ih$ into the polarization of \eqref{eq:hessian-scalar-form} and using $B_i^*=B_i$, we obtain \eqref{eq:hessian-two-corner-cross-term}.
By \eqref{eq:hessian-corner-quadratic-form} and
\eqref{eq:hessian-diagonal-degeneracy}, the two subtracted terms vanish. If $W_{rst}\ne0$, let $W_{rst}=v|W_{rst}|$ be its polar decomposition. We choose $\kappa>0$ and a nonzero spectral projection $e$ of $|W_{rst}|$ with $|W_{rst}|e\ge\kappa e$, and a nonzero projection $f\le p_t$. For $x\in eMf$,
\begin{equation}
 Q_C\bigl(x+x^*\pm(vx+x^*v^*)\bigr)
 =\pm2\tau(xx^*|W_{rst}|).
 \label{eq:hessian-opposite-sign-directions}
\end{equation}
For either sign, the displayed vector has squared $L^2$-norm
$4\|x\|_2^2$. Thus each of the two real-linear maps $x\mapsto x+x^*\pm(vx+x^*v^*)$ on $eMf$ is injective, and its range is infinite dimensional because $M$ is a ${\rm II}_1$ factor. On the ranges corresponding to
the plus and minus signs, respectively, we have
$Q_C(h)\ge(\kappa/2)\norm h_2^2$ and
$Q_C(h)\le-(\kappa/2)\norm h_2^2$.

Suppose now that all $W_{rst}$ vanish. By the second relation in \eqref{eq:hessian-cokernel-relations}, with $r$ and $s$ interchanged, we have
\begin{equation}
 p_sB_rp_r=p_sB_sp_r\quad(r\ne s).
 \label{eq:hessian-off-diagonal-identity}
\end{equation}
The blocks in \eqref{eq:hessian-off-diagonal-identity} vanish when
either index is $a$. When neither is $a$, the same conclusion follows
by taking $t=a$ in $W_{rst}=0$. Taking $t=i$ also gives
$p_sB_ip_r=0$ for $i\notin\{r,s\}$. Thus all off-diagonal blocks of each
$B_i$ vanish, and by \eqref{eq:hessian-diagonal-block-identity}, $B_i$ is
scalar.
Since $B_i$ has zero coefficient at the nonzero projection $q_b$, we have $B_i=0$,
and hence $C=0$, a contradiction. The argument also covers $a=2$,
when there are no triples of distinct indices.

Finally, let $K_+$ and $K_-$ be the intersections of the subspaces obtained above on which $Q_C$ is positive definite and negative definite, respectively, with $\ker A$. Since $\ker A$ is defined by finitely many real-linear conditions, both spaces remain infinite dimensional.
The inequalities $Q_C(h)\ge\gamma\norm h_2^2$ on $K_+$ and
$Q_C(h)\le-\gamma\norm h_2^2$ on $K_-$ still hold for a common
$\gamma>0$.
\end{proof}

\begin{lemma}
\label{lem:partition-quadratic-loop}
Under the preceding assumptions, suppose $R\ne0$ and set
\[
\mathcal B(h,k):=P_RD^2F(0)[h,k]
\qquad(h,k\in K).
\]
There are $L^2$-orthonormal elements $k_1,\ldots,k_N\in K$ and numbers
$\lambda_1,\ldots,\lambda_N>0$ such that, with
\[
z_m:=\mathcal B(k_m,k_m),
\]
we have \begin{equation} \mathcal B(k_m,k_\ell)=0\quad(m\ne\ell),
\qquad
\operatorname{span}\{z_1,\ldots,z_N\}=R,
\qquad
\sum_m\lambda_mz_m=0.
\label{eq:partition-positive-quadratic-loop}
\end{equation}
Consequently, for $k_0:=\sum_m\sqrt{\lambda_m}\,k_m$,
\begin{equation}
\mathcal B(k_0,k_0)=0,
\qquad
\mathcal B(k_0,k_m)=\sqrt{\lambda_m}\,z_m,
\label{eq:partition-loop-surjectivity}
\end{equation}
and the vectors $\mathcal B(k_0,k_m)$ span $R$.
\end{lemma}

\begin{proof}
Let $d=(a-1)(b-1)$. Applying H\"older's inequality to the double commutator, we obtain
\begin{equation}
|\tau([k,[k,p_i]]q_j)|\le4\norm k_2^2,
\qquad \norm{\mathcal B(k,k)}\le4\sqrt d\,\norm k_2^2.
\label{eq:partition-quadratic-bound}
\end{equation}
For each unit vector $c\in R$, we have $\langle c,\mathcal B(k,k)\rangle=Q_c(k)$ for $k\in K$, since $P_R$ is the orthogonal projection onto $R$. Hence, by Lemma~\ref{lem:partition-hessian-indefiniteness} with $C=c$,
there exist an infinite-dimensional real subspace $L_c\subseteq K$ and $\gamma_c>0$
with
\begin{equation}
\langle c,\mathcal B(k,k)\rangle\le-\gamma_c\norm k_2^2
\quad(k\in L_c).
\label{eq:partition-negative-sign-space}
\end{equation}
By \eqref{eq:partition-quadratic-bound}, there is a neighborhood
$U_c$ of $c$ such that $\langle c',\mathcal B(k,k)\rangle\le-(\gamma_c/2)\norm k_2^2$ for all $c'\in U_c$ and $k\in L_c$.
We choose a finite cover of the unit sphere of $R$ by $U_{c_1},\ldots,U_{c_N}$ and recursively choose $L^2$-unit vectors $k_\ell\in L_{c_\ell}$ satisfying
\begin{equation}
\tau(k_\ell k_m)=0,\qquad \mathcal B(k_\ell,k_m)=0\quad(m<\ell).
\label{eq:partition-loop-orthogonality}
\end{equation}
Indeed, at each step, \eqref{eq:partition-loop-orthogonality} imposes only finitely many real-linear conditions on $k_\ell$, so the elements of $L_{c_\ell}$ satisfying them form an infinite-dimensional subspace.
For $z_\ell=\mathcal B(k_\ell,k_\ell)$, the cover gives
\begin{equation}
\text{for every }0\ne c\in R,\quad
\langle c,z_\ell\rangle<0\text{ for some }\ell.
\label{eq:partition-loop-separation}
\end{equation}
The cone generated by $z_1,\ldots,z_N$ is closed and convex. If it were not all of $R$, some $c\ne0$ would satisfy $\langle c,z_\ell\rangle\ge0$ for all $\ell$, contradicting \eqref{eq:partition-loop-separation}. Hence
\begin{equation}
\operatorname{cone}(z_1,\ldots,z_N)=R.
\label{eq:partition-loop-full-cone}
\end{equation}
Write $-z_\ell=\sum_m\theta_{\ell m}z_m$ with $\theta_{\ell m}\ge0$.
Then $\lambda_m=1+\sum_\ell\theta_{\ell m}>0$ satisfies
$\sum_m\lambda_mz_m=0$, proving
\eqref{eq:partition-positive-quadratic-loop}.
By bilinearity and the relations $\mathcal B(k_m,k_\ell)=0$ for $m\ne\ell$, we
obtain \eqref{eq:partition-loop-surjectivity}.
\end{proof}

\begin{proof}[Second-order proof of Corollary~\ref{thm:faithful-relative-partition-stability}]
Suppose that the corollary fails for fixed $a,b,\delta,\rho$. Then we can choose
${\rm II}_1$ factors $(M_n,\tau_n)$ and pairs $(\mathbf p_n,\mathbf q_n)$ of partitions of $1$ with all
traces at least $\delta$, satisfying
\begin{equation}
 \max_{i,j}|\Delta_{ij}(\mathbf p_n,\mathbf q_n)|<1/n,
 \label{eq:faithful-bad-sequence-defect}
\end{equation}
but for which no $h=h^*\in M_n$ with $\|h\|_\infty<\rho$ makes
$(e^{ih}p_{i,n}e^{-ih})_i$ independent of $(q_{j,n})_j$.
In a nonprincipal tracial ultraproduct $M=\prod_\omega(M_n,\tau_n)$,
a ${\rm II}_1$ factor, the partitions $\mathbf p$ and $\mathbf q$ with $p_i=(p_{i,n})_\omega$ and $q_j=(q_{j,n})_\omega$ are independent and have positive traces. Let $F$ be the map \eqref{eq:faithful-covariance-map} for these partitions.

We choose a basis $(e_\nu)$ of a complement $E\subset M_{\rm sa}$ to
$\ker DF(0)$. If $R\ne0$, let $k_1,\ldots,k_N$ be as in Lemma~\ref{lem:partition-quadratic-loop}. If $R=0$, we take no $k_m$. We fix
uniformly bounded self-adjoint representatives
$(e_{\nu,n})_n$ of $e_\nu$ and $(k_{m,n})_n$ of $k_m$, and set
\[
G_n(s,t)=\sum_\nu s_\nu e_{\nu,n}+\sum_m t_mk_{m,n},
\qquad f=F\circ G,
\]
where $G(s,t)=\sum_\nu s_\nu e_\nu+\sum_m t_mk_m$. Let
\begin{equation}
 f_n(s,t)_{ij}
 =\tau_n(e^{iG_n(s,t)}p_{i,n}e^{-iG_n(s,t)}q_{j,n})
 -\tau_n(p_{i,n})\tau_n(q_{j,n})\quad(i<a,\ j<b).
 \label{eq:faithful-sequence-jet-convergence}
\end{equation}
The values and first two derivatives of $f_n$ at zero converge to
those of $f$ along $\omega$. Passing to a subsequence, we may assume that these convergences
hold in the ordinary sense. By differentiating the exponentials three times and using the
uniformly bounded representatives, we obtain uniform bounds for the third derivatives on a fixed
neighborhood of zero.
Since $G$ is linear and $DF(0)$ maps $E$ onto $Y_0$, the map $Df(0)=DF(0)\circ G$ has range $Y_0$, and $D^2f(0)[x,x']=D^2F(0)[G(x),G(x')]$. If $R\ne0$, then $G$ maps the coordinate vectors of $t_1,\ldots,t_N$ to $k_1,\ldots,k_N\in K$, so these coordinate vectors lie in $\ker Df(0)$ and, by \eqref{eq:partition-positive-quadratic-loop}, satisfy \eqref{eq:second-order-positive-loop-hypotheses} for the decomposition $Y=Y_0\oplus R$.
Hence, by Lemma~\ref{lem:second-order-inverse}, there are $x_n\to0$ with $f_n(x_n)=0$.
Then $\norm{G_n(x_n)}_\infty\to0$. For $h=G_n(x_n)$, since $f_n(x_n)=0$, the covariances of $(e^{ih}p_{i,n}e^{-ih})_i$ and $(q_{j,n})_j$ with $i<a$ and $j<b$ vanish. Since the covariance matrix has zero row and column sums by \eqref{eq:covariance-zero-row-column-sums}, the covariances with $i=a$ or $j=b$ also vanish.
Thus, for large $n$, we have $\norm h_\infty<\rho$, and $h$ makes $(e^{ih}p_{i,n}e^{-ih})_i$ independent of $(q_{j,n})_j$, contradicting
the choice of the sequence.
\end{proof}

\end{document}